# APPLICATION OF BIMATRICES TO SOME FUZZY AND NEUTROSOPHIC MODELS

**W. B. Vasantha Kandasamy**
**Florentin Smarandache**
**K. Ilanthenral**

# Application of Bimatrices to some Fuzzy and Neutrosophic Models

**W. B. Vasantha Kandasamy**
Department of Mathematics
Indian Institute of Technology, Madras
Chennai – 600036, India
e-mail: vasantha@iitm.ac.in
web: http://mat.iitm.ac.in/~wbv

**Florentin Smarandache**
Department of Mathematics
University of New Mexico
Gallup, NM 87301, USA
e-mail: smarand@unm.edu

**K. Ilanthenral**
Editor, Maths Tiger, Quarterly Journal
Flat No.11, Mayura Park,
16, Kazhikundram Main Road, Tharamani,
Chennai – 600 113, India
e-mail: ilanthenral@gmail.com

**2005**

# CONTENTS

## Chapter Four

## NEUTROSOPHIC BIGRAPHS THEIR GENERALIZATIONS AND APPLICATIONS TO NEUTROSOPHIC MODELS

## Preface

Graphs and matrices play a vital role in the analysis and study of several of the real world problems which are based only on unsupervised data. The fuzzy and neutrosophic tools like fuzzy cognitive maps invented by Kosko and neutrosophic cognitive maps introduced by us help in the analysis of such real world problems and they happen to be mathematical tools which can give the hidden pattern of the problem under investigation. This book, in order to generalize the two models, has systematically invented mathematical tools like bimatrices, trimatrices, n-matrices, bigraphs, trigraphs and n-graphs and describe some of its properties. These concepts are also extended neutrosophically in this book.

Using these new tools we define fuzzy cognitive bimaps, fuzzy cognitive trimaps, fuzzy relational bimaps and fuzzy relational trimaps which exploit the new notions of bimatrices, bigraphs, trimatrices and trigraphs. It is worth mentioning that these can be extended to n-array models as a simple exercise. The main advantage of these models are:

1. Comparative study of the views stage-by-stage is possible.
2. Several experts' opinion can be compared not only at each stage but also the final result can easily be compared.
3. Saves time and economy.

This book is organized into four chapters. The first chapter recalls the basic concepts of bimatrices and neutrosophic bimatrices. Second chapter introduces several new notions of graphs like bigraphs, trigraphs and their properties. Chapter three illustrates how these new tools are used in the construction of fuzzy cognitive bimaps, trimaps and n-maps and fuzzy relational bimaps, trimaps and n-maps. The neutrosophic analogues of chapter 3 is carried out in the fourth chapter.

In this book, we have given nearly 95 examples to make the reader easily follow the definitions. We have also provided 125 figures for help in easily understanding the definition and examples. We have also given 25 real world problems as applications of Bimatrices to Fuzzy and Neutrosophic models.

This book contains generalization of some of the results given in our earlier book Fuzzy Cognitive Maps and Neutrosophic Cognitive Maps (2003). Those who wish to simultaneously study several models that are at times time dependent or when a comparative analysis is needed can use these models.

On a personal note, we thank Dr.K.Kandasamy for proof reading this book.

Finally, we dedicate this book to Prof. Dr. Bart Kosko, the father of Fuzzy Cognitive Maps.

W.B.VASANTHA KANDASAMY
FLORENTIN SMARANDACHE
K.ILANTHENRAL

**Chapter One**

# BASIC CONCEPTS ON BIMATRICES

In this chapter we just briefly recall the definition of bimatrices and neutrosophic bimatrices, illustrate them with examples. This chapter has two sections. In the first section we introduce the definition of bimatrices and in section two the new notion of neutrosophic bimatrices are introduced and illustrated with examples.

## 1.1 Definition and Basic Operations on Bimatrices

In this section we recall the definition of bimatrix and illustrate them with examples and define some of the basic operations on it.

**DEFINITION 1.1.1:** *A bimatrix $A_B$ is defined as the union of two rectangular array of numbers $A_1$ and $A_2$ arranged into rows and columns. It is written as follows $A_B = A_1 \cup A_2$ where $A_1 \neq A_2$ with*

$$A_1 = \begin{bmatrix} a_{11}^1 & a_{12}^1 & \cdots & a_{1n}^1 \\ a_{21}^1 & a_{22}^1 & \cdots & a_{2n}^1 \\ \vdots & & & \\ a_{m1}^1 & a_{m2}^1 & \cdots & a_{mn}^1 \end{bmatrix}$$

*and*

$$A_2 = \begin{bmatrix} a_{11}^2 & a_{12}^2 & \cdots & a_{1n}^2 \\ a_{21}^2 & a_{22}^2 & \cdots & a_{2n}^2 \\ \vdots & & & \\ a_{m1}^2 & a_{m2}^2 & \cdots & a_{mn}^2 \end{bmatrix}$$

*'$\cup$' is just the notational convenience (symbol) only.*

The above array is called a m by n bimatrix (written as B(m × n) since each of $A_i$ (i = 1, 2) has m rows and n columns. It is to be noted a bimatrix has no numerical value associated with it. It is only a convenient way of representing a pair of arrays of numbers.

*Note:* If $A_1 = A_2$ then $A_B = A_1 \cup A_2$ is not a bimatrix. A bimatrix $A_B$ is denoted by $\left(a_{ij}^1\right) \cup \left(a_{ij}^2\right)$. If both $A_1$ and $A_2$ are m × n matrices then the bimatrix $A_B$ is called the m × n rectangular bimatrix. But we make an assumption the zero bimatrix is a union of two zero matrices even if $A_1$ and $A_2$ are one and the same; i.e., $A_1 = A_2 = (0)$.

***Example 1.1.1:*** The following are bimatrices:

i. $$A_B = \begin{bmatrix} 3 & 0 & 1 \\ -1 & 2 & 1 \end{bmatrix} \cup \begin{bmatrix} 0 & 2 & -1 \\ 1 & 1 & 0 \end{bmatrix}$$

is a rectangular 2 × 3 bimatrix.

ii. $$A'_B = \begin{bmatrix} 3 \\ 1 \\ 2 \end{bmatrix} \cup \begin{bmatrix} 0 \\ -1 \\ 0 \end{bmatrix}$$

is a column bimatrix.

iii. $A''_B = (3, -2, 0, 1, 1) \cup (1, 1, -1, 1, 2)$

is a row bimatrix.

In a bimatrix $A_B = A_1 \cup A_2$ if both $A_1$ and $A_2$ are $m \times n$ rectangular matrices then the bimatrix $A_B$ is called the rectangular $m \times n$ bimatrix.

**DEFINITION 1.1.2:** *Let $A_B = A_1 \cup A_2$ be a bimatrix. If both $A_1$ and $A_2$ are square matrices then $A_B$ is called the square bimatrix.*

*If one of the matrices in the bimatrix $A_B = A_1 \cup A_2$ is square and other is rectangular or if both $A_1$ and $A_2$ are rectangular matrices say $m_1 \times n_1$ and $m_2 \times n_2$ with $m_1 \neq m_2$ or $n_1 \neq n_2$ then we say $A_B$ is a mixed bimatrix.*

The following are examples of square bimatrices and mixed bimatrices.

***Example 1.1.2:*** Given

$$A_B = \begin{bmatrix} 3 & 0 & 1 \\ 2 & 1 & 1 \\ -1 & 1 & 0 \end{bmatrix} \cup \begin{bmatrix} 4 & 1 & 1 \\ 2 & 1 & 0 \\ 0 & 0 & 1 \end{bmatrix}$$

is a $3 \times 3$ square bimatrix.

$$A'_B = \begin{bmatrix} 1 & 1 & 0 & 0 \\ 2 & 0 & 0 & 1 \\ 0 & 0 & 0 & 3 \\ 1 & 0 & 1 & 2 \end{bmatrix} \cup \begin{bmatrix} 2 & 0 & 0 & -1 \\ -1 & 0 & 1 & 0 \\ 0 & -1 & 0 & 3 \\ -3 & -2 & 0 & 0 \end{bmatrix}$$

is a $4 \times 4$ square bimatrix.

***Example 1.1.3:*** Let

$$A_B = \begin{bmatrix} 3 & 0 & 1 & 2 \\ 0 & 0 & 1 & 1 \\ 2 & 1 & 0 & 0 \\ 1 & 0 & 1 & 0 \end{bmatrix} \cup \begin{bmatrix} 1 & 1 & 2 \\ 0 & 2 & 1 \\ 0 & 0 & 4 \end{bmatrix}$$

then $A_B$ is a mixed square bimatrix.
Let

$$A'_B = \begin{bmatrix} 2 & 0 & 1 & 1 \\ 0 & 1 & 0 & 1 \\ -1 & 0 & 2 & 1 \end{bmatrix} \cup \begin{bmatrix} 2 & 0 \\ 4 & -3 \end{bmatrix}$$

$A'_B$ is a mixed bimatrix.

Now we proceed on to give the bimatrix operations.

Let $A_B = A_1 \cup A_2$ and $C_B = C_1 \cup C_2$ be two bimatrices we say $A_B$ and $C_B$ are equal written as $A_B = C_B$ if and only if $A_1$ and $C_1$ are identical and $A_2$ and $C_2$ are identical i.e., $A_1 = C_1$ and $A_2 = C_2$.

If $A_B = A_1 \cup A_2$ and $C_B = C_1 \cup C_2$, we say $A_B$ is not equal to $C_B$ we write $A_B \neq C_B$ if and only if $A_1 \neq C_1$ or $A_2 \neq C_2$.

***Example 1.1.4:*** Let

$$A_B = \begin{bmatrix} 3 & 2 & 0 \\ 2 & 1 & 1 \end{bmatrix} \cup \begin{bmatrix} 0 & -1 & 2 \\ 0 & 0 & 1 \end{bmatrix}$$

and

$$C_B = \begin{bmatrix} 1 & 1 & 1 \\ 0 & 0 & 0 \end{bmatrix} \cup \begin{bmatrix} 2 & 0 & 1 \\ 1 & 0 & 2 \end{bmatrix}$$

clearly $A_B \neq C_B$. Let

$$A_B = \begin{bmatrix} 0 & 0 & 1 \\ 1 & 1 & 2 \end{bmatrix} \cup \begin{bmatrix} 0 & 4 & -2 \\ -3 & 0 & 0 \end{bmatrix}$$

$$C_B = \begin{bmatrix} 0 & 0 & 1 \\ 1 & 1 & 2 \end{bmatrix} \cup \begin{bmatrix} 0 & 0 & 0 \\ 1 & 0 & 1 \end{bmatrix}$$

clearly $A_B \neq C_B$.

If $A_B = C_B$ then we have $C_B = A_B$.

We now proceed on to define multiplication by a scalar. Given a bimatrix $A_B = A_1 \cup B_1$ and a scalar $\lambda$, the product of $\lambda$ and $A_B$ written $\lambda A_B$ is defined to be

$$\lambda A_B = \begin{bmatrix} \lambda a_{11} & \cdots & \lambda a_{1n} \\ \vdots & & \vdots \\ \lambda a_{m1} & \cdots & \lambda a_{mn} \end{bmatrix} \cup \begin{bmatrix} \lambda b_{11} & \cdots & \lambda b_{1n} \\ \vdots & & \vdots \\ \lambda b_{m1} & \cdots & \lambda b_{mn} \end{bmatrix}$$

each element of $A_1$ and $B_1$ are multiplied by $\lambda$. The product $\lambda A_B$ is then another bimatrix having m rows and n columns if $A_B$ has m rows and n columns.

We write

$$\begin{aligned} \lambda A_B &= \left[\lambda a_{ij}\right] \cup \left[\lambda b_{ij}\right] \\ &= \left[a_{ij}\lambda\right] \cup \left[b_{ij}\lambda\right] \\ &= A_B \lambda. \end{aligned}$$

***Example 1.1.5:*** Let

$$A_B = \begin{bmatrix} 2 & 0 & 1 \\ 3 & 3 & -1 \end{bmatrix} \cup \begin{bmatrix} 0 & 1 & -1 \\ 2 & 1 & 0 \end{bmatrix}$$

and $\lambda = 3$ then

$$3A_B = \begin{bmatrix} 6 & 0 & 3 \\ 9 & 9 & -3 \end{bmatrix} \cup \begin{bmatrix} 0 & 3 & -3 \\ 6 & 3 & 0 \end{bmatrix}.$$

If $\lambda = -2$ then, for

$$A_B = [3\ 1\ 2\ -4] \cup [0\ 1\ -1\ 0],$$

$$\lambda A_B = [-6\ -2\ -4\ 8] \cup [0\ -2\ 2\ 0].$$

Let $A_B = A_1 \cup B_1$ and $C_B = A_2 \cup B_2$ be any two $m \times n$ bimatrices. The sum $D_B$ of the bimatrices $A_B$ and $C_B$ is defined as $D_B = A_B + C_B = [A_1 \cup B_1] + [A_2 \cup B_2] = (A_1 + A_2) \cup [B_2 + B_2]$; where $A_1 + A_2$ and $B_1 + B_2$ are the usual addition of matrices i.e., if

$$A_B = \left(a^1_{ij}\right) \cup \left(b^1_{ij}\right)$$

and

$$C_B = \left(a^2_{ij}\right) \cup \left(b^2_{ij}\right)$$

then

$$A_B + C_B = D_B = \left(a^1_{ij} + b^2_{ij}\right) \cup \left(b^1_{ij} + b^2_{ij}\right)\ \left(\forall\, ij\right).$$

If we write in detail

$$A_B = \begin{bmatrix} a^1_{11} & \cdots & a^1_{1n} \\ \vdots & & \\ a^1_{m1} & \cdots & a^1_{mn} \end{bmatrix} \cup \begin{bmatrix} b^1_{11} & \cdots & b^1_{1n} \\ \vdots & & \\ b^1_{m1} & \cdots & b^1_{mn} \end{bmatrix}$$

$$C_B = \begin{bmatrix} a^2_{11} & \cdots & a^2_{1n} \\ \vdots & & \\ a^2_{m1} & \cdots & a^2_{mn} \end{bmatrix} \cup \begin{bmatrix} b^2_{11} & \cdots & b^2_{1n} \\ \vdots & & \\ b^2_{m1} & \cdots & b^2_{mn} \end{bmatrix}$$

$A_B + C_B =$

$$\begin{bmatrix} a^1_{11} + a^2_{11} & \dots & a^1_{1n} a^2_{1n} \\ \vdots & & \vdots \\ a^1_{m1} + a^2_{m1} & \dots & a^1_{mn} + a^2_{mn} \end{bmatrix} \cup \begin{bmatrix} b^1_{11} + b^2_{11} & \dots & b^1_{1n} + b^2_{1n} \\ \vdots & & \vdots \\ b^1_{m1} + b^2_{m1} & \dots & b^1_{mn} + b^2_{mn} \end{bmatrix}.$$

The expression is abbreviated to

$$\begin{aligned} D_B &= A_B + C_B \\ &= (A_1 \cup B_1) + (A_2 \cup B_2) \end{aligned}$$

$= \quad (A_1 + A_2) \cup (B_1 + B_2).$

Thus two bimatrices are added by adding the corresponding elements only when compatibility of usual matrix addition exists.

*Note*: If $A_B = A_1 \cup A_2$ be a bimatrix we call $A_1$ and $A_2$ as the components matrices of the bimatrix $A_B$.

## 1.2. Neutrosophic Bimatrices

Here for the first time we define the notion of neutrosophic bimatrix and illustrate them with examples. Also we define fuzzy neutrosophic bimatrices.

**DEFINITION 1.2.1:** *Let $A = A_1 \cup A_2$ where $A_1$ and $A_2$ are two distinct neutrosophic matrices with entries from a neutrosophic field. Then $A = A_1 \cup A_2$ is called the neutrosophic bimatrix.*

*It is important to note the following:*

(1) *If both $A_1$and $A_2$ are neutrosophic matrices we call A a neutrosophic bimatrix.*

(2) *If only one of $A_1$ or $A_2$ is a neutrosophic matrix and other is not a neutrosophic matrix then we all $A = A_1 \cup A_2$ as the semi neutrosophic bimatrix. (It is clear all neutrosophic bimatrices are trivially semi neutrosophic bimatrices).*

*It both $A_1$ and $A_2$ are $m \times n$ neutrosophic matrices then we call $A = A_1 \cup A_2$ a $m \times n$ neutrosophic bimatrix or a rectangular neutrosophic bimatrix.*

*If $A = A_1 \cup A_2$ be such that $A_1$ and $A_2$ are both $n \times n$ neutrosophic matrices then we call $A = A_1 \cup A_2$ a square or a $n \times n$ neutrosophic bimatrix. If in the neutrosophic* *bimatrix $A = A_1 \cup A_2$ both $A_1$ and $A_2$ are square matrices but of different order say $A_1$ is a $n \times n$ matrix and $A_2$ a $s \times s$ matrix then we call $A = A_1 \cup A_2$ a mixed neutrosophic square bimatrix. (Similarly one can define mixed square semi neutrosophic bimatrix).*

*Likewise in $A = A_1 \cup A_2$ if both $A_1$ and $A_2$ are rectangular matrices say $A_1$ is a $m \times n$ matrix and $A_2$ is a $p \times q$ matrix then we call $A = A_1 \cup A_2$ a mixed neutrosophic rectangular bimatrix. (If $A = A_1 \cup A_2$ is a semi neutrosophic bimatrix then we call A the mixed rectangular semi neutrosophic bimatrix).*

Just for the sake of clarity we give some illustrations.

**Notation:** We in this book denote a neutrosophic bimatrix by $A_N = A_1 \cup A_2$.

***Example 1.2.1:*** Let

$$A_N = \begin{bmatrix} 0 & I & 0 \\ 1 & 2 & -1 \\ 3 & 2 & I \end{bmatrix} \cup \begin{bmatrix} 2 & I & 1 \\ I & 0 & I \\ 1 & 1 & 2 \end{bmatrix}$$

$A_N$ is the $3 \times 3$ square neutrosophic bimatrix.

***Example 1.2.2:*** Let

$$A_N = \begin{bmatrix} 2 & 0 & I \\ 4 & I & 1 \\ 1 & 1 & 2 \end{bmatrix} \cup \begin{bmatrix} 3 & I & 0 & 1 & 5 \\ 0 & 0 & I & 3 & 1 \\ I & 0 & 0 & I & 2 \\ 1 & 3 & 3 & 5 & 4 \\ 2 & 1 & 3 & 0 & I \end{bmatrix}$$

$A_N$ is a mixed square neutrosophic bimatrix.

***Example 1.2.3:*** Let

$$A_N = \begin{bmatrix} 3 & 1 & 1 & 1 & I \\ I & 0 & 2 & 3 & 4 \end{bmatrix} \cup \begin{bmatrix} I & 2 & 0 & I \\ 3 & 1 & 2 & 1 \\ 4 & 1 & 0 & 0 \\ 3 & 3 & 1 & 1 \\ 1 & I & 0 & I \end{bmatrix}$$

$A_N$ is a mixed rectangular neutrosophic bimatrix. We denote $A_N$ by

$$A_N = A_1 \cup A_2 = A_1^{2\times 5} \cup A_2^{5\times 4}.$$

***Example 1.2.4:*** Let

$$A_N = \begin{bmatrix} 3 & 1 \\ 1 & 2 \\ I & 0 \\ 3 & I \end{bmatrix} \cup \begin{bmatrix} 3 & 3 \\ 1 & I \\ 4 & I \\ -I & 0 \end{bmatrix}$$

$A_N$ is $4 \times 2$ rectangular neutrosophic bimatrix.

***Example 1.2.5*****:** Let

$$A_N = \begin{bmatrix} 3 & 1 & 1 \\ 2 & 2 & 2 \end{bmatrix} \cup \begin{bmatrix} -I & 1 & 2 \\ 0 & I & 3 \end{bmatrix}$$

$A_N$ is a rectangular semi neutrosophic bimatrix for $A = A_1 \cup A_2$ with

$$A_1 = \begin{bmatrix} 3 & 1 & 1 \\ 2 & 2 & 2 \end{bmatrix}$$

is not a neutrosophic matrix, only

$$A_2 = \begin{bmatrix} -I & 1 & 2 \\ 0 & I & 3 \end{bmatrix}$$

is a neutrosophic matrix.

***Example 1.2.6:*** Let

$$A_N = \begin{bmatrix} 3 & 1 & 1 & 1 \\ 0 & I & 1 & 2 \\ 0 & 0 & 0 & 3 \\ I & 1 & 1 & 1 \end{bmatrix} \cup \begin{bmatrix} 0 & 2 & 2 & 2 \\ 1 & 0 & 0 & 0 \\ 2 & 0 & 0 & 1 \\ 5 & 0 & -1 & 2 \end{bmatrix}$$

$A_N$ is a square semi neutrosophic bimatrix.

***Example 1.2.7:*** Let

$$A_N = \begin{bmatrix} 1 & 1 & 1 & 1 \\ 0 & 0 & 0 & 0 \end{bmatrix} \cup \begin{bmatrix} 0 \\ I \\ 1 \\ 2 \\ 3 \end{bmatrix}.$$

$A_N$ is a rectangular mixed semi neutrosophic bimatrix.

Thus as in case of bimatrices we may have square, mixed square, rectangular or mixed rectangular neutrosophic (semi neutrosophic) bimatrices.

Now we can also define the neutrosophic bimatrices or semi neutrosophic bimatrices over different fields. When both $A_1$ and $A_2$ in the bimatrix $A = A_1 \cup A_2$ take its values from the same neutrosophic field K we call it a neutrosophic (semi neutrosophic) bimatrix. If in the bimatrix $A = A_1 \cup A_2$, $A_1$ is defined over a neutrosophic field F and $A_2$ over some other neutrosophic field $F^1$ then we call the neutrosophic bimatrix as strong neutrosophic bimatrix $F^1 \not\subset F$ or $F \not\subset F^1$. If on the other hand the neutrosophic bimatrix $A = A_1 \cup A_2$ is such that $A_1$ takes entries from the neutrosophic field F and $A_2$ takes its entries from a proper subfield of a neutrosophic field then we call A the weak neutrosophic bimatrix. All properties of bimatrices can be carried on to neutrosophic bimatrices and semi neutrosophic bimatrices.

Now we proceed on to define fuzzy bimatrix, fuzzy neutrosophic bimatrix, semi-fuzzy bimatrix, semi fuzzy neutrosophic bimatrix and illustrate them with examples.

**DEFINITION 1.2.2:** *Let $A = A_1 \cup A_2$ where $A_1$ and $A_2$ are two distinct fuzzy matrices with entries from the interval [0, 1]. Then $A = A_1 \cup A_2$ is called the fuzzy bimatrix.*

*It is important to note the following:*

1. *If both $A_1$ and $A_2$ are fuzzy matrices we call A a fuzzy bimatrix.*
2. *If only one of $A_1$ or $A_2$ is a fuzzy matrix and other is not a fuzzy matrix then we all $A = A_1 \cup A_2$ as the semi fuzzy bimatrix. (It is clear all fuzzy matrices are trivially semi fuzzy matrices).*

*It both $A_1$ and $A_2$ are $m \times n$ fuzzy matrices then we call $A = A_1 \cup A_2$ a $m \times n$ fuzzy bimatrix or a rectangular fuzzy bimatrix.*

*If $A = A_1 \cup A_2$ is such that $A_1$ and $A_2$ are both $n \times n$ fuzzy matrices then we call $A = A_1 \cup A_2$ a square or a $n \times n$ fuzzy bimatrix. If in the fuzzy bimatrix $A = A_1 \cup A_2$ both $A_1$ and $A_2$ are square matrices but of different order say $A_1$ is a $n \times n$ matrix and $A_2$ a $s \times s$ matrix then we call $A = A_1 \cup A_2$ a mixed fuzzy square bimatrix. (Similarly one can define mixed square semi fuzzy bimatrix).*

*Likewise in $A = A_1 \cup A_2$ if both $A_1$ and $A_2$ are rectangular matrices say $A_1$ is a $m \times n$ matrix and $A_2$ is a $p \times q$ matrix then we call $A = A_1 \cup A_2$ a mixed fuzzy rectangular bimatrix. (If $A = A_1 \cup A_2$ is a semi fuzzy bimatrix then we call A the mixed rectangular semi fuzzy bimatrix).*

Just for the sake of clarity we give some illustration.

**Notation:** We denote a fuzzy bimatrix by $A_F = A_1 \cup A_2$.

***Example 1.2.8:*** Let

$$A_F = \begin{bmatrix} 0 & .1 & 0 \\ .1 & .2 & .1 \\ .3 & .2 & .1 \end{bmatrix} \cup \begin{bmatrix} .2 & .1 & .1 \\ .1 & 0 & .1 \\ .2 & .1 & .2 \end{bmatrix}$$

$A_F$ is the 3 × 3 square fuzzy bimatrix.

***Example 1.2.9:*** Let

$$A_F = \begin{bmatrix} .2 & 0 & 1 \\ .4 & .2 & 1 \\ .3 & 1 & .2 \end{bmatrix} \cup \begin{bmatrix} .3 & 1 & 0 & .4 & .5 \\ 0 & 0 & 1 & .8 & .2 \\ 1 & 0 & 0 & .1 & .2 \\ .1 & .3 & .3 & .5 & .4 \\ .2 & .1 & .3 & 0 & 1 \end{bmatrix}$$

$A_F$ is a mixed square fuzzy bimatrix.

***Example 1.2.10:*** Let

$$A_F = \begin{bmatrix} .3 & 1 & .5 & 1 & .9 \\ .6 & 0 & .2 & .3 & .4 \end{bmatrix} \cup \begin{bmatrix} 1 & .2 & 0 & 0 \\ .3 & 1 & .2 & 1 \\ .4 & 1 & 0 & 0 \\ .3 & .3 & .2 & 1 \\ 1 & .5 & .7 & .6 \end{bmatrix}$$

$A_F$ is a mixed rectangular fuzzy bimatrix. We denote $A_F$ by

$$A_F = A_1 \cup A_2 = A_1^{2\times 5} \cup A_2^{5\times 4}.$$

***Example 1.2.11:*** Let

$$A_F = \begin{bmatrix} .3 & 1 \\ 1 & .2 \\ .5 & 0 \\ .3 & .6 \end{bmatrix} \cup \begin{bmatrix} .3 & .7 \\ 1 & 1 \\ .4 & 1 \\ .2 & 0 \end{bmatrix}$$

$A_F$ is 4 × 2 rectangular fuzzy bimatrix.

***Example 1.2.12*:** Let

$$A_F = \begin{bmatrix} 3 & 1 & 1 \\ 2 & 2 & 2 \end{bmatrix} \cup \begin{bmatrix} .5 & .7 & .2 \\ 0 & .1 & .3 \end{bmatrix}$$

$A_F$ is a rectangular semi fuzzy bimatrix for

$$A = A_1 \cup A_2$$

with

$$A_1 = \begin{bmatrix} 3 & 1 & 1 \\ 2 & 2 & 2 \end{bmatrix}$$

which is not a fuzzy matrix, only

$$A_2 = \begin{bmatrix} .5 & .7 & .2 \\ 0 & .1 & .3 \end{bmatrix}$$

is a fuzzy matrix.

***Example 1.2.13:*** Let

$$A_F = \begin{bmatrix} .3 & 1 & 1 & 1 \\ 0 & 0 & .1 & .2 \\ 0 & 0 & 0 & .3 \\ .3 & 1 & 1 & 1 \end{bmatrix} \cup \begin{bmatrix} 0 & 2 & 2 & 2 \\ 1 & 0 & 0 & 0 \\ 2 & 0 & 0 & 1 \\ 5 & 0 & -1 & 2 \end{bmatrix}$$

$A_N$ is a square semi fuzzy bimatrix.

***Example 1.2.14:*** Let

$$A_F = \begin{bmatrix} 1 & 1 & 1 & 1 \\ 0 & 0 & 0 & 0 \end{bmatrix} \cup \begin{bmatrix} 0 \\ 1 \\ 1 \\ 2 \\ 3 \end{bmatrix}.$$

$A_F$ is a rectangular mixed semi fuzzy bimatrix.

Thus as in case of bimatrices we may have square, mixed square, rectangular or mixed rectangular fuzzy (semi fuzzy) bimatrices.

Now we proceed on to define fuzzy integral neutrosophic bimatrix.

**DEFINITION 1.2.3:** *Let $A_{FN} = A_1 \cup A_2$ where $A_1$ and $A_2$ are distinct integral fuzzy neutrosophic matrices. Then $A_{FN}$ is called the integral fuzzy neutrosophic bimatrix. If both $A_1$ and $A_2$ are $m \times m$ distinct integral fuzzy neutrosophic matrices then $A_{FN} = A_1 \cup A_2$ is called the square integral fuzzy neutrosophic bimatrix.*

As in case of neutrosophic bimatrices we can define rectangular integral fuzzy neutrosophic bimatrix, mixed square integral fuzzy neutrosophic bimatrix and so on.

If in $A_{FN} = A_1 \cup A_2$ one of $A_1$ or $A_2$ is a fuzzy neutrosophic matrix and the other is just a fuzzy matrix or a neutrosophic matrix we call $A_{FN}$ the semi integral fuzzy neutrosophic bimatrix.

Now we will illustrate them with examples.

***Examples 1.2.15:*** Let $A_{FN} = A_1 \cup A_2$ where

$$A_{FN} = \begin{bmatrix} 0 & I & .3 \\ .2I & .4 & 1 \\ 0 & .3 & -.6 \end{bmatrix} \cup \begin{bmatrix} 1 & I & 0 \\ I & 1 & .8 \\ .6 & 1 & .7I \end{bmatrix}$$

then $A_{FN}$ is a square fuzzy neutrosophic bimatrix.

***Example 1.2.16:*** Consider $A_{FN} = A_1 \cup A_2$ where

$$A_1 = \begin{bmatrix} 0 & .2I & 1 \\ I & .7 & 0 \\ 1 & 0 & .1 \end{bmatrix}$$

and

$$A_2 = \begin{bmatrix} .1 & 0 \\ .2 & I \\ 1 & .2I \end{bmatrix}.$$

Clearly $A_{FN}$ is a mixed fuzzy neutrosophic bimatrix.

***Example 1.2.17:*** Let $A_{FN} = A_1 \cup A_2$ where

$$A_1 = \begin{bmatrix} 2 & 0 & 1 \\ 1 & 2 & 3 \\ I & 0 & I \end{bmatrix}$$

and

$$A_2 = \begin{bmatrix} I & .2I & .6 & .1 \\ .3 & 1 & 0 & I \\ 0 & 0 & .2 & 1 \end{bmatrix}.$$

Clearly $A_{FN}$ is a mixed semi fuzzy neutrosophic bimatrix.

***Example 1.2.18:*** Let

$$A_{FN} = \begin{bmatrix} I & .3I & 0 \\ 1 & .2 & .6 \end{bmatrix} \cup \begin{bmatrix} .6 & 0 & .3 & 1 \\ 1 & 1 & .6 & .2 \\ .3 & 0 & 0 & .5 \end{bmatrix},$$

$A_{FN}$ is a mixed semi fuzzy neutrosophic bimatrix.

***Example 1.2.19:*** Let

$$A_{FN} = \begin{bmatrix} .3 \\ I \\ .2 \\ 0 \end{bmatrix} \cup \begin{bmatrix} I \\ 7 \\ 2 \\ 1 \end{bmatrix}.$$

$A_{FN}$ is a column semi fuzzy neutrosophic bimatrix.

For more about bimatrices and its properties please refer [154, 155].

**Chapter Two**

# BIGRAPHS AND THEIR APPLICATION TO BIMATRICES

Before we proceed onto give some of the applications of bimatrices in the fuzzy models in general and fuzzy cognitive maps in particular. To have such models we have to introduce bigraphs and their related bimatrices. FCMs work on opinion of the experts given on the unsupervised data. The main advantage of this system is it can give the hidden pattern of the problem. To the best of our knowledge this is the model which gives the hidden pattern.

Bart Kosko had given FCM's and talked about combined FCMs; several other researchers have worked on several new types of modified or super imposed models of FCMs [72 to 76]. The book on FCM and NCM by [151] gives a brief description of working of models. One of the modified models are disjoint overlap FCM and overlap FCM when some of the attributes overlap i.e. in the study or analysis of a model when we have some common attributes given by experts i.e. the concepts / attributes are not totally disjoint but have certain concepts / attributes to be common. Such study has already been analyzed. Here we wish to describe the problem when the FCM model is to analyze two sets of attributes and in that pair of concepts / attributes only when a point or a edge or some points and some edges are common in the directed graphs obtained between those pairs of models. We now proceed on to describe these models for this we need some properties of bigraphs and bimatrices. The chapter has only one section in which we introduce bigraphs and give some of its properties.

Recently the notion of bimatrices have been introduced, now the concept of bigraphs will be introduced as in many cases bigraphs have an association with bimatrices.

We can have three types of bimatrices associated with bigraphs.

1. Weighted bigraph's bimatrices
2. Incidence bimatrices.
3. Kirchloff bimatrix K(G)

1. The weighted bigraph's bimatrices will always be square bimatrices or mixed square bimatrices. It is always a square bimatrix.
2. In case of Incidence bimatrices we see that the bimatrix can be square bimatrix or mixed square bimatrix or rectangular bimatrix or mixed rectangular bimatrix depending on the number of edges and vertices.
   We have defined and introduced the notion of bimatrices in [154, 155]. Here we only recall just the definition merely for the sake of completeness.
3. Kirchhoff bimatrix will always to be a square or a mixed square bimatrix denoted by $K(G) = K_1(G_1) \cup K_2(G_2)$
   $K^1_{ii} = d^-(v_i)$ in-degree of the $i^{th}$ vertex of $G_1$
   $K^2_{ii} = d^{-1}(v_i)$ in-degree of the $i^{th}$ vertex of $G_2$
   $K^1_{ij} = -x^1_{ij}$ $(i, j)^{th}$ entry of the adjacency matrix with a negative sign in $G_1$
   $K^2_{ij} = -x^2_{ij}$ $(i, j)^{th}$ entry of the adjacency matrix with a negative sign in $G_2$.

As directed graphs play a role in the study of FCMs we see the major role played by bimatrices in Fuzzy cognitive bimaps. Now we proceed on to define certain new notions like bigraph, point wise glued bigraph, edge / curve glued bigraph points and curve / edge wise glued bigraph. Here we mention the notion of bigraph, and this concept is different from the concept of bipartite graph. Infact bipartite graphs will be shown to be bigraphs only in very special cases but bigraphs in general are not bipartite.

**DEFINITION 2.1.1:** *$G = G_1 \cup G_2$ is said to be a bigraph if $G_1$ and $G_2$ are two graphs such that $G_1$ is not a subgraph of $G_2$ or $G_2$ is not a subgraph of $G_1$, i.e., they have either distinct vertices or edges.*

***Example 2.1.1:*** Let $G = G_1 \cup G_2$ where

$G_1 =$

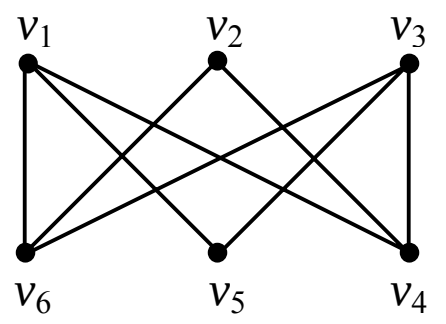

FIGURE: 2.1.1a

and

$G_2 =$

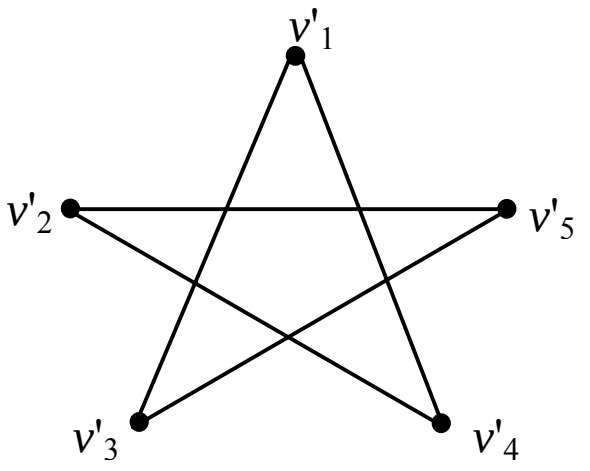

FIGURE: 2.1.1b

$G = G_1 \cup G_2$ is a bigraph. G can also be represented as

$$G = \{v_1, v_2, v_3, v_4, v_5, v_6\} \cup \{v'_1, v'_2, v'_3, v'_4, v'_5\};$$

the vertices of the two graphs $G_1$ and $G_2$ respectively.

*Note:* It is very important to note that the vertices of the bigraphs i.e. $G = G_1 \cup G_2$ in general need not form disjoint subsets of G such that $G_1 \cap G_2 = \phi$.

In the above example we have $G_1 \cap G_2$ is empty.
Now we can have bigraphs given by the following example.

***Example 2.1.2:*** Let $G = G_1 \cup G_2$ be the bigraph given by the following figure.

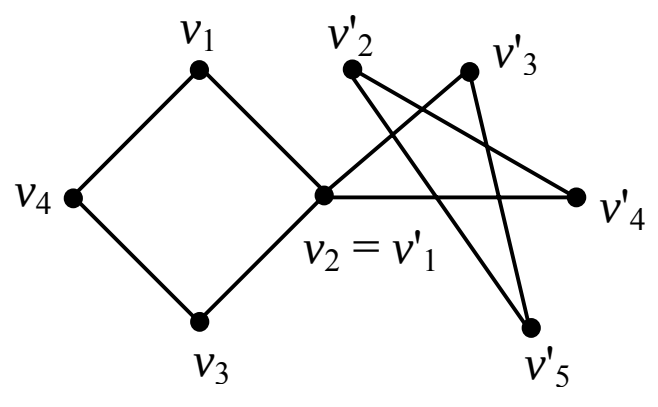

FIGURE: 2.1.2

Here $G = G_1 \cup G_2$ , $\{ v_1, v_2 = v'_1, v_3, v_4\} \cup \{ v'_1, v'_2, v'_3, v'_4, v'_5\}$ we see that this bigraph is very special in its own way for it has only one point in common, viz. $v_2 = v'_1$.

***Example 2.1.3:*** Let $G = G_1 \cup G_2$ be a bigraph given by the following figure.

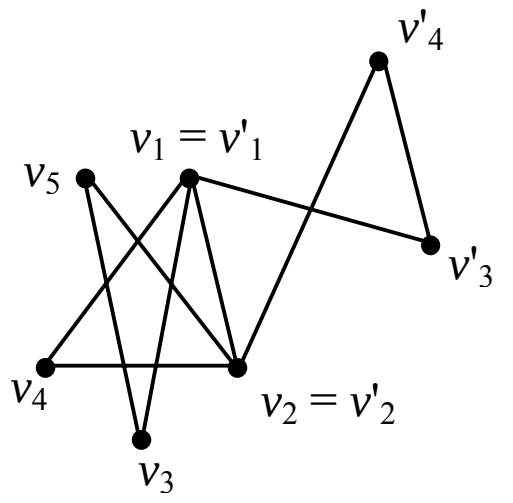

FIGURE: 2.1.3

It is interesting to observe that the graphs have only two points in common viz. $v_1 = v'_1$ and $v_2$ and $v'_2$ .
Thus

$$G = \{v_1, v_2, v_3, v_4, v_5\} \cup \{ v'_1, v'_2, v'_3, v'_4\}.$$

Here $v_1 = v'_1$ and $v_2 = v'_2$.

Now we proceed onto see a graph which has only one edge and two points in common.

***Example 2.1.4:*** Let $G = G_1 \cup G_2$ where $G_1 = \{v_1, v_2, v_3, v_4, v_5, v_6\}$ and $G_2 = \{v'_1, v'_2, v'_3, v'_4, v'_5, v'_6\}$.
$G_1 \cup G_2$ be the bigraph given by the following figure.

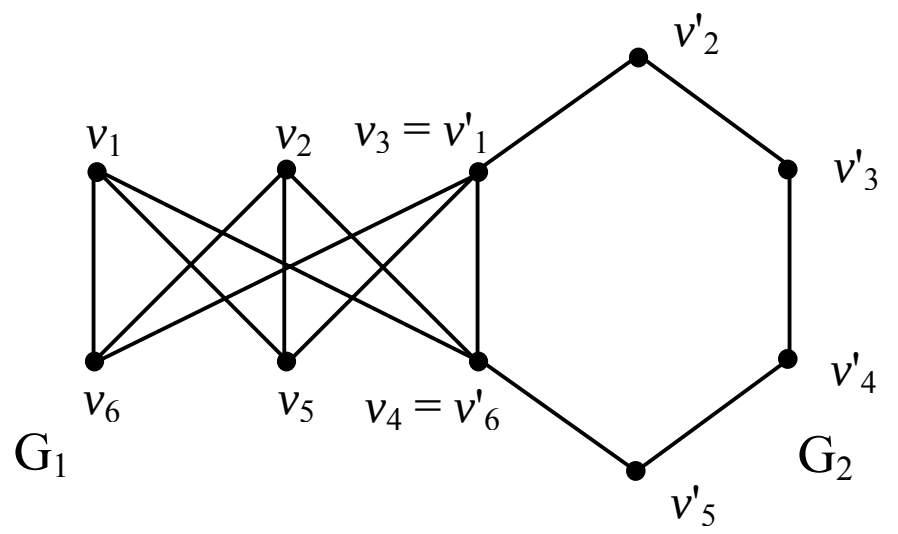

FIGURE: 2.1.4

Clearly this bigraph has two points $v_3 = v'_1$ and $v_4 = v'_6$ common, also the edges connecting $v'_1$ and $v'_6$ in $G_2$ and $v_3$ and $v_4$ in $G_1$ is the only common edge.

Now we proceed on to define bigraphs which can have more than one edge and vertices to be common given by the following example.

***Example 2.1.5:*** Consider the bigraph given by the following figure, $G = G_1 \cup G_2$

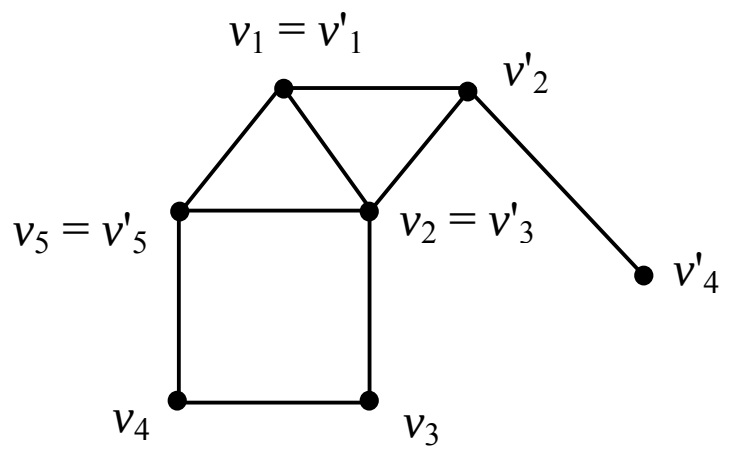

FIGURE: 2.1.5

The vertices common with these graphs are $v_1 = v'_1$, $v_2 = v'_3$ and $v_5 = v'_5$ and also 3 edges or to be more precise we see a subgraph itself is common for both the bigraphs.

Based on these examples we define some new notions.

**DEFINITION 2.1.2:** *Let G = $G_1 \cup G_2$ be a bigraph we say G is said to be a disjoint bigraph if G = $G_1 \cup G_2$ are such that* $G_1 \cap G_2 = \phi$.

We see the bigraph given by the following example is the disjoint bigraph.

***Example 2.1.6:*** Let G = $G_1 \cup G_2$ be given by the following figure.

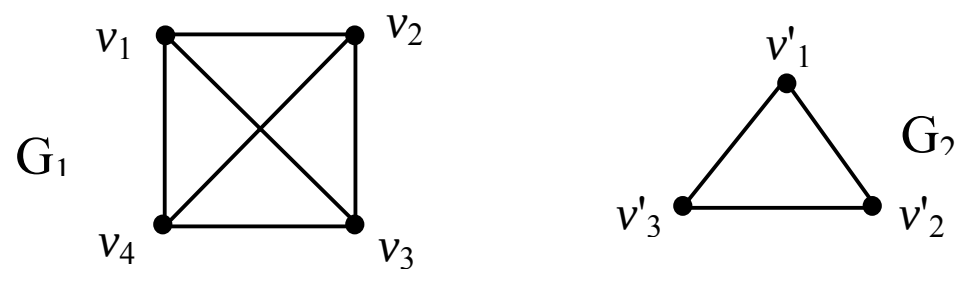

FIGURE: 2.1.6

The bigraph is such that $G_1 \cap G_2 = \phi$. This is a disjoint bigraph.

**DEFINITION 2.1.3:** *Let G = $G_1 \cup G_2$ be the bigraph. If $G_1$ and $G_2$ are graphs such that they have a single point in common i.e. $G_1 \cap G_2$ = {single vertex}, then we say the bigraph is a pair of graphs glued at a point i.e. single point glued bigraph.*

***Example 2.1.7:*** The graph G = $G_1 \cup G_2$ where

$$G = \{v_1, v_2, v_3, v_4\} \cup \{v'_1, v'_2, v'_3, v'_4, v'_5, v'_6, v'_7\}$$

is given by the following figure

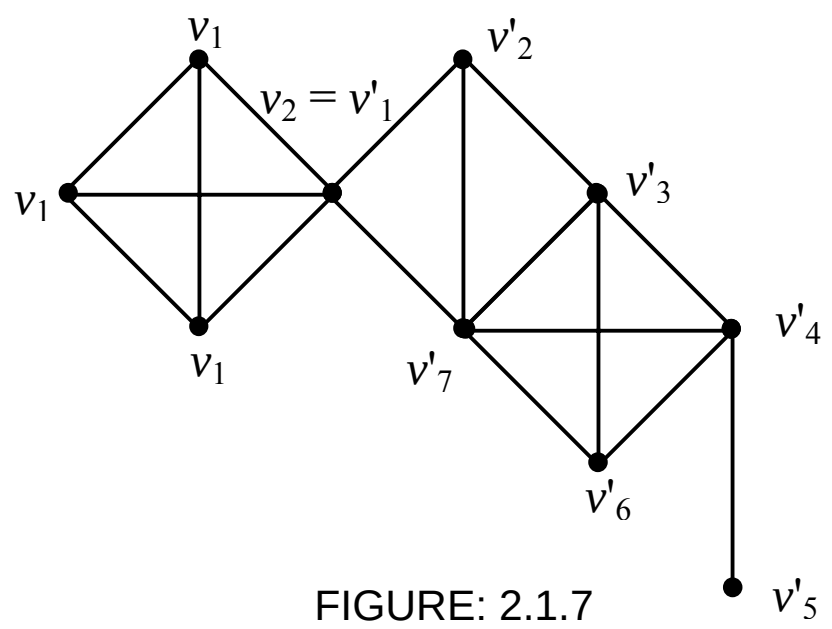

FIGURE: 2.1.7

This is a single point glued bigraph glued at the vertex $v_2 = v'_1$.

Now be proceed on to define bigraphs glued by a line / edge.

**DEFINITION 2.1.4:** *$G = G_1 \cup G_2$ be a bigraph. If the bigraph is such that the graphs $G_1$ and $G_2$ have a common edge then we call the bigraph to be a edge glued bigraph or a line glued bigraph.*

The following example will show how a line glued bigraph looks.

***Example 2.1.8:*** Let $G = G_1 \cup G_2$ be a bigraph given by the following figure.

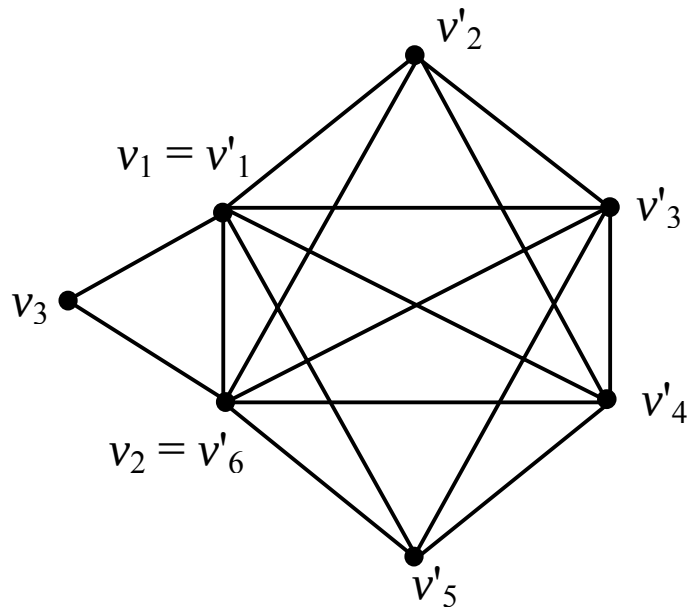

FIGURE: 2.1.8

The graphs $G_1 = \{v_1, v_2, v_3\}$

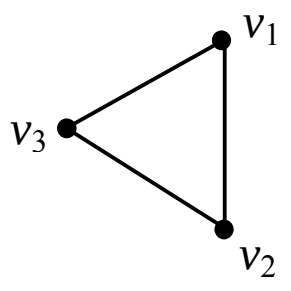

FIGURE: 2.1.8a

and the graph $G_2$

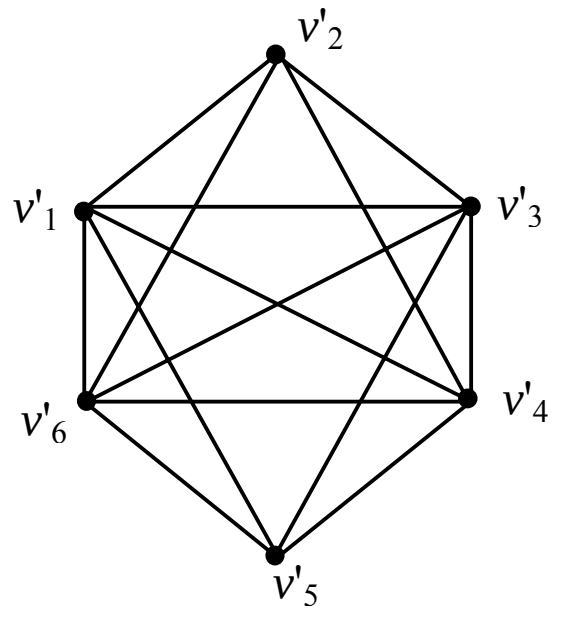

FIGURE: 2.1.8b

Clearly the edge connection $v_1$ and $v_2$ of $G_1$ is common with the edge connecting $v'_1$ and $v'_6$. Thus this is an example of a edge glued bigraph.

Now it is still interesting to see bigraphs which have edge glued or edge glued bigraphs are also point or vertex glued bigraphs.

Thus we can have the following Theorem.

**THEOREM 2.1.1:** *Let $G = G_1 \cup G_2$ be a bigraph if G is a edge glued bigraph then G is a vertex glued bigraph.*

*A vertex glued bigraph in general need not be edge glued bigraph.*

*Proof:* Let G = $G_1 \cup G_2$ be a bigraph. Suppose G is a edge glued bigraph then certainly G is a vertex glued bigraph for edge will certainly include at least two vertices. So always a edge glued bigraph will be a vertex glued bigraph.

However a vertex glued bigraph in general need not be a edge glued bigraph even if the two graphs $G_1$ and $G_2$ have more than two vertices in common. This will be proved by the following example.

The bigraph given in example 2.1.9 is such that it has 3 vertices in common but have no edge in common.

***Example 2.1.9:*** $G = G_1 \cup G_2 = \{v_1, v_2, v_3, \ldots, v_9, v_{10}\} \cup \{v'_1, v'_2, v'_3, v'_4\}$, given in figure 2.1.9.

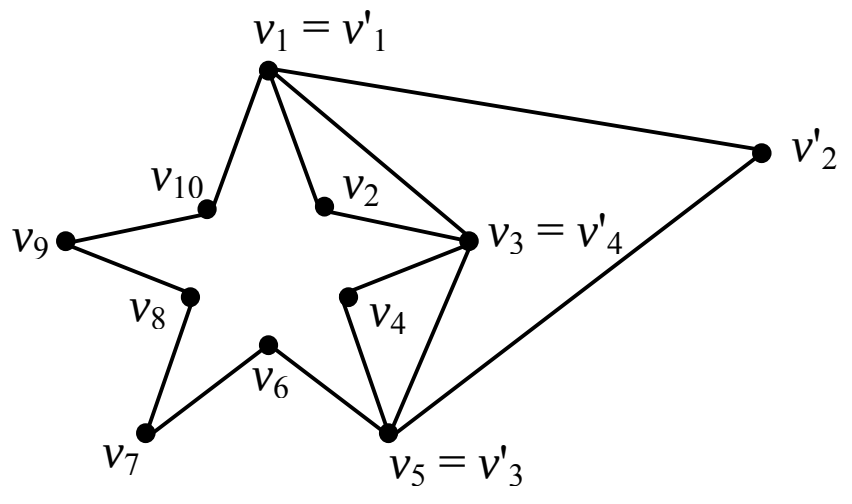

FIGURE: 2.1.9

The two separate graphs are as follows.

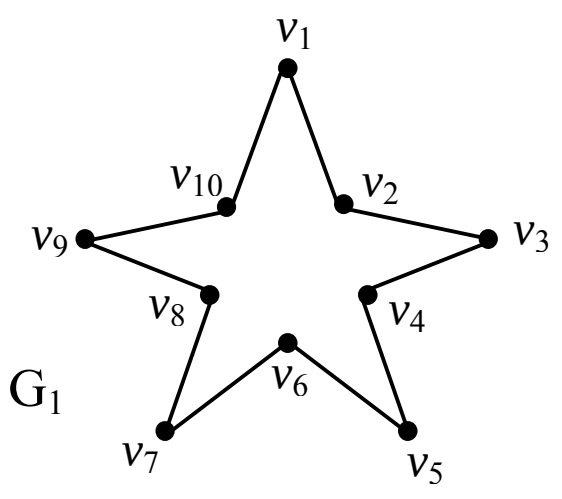

FIGURE: 2.1.9a

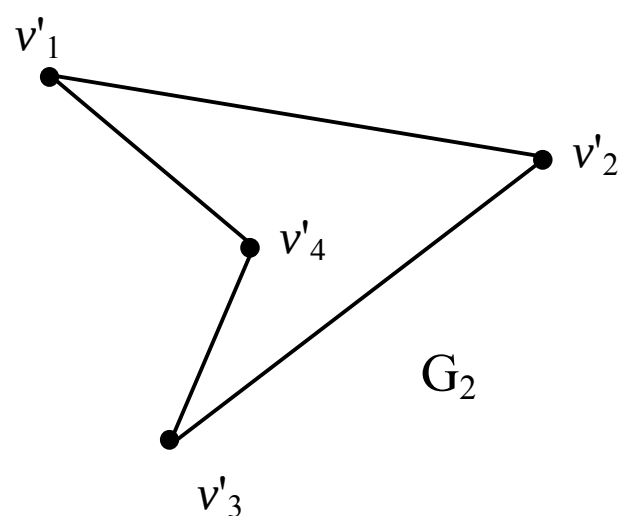

FIGURE: 2.1.9b

This bigraph has three vertices in common. i.e. $v_1 = v'_1$, $v_3 = v'_4$ and $v_5 = v'_3$ but this bigraph has no edge in common.

Hence the claim of the theorem.

Now it is still interesting to note the following result.

**THEOREM 2.1.2:** *A bigraph glued by even more than two vertices need not in general be glued by an edge.*

*Proof***:** This result is proved by the following example.

***Example 2.1.10:*** Consider the bigraph $G = G_1 \cup G_2$ given by the following figure.

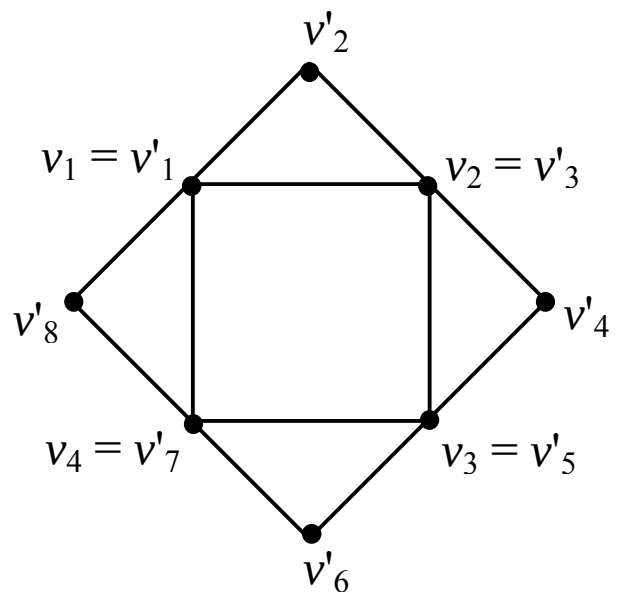

FIGURE: 2.1.10

That is

$G_1 = \{v_1, v_2, v_3, v_4\}$ and $G_2 = \{v'_1, v'_2, v'_3, v'_4, v'_5, v'_6, v'_7, v'_8\}$.

This bigraph has four vertex in common but no edge. The graphs related with $G_1$ and $G_2$ are given by the following figures 2.1.10a and 2.1.10b.

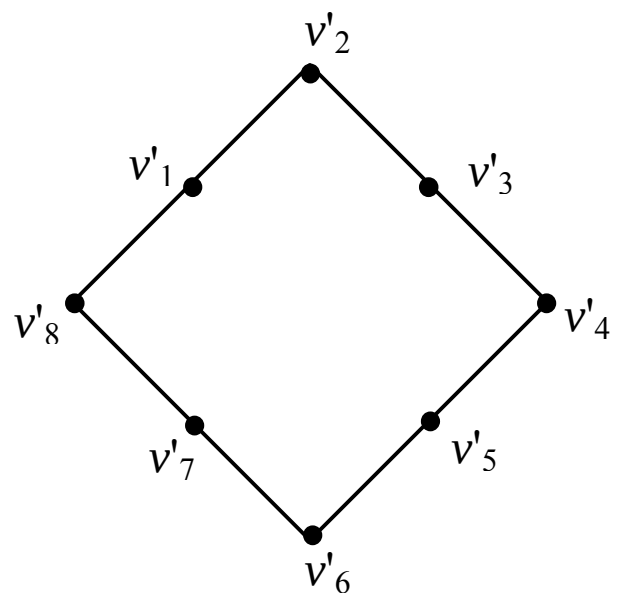

FIGURE: 2.1.10a

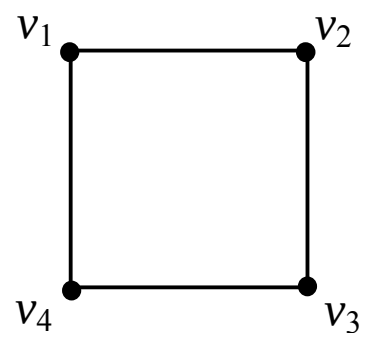

FIGURE: 2.1.10b

Hence the claim.

We can have bigraphs glued by more than one edge or by a subgraph which happen to be the subgraph of both the graphs. We define it as follows.

**DEFINITION 2.1.5:** *Let $G = G_1 \cup G_2$ be a bigraph. Suppose the bigraphs is glued such that it has a subgraph (with more than one vertex and more than one edge) in common then we define this bigraph as a bigraph glued by a strong subgraph or a strong subgraph glued bigraph.*

We illustrate this by the following example.

***Example 2.1.11:*** $G = G_1 \cup G_2$ be a bigraph given by the following figure.

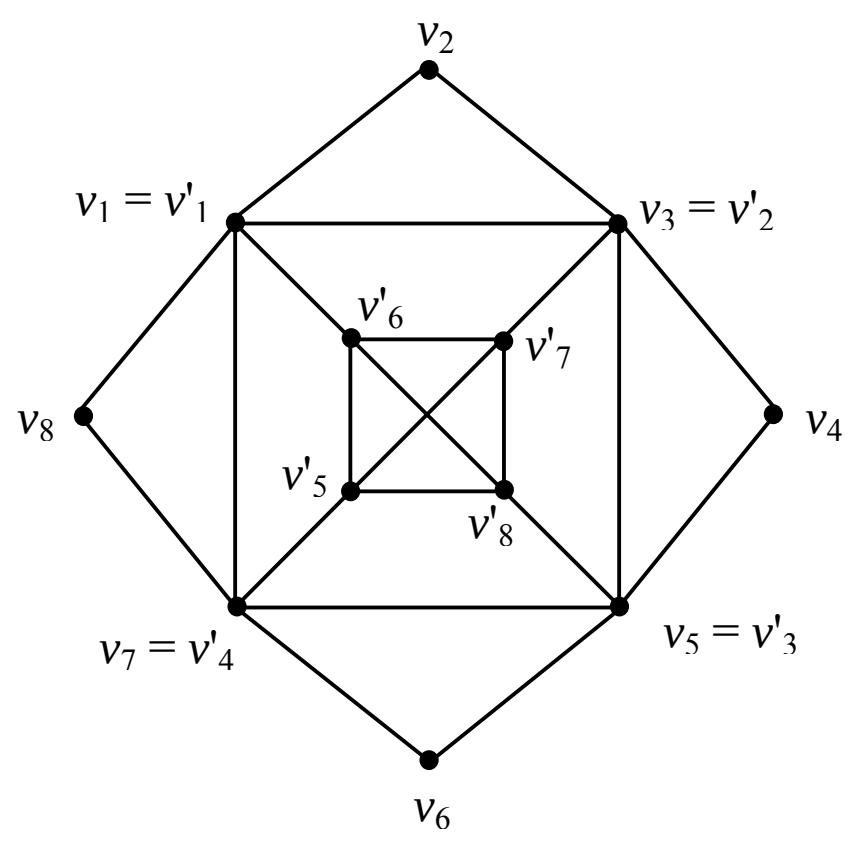

FIGURE: 2.1.11

Here $G = G_1 \cup G_2 = \{v_1, v_2, v_3, \ldots, v_8\} \cup \{v'_1, v'_2, v'_3, \ldots, v'_8\}$.

The graphs associated with $G_1$ and $G_2$ are given by the following figure.

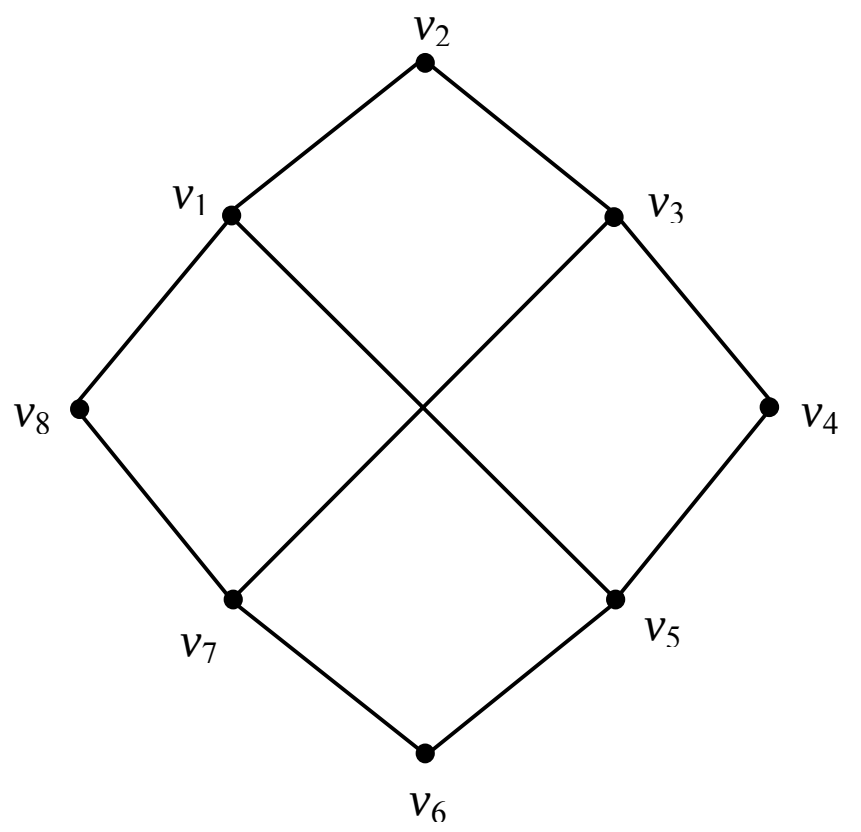

FIGURE: 2.1.11a

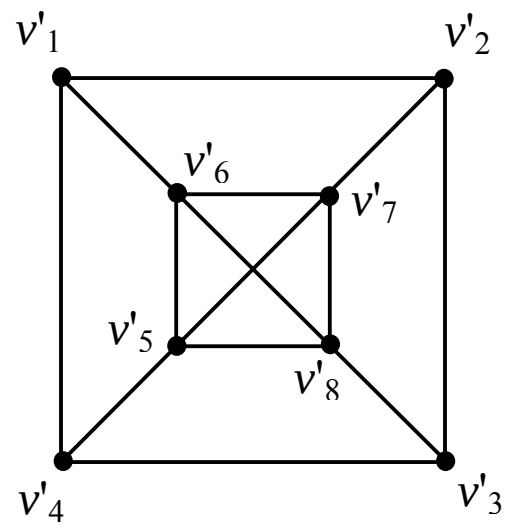

FIGURE: 2.1.11b

This bigraph has a subgraph in common given by the following figure.

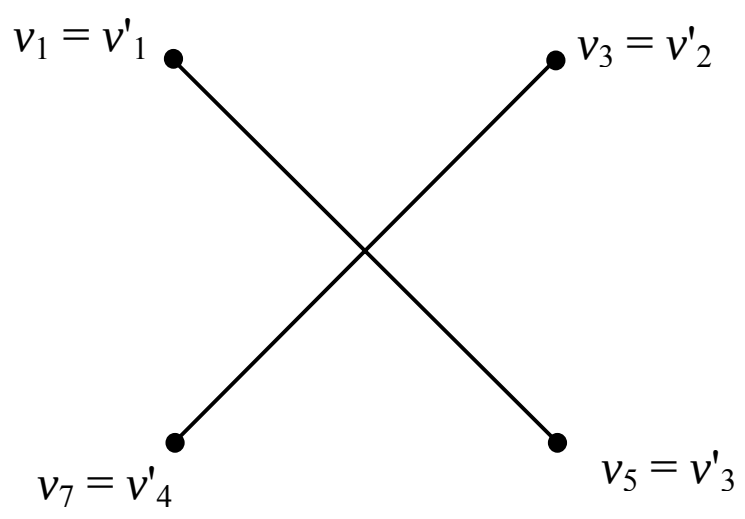

FIGURE: 2.1.11c

Thus we have the following interesting result.

**THEOREM 2.1.3:** *Let $G = G_1 \cup G_2$ be a bigraph which is strong subgraph glued bigraph, then G is a vertex glued graph and a edge glued graph.*

*Proof:* Since single point or an edge is subgraph we see vertex glued bigraph and a edge glued bigraph are also subgraph glued bigraphs, but cannot be called as strong subgraph glued bigraph. But clearly in case of strong subgraph glued bigraphs we see it is both a vertex glued bigraph and edge glued bigraph.

*Note:* We say a vertex glued bigraph or an edge glued bigraph will be called as just a subgraph glued bigraphs.

Now we proceed on to define the notion of subbigraph.

**DEFINITION 2.1.6:** *Let $G = G_1 \cup G_2$ be a bigraph. A non empty subset H of G is said to be a subbigraph of G if H itself is a bigraph of G.*

***Example 2.1.12:*** Let $G = G_1 \cup G_2$ be a bigraph given by the following figure.

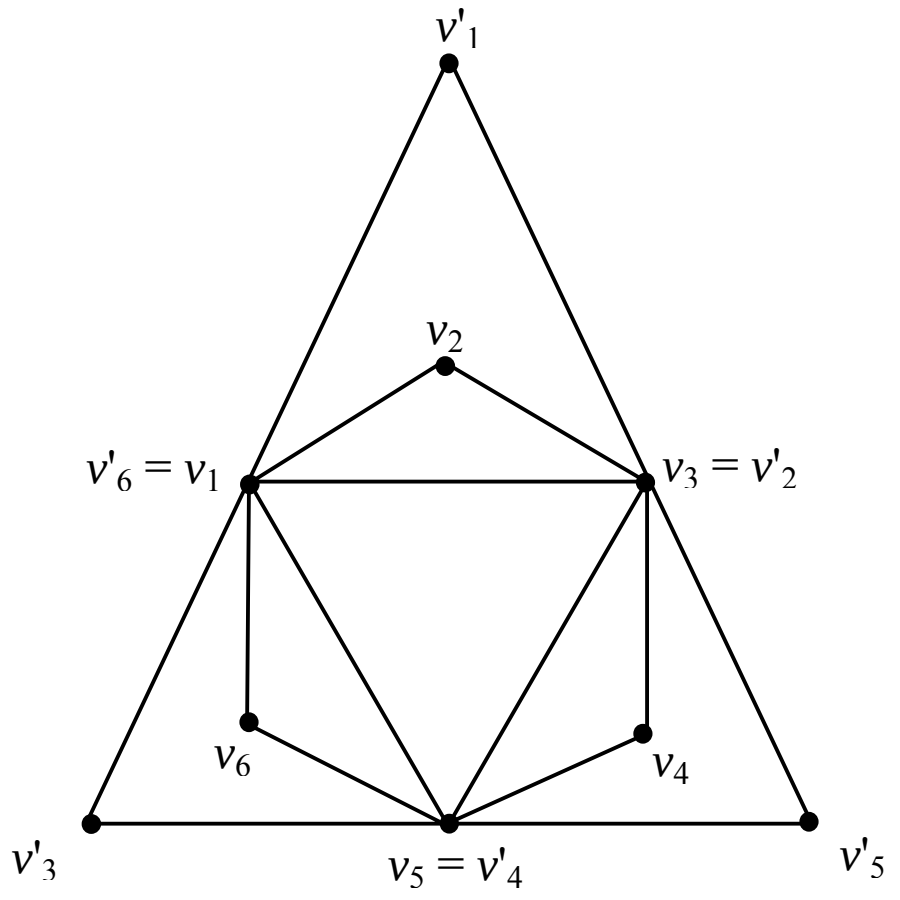

FIGURE: 2.1.12

where $G_1 = \{v_1, v_2, v_3, \ldots, v_6\}$ and $G_2 = \{v'_1, v'_2, v'_3, \ldots, v'_6\}$. The graphs of $G_1$ and $G_2$ are given by following figures.

Graph of $G_1$

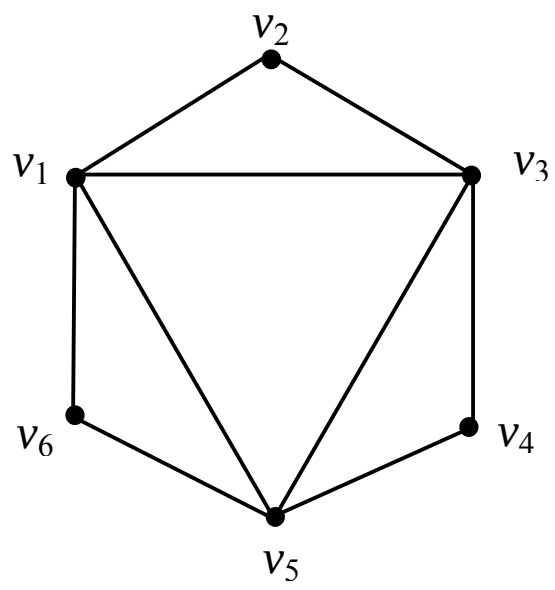

FIGURE: 2.1.12a

Graph of $G_2$

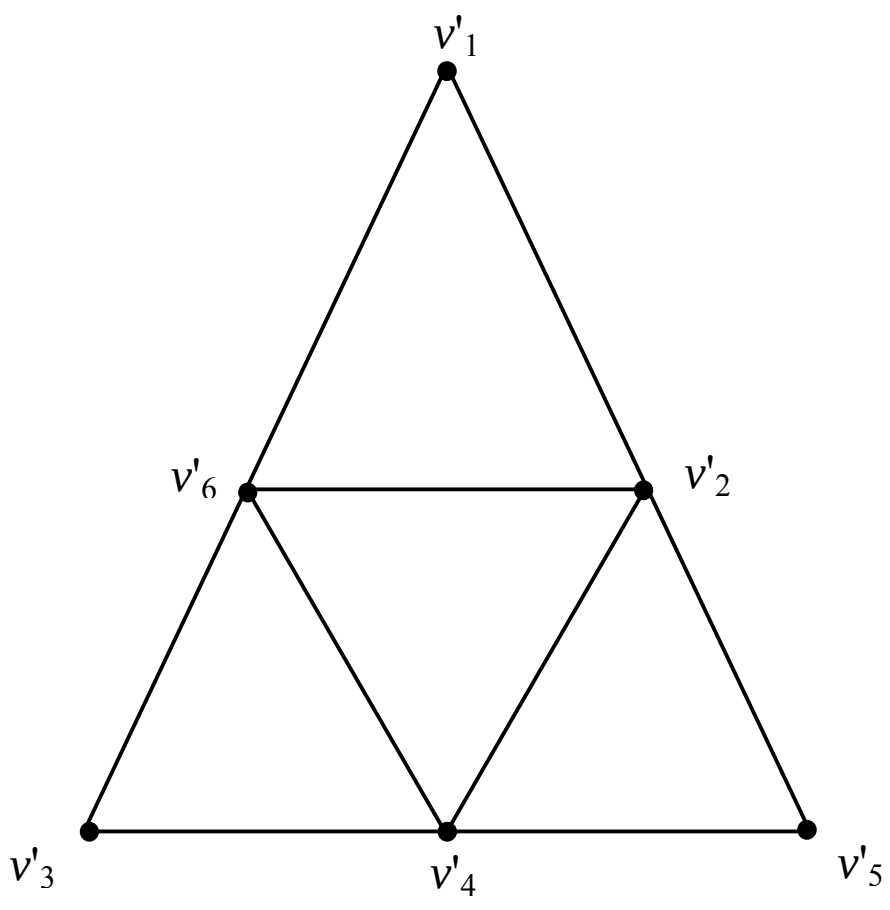

FIGURE: 2.1.12b

The subbigraph H is given by the following figure 2.1.12c.

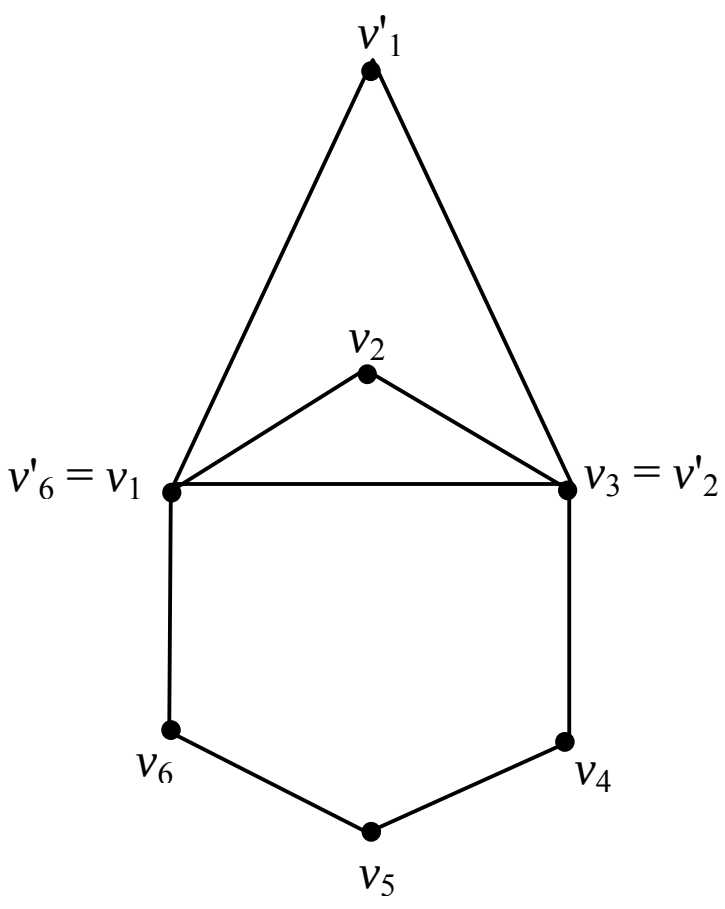

FIGURE: 2.1.12c

Here $H = \{v_1, \ldots, v_6\} \cup \{v'_1, v'_2, \ldots, v'_6\}$

$H_1 = \{v_1, \ldots, v_6\}$ and $H_2 = \{v'_1, v'_2, \ldots, v'_6\}$.

Now we proceed on to give the three types of bimatrices associated with the bigraphs.

We now show by the following example that a bipartite graph in general is not a bigraph.

***Example 2.1.13:*** Consider the bipartite graph B (G) of G given by the following figure.

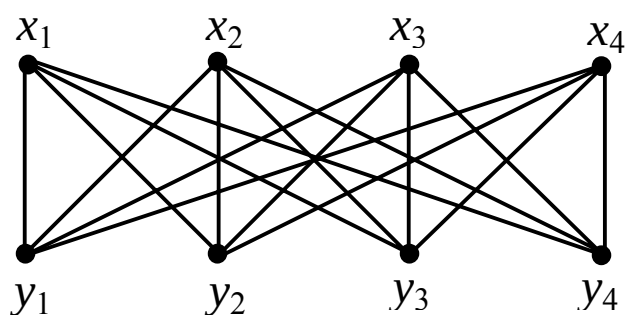

FIGURE: 2.1.13

Clearly this is not a Bigraph.

***Example 2.1.14:*** Now consider the bipartite graph given by the following figure 2.1.14.

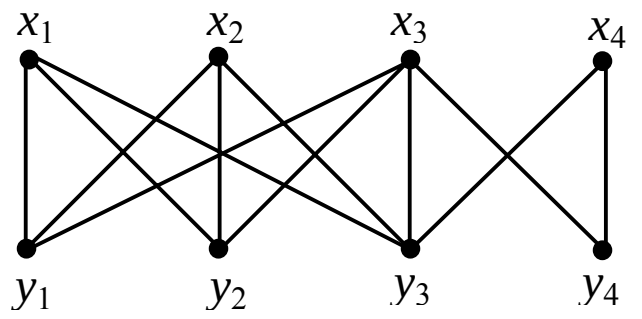

FIGURE: 2.1.14

Take G = $G_1 \cup G_2$

= $\{x_1, x_2, x_3, y_1, y_2, y_3\} \cup \{x_3\ y_3\ x_4\ y_4\}$.

Clearly G is a bigraph with $G_1$ given by the following figure 2.1.14a.

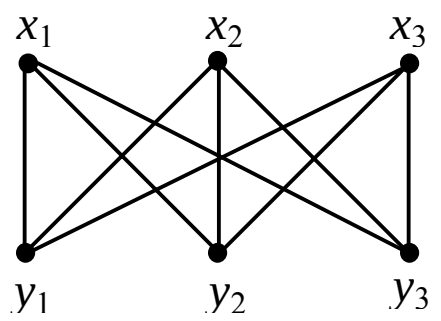

FIGURE: 2.1.14a

and $G_2$ is given by the following figure 2.1.14b.

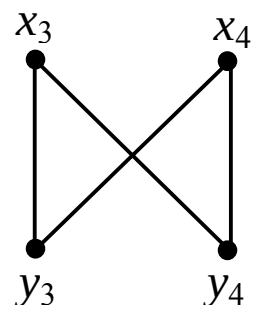

FIGURE: 2.1.14b

The study of conditions making a bigraph into a bipartite graph is an interesting one. We say as in case of graphs in bigraphs also the following.

Since a bigraph can be realized as the 'union' of two graphs here the 'union' is distinctly different from terminology union of graphs used. The symbol just denotes only connection or union as subsets. Thus a bigraph G can be realized as G = (V ($G_1$), E ($G_1$), $I_{G_1}$) $\cup$ (V ($G_2$), E ($G_2$), $I_{G_2}$) where $V_i$ (G) is a nonempty set (i = 1, 2) and $E_i$ (G) is a set disjoint from $V_i$ (G) for i = 1, 2 and $I_{G_1}$ and $I_{G_2}$ are incidence maps that associates with each element of $E_i$ (G) an unordered pair of elements of $V_i$ (G) i = 1, 2.

It is to be noted that
$E_1 (G_1) \not\subset E_2 (G_2)$ or $V_1 (G_1) \not\subseteq V_2 (G_2)$, $I_{G_1} \not\subseteq I_{G_2}$ and
$E_2 (G_2) \not\subset E_1 (G_1)$ or $V_2 (G_2) \not\subset V_1 (G_1)$, $I_{G_1} \not\subseteq I_{G_2}$.

Further we can have common edges between $G_1$ and $G_2$ i.e. in general

$$\left.\begin{array}{l} E_1(G_1) \cap E_2(G_2) \neq \phi \\ I_1(G_1) \cap I_2(G_2) \neq \phi \text{ and} \\ V_1(G_1) \cap V_2(G_2) \neq \phi \end{array}\right\} \qquad \text{(I)}$$

But nothing prevents us from having

$$\left.\begin{array}{l} V_1(G_1) \cap V_2(G_2) = \phi \text{ or} \\ E_1(G_1) \cap E_2(G_2) = \phi \text{ and} \\ I_{G_1} \cap I_{G_2} = \phi. \end{array}\right\} \qquad \text{(II)}$$

In case of disjoint bigraphs we see all the equation described in II are satisfied i.e. they are union of disjoint graphs. Thus union of any two graphs is a bigraph provided the union does not take place between the graph and its subgraph.

We have just now seen the diagrammatic representation of bigraph.

Almost all terms associated with graphs can also be easily extended in case of bigraphs like and vertices, incident, multiple incident, neighbor of V in G etc in an appropriate way with suitable or appropriate modifications.

**DEFINITION 2.1.7:** *Let $G = G_1 \cup G_2$ be a bigraph. A vertex $u_1 \in G_1$ is a bineighbour of $u_2$ in $G_2$ if $u_1$, $u_2$ is an edge in G such edges are called as biedges of the bigraph.*

The set of all bineighbours of $u_2$ in $G_2$ ($u_1$ in $G_1$) is the open bineighbourhood of $u_2$ (or $u_1$) or the bineighbour set of $u_2$ (or $u_1$) and will be denoted by

$$BN(u_2) \text{ (or } BN(u_1))$$

$$BN(u_2) = BN(u_2) \cup \{u_2\}.$$

$$\text{Similarly } BN(u_1) = BN(u_1) \cup \{u_1\} \text{ is}$$

called the closed neighbourhood of $u_2$ (or $u_1$ in $G_1$). We may not have concept of bineighbour or open bineighbourhood if $G_1 \cap G_2 = \phi$ i.e. when the bigraph is a disjoint bigraph.

A bigraph $G = G_1 \cup G_2$ is simple if both $G_1$ and $G_2$ are simple. Two biedges are biadjacent if and only if they have a common end vertex which is in the intersection of $G_1$ and $G_2$. A bigraph $G = G_1 \cup G_2$ is called finite if both $G_1$ and $G_2$ are finite graphs even if one of $G_1$ or $G_2$ is not finite then the bigraph G is infinite. The number $m(G) = m_1(G) \cup m_2(G)$ is called the bisize of the bigraph G where $m_i$ (G) is the number of biedges of G, i = 1, 2.

A bigraph $G = G_1 \cup G_2$ is labeled if both $G_1$ and $G_2$ are labeled. Let $G = G_1 \cup G_2$ and $H = H_1 \cup H_2$ be any two bigraphs. The bigraph biisomorphism is one of graph isomorphisms from $G_1$ to $H_1$ and $G_2$ to $H_2$ 'or' $G_2$ $\underrightarrow{\text{to}}$ $H_1$, $G_1$ $\underrightarrow{\text{to}}$ $H_2$.

A bigraph $G = G_1 \cup G_2$ is simple if both $G_1$ and $G_2$ are simple.

A bigraph $G = G_1 \cup G_2$ is bisimple if $|G| = n$ that is number of elements in $G_1$ and $G_2$ without repetition i.e. $o(G_1 \cup G_2) = o(G) = o(G_1) + o(G_2) - o(G_1 \cap G_2)$ and if Bm(G) denotes the biedges then $0 \leq Bm(G) \leq n(n-1)/2$, if G is a disjoint bigraph then Bm(G) = 0.

It is an interesting problem to study Bm (G) when G is a vertex glued bigraph or a edge glued bigraph or a strong subgraph glued bigraph. Study in that direction may throw light on the properties of the bigraph.

**DEFINITION 2.1.8:** *Let $G = G_1 \cup G_2$ be a bigraph we say the bigraph is simple if G is simple, pseudo simple if both $G_1$ and $G_2$ is simple but G is not simple. The simple bigraph G is said to be bicomplete if every pair of distinct vertices of G are adjacent in G.*

**THEOREM 2.1.4:** *If $G_1$ and $G_2$ are complete the bigraph $G = G_1 \cup G_2$ need not in general be complete.*

*Proof:* If G is a disjoint bigraph even if $G_1$ and $G_2$ are complete G need not be complete.
This clear from the following example.

***Example 2.1.15:*** Let the bigraph $G = G_1 \cup G_2$ given by the following figure 2.1.15.

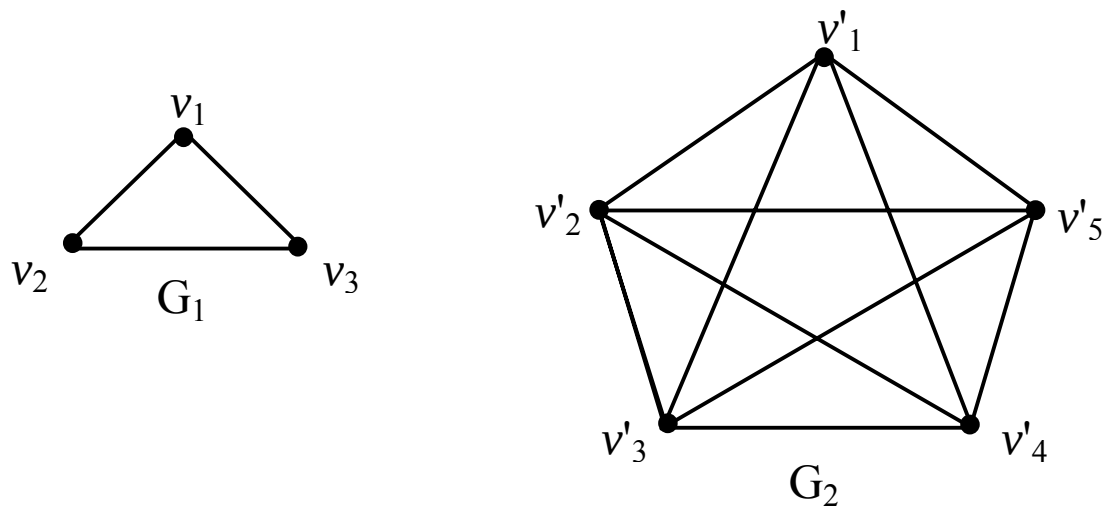

FIGURE: 2.1.15

Clearly $G_1$ and $G_2$ are complete but $G = G_1 \cup G_2$ is not complete for G is a disjoint bigraph.
Consider the vertex glued bigraph $G = G_1 \cup G_2$ where both $G_1$ and $G_2$ are complete but G is not complete given by the following figure 2.1.15a.

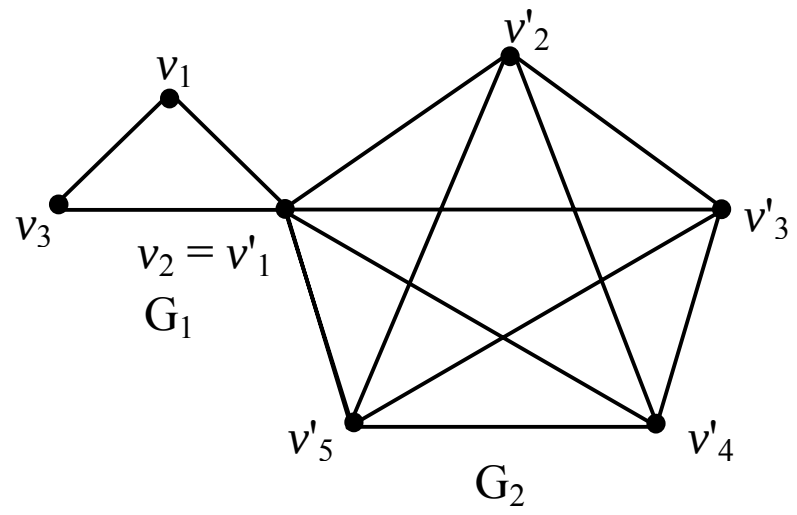

FIGURE: 2.1.15a

Let $G = G_1 \cup G_2$ be a edge glued bigraph where both $G_1$ and $G_2$ are complete. Clearly $G = G_1 \cup G_2$ is not a complete bigraph. The following figure 2.1.15b of the edge glued bigraph of two complete graphs which is not complete.

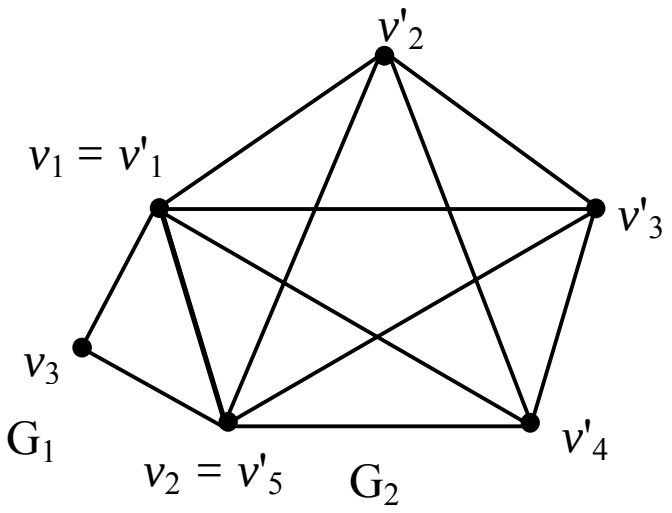

FIGURE: 2.1.15b

The edge $v_1$, $v_2$ of $G_1$ is glued with the edge $v'_1$, $v'_5$ of $G_2$.
A bigraph can also be a bipartite bigraph.

**DEFINITION 2.1.9:** *A bigraph $G = G_1 \cup G_2$ has complement $G^C$ defined by $V(G^C) = V(G)$ i.e. $V\left(G_1^C\right) = V(G)$ and $V\left(G_2^C\right) = V(G)$ and making two adjacent vertices u and v adjacent in $G^C$ if and only if they are non adjacent in G.*

***Example 2.1.16:*** Let $G = G_1 \cup G_2$ be a bigraph where $G = G_1 \cup G_2$ is given by the following figure 2.1.16a.

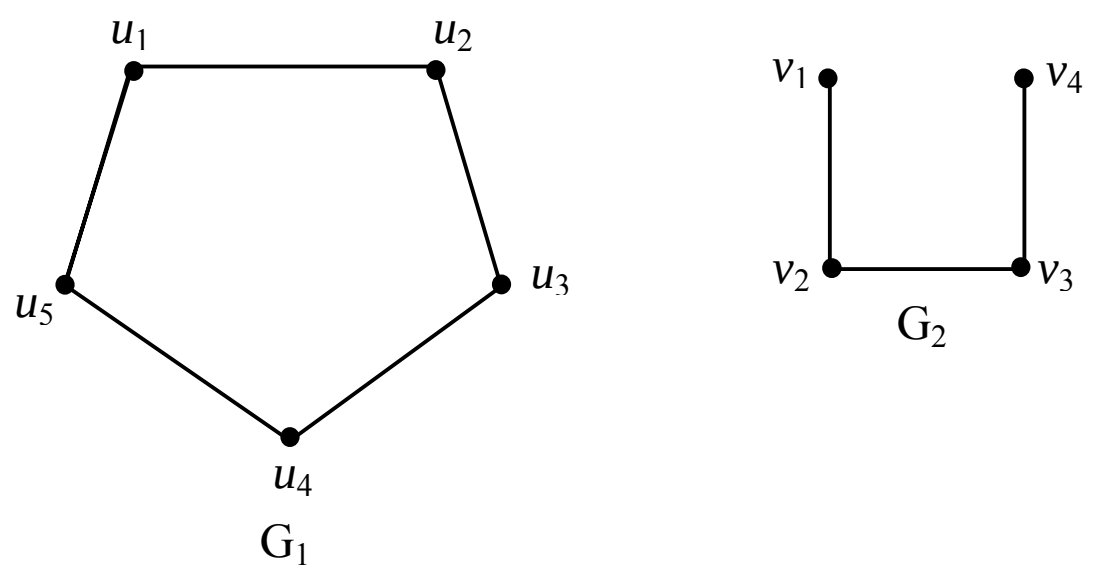

FIGURE: 2.1.16a

The complement $G^C$ of the bigraph G is given by the following figure 2.1.16b.

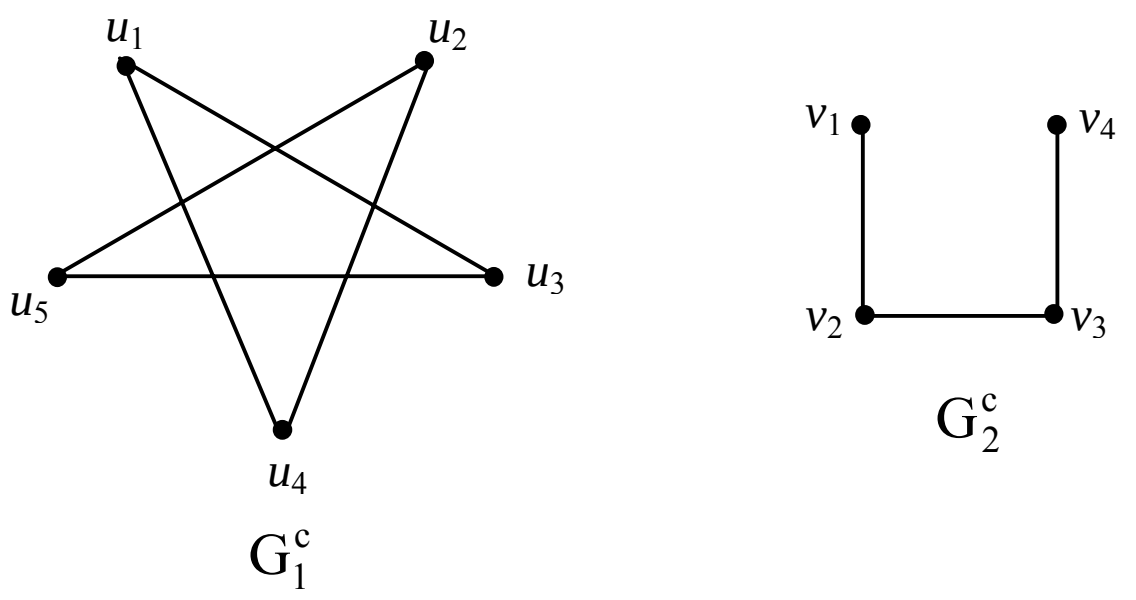

FIGURE: 2.1.16b

We see if G is a disjoint bigraph so is its complement.

**THEOREM 2.1.5:** *If G = $G_1 \cup G_2$ be a disjoint bigraph such that $G_1 = G_2^C$ then G is a self complementary bigraph.*

*Proof:* Follows directly by the definition and the fact the bigraph is disjoint.

We need to work more only when the bigraphs are not disjoint.

**DEFINITION 2.1.10:** *Let G = $G_1 \cup G_2$ be a bigraph and v ∈ V ($G_1$) ∪ V ($G_2$). The number of edges incident at v in G = $G_1 \cup G_2$ is called the degree (or valency) of the vertex v in G and is denoted by $d_G$ (v) or simply by d(v) when G requires no explicit reference. A graph G is called K regular if every vertex has degree K. Now if v ∈ $G_1 \cap G_2$ then the degree associated with v is called as bidegree. A bigraph G is called $K_1 + K_2$ biregular if every bivertex is of degree $K_1 + K_2$ where $G_1$ has degree $K_1$ and $G_2$ has degree $K_2$.*

***Example 2.1.17:*** Let G = $G_1 \cup G_2$ be a bigraph given by the following figure 2.1.17.

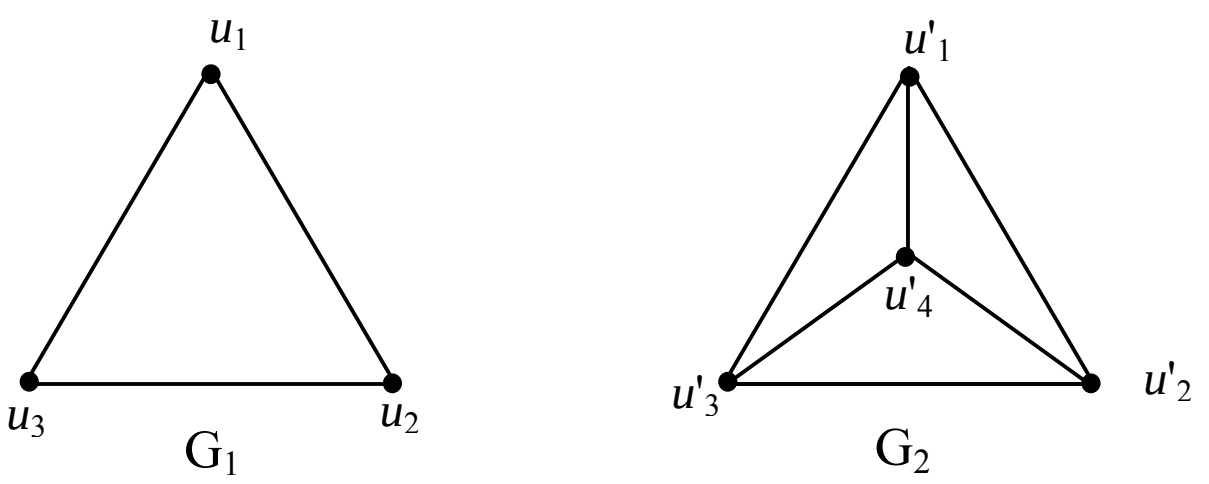

FIGURE: 2.1.17

The degree of any vertex in G is either 1 or 2. This bigraph has no bidegree associated with it. Thus if G = $G_1 \cup G_2$ is a disjoint bigraph it has no bidegree associated with it.

Now we give an example of a bigraph $G = G_1 \cup G_2$ which has a vertex with bidegree associated with it.

***Example 2.1.18:*** Let $G = G_1 \cup G_2$ be a bigraph given by the following figure 2.1.18.

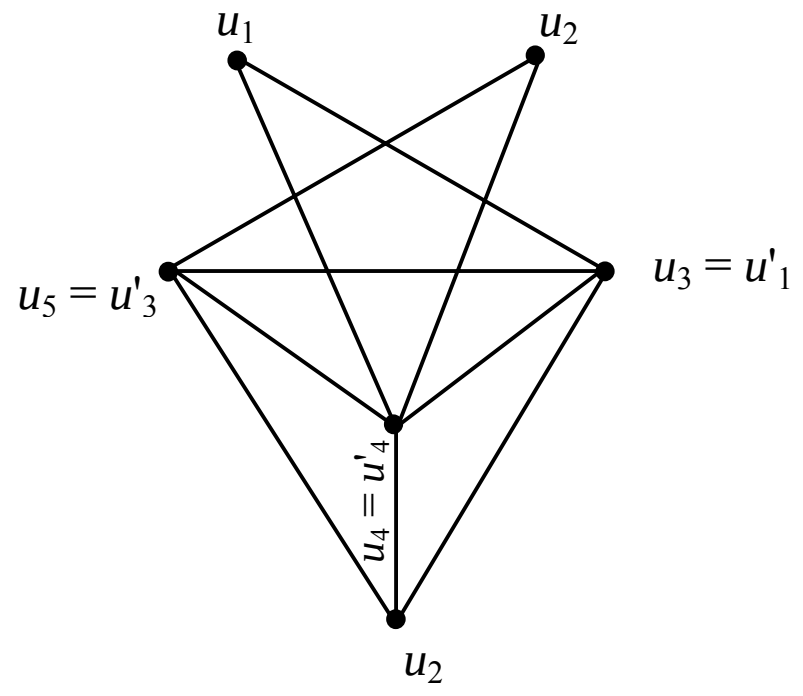

FIGURE: 2.1.18

Here

$$
\begin{aligned}
G &= G_1 \cup G_2 \\
&= \{u_1, u_2, u_3, u_4, u_5\} \cup \{u'_1, u'_2, u'_3, u'_4\}.
\end{aligned}
$$

The graphs of $G_1$ and $G_2$ are as follows; given in figure 2.1.18a and figure 2.1.18b.

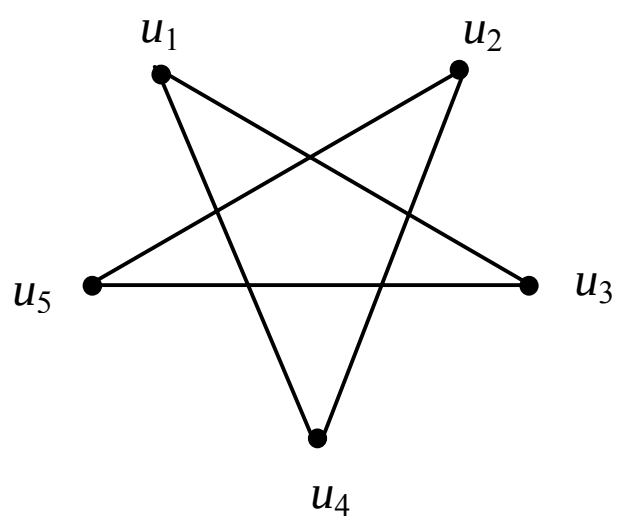

FIGURE: 2.1.18a

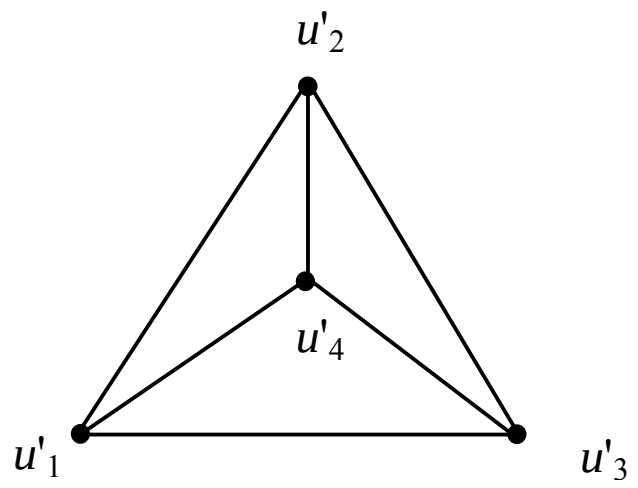

FIGURE: 2.1.18b

Clearly $G_1$ is 2-regular and $G_2$ is 3-regular but the bigraph $G = G_1 \cup G_2$ is not regular. Further G is not $K_1 + K_2$ biregular for we have 3 bivertices $u_3 = u'_1$, $u_4 = u'_4$ and $u_5 = u'_3$ and only the vertex $u_4 = u'_4$ has degree 2 + 3 and the other vertices are just of degree 4.

Now we proceed onto give an example of a bigraph which is biregular.

***Example 2.1.19:*** Let $G = G_1 \cup G_2$ be a bigraph given by the following figure 2.1.19. $G = G_1 \cup G_2$.

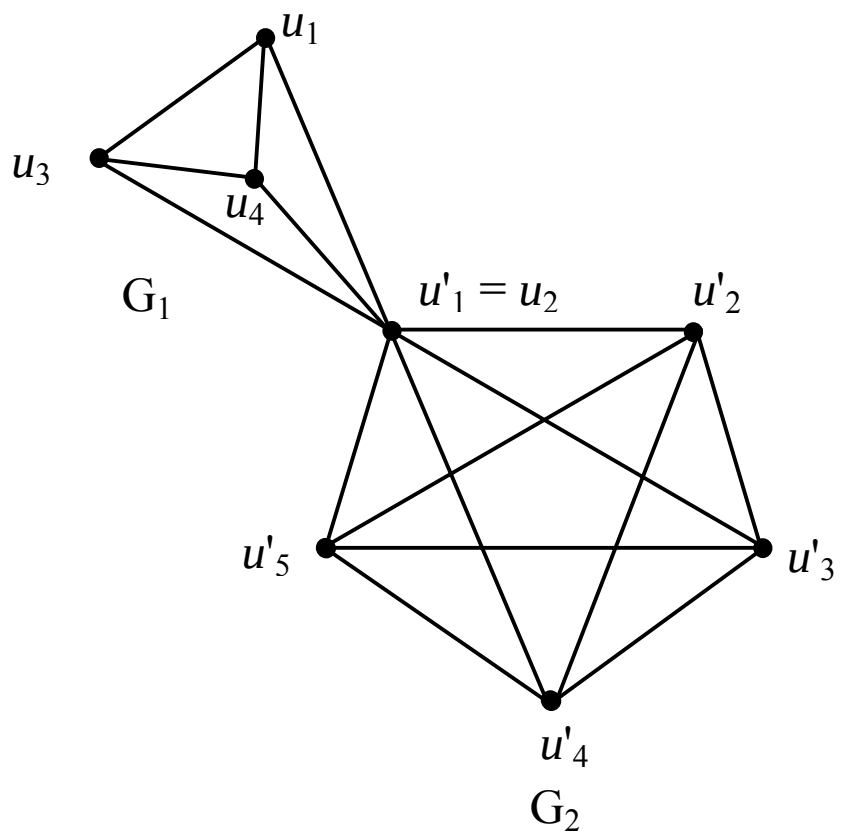

FIGURE: 2.1.19

$G = G_1 \cup G_2$ is a vertex glued bigraph. The vertex $u_2$ and $u'_1$ are glued. Clearly $G_1$ is 3-regular and $G_2$ is 4 regular. G =

$G_1 \cup G_2$ is such that it has only one vertex $u'_1 = u_2$ in common and degree of $u'_1 = u_2$ is 5. So G is a 5-biregular bigraph.

The notion of isolated vertex of the bigraph $G = G_1 \cup G_2$ and pendent vertex of G can be defined as in case of graphs.

The following example will illustrate the isolated vertex and pendent vertex of the bigraph $G = G_1 \cup G_2$.

***Example 2.1.20:*** Let $G = G_1 \cup G_2$ be a bigraph given by the following figure 2.1.20.

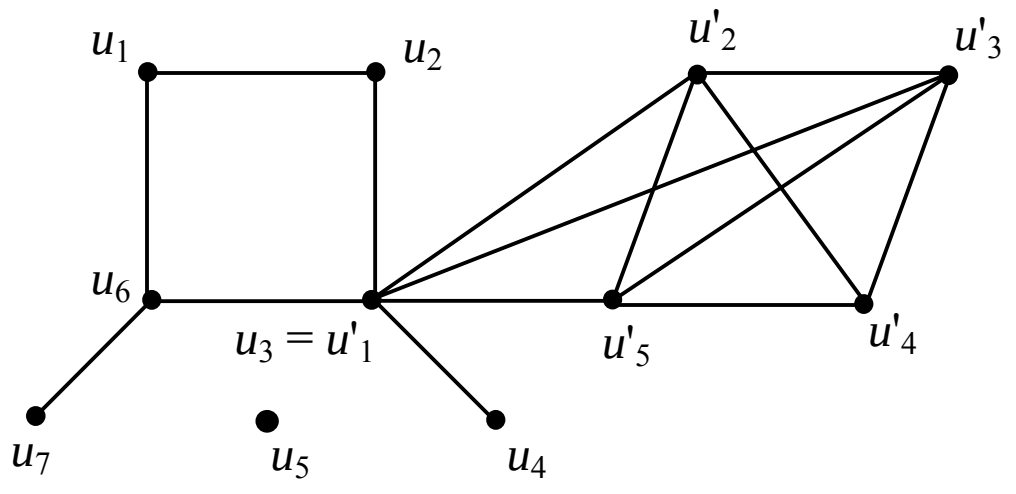

FIGURE: 2.1.20

$G_1 = \{u_1, u_2, \ldots, u_7\}$.

This is a vertex glued bigraph which has pendent vertices $u_4$ and $u_7$ and an isolated vertex $u_5$.

However the graph $G_2$ has no pendent vertex or an isolated vertex.

Now we give an example of a bigraph $G = G_1 \cup G_2$ where $G_1$ has both isolated vertex and pendent vertex but $G_2$, the graph has no pendent vertex or isolated vertex. This sort of bigraph $G = G_1 \cup G_2$ is said to have pseudo pendent vertex and pseudo isolated vertex. The following bigraph is an example of a bigraph with pseudo isolated vertex and pseudo pendent vertex.

***Example 2.1.21:*** Let $G = G_1 \cup G_2$ be a bigraph given by the following figure 2.1.21.

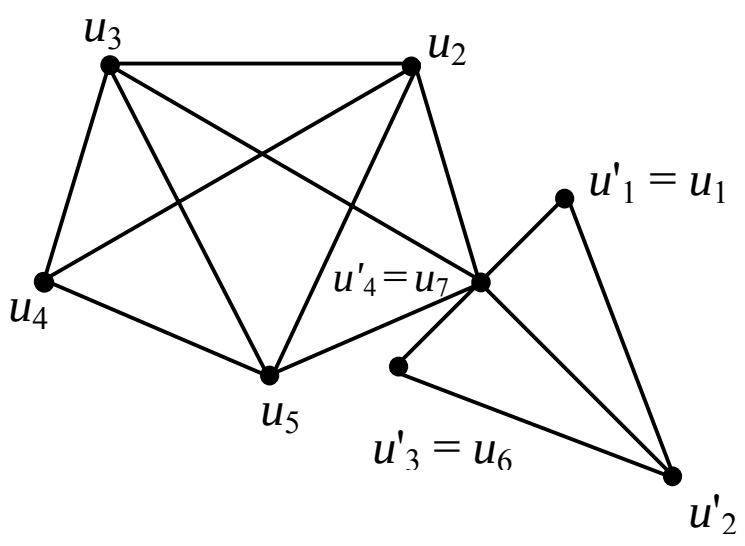

FIGURE: 2.1.21

with $G_1 = \{u_1, \ldots, u_7\}$ and $G_2 = \{u'_1, u'_2, u'_3, u'_4\}$.
The individual graphs of $G_1$ and $G_2$ is given in figure 2.1.21a and figure2.1.21b.

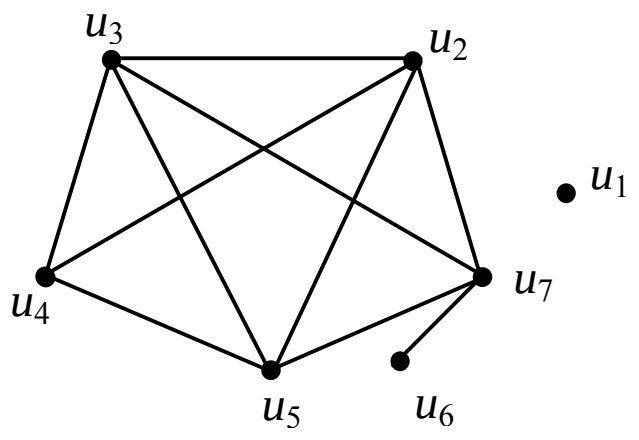

FIGURE: 2.1.21a

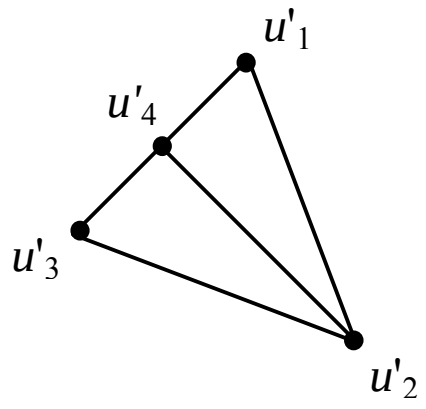

FIGURE: 2.1.21b

Clearly $G_1$ has both an isolated vertex $u_1$ and a pendent vertex $u_6$ but the bigraph $G = G_1 \cup G_2$ has no pendent vertex or isolated vertex.

It is easily seen that the Euler's theorem is true for all bigraphs. Let $G = G_1 \cup G_2$ be a bigraph. The join of $G_1 \vee G_2 = G$ is different from the bigraph G.

This is illustrated by the following example.

***Example 2.1.22:*** $G = G_1 \cup G_2$ is a disjoint bigraph given in figure 2.1.22.

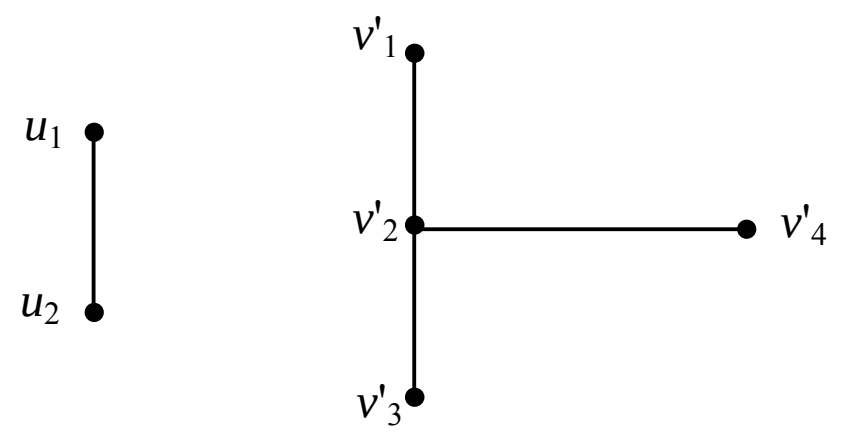

FIGURE: 2.1.22

But the join of $G_1 \vee G_2 = G^1$ is given by figure 2.1.22a.

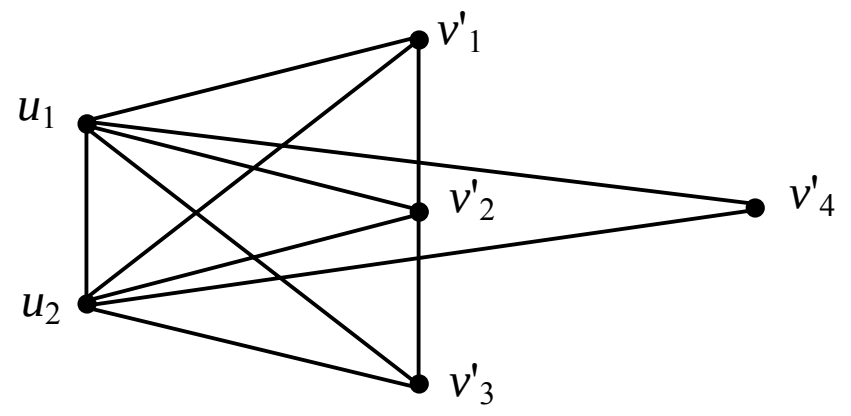

FIGURE: 2.1.22a

Clearly G and $G^1$ are distinct.
Suppose $G = G_1 \cup G_2$ got as a bigraph which is glued by an edge, given in the following figure 2.1.22b.

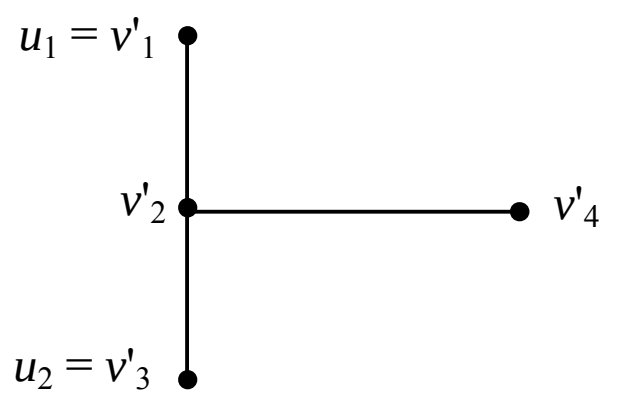

FIGURE: 2.1.22b

Thus $G = G_1 \cup G_2 \neq G^1$.
Now consider the bigraph $G = G_1 \cup G_2$ glued by a vertex given by the following figure 2.1.22c.

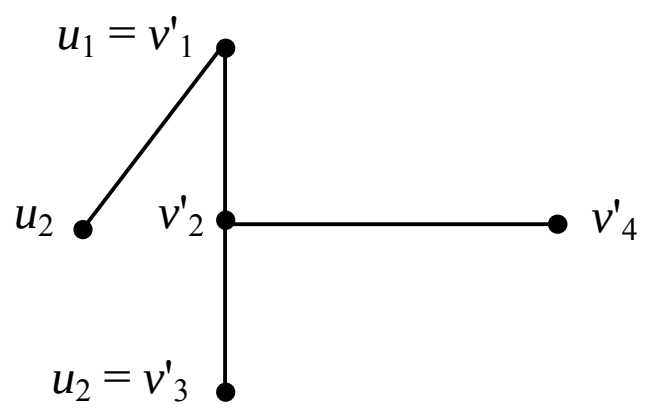

FIGURE: 2.1.22c

Clearly $G = G_1 \cup G_2 \neq G^1$. Thus the join of two graphs $G_1 \vee G_2$ is not the same as the bigraph given by $G = G_1 \cup G_2$.
Also we can show that in general the direct product two graphs $G_1$ and $G_2$ cannot be got as a bigraph $G_1 \cup G_2$. i.e. $G_1 \times G_2 \neq G_1 \cup G_2$.

For this is clear from the following example.

***Example 2.1.23:***

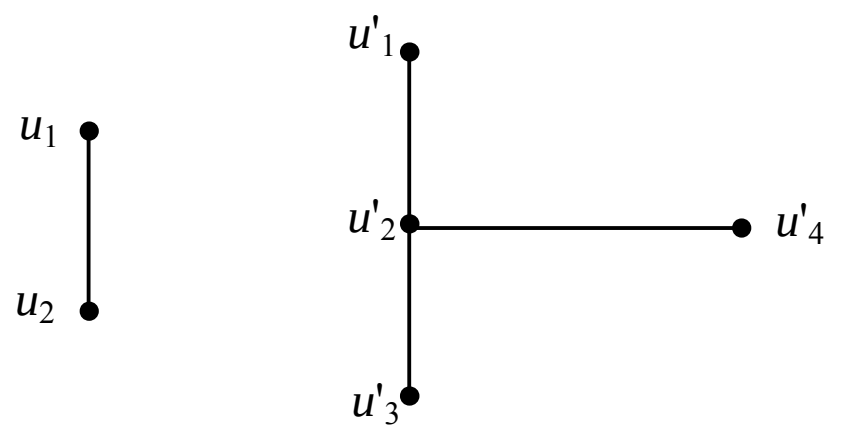

FIGURE: 2.1.23

$G = G_1 \cup G_2$ is the disjoint bigraphs of $G_1$ and $G_2$ given in figure 2.1.23.

Consider the direct product $G^1 = G_1 \times G_2$ is given by the following figure 2.1.23a.

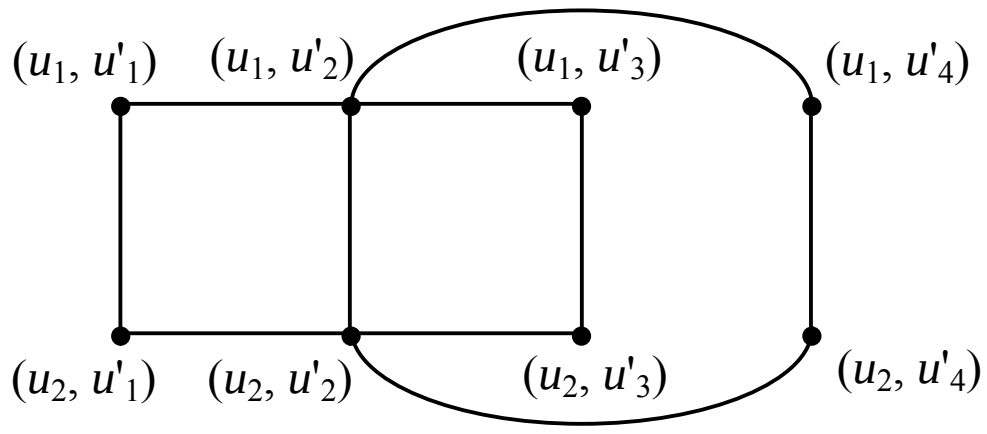

$G_1 \times G_2$

FIGURE: 2.1.23a

Now we proceed on to define the notion of directed bigraph.

***DEFINITION 2.1.11:*** *A directed bigraph $G = G_1 \cup G_2$ is a pair of ordered triple $\{(V(G_1), A(G_1), I_{G_1}), (V(G_2), A(G_2), I_{G_2})\}$ where $V(G_1)$ and $V(G_2)$ are non empty proper sets of $V(G)$ called the set of vertices of $G = G_1 \cup G_2$. $A(G_1)$ and $A(G_2)$ is a set disjoint from $V(G_1)$ and $V(G_2)$ respectively called the set of arcs of $G_1$ and $G_2$ and $I_{G_1}$ and $I_{G_2}$ are incidence map that associates with each arc of $G_1$ and $G_2$ an ordered pair of vertices of $G_1$ and $G_2$ respectively. A directed bigraph is called the dibigraph.*

All concepts incident into, incident out, outneighbour and inneighbour are defined as in case of graphs.

***Example 2.1.24:*** Let $G = G_1 \cup G_2$ be a bigraph given in figure 2.1.24.

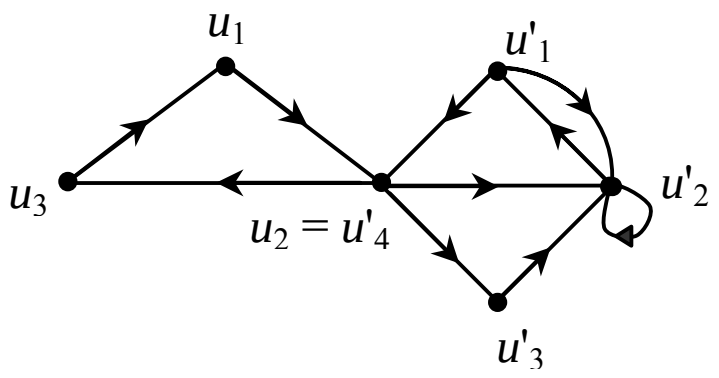

FIGURE: 2.1.24

The concept of dibigraph will prove a major role in the applications of bimodels both in fuzzy theory and neutrosophic theory.

Finally we give the matrix representations of bigraphs. First we give the simple bigraph and the related adjacency bimatrix.

***Example 2.1.25:*** Let $G = G_1 \cup G_2$ be the bigraph $G = G_1 \cup G_2$ in which both $G_1$ and $G_2$ are simple. The bigraph $G = G_1 \cup G_2$ is given by the figure 2.1.25.

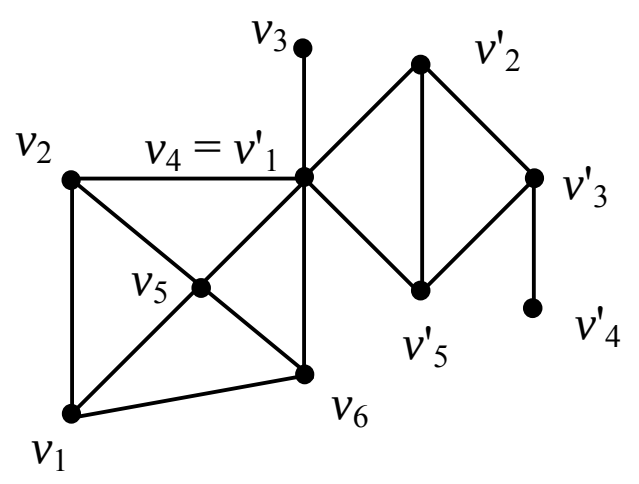

FIGURE: 2.1.25

The adjacency bimatrix of the bigraph is a mixed square bimatrix given as $X = X_1 \cup X_2$ where

$$X_1 = \begin{array}{c} \\ v_1 \\ v_2 \\ v_3 \\ v_4 \\ v_5 \\ v_6 \end{array} \begin{array}{c} \begin{array}{cccccc} v_1 & v_2 & v_3 & v_4 & v_5 & v_6 \end{array} \\ \left[\begin{array}{cccccc} 0 & 1 & 0 & 0 & 1 & 1 \\ 1 & 0 & 0 & 1 & 1 & 0 \\ 0 & 0 & 0 & 1 & 0 & 0 \\ 0 & 1 & 1 & 0 & 1 & 1 \\ 1 & 1 & 0 & 1 & 0 & 1 \\ 1 & 0 & 0 & 1 & 1 & 0 \end{array}\right] \end{array}$$

$$X_2 = \begin{array}{c} \\ v_1 \\ v_2 \\ v_3 \\ v_4 \\ v_5 \end{array} \begin{array}{c} \begin{array}{ccccc} v'_1 & v'_2 & v'_3 & v'_4 & v'_5 \end{array} \\ \left[\begin{array}{ccccc} 0 & 1 & 0 & 0 & 1 \\ 1 & 0 & 1 & 0 & 1 \\ 0 & 1 & 0 & 1 & 1 \\ 0 & 0 & 1 & 0 & 0 \\ 1 & 1 & 1 & 0 & 0 \end{array}\right] \end{array}$$

Thus $X = X_1 \cup X_2$ is the adjacency bimatrix.

Always the adjacency bimatrix of a bigraph will be the square bimatrix.

Now we proceed on to illustrate the weighted bigraph of the bigraph $G = G_1 \cup G_2$ by the following example.

***Example 2.1.26*****:** Let $G = G_1 \cup G_2$ be a bigraph which is a weighted bigraph given by the following figure 2.1.26.

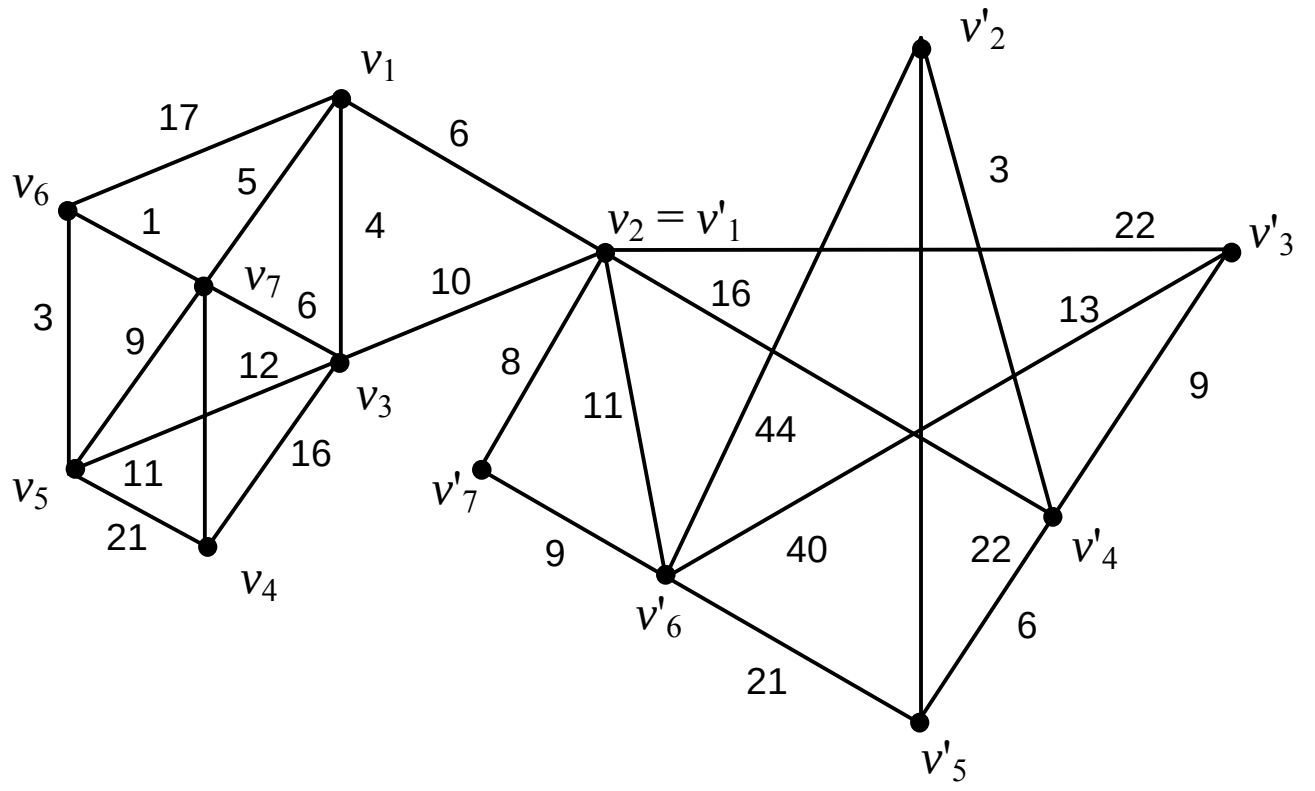

FIGURE: 2.1.26

We give below the weighted matrix related with the bigraph $G = G_1 \cup G_2$. The bimatrix of a weighted bigraph is always a square bimatrix. Let $W = W_1 \cup W_2$ be the weighted bimatrix of $G = G_1 \cup G_2$.

$$W_1 = \begin{array}{c} \\ v_1 \\ v_2 \\ v_3 \\ v_4 \\ v_5 \\ v_6 \\ v_7 \end{array} \begin{array}{c} \begin{array}{ccccccc} v_1 & v_2 & v_3 & v_4 & v_5 & v_6 & v_7 \end{array} \\ \begin{bmatrix} \infty & 6 & 4 & \infty & \infty & 17 & 5 \\ 6 & \infty & 10 & \infty & \infty & \infty & \infty \\ 4 & 10 & \infty & 16 & 11 & \infty & 6 \\ \infty & \infty & 16 & \infty & 21 & \infty & 12 \\ \infty & \infty & 11 & 21 & \infty & 3 & 9 \\ 17 & \infty & \infty & \infty & 3 & \infty & 1 \\ 5 & \infty & 6 & 12 & 9 & 1 & \infty \end{bmatrix} \end{array} \cup$$

$$W_2 = \begin{array}{c} \\ v_1' \\ v_2' \\ v_3' \\ v_4' \\ v_5' \\ v_6' \\ v_7' \end{array} \begin{array}{c} \begin{array}{ccccccc} v_1' & v_2' & v_3' & v_4' & v_5' & v_6' & v_7' \end{array} \\ \begin{bmatrix} \infty & \infty & 22 & 16 & \infty & 11 & 8 \\ \infty & \infty & \infty & 3 & 22 & 4 & \infty \\ 22 & \infty & \infty & 9 & \infty & 13 & \infty \\ 16 & 3 & 9 & \infty & 6 & \infty & \infty \\ \infty & 22 & \infty & 6 & \infty & 21 & \infty \\ 11 & 4 & 40 & \infty & 21 & \infty & 9 \\ 8 & \infty & \infty & \infty & \infty & 9 & \infty \end{bmatrix} \end{array}.$$

'$\infty$' symbol denotes when the vertices are non adjacent. It is to be noted that given any weighted bimatrix which is symmetric one can get back to the bigraph and given the bigraph we always have a bimatrix which weighed can be associated with it. Suppose $G = G_1 \cup G_2$ is a dibigraph then we can have the incident bimatrix associated with it. The incidence matrix of a digraph with n vertices and e edges and no self loops.

We show by the following example the incidence bimatrix associated with the dibigraph.

***Example 2.1.27:*** Let $G = G_1 \cup G_2$ be a bigraph given by the following figure 2.1.27.

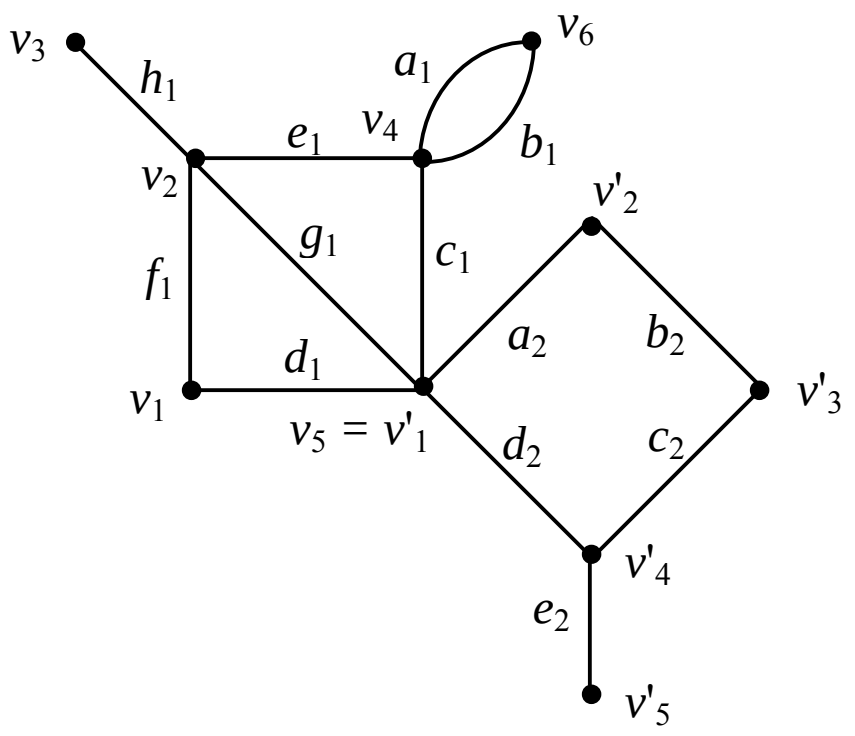

FIGURE: 2.1.27

The incidence bimatrix of the bigraph.

$$\begin{array}{c} \\ v_1 \\ v_2 \\ v_3 \\ v_4 \\ v_5 \\ v_6 \end{array} \begin{array}{c} \begin{array}{cccccccc} a_1 & b_1 & c_1 & d_1 & e_1 & f_1 & g_1 & h_1 \end{array} \\ \begin{bmatrix} 0 & 0 & 0 & 1 & 0 & 1 & 0 & 0 \\ 0 & 0 & 0 & 0 & 1 & 1 & 1 & 1 \\ 0 & 0 & 0 & 0 & 0 & 0 & 0 & 1 \\ 1 & 1 & 1 & 0 & 1 & 0 & 0 & 0 \\ 0 & 0 & 1 & 1 & 0 & 0 & 1 & 0 \\ 1 & 1 & 0 & 0 & 0 & 0 & 0 & 0 \end{bmatrix} \end{array} \cup \begin{array}{c} \\ v'_1 \\ v'_2 \\ v'_3 \\ v'_4 \\ v'_5 \end{array} \begin{array}{c} \begin{array}{ccccc} a_2 & b_2 & c_2 & d_2 & e_2 \end{array} \\ \begin{bmatrix} 1 & 0 & 0 & 1 & 0 \\ 1 & 1 & 0 & 0 & 0 \\ 0 & 1 & 1 & 0 & 0 \\ 0 & 0 & 1 & 1 & 1 \\ 0 & 0 & 0 & 1 & 0 \end{bmatrix} \end{array}$$

The incidence bimatrix associated with the bigraph is a mixed rectangular bimatrix. The incidence bimatrix can also be written with directions. It is very interesting to note all bigraphs which are disjoint bigraphs are trivially separable bigraphs. Further all bigraphs glued by a vertex or single vertex glued bigraphs are separable.

**THEOREM 2.1.6:** Let $G = G_1 \cup G_2$ which is a single vertex glued bigraph G is separable.

*Proof:* Given $G = G_1 \cup G_2$ is a bigraph which is a single vertex glued bigraph say let them be glued by the vertex $v_j = v'_1$ by removing that vertex, the bigraph becomes the separable bigraph.

**DEFINITION 2.1.12:** *A bigraph $G = G_1 \cup G_2$ is connected if there is at least one path between every pair of vertices in G other wise G is disconnected. The disjoint bigraph $G = G_1 \cup G_2$ is disconnected. A bigraph $G = G_1 \cup G_2$ is connected if both $G_1$ and $G_2$ are connected and the bigraph is vertex glued bigraph or edge glued bigraph or a subgraph glued bigraph. A bitree is a connected bigraph without any circuits.*

Now we proceed on to define cut-set in a connected bigraph.

**DEFINITION 2.1.13:** *In a connected bigraph $G = G_1 \cup G_2$ a bicut set is a set of edges (together with their end vertices) whose removal from $G = G_1 \cup G_2$ leaves both the graphs $G_1$ and $G_2$ disconnected, provided removal of no proper subset of these edges disconnects G.*
*Further we see the bicut set of a connected bigraph $G = G_1 \cup G_2$ need not be unique.*

We can have more than one bicut set of a bigraph.

We illustrate this by an example.

***Example 2.1.28:*** Let $G = G_1 \cup G_2$ be a connected bigraph given by the following figure 2.1.28.

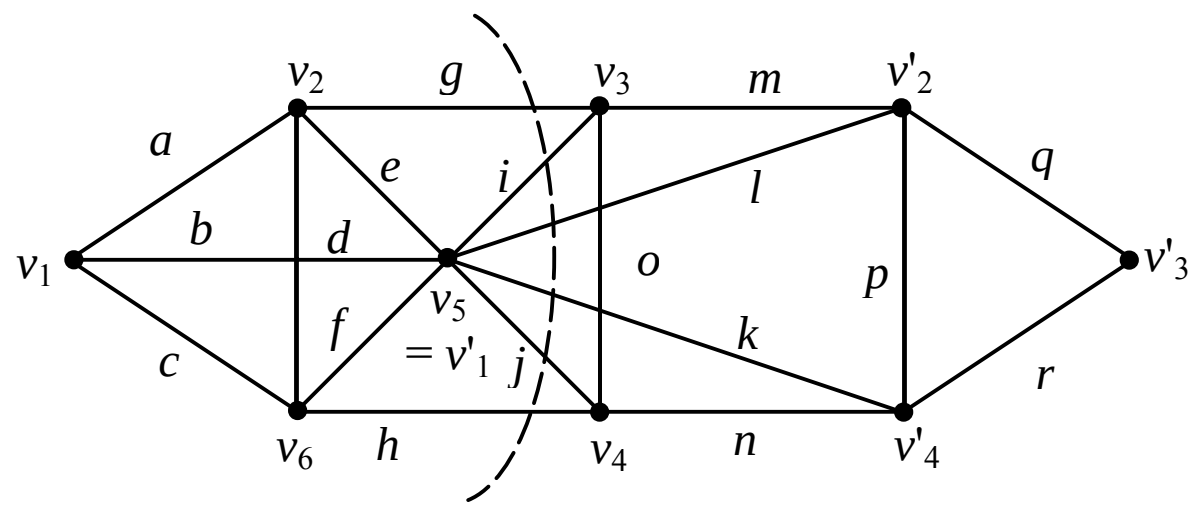

FIGURE: 2.1.28

$G = G_1 \cup G_2 = \{v_1, v_2, v_3, \ldots, v_6\} \cup \{v'_1, v'_2, v'_3, v'_4\}$.

The graph of $G_1$ and $G_2$ are separately given by the following figures.

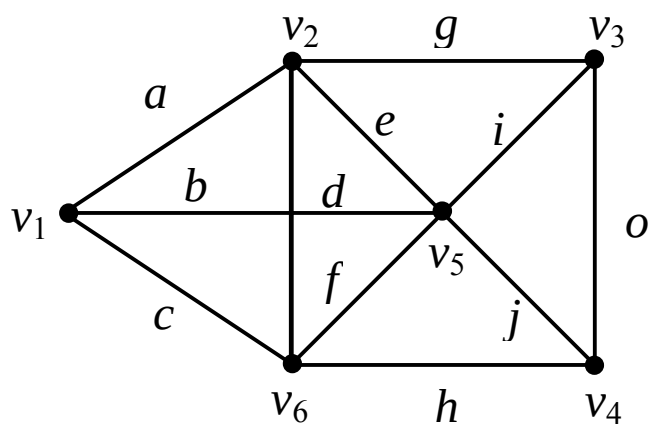

FIGURE: 2.1.28a

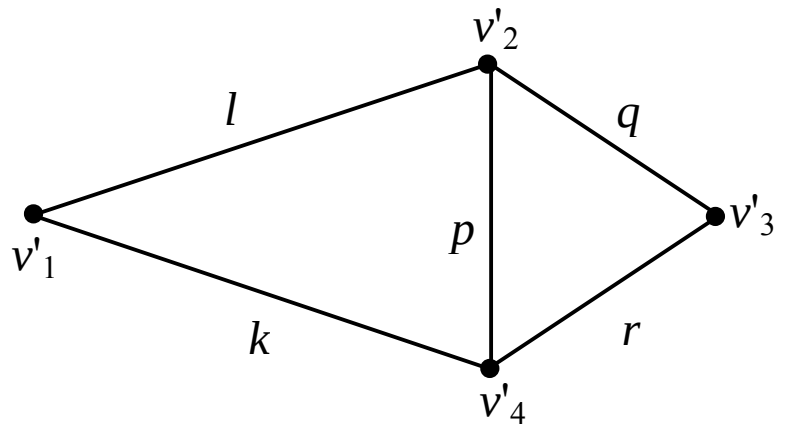

FIGURE: 2.1.28b

Removal of the bicut set {g, i, l, k, j, h} from the bigraph is given by the following figures.

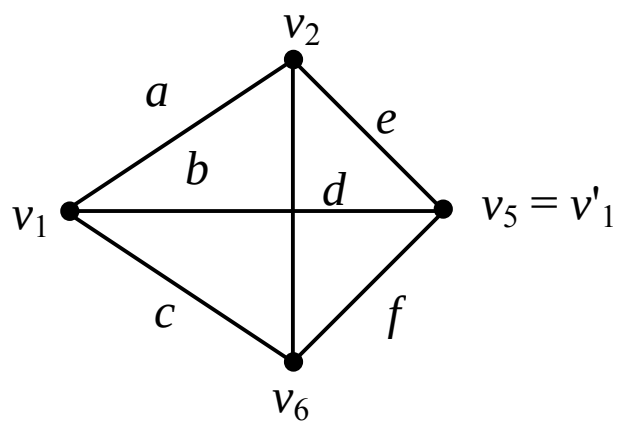

FIGURE: 2.1.28c

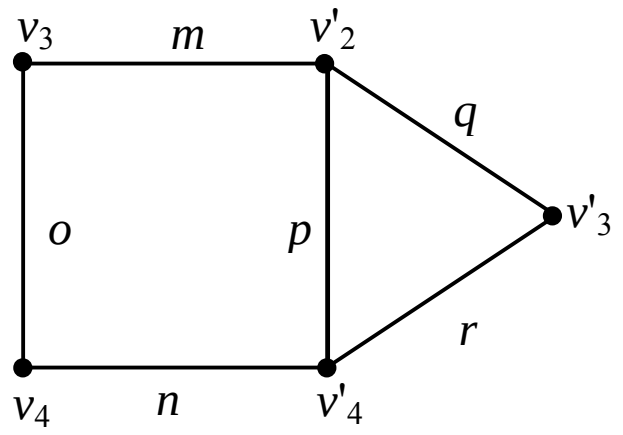

FIGURE: 2.1.28d

Clearly H and K are bigraphs where

$$H = H_1 \cup H_2$$
$$= \{v_1, v_2, v_5, v_6\} \cup \{v'_1\}.$$

and

$$K = K_1 \cup K_2$$
$$= \{v_3, v_4\} \cup \{v'_2, v'_4, v'_3\}.$$

Several related results can be derived; however we define the edge connectivity of a bigraph.

**DEFINITION 2.1.14:** *Each bicut set of a connected bigraph $G = G_1 \cup G_2$ consists of a certain number of edges. The number of edges in the smallest bicut set defined as the edge connectivity of the bigraph $G = G_1 \cup G_2$.*

We are interested in applying bigraphs to network flows. We now proceed on to define separable bigraphs.

**DEFINITION 2.1.15:** *A separable bigraph consists of two or more non separable bisubgraphs. Each of the largest non separable sub bigraphs is called a biblock or bicomponent.*

We just illustrate this by the following example.

***Example 2.1.29:*** Let $G = G_1 \cup G_2$ be a separable bigraph given by the following figure.

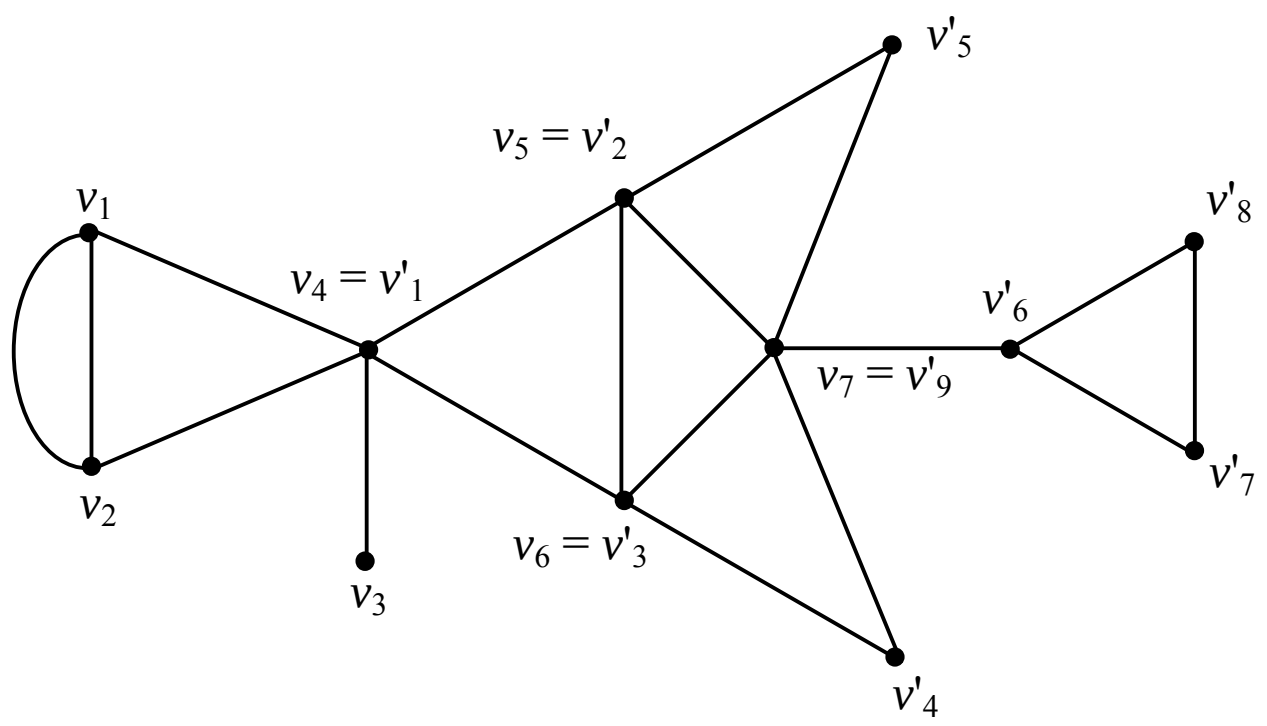

FIGURE: 2.1.29

$G = G_1 \cup G_2$

where

$$V(G_1) = \{v_1, v_2, v_3, \ldots, v_7\}$$

and

$$V(G_2) = \{v'_1, v'_2, v'_3, \ldots, v'_9\}.$$

The separate graphs of $G_1$ and $G_2$ are given by the following figure. H and K are subbigraphs of $G = G_1 \cup G_2$.

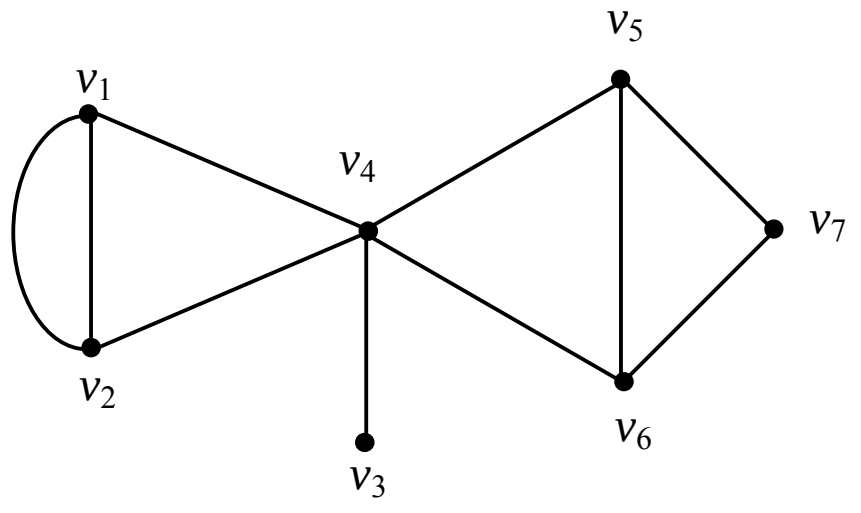

FIGURE: 2.1.29a

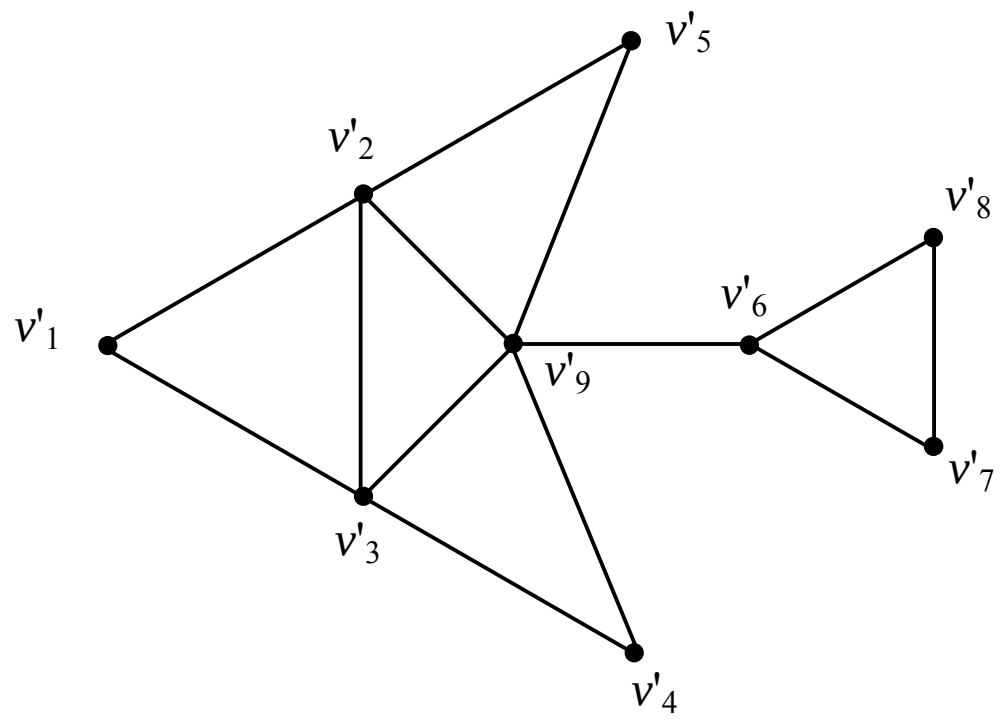

FIGURE: 2.1.29b

Now we see the division of the bigraph $G = G_1 \cup G_2$ into subbigraphs given in figures 2.1.29c and 2.1.29d.

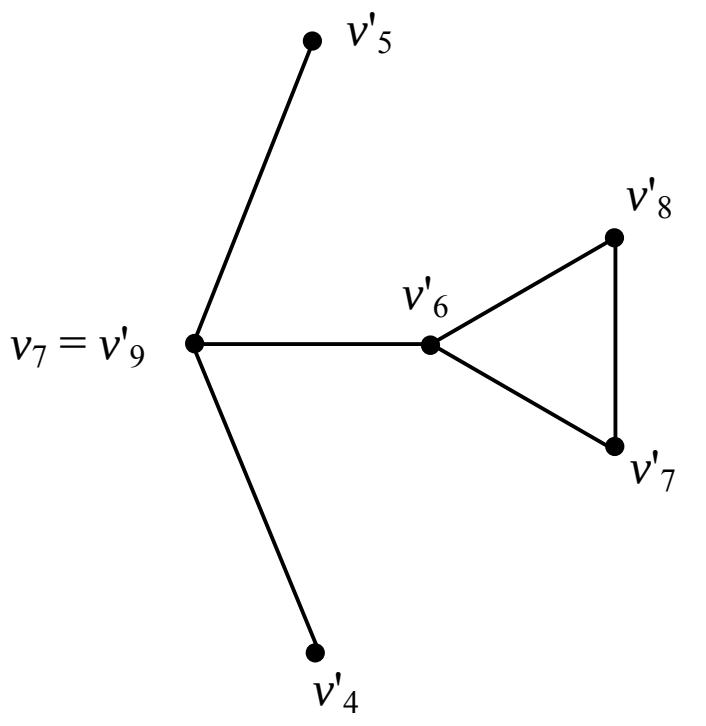

FIGURE: 2.1.29c

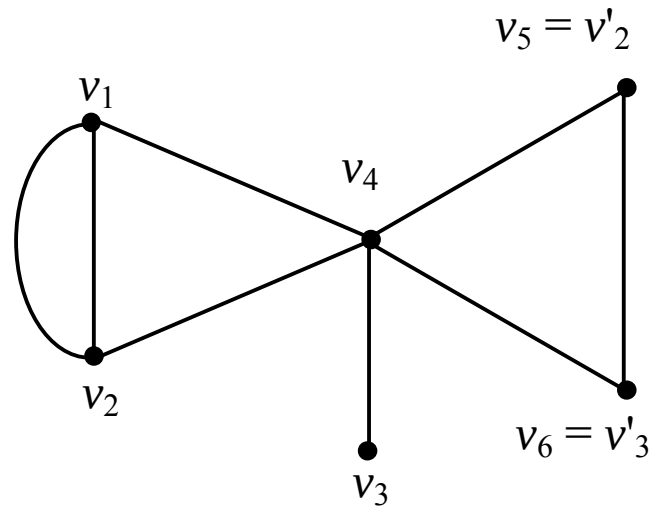

FIGURE: 2.1.29d

Now we can extend the notion of bigraphs to trigraphs, quadruple graphs and so on say to n-graphs $n \geq 2$.

**DEFINITION 2.1.16:** *$G = G_1 \cup G_2 \cup G_3$ is said to be a trigraph if G is the union of three distinct graphs i.e. $G_i$ is not a subgraph of $G_j$ $1 \leq i \leq 3$ and $1 \leq j \leq 3$. When the vertex set of $G_1$, $G_2$ and $G_3$ are disjoint we call G as the strongly disjoint trigraph. If a pair $G_i \cup G_j$, $i \neq j$ is not a disjoint bigraph but $G_k$ is disjoint from $G_i \cup G_j$ we call this disjoint trigraph. A trigraph can be realized as the union of a bigraph and a graph.*

Now we illustrate them by the following examples.

***Example 2.1.30:*** Let G = $G_1 \cup G_2 \cup G_3$ be the trigraph given by the following figure 2.1.30.

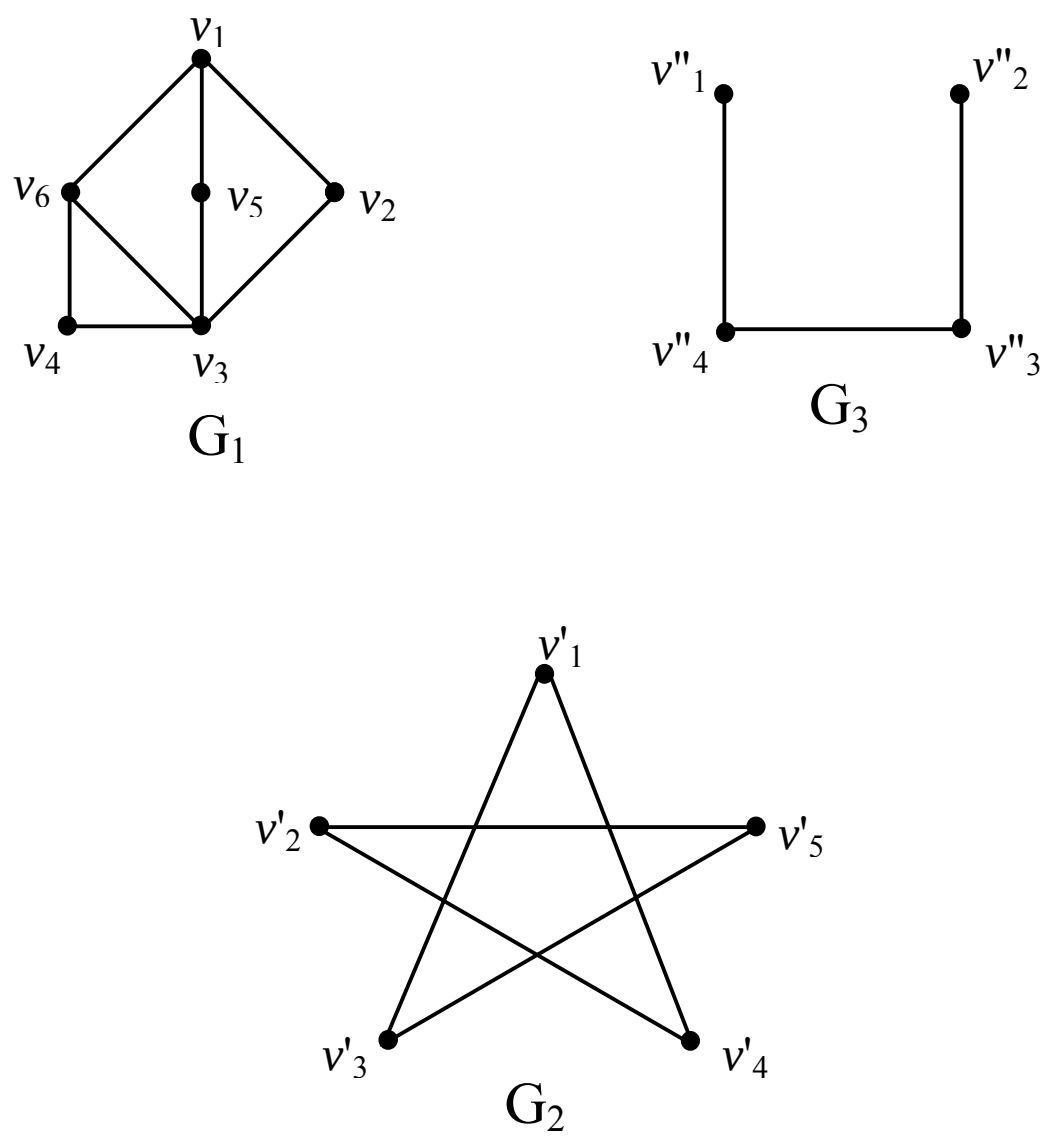

FIGURE: 2.1.30

Clearly G is a strongly disjoint trigraph.

***Example 2.1.31:*** Let G = $G_1 \cup G_2 \cup G_3$ be a trigraph given by the following figure 2.1.31.

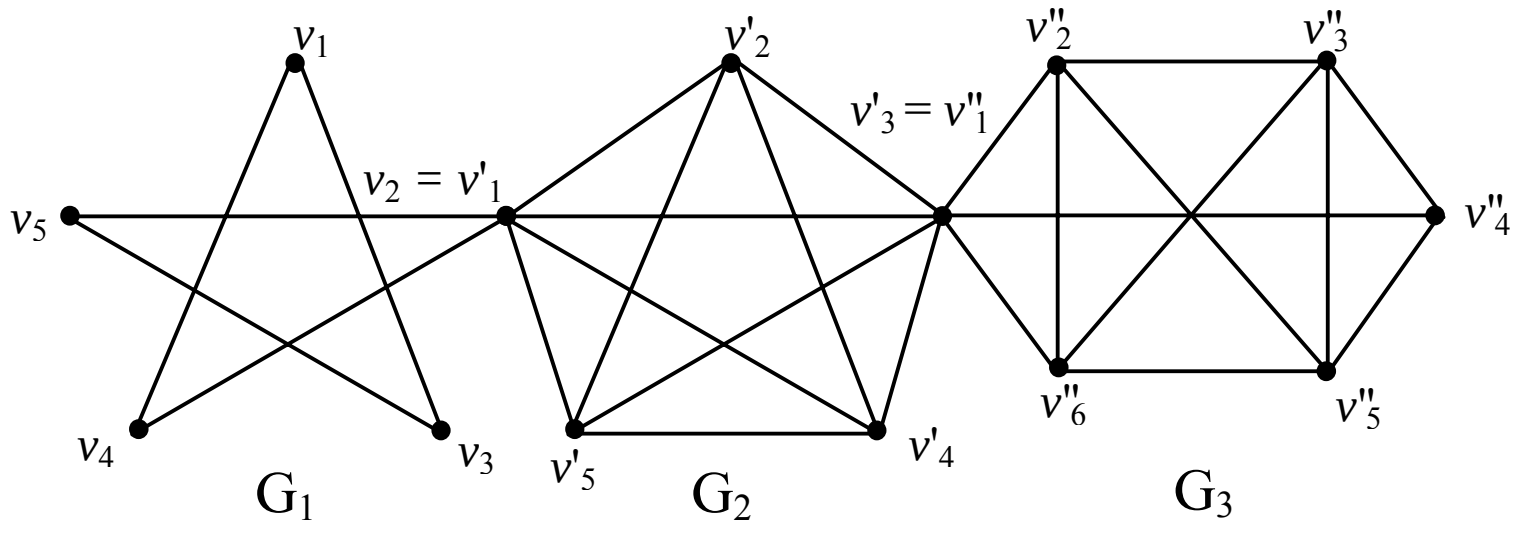

FIGURE 2.1.31

The vertex set of G has only 14 points and the individual graphs of the trigraph is given by the following figures.

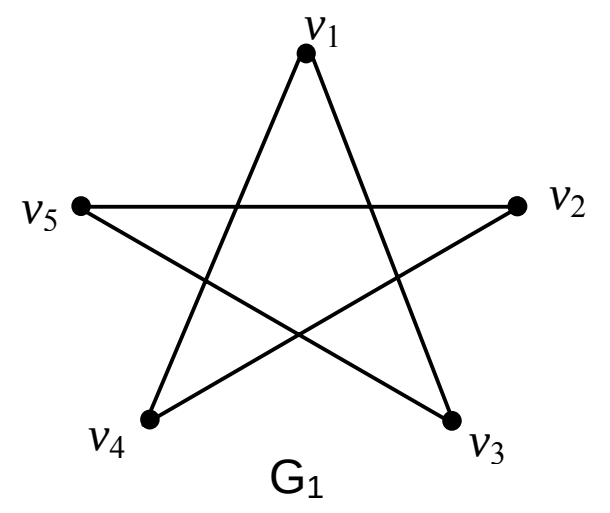

FIGURE 2.1.31a

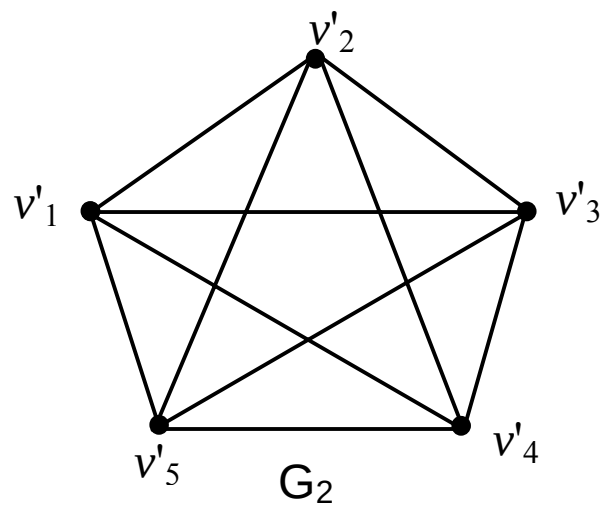

FIGURE 2.1.31b

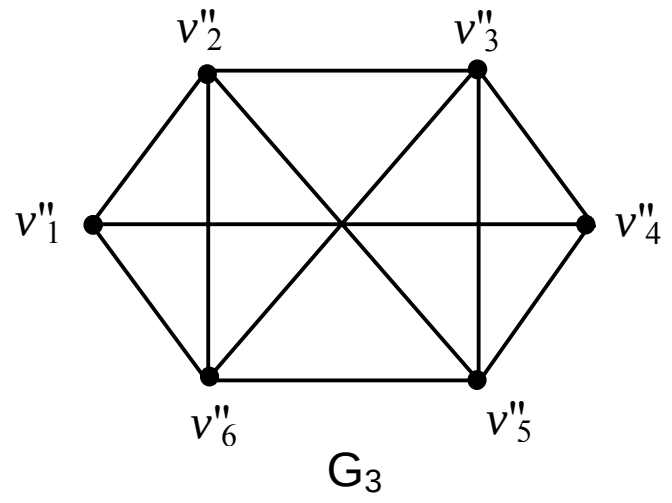

FIGURE 2.1.31c

Infact $G = G_1 \cup G_2 \cup G_3$ is a trigraph which is connected. Infact G is not a disjoint trigraph.

Now we proceed on to give an example of just a disjoint trigraph.

***Example 2.1.32:*** Let G = $G_1 \cup G_2 \cup G_3$ be the trigraph given by the following figure 2.1.32.

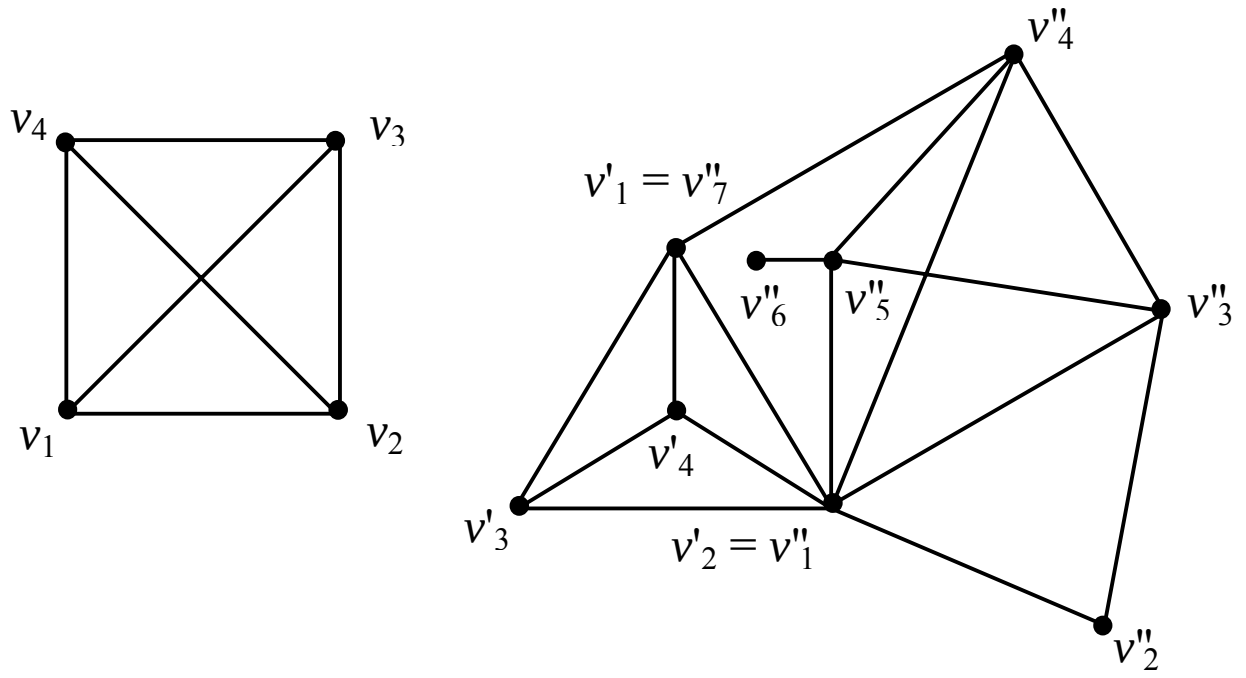

FIGURE 2.1.32

Clearly the individual graphs of the trigraph are given as

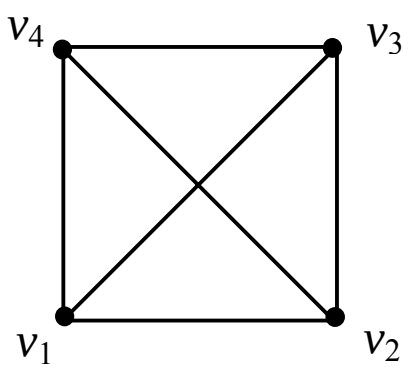

FIGURE 2.1.32a

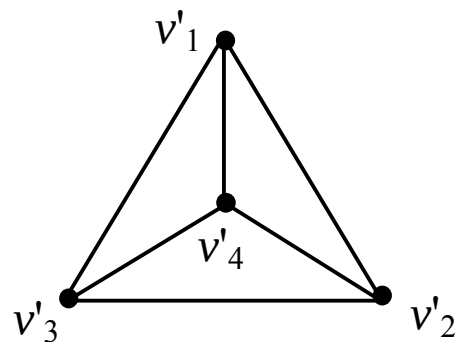

FIGURE 2.1.32b

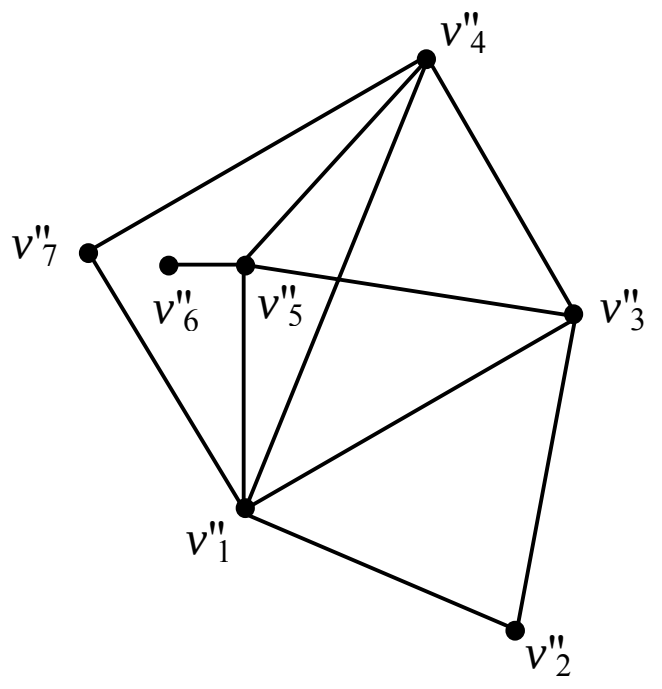

FIGURE 2.1.32c

Clearly this trigraph is a disjoint trigraph and the bigraph in it is a edge glued bigraph. Thus we have seen three types of trigraphs. In fact all properties of graphs can be proved with suitable and necessary modification for every trigraph is also a graph.

Now we define a quadgraph $G = G_1 \cup G_2 \cup G_3 \cup G_4$ as the union of four graphs $G_1$, $G_2$, $G_3$ and $G_4$ such that the vertex set of any one of them is not a proper subset of the other. If the four graphs happen to be disjoint we call them disjoint quadgraph. If all the graphs are connected by vertex we call the quadgraphs the vertex glued graph.

If they are glued by an edge they will be known as edge glued quadgraph. A quadgraph can also be realized as the union of a pair of bigraph or as the union of a graph and a trigraph.

We give some examples to illustrate the quadgraphs.

***Example 2.1.33:*** Let $G = G_1 \cup G_2 \cup G_3 \cup G_4$ be a quadgraph given by the following figure.

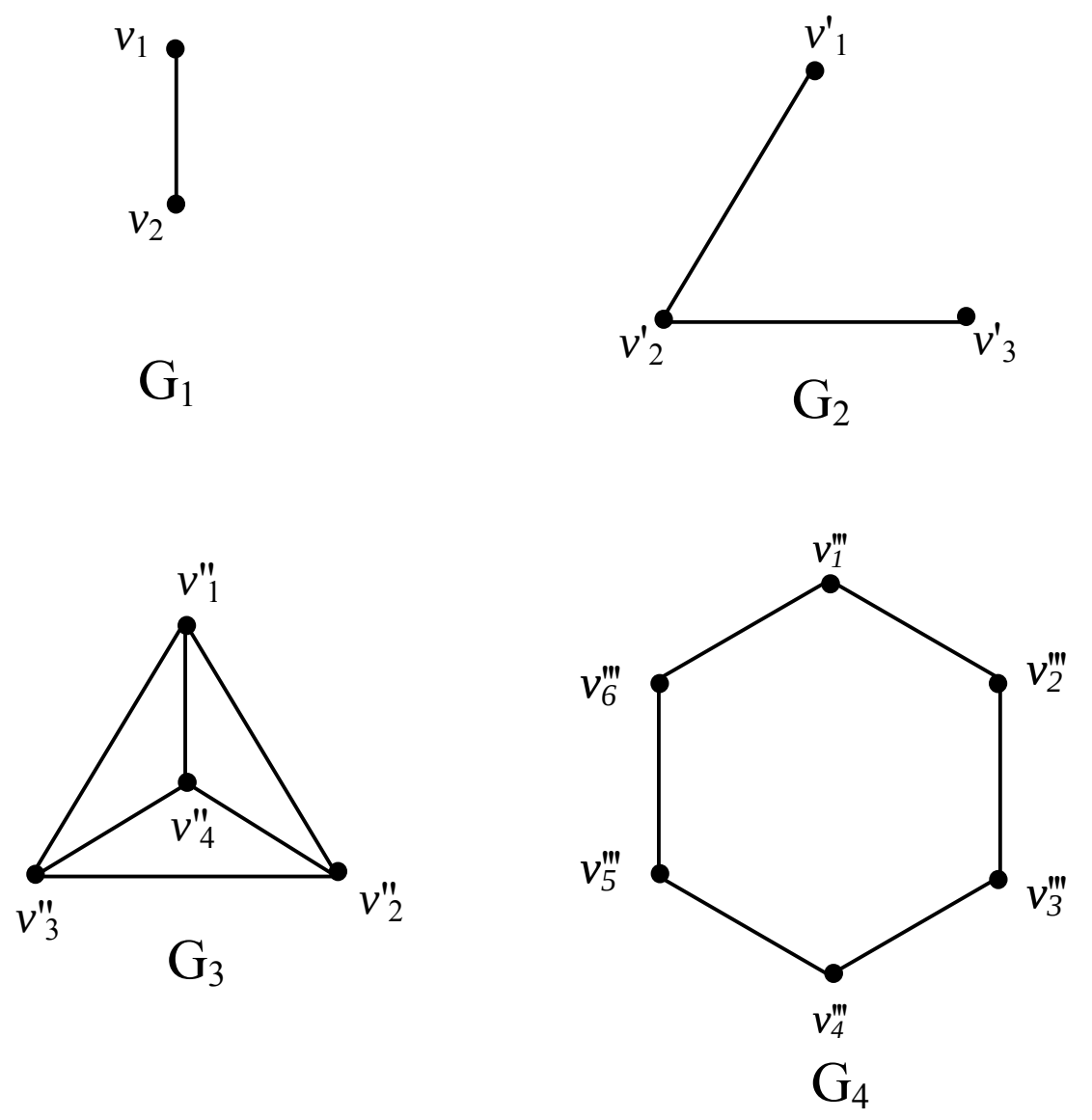

FIGURE 2.1.33

The vertex set of $G = G_1 \cup G_2 \cup G_3 \cup G_4$ is

$$V(G) = \{v_1, v_2, v'_1, v'_2, v'_3, v''_1, v''_2, v''_3, v'''_1, v'''_2, v'''_3, v'''_4, v'''_5\}.$$

This is a strongly disjoint quadgraph.

Now we proceed on to give an illustration of a disjoint quadgraph.

***Example 2.1.34:***

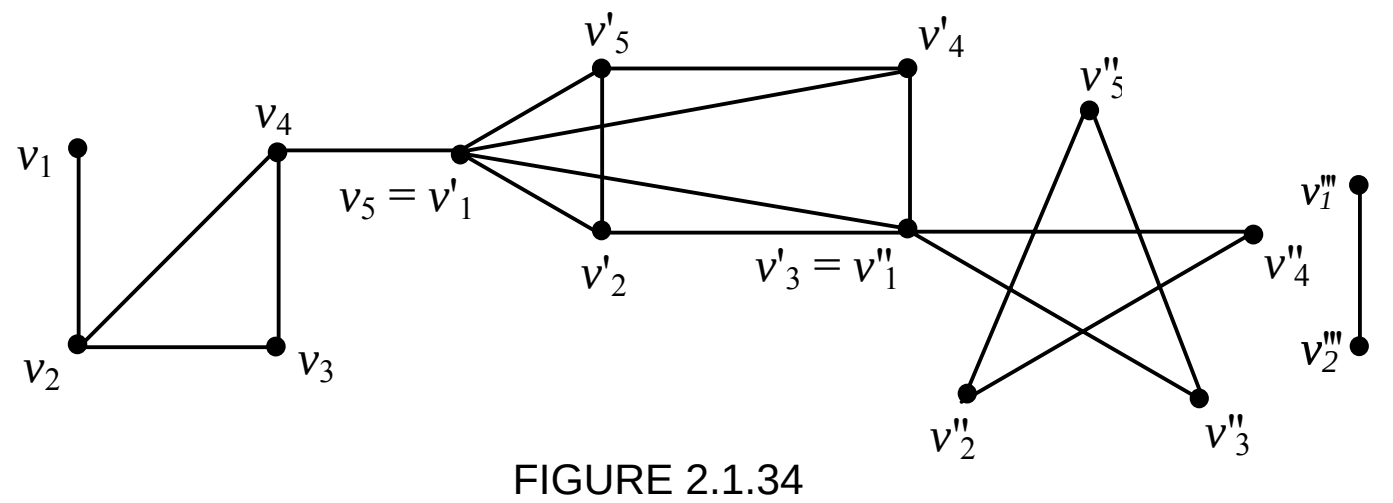

FIGURE 2.1.34

$$G = G_1 \cup G_2 \cup G_3 \cup G_4$$

We give the figures of each of the graphs $G_1$, $G_2$, $G_3$ and $G_4$. (Figures 2.1.34a, 2.1.34b, 2.1.34c and 2.1.34d).

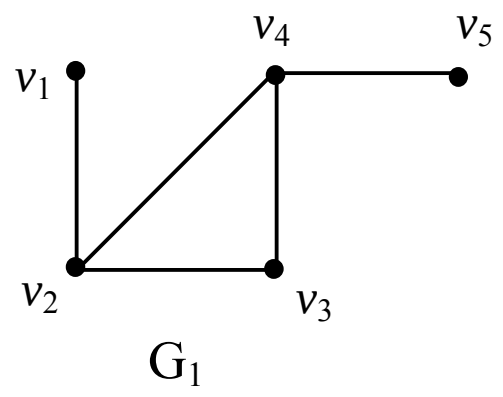

FIGURE 2.1.34a

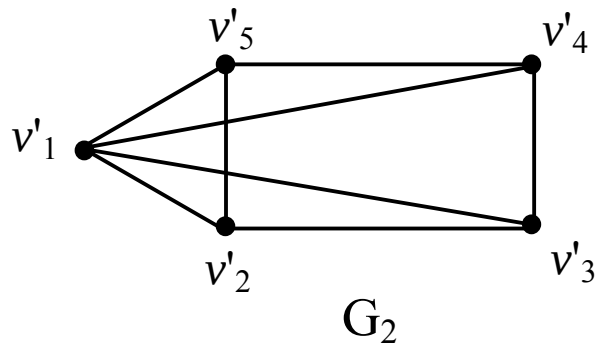

FIGURE 2.1.34b

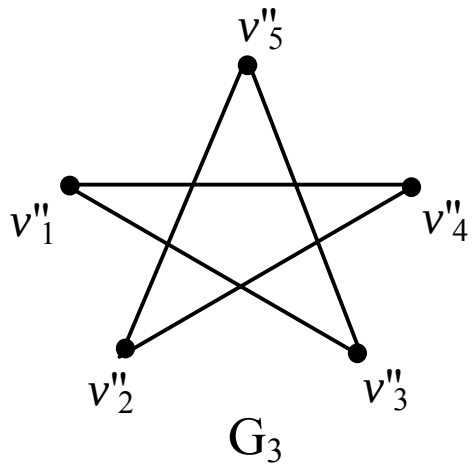

FIGURE 2.1.34c

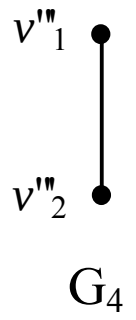

FIGURE 2.1.34d

Now we proceed on to define when a quadgraph is said to be disjoint bigraph of a quadgraph.

This is illustrated by the following example.

***Example 2.1.35:*** Let $G = G_1 \cup G_2 \cup G_3 \cup G_4$ be the quadgraph given by the following diagram. The quadgraph is a disjoint union of two connected bigraphs.

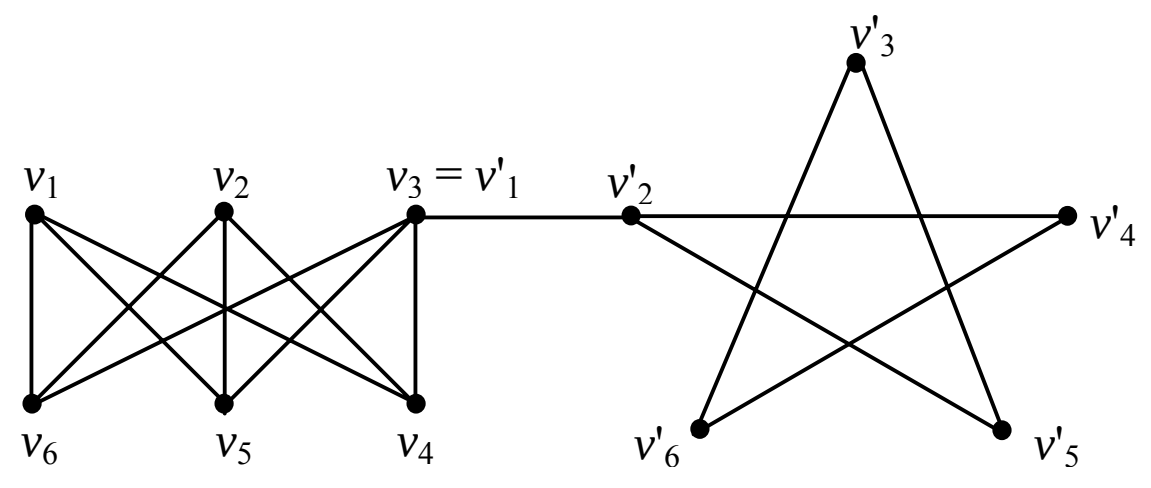

$G_1 \cup G_2$

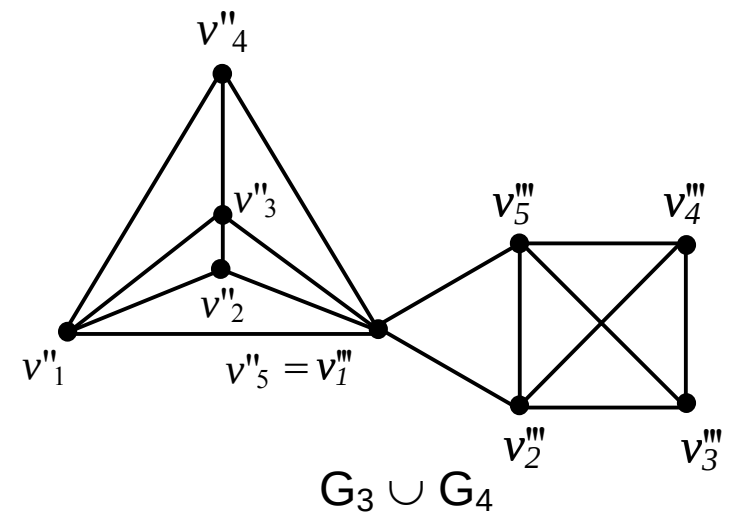

$G_3 \cup G_4$

FIGURE 2.1.35

Thus $G = G_1 \cup G_2 \cup G_3 \cup G_4$ is a quadgraph which is the union of two bigraphs $G_1 \cup G_2$ and $G_3 \cup G_4$. Thus a quadgraph can also be defined as $G = H_1 \cup H_2$ where $H_1$ and $H_2$ are bigraphs.

Thus we can have n-graphs $G = G_1 \cup G_2 \cup \ldots \cup G_n$ all properties studied and defined for bigraphs can be extended easily to n-graphs as even a n-graph is also a special type of graph.

Here it is very pertinent to mention that all bistructures behave very differently from a bigraph for all n-graphs are also bigraphs but while giving it the adjacency matrix representation a graph and an n graphs are differentiated. That is only the main place a vast difference occurs for every graph is associated with a matrix but a bigraph is represented by a bimatrix and a bimatrix is not a matrix.

Thus by a trimatrix $A = A_1 \cup A_2 \cup A_3$ we mean A is just identified with three matrices the only conditions is that all the three matrices must be different, can be of any rank or order.

***Example 2.1.36:*** Thus we show the following is a trimatrix $A = A_1 \cup A_2 \cup A_3$.

$$A = \begin{bmatrix} 3 & 1 \\ 5 & 6 \\ 0 & 1 \end{bmatrix} \cup \begin{bmatrix} 6 & 0 & 0 & 0 \\ 1 & 1 & 0 & 0 \\ 0 & 1 & 1 & 0 \\ -4 & 2 & 0 & -1 \end{bmatrix} \cup \begin{bmatrix} 5 \\ 4 \\ 3 \\ 2 \\ 1 \\ 0 \end{bmatrix}$$

This type of trimatrix will be called as the mixed trimatrix.

***Example 2.1.37:*** The trimatrix of the form

$$A = \begin{bmatrix} 3 & 0 \\ 1 & 1 \end{bmatrix} \cup \begin{bmatrix} 0 & 1 \\ 0 & 1 \end{bmatrix} \cup \begin{bmatrix} -2 & 0 \\ 0 & 0 \end{bmatrix}$$

will be called as a $2 \times 2$ square trimatrix.

***Example 2.1.38:*** A trimatrix of the form

$$A = \begin{bmatrix} 3 & 1 \\ 0 & 1 \\ 1 & 1 \end{bmatrix} \cup \begin{bmatrix} 1 & 1 \\ 0 & 0 \\ 0 & 0 \end{bmatrix} \cup \begin{bmatrix} 1 & 0 \\ 1 & 0 \\ 1 & 0 \end{bmatrix}$$

will be called as a 3 × 2 rectangular trimatrix.

***Example 2.1.39:*** Let $A = A_1 \cup A_2 \cup A_3$

$$= \begin{bmatrix} 2 & 0 \\ 1 & 1 \end{bmatrix} \cup \begin{bmatrix} -1 & 0 & 0 \\ 0 & 1 & 0 \\ 2 & 2 & -1 \end{bmatrix} \cup \begin{bmatrix} 0 & 1 & 0 & 1 \\ 0 & 2 & 0 & -4 \\ 1 & 1 & 0 & 0 \\ 1 & 0 & 1 & 0 \end{bmatrix}$$

is a mixed square trimatrix.

On similar lines we define n-matrix A which is the 'union' of n distinct matrices $A_1, A_2, \ldots, A_n$ denoted by $A = A_1 \cup A_2 \cup \ldots \cup A_n$ where '$\cup$' is just the symbol.
Now we give the adjacency matrix associated with the trigraph $G = G_1 \cup G_2 \cup G_3$. Clearly the adjacency connection trimatrix is just a mixed square trimatrix.

***Example 2.1.40:*** Let $G = G_1 \cup G_2 \cup G_3$ be the trigraph given by the following figure.

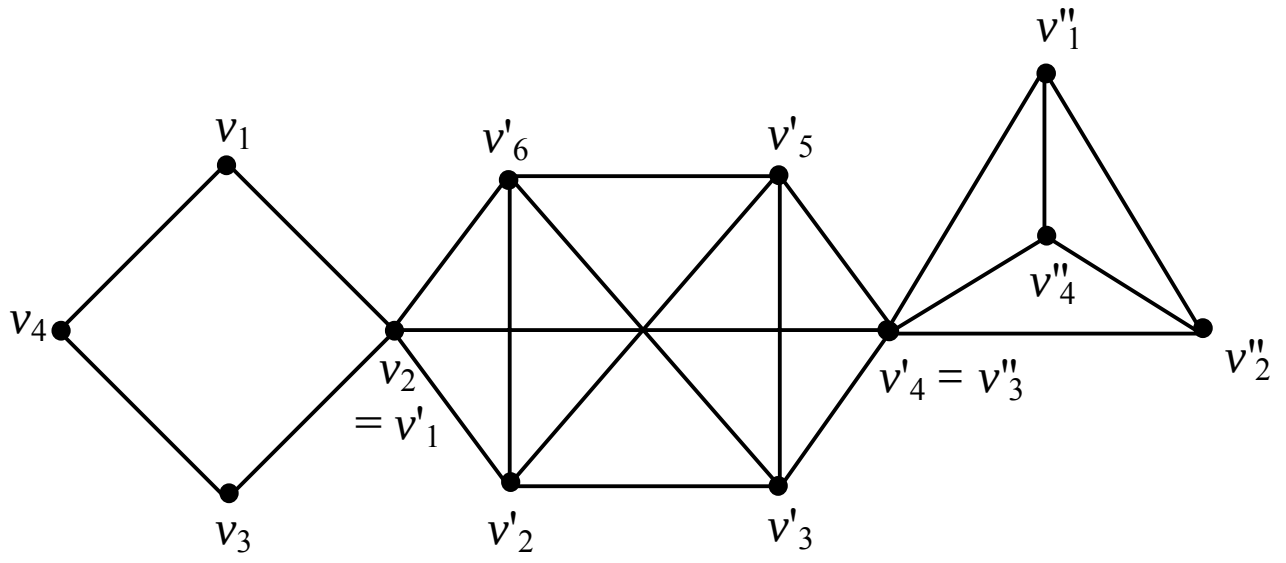

FIGURE 2.1.36

$$G = G_1 \cup G_2 \cup G_3$$

The separate graphs of $G_1$, $G_2$, $G_3$ is given by the following figures.

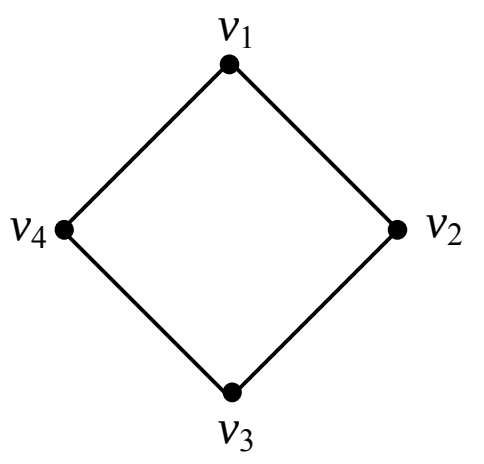

G₁

FIGURE 2.1.36a

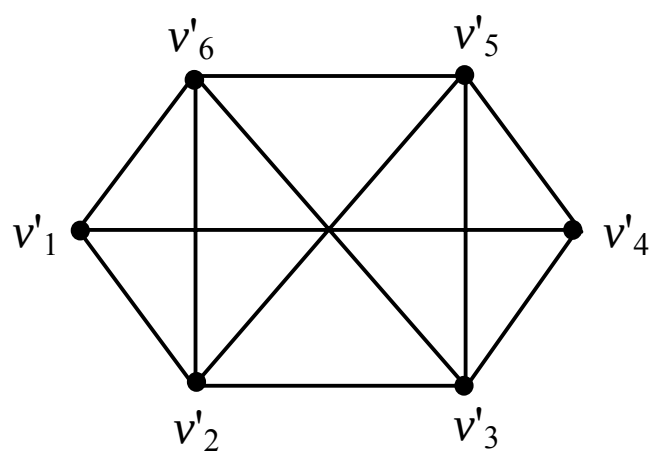

G₂

FIGURE 2.1.36b

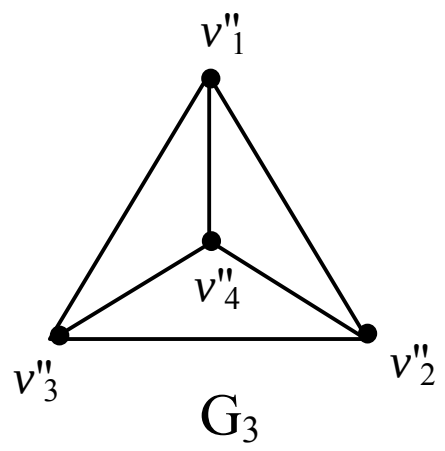

G₃

FIGURE: 2.1.36c

The related trimatrix for the trigraph is $C = C_1 \cup C_2 \cup C_3$

$$C = \begin{array}{c} \\ v_1 \\ v_2 \\ v_3 \\ v_4 \end{array}\begin{array}{c} \begin{array}{cccc} v_1 & v_2 & v_3 & v_4 \end{array} \\ \begin{bmatrix} 0 & 1 & 0 & 1 \\ 1 & 0 & 1 & 0 \\ 0 & 1 & 0 & 1 \\ 1 & 0 & 1 & 0 \end{bmatrix} \end{array} \cup \begin{array}{c} \\ v'_1 \\ v'_2 \\ v'_3 \\ v'_4 \\ v'_5 \\ v'_6 \end{array}\begin{array}{c} \begin{array}{cccccc} v'_1 & v'_2 & v'_3 & v'_4 & v'_5 & v'_6 \end{array} \\ \begin{bmatrix} 0 & 1 & 0 & 1 & 0 & 1 \\ 1 & 0 & 1 & 0 & 1 & 1 \\ 0 & 1 & 0 & 1 & 1 & 1 \\ 1 & 0 & 1 & 0 & 1 & 1 \\ 0 & 1 & 1 & 1 & 0 & 1 \\ 1 & 1 & 1 & 0 & 1 & 0 \end{bmatrix} \end{array} \cup$$

$$\begin{array}{c} \\ v''_1 \\ v''_2 \\ v''_3 \\ v''_4 \end{array}\begin{array}{c} \begin{array}{cccc} v''_1 & v''_2 & v''_3 & v''_4 \end{array} \\ \begin{bmatrix} 0 & 1 & 1 & 1 \\ 1 & 0 & 1 & 1 \\ 1 & 1 & 0 & 1 \\ 1 & 1 & 1 & 0 \end{bmatrix} \end{array}$$

On similar lines we can also get the weighted trimatrix associated with a trigraph. Likewise the rectangular path matrix which will be a trimatrix can also be obtained.

**Chapter Three**

# APPLICATION OF BIMATRICES TO NEW FUZZY MODELS

In this chapter we give applications of bimatrices to specially constructed fuzzy models. This chapter has six sections. In the first section we just recall the definition of FCMs. In section two for the first time we introduce the notion of fuzzy cognitive bimaps and give its applications. Section three gives the extension of applications. Section three gives the extension of bimaps to fuzzy cognitive trimaps and illustrates with applications. The concept of fuzzy relational maps are recalled in section four. Fuzzy relational bimaps and its application are introduced for the first time in this section five. The final section introduces the notion of fuzzy relational trimaps and illustrates it with examples.

## 3.1 Definition of Fuzzy Cognitive Maps

In this section we recall the notion of Fuzzy Cognitive Maps (FCMs), which was introduced by Bart Kosko [72 to 76]. We also give several of its interrelated definitions. FCMs have a major role to play mainly when the data concerned is an unsupervised one. Further this method is most simple and an effective one as it can analyse the data by directed graphs and connection matrices.

**DEFINITION 3.1.1:** *An FCM is a directed graph with concepts like policies, events etc. as nodes and causalities as edges. It represents causal relationship between concepts.*

We illustrate this by the example 3.1.1:

***Example 3.1.1:*** In Tamil Nadu (a southern state in India) in the last decade several new engineering colleges have been approved and started. The resultant increase in the production of engineering graduates in these years is disproportionate with the need of engineering graduates. This has resulted in thousands of unemployed and underemployed graduate engineers. Using an expert's opinion we study the effect of such unemployed people on the society. An expert spells out the five major concepts relating to the unemployed graduated engineers as

$E_1$ – Frustration
$E_2$ – Unemployment
$E_3$ – Increase of educated criminals
$E_4$ – Under employment
$E_5$ – Taking up drugs etc.

The directed graph where $E_1, \ldots, E_5$ are taken as the nodes and causalities as edges as given by an expert is given in the following Figure 3.1.1:

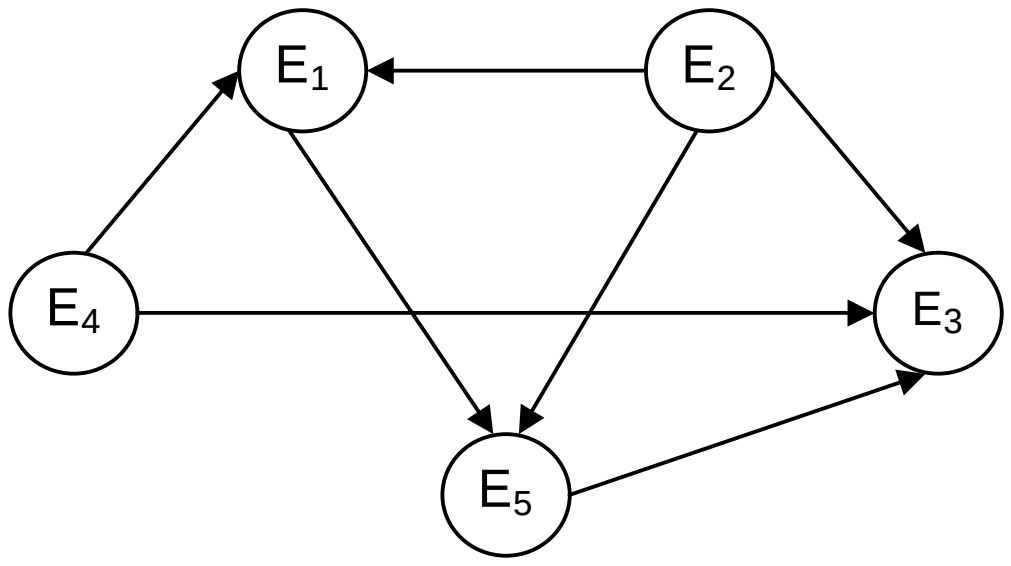

FIGURE: 3.1.1

According to this expert, increase in unemployment increases frustration. Increase in unemployment, increases the educated criminals. Frustration increases the graduates to take up to evils like drugs etc. Unemployment also leads to the increase in number of persons who take up to drugs, drinks etc. to forget their worries and unoccupied time. Under-employment forces then to do criminal acts like theft (leading to murder) for want of more money and so on. Thus one cannot actually get data for this but can use the expert's opinion for this unsupervised data to obtain some idea about the real plight of the situation. This is just an illustration to show how FCM is described by a directed graph.

{If increase (or decrease) in one concept leads to increase (or decrease) in another, then we give the value 1. If there exists no relation between two concepts the value 0 is given. If increase (or decrease) in one concept decreases (or increases) another, then we give the value –1. Thus FCMs are described in this way.}

**DEFINITION 3.1.2:** *When the nodes of the FCM are fuzzy sets then they are called as fuzzy nodes.*

**DEFINITION 3.1.3:** *FCMs with edge weights or causalities from the set {–1, 0, 1} are called simple FCMs.*

**DEFINITION 3.1.4:** *Consider the nodes / concepts $C_1, \ldots, C_n$ of the FCM. Suppose the directed graph is drawn using edge weight $e_{ij} \in \{0, 1, -1\}$. The matrix E be defined by $E = (e_{ij})$ where $e_{ij}$ is the weight of the directed edge $C_i C_j$. E is called the adjacency matrix of the FCM, also known as the connection matrix of the FCM.*

It is important to note that all matrices associated with an FCM are always square matrices with diagonal entries as zero.

**DEFINITION 3.1.5:** *Let $C_1$, $C_2$, … , $C_n$ be the nodes of an FCM. $A = (a_1, a_2, \ldots, a_n)$ where $a_i \in \{0, 1\}$. A is called the instantaneous state vector and it denotes the on-off position of the node at an instant.*

$$a_i = 0 \text{ if } a_i \text{ is off and}$$
$$a_i = 1 \text{ if } a_i \text{ is on}$$

*for $i = 1, 2, \ldots, n$.*

**DEFINITION 3.1.6:** *Let $C_1$, $C_2$, … , $C_n$ be the nodes of an FCM. Let $\overrightarrow{C_1C_2}$, $\overrightarrow{C_2C_3}$, $\overrightarrow{C_3C_4}$, … , $\overrightarrow{C_iC_j}$ be the edges of the FCM ($i \neq j$). Then the edges form a directed cycle. An FCM is said to be cyclic if it possesses a directed cycle. An FCM is said to be acyclic if it does not possess any directed cycle.*

**DEFINITION 3.1.7:** *An FCM with cycles is said to have a feedback.*

**DEFINITION 3.1.8:** *When there is a feedback in an FCM, i.e., when the causal relations flow through a cycle in a revolutionary way, the FCM is called a dynamical system.*

**DEFINITION 3.1.9:** *Let $\overline{C_1C_2}$, $\overline{C_2C_3}$, … , $\overrightarrow{C_{n-1}C_n}$ be a cycle. When $C_i$ is switched on and if the causality flows through the edges of a cycle and if it again causes $C_i$ , we say that the dynamical system goes round and round. This is true for any node $C_i$ , for $i = 1, 2, \ldots, n$. The equilibrium state for this dynamical system is called the hidden pattern.*

**DEFINITION 3.1.10:** *If the equilibrium state of a dynamical system is a unique state vector, then it is called a fixed point.*

***Example 3.1.2:*** Consider a FCM with $C_1$, $C_2$, …, $C_n$ as nodes. For example let us start the dynamical system by switching on $C_1$. Let us assume that the FCM settles down with $C_1$ and $C_n$ on i.e. the state vector remains as (1, 0, 0,

…, 0, 1) this state vector (1, 0, 0, …, 0, 1) is called the fixed point.

**DEFINITION 3.1.11:** *If the FCM settles down with a state vector repeating in the form*

$$A_1 \to A_2 \to \ldots \to A_i \to A_1$$

*then this equilibrium is called a limit cycle.*

Methods of finding the hidden pattern are discussed in the following.

**DEFINITION 3.1.12:** *Finite number of FCMs can be combined together to produce the joint effect of all the FCMs. Let $E_1, E_2, \ldots, E_p$ be the adjacency matrices of the FCMs with nodes $C_1, C_2, \ldots, C_n$ then the combined FCM is got by adding all the adjacency matrices $E_1, E_2, \ldots, E_p$.*

*We denote the combined FCM adjacency matrix by $E = E_1 + E_2 + \ldots + E_p$.*

**NOTATION:** Suppose $A = (a_1, \ldots, a_n)$ is a vector which is passed into a dynamical system E. Then $AE = (a'_1, \ldots, a'_n)$ after thresholding and updating the vector suppose we get $(b_1, \ldots, b_n)$ we denote that by

$$(a'_1, a'_2, \ldots, a'_n) \to (b_1, b_2, \ldots, b_n).$$

Thus the symbol '→' means the resultant vector has been thresholded and updated.

FCMs have several advantages as well as some disadvantages. The main advantage of this method it is simple. It functions on expert's opinion. When the data happens to be an unsupervised one the FCM comes handy. This is the only known fuzzy technique that gives the hidden pattern of the situation. As we have a very well known theory, which states that the strength of the data depends on, the number of experts opinion we can use combined FCMs with several experts opinions.

At the same time the disadvantage of the combined FCM is when the weightages are 1 and –1 for the same $C_i$ $C_j$, we have the sum adding to zero thus at all times the connection matrices $E_1, \ldots, E_k$ may not be conformable for addition.

Combined conflicting opinions tend to cancel out and assisted by the strong law of large numbers, a consensus emerges as the sample opinion approximates the underlying population opinion. This problem will be easily overcome when the FCM model has entries only 0 and 1.

We have just briefly recalled the definitions. For more about FCMs please refer Kosko [72 to 76].

## 3.2 Fuzzy Cognitive Bimaps and their Applications

Just for the sake of completeness we have given the basic concepts of FCM in the earlier section. Now we proceed on to define the notion of Fuzzy cognitive bimaps.

Suppose we have some unsupervised data to be analyzed and suppose it has two sets of disjoint attributes to be analyzed. The two sets of attributes are unrelated but are to be analyzed using FCMs they would have 2 directed graphs or a bigraph related with it so a bimatrix is a dynamical bisystem which will give us the bihidden pattern i.e. the limit bicycle or the fixed bipoint. How to work for such types of models and construct such models so that we can analyze the problem. Thus we now proceed on to define Fuzzy Cognitive bimaps.

**DEFINITION 3.2.1:** *Fuzzy Cognitive bimaps (FCBMs) are fuzzy signed directed bigraphs with feed back. The directed edge $e_{ij}^p$ from causal concept $c_i^p$ to concept $c_j^p$ measures how much $c_i^p$ causes $c_j^p$, (p = 1, 2). The time varying* *concept function $c_i^p(t)$ measures the non negative occurrence of some fuzzy event, perhaps the strength of a political statement, in medical analysis or so on.*
*The edge values $e_{ij}^p$ takes values in the fuzzy causal interval [-1, 1], $e_{ij}^p = 0$ indicates no causality, $e_{ij}^p > 0$ indicates causal increase, $c_j^p$ increases as $c_i^p$ increases (or $c_j^p$ decreases as $c_i^p$ decreases); $e_{ij}^p < 0$ indicates causal decrease or negative causality $c_i^p$ decreases as $c_j^p$ increases (and or $c_i^p$ increases or $c_j^p$ decreases) (p = 1, 2).*

Simple FCBMs have edge values $e_{ij}^p \in \{-1, 0, 1\}$, (p = 1, 2). We give a simple illustration of FCBMs.

***Example 3.2.1:*** Let us consider the model which is to analyze the problem faced by the industry. All problems faced while running a industry or a factory cannot be put as a statistical data. Several of them are feelings involving a great deal of uncertainty and impreciseness. In order to run the industry smoothly and with atleast some profit one should know and try to analyze the problem. To get some sort of frictionless feelings among workers, among the financers and above all the relation in between the workers and boss i.e. what we mean the relation between the employee and the employers.

Thus to have good profit the sales should be good which indirectly means the impact of their products in the public has a good standing and rapport. So the problem involved is multi dimensional.

So if we wish to divide this into two sets of problems and want to make use of the fuzzy cognitive bimaps (FCBMs). Suppose the industry on one side analyzes the factors promoting business and considers the following five attributes.

| | | |
|---|---|---|
| $C_1$ | – | Good business |
| $C_2$ | – | Good investment |
| $C_3$ | – | Customer Satisfaction |
| $C_4$ | – | Establishment |
| $C_5$ | – | Marketing strategies |

and also at the same time wishes to study the employee problems .

Employer relationship with employee so that the factory runs smoothly. The attributes given by an expert in the analysis of the employee-employer model is as follows:

| | | |
|---|---|---|
| $E_1$ | – | Maximum profit to the employer |
| $E_2$ | – | Just profit to the employer |
| $E_3$ | – | Neither profit nor loss to the employer |
| $E_4$ | – | Loss to the employer |
| $E_5$ | – | Best performance by the employee |
| $E_6$ | – | Only pay to employee |
| $E_7$ | – | Employee workers more number of hours. |
| $E_8$ | – | Average performance by few employee |
| $E_9$ | – | Poor performance by some employee |

Now we give the Fuzzy Cognitive bimaps. The directed bigraphs related with the model is as follows.

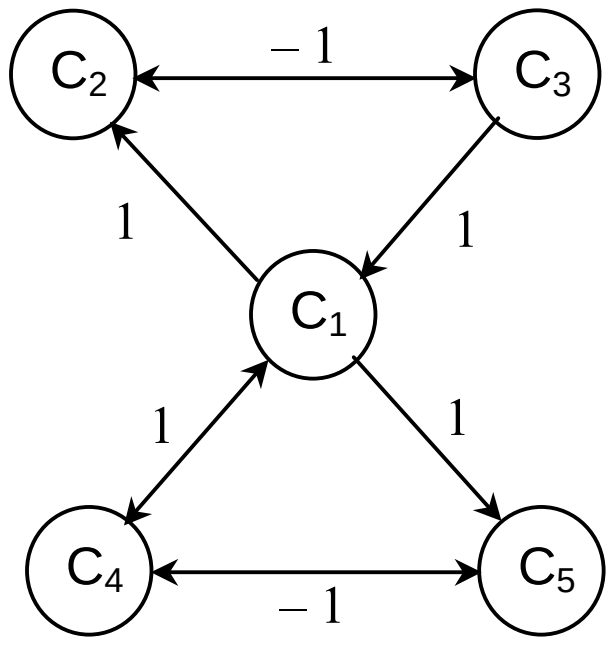

FIGURE: 3.2.1

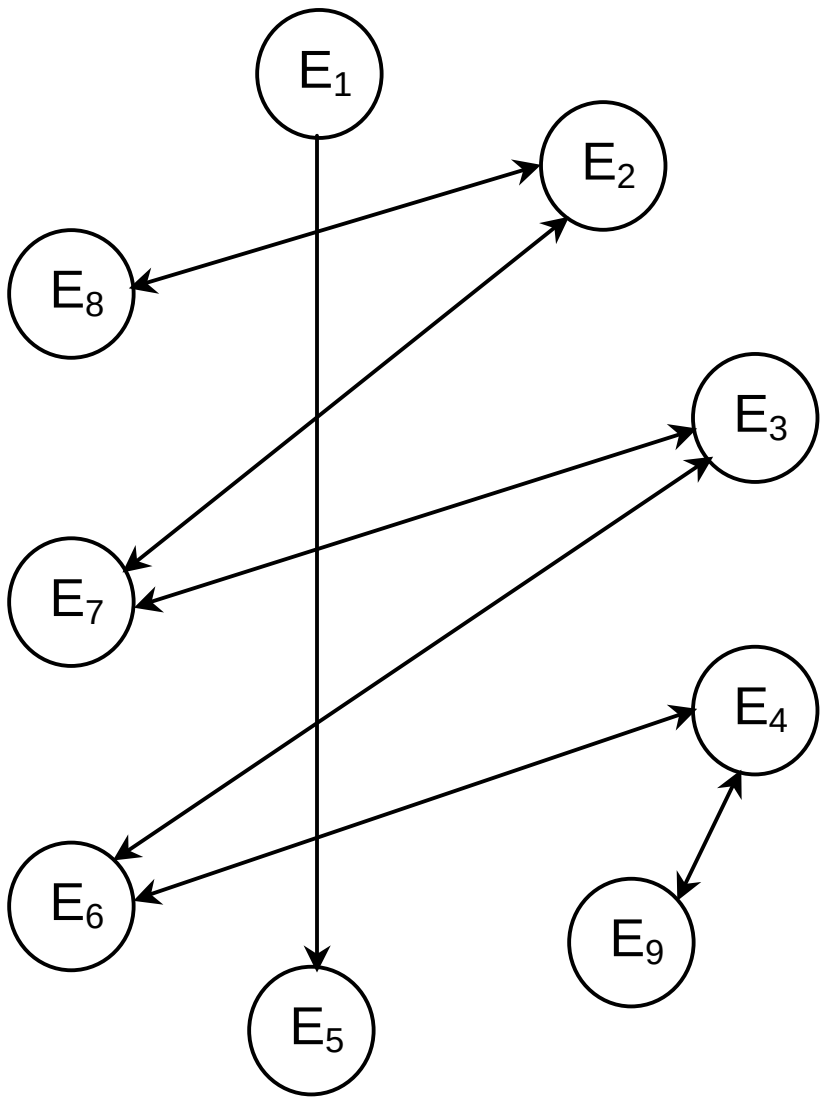

FIGURE: 3.2.2

We see the attributes $\{C_1, C_2, C_3, C_4, C_5\}$ are disjoint with the attributes $\{E_1, E_2, \ldots, E_9\}$.

This bigraph is a disjoint bigraph as we have no common attributes with the given set of concepts.
The related bimatrix or the connection bimatrix is a mixed square bimatrix B.

$$B = B_1 \cup B_2.$$

$$B = \begin{array}{c} \\ C_1 \\ C_2 \\ C_3 \\ C_4 \\ C_5 \end{array} \begin{array}{c} \begin{array}{ccccc} C_1 & C_2 & C_3 & C_4 & C_5 \end{array} \\ \begin{bmatrix} 0 & 1 & 0 & 1 & 1 \\ 0 & 0 & -1 & 0 & 0 \\ 1 & -1 & 0 & 0 & 0 \\ 1 & 0 & 0 & 0 & -1 \\ 0 & 0 & 0 & -1 & 0 \end{bmatrix} \end{array} \cup$$

$$\begin{array}{c} \\ E_1 \\ E_2 \\ E_3 \\ E_4 \\ E_5 \\ E_6 \\ E_7 \\ E_8 \\ E_9 \end{array} \begin{array}{c} \begin{array}{ccccccccc} E_1 & E_2 & E_3 & E_4 & E_5 & E_6 & E_7 & E_8 & E_9 \end{array} \\ \begin{bmatrix} 0 & 0 & 0 & 0 & 1 & 0 & 0 & 0 & 0 \\ 0 & 0 & 0 & 0 & 0 & 0 & 1 & 1 & 0 \\ 0 & 0 & 0 & 0 & 0 & 1 & 1 & 0 & 0 \\ 0 & 0 & 0 & 0 & 0 & 1 & 0 & 0 & 1 \\ 0 & 0 & 0 & 0 & 0 & 0 & 0 & 0 & 0 \\ 0 & 0 & 1 & 1 & 0 & 0 & 0 & 0 & 0 \\ 0 & 1 & 1 & 0 & 0 & 0 & 0 & 0 & 0 \\ 0 & 1 & 0 & 0 & 0 & 0 & 0 & 0 & 0 \\ 0 & 0 & 0 & 1 & 0 & 0 & 0 & 0 & 0 \end{bmatrix} \end{array}$$

Now these two bimatrices can be used to find the effect of any state vector. Here it is important to note that a state vector for any bimatrix which will only be a row bivector. If we wish to make only one system to work then in the initial row bivector we use one of the row vectors to be a zero vector. Now we explain how this bimodel functions.

Any initial state vector would be a pair of row bivector in this case any state bivector

$$X = \left\{\left(a_1^1, a_2^1, a_3^1, a_4^1, a_5^1\right) \cup \left(a_1^2, a_2^2, a_3^2, a_4^2, a_5^2, \ldots, a_9^2\right)\right\}$$

$$a_j^1 = \begin{cases} 0 \text{ if the } j^{th} \text{ state is off} \\ 1 \text{ if the } j^{th} \text{ state is on} \end{cases}$$

$$a_j^2 = \begin{cases} 0 \text{ if the } j^{th} \text{ state is off} \\ 1 \text{ if the } j^{th} \text{ state is on} \end{cases}.$$

This bivector X will be known as the instantaneous state bivector.

Now we study the effect of any state bivector on the dynamical bisystem $B = B_1 \cup B_2$.

Now we find the effect of the initial state bivector

$$\begin{aligned} X &= (1\,0\,0\,0\,0) \cup (0\,1\,0\,0\,0\,0\,0\,0\,0) \\ &= X_1 \cup X_2 \end{aligned}$$

on $B = B_1 \cup B_2$.

$$\begin{aligned} XB &= (X_1 \cup X_2)\, B \\ &= (X_1 \cup X_2)\, (B_1 \cup B_2) \\ &= X_1 B_1 \cup X_2 B_2 \\ &= [(0\ 1\ 0\ 1\ 1) \cup (0\ 0\ 0\ 0\ 0\ 0\ 1\ 1\ 0)]. \end{aligned}$$

Now updating and thresholding the bivector we get

$$\begin{aligned} XB &= [(1\ 1\ 0\ 1\ 1) \cup (0\ 1\ 0\ 0\ 0\ 0\ 1\ 1\ 0)] \\ &= Y = Y_1 \cup Y_2. \end{aligned}$$

$$\begin{aligned} YB &= (Y_1 \cup Y_2)\, (B) \\ &= (Y_1 \cup Y_2)\, (B_1 \cup B_2) \\ &= Y_1 B_1 \cup Y_2 B_2 \\ &= [(1\ 1\ {-1}\ 0\ 0)] \cup [(0\ 2\ 1\ 0\ 0\ 0\ 1\ 1\ 0)]; \end{aligned}$$

after updating and thresholding we get the resultant as

$$\begin{aligned} &= Z \\ &= [(1\ 1\ 0\ 0\ 0) \cup (0\ 1\ 1\ 0\ 0\ 0\ 1\ 1\ 0)] \\ &= Z_1 \cup Z_2. \end{aligned}$$

Consider

$$\begin{aligned} ZB &= (Z_1 \cup Z_2)\, B \\ &= (Z_1 \cup Z_2)\, (B_1 \cup B_2) \\ &= Z_1 B_1 \cup Z_2 B_2 \\ &= [(0\ 1\ {-1}\ 1\ 1)] \cup [(0\ 2\ 1\ 0\ 0\ 1\ 1\ 1\ 0)] \end{aligned}$$

after updating and thresholding we get the resultant bivector as S where

$$\begin{aligned} S &= [(1\ 1\ 0\ 1\ 1) \cup (0\ 1\ 1\ 0\ 0\ 1\ 1\ 1\ 0)] \\ &= S_1 \cup S_2. \end{aligned}$$

Now the effect of S on the dynamical system B is given by

$$BS = (1\ 1\ {-1}\ 0\ 0) \cup (0\ 2\ 2\ 1\ 0\ 2\ 2\ 1\ 1)$$

after updating and thresholding the state bivector we get the resultant bivector as

$$R = (1\ 1\ 0\ 0\ 0) \cup (0\ 1\ 1\ 1\ 0\ 1\ 1\ 1\ 1)$$
$$= R_1 \cup R_2\ .$$

Now we study the effect of R on the dynamical system B,

$$RB = (0\ 1\ {-1}\ 1\ 1) \cup (0\ 2\ 2\ 1\ 0\ 2\ 2\ 1\ 1)$$

after updating and thresholding resultant state bivector we get the resultant bivector; we get the resultant vector T as

$$T = (1\ 1\ 0\ 1\ 1) \cup (0\ 1\ 1\ 1\ 0\ 1\ 1\ 1\ 1)$$
$$= T_1 \cup T_2.$$

Now the effect of T on the bimatrix B is given by

$$TB = (T_1 \cup T_2)\ (B_1 \cup B_2)$$
$$= T_1\ B_1 \cup T_2\ B_2$$
$$= (1\ 1\ {-1}\ 0\ 0) \cup (0\ 2\ 2\ 2\ 0\ 2\ 2\ 1\ 1)$$

after thresholding and updating the resultant bivector we get the bivector as

$$A = [(1\ 1\ 0\ 0\ 0\ ) \cup (0\ 1\ 1\ 1\ 0\ 1\ 1\ 1\ 1)]$$
$$= A_1 \cup A_2$$

which is a hidden pattern. The hidden pattern is a limit cycle combined with the fixed point given by

$$\{(1\ 1\ 0\ 1\ 1) \cup (0\ 1\ 1\ 1\ 0\ 1\ 1\ 1\ 1)\} \text{ or}$$
$$\{(1\ 1\ 0\ 0\ 0) \cup (0\ 1\ 1\ 1\ 0\ 1\ 1\ 1\ 1)\}.$$

So in the system B when we consider a initial state bivector $(1\ 0\ 0\ 0\ 0\ ) \cup (0\ 1\ 0\ 0\ 0\ 0\ 0\ 0\ 0\ )$ i.e. Good business coupled with just profit to the employer we see it gives good investment with other states like neither profit nor loss to the employer is on, when he contemplates on good business.

One can expect to give only pay to employee and for good business with just profit, the employee works for more number of hours, it still gives only an average performance by few employee and poor performance by some employee.

Thus this state bivector (1 0 0 0 0) $\cup$ (0 1 0 0 0 0 0 0 0 ) gives a unique form of the hidden pattern i.e. one row vector in the birow vector is a limit cycle and the other happens to be a fixed point.

It is important to note here that at times, the hidden pattern is such that both the bivectors are fixed point or both of them are limit cycles. In case when both are limit cycle it is still interesting to note that the interpretation of the state bivectors vary from stage to stage, which is pertinent for in practical situation such types of solutions are possible in the real world problems that too when we use an unsupervised data. Thus we can get four possibilities as the hidden pattern which we may from now onwards call as bihidden pattern. It may be a fixed bipoint or limit bicycle or fixed point and limit cycle.

Now we consider the state bivector

$$Y = (0\ 0\ 1\ 0\ 0) \cup (0\ 0\ 0\ 0\ 1\ 0\ 0\ 0\ 0) = Y_1 \cup Y_2$$

i.e. only the state $C_3$ i.e. the customer satisfaction and $E_5$ – Best performance by the employees are in the on state all other state vectors in the bivector is in the off state. Now we analyse the effect of the state bivector Y on the dynamical system B.

$$\begin{aligned} YB &= (Y_1 \cup Y_2)\,(B_1 \cup B_2) \\ &= Y_1 B_1 \cup Y_2 B_2 \\ &= (1\ {-1}\ 0\ 0\ 0) \cup (0\ 0\ 0\ 0\ 1\ 0\ 0\ 0\ 0) \end{aligned}$$

(after thresholding and updating the resultant bivector) we get Z

$$Z = (1\ 0\ 1\ 0\ 0) \cup (0\ 0\ 0\ 0\ 1\ 0\ 0\ 0\ 0) = Z_1 \cup Z_2\ .$$

ZB = $Z_1 B_1 \cup Z_2 B_2$
= $(1\ 0\ 0\ 1\ 1) \cup (0\ 0\ 0\ 0\ 1\ 0\ 0\ 0\ 0)$ .

after updating we get the resultant state bivector

X = $(1\ 0\ 1\ 1\ 1) \cup (0\ 0\ 0\ 0\ 1\ 0\ 0\ 0\ 0)$

The effect of the bivector X on the dynamical system B gives

XB = $(2\ 0\ 0\ 0\ 0) \cup (0\ 0\ 0\ 0\ 1\ 0\ 0\ 0\ 0)$

After updating and thresholding we get the resultant bivector as

T = $T_1 \cup T_2$
= $(1\ 0\ 1\ 0\ 0) \cup (0\ 0\ 0\ 0\ 1\ 0\ 0\ 0\ 0)$

The effect the birow vector T on B gives

TB = $T_1 B_1 \cup T_2 B_2$
= $(1\ 0\ 0\ 1\ 1) \cup (0\ 0\ 0\ 0\ 1\ 0\ 0\ 0\ 0)$

after thresholding and updating the resultant bivector is

U = $(1\ 0\ 1\ 1\ 1) \cup (0\ 0\ 0\ 0\ 1\ 0\ 0\ 0\ 0)$
= $U_1 \cup U_2$ .

Thus the bihidden pattern is a fixed point and a limit cycle, we see when the concept / attribute, Best Performance by the employee is in the on state the system remains static on other attribute ever becomes on. But on the other hand when the customer satisfaction is in the on state we see the bihidden pattern is a limit cycle which at one point makes all the states to be or $C_2$ alone is in the off state. Thus it fluctuates from (1 0 1 0 0) and (1 0 1 1 1).

It is still important to note that when we analyze a problem with a FCBMs both the sets of attributes need not always be disjoint. It can also be a over lapping set. For even to analyze the same problem one can use FCBMs.

The application of FCBMs to analyze the factors promoting business. After discussion with several experts we have taken up the following attributes. We have conditioned them to give only 5 attributes related with the models. We have simultaneously taken the opinion of two experts. The attributes given by the first expert is

| | | |
|---|---|---|
| $E_1$ | – | Good business |
| $E_2$ | – | Good investment |
| $E_3$ | – | Customer satisfaction |
| $E_4$ | – | Establishment |
| $E_5$ | – | Good Marketing Strategies. |

The attributes given by the second expert is as follows

| | | |
|---|---|---|
| $E'_1$ | – | Good business |
| $E'_2$ | – | Appropriate locality |
| $E'_3$ | – | Selling quality products |
| $E'_4$ | – | Updation of techniques |
| $E'_5$ | – | Knowledge about the policies of the government. |

Now the related directed bigraph is as follows.

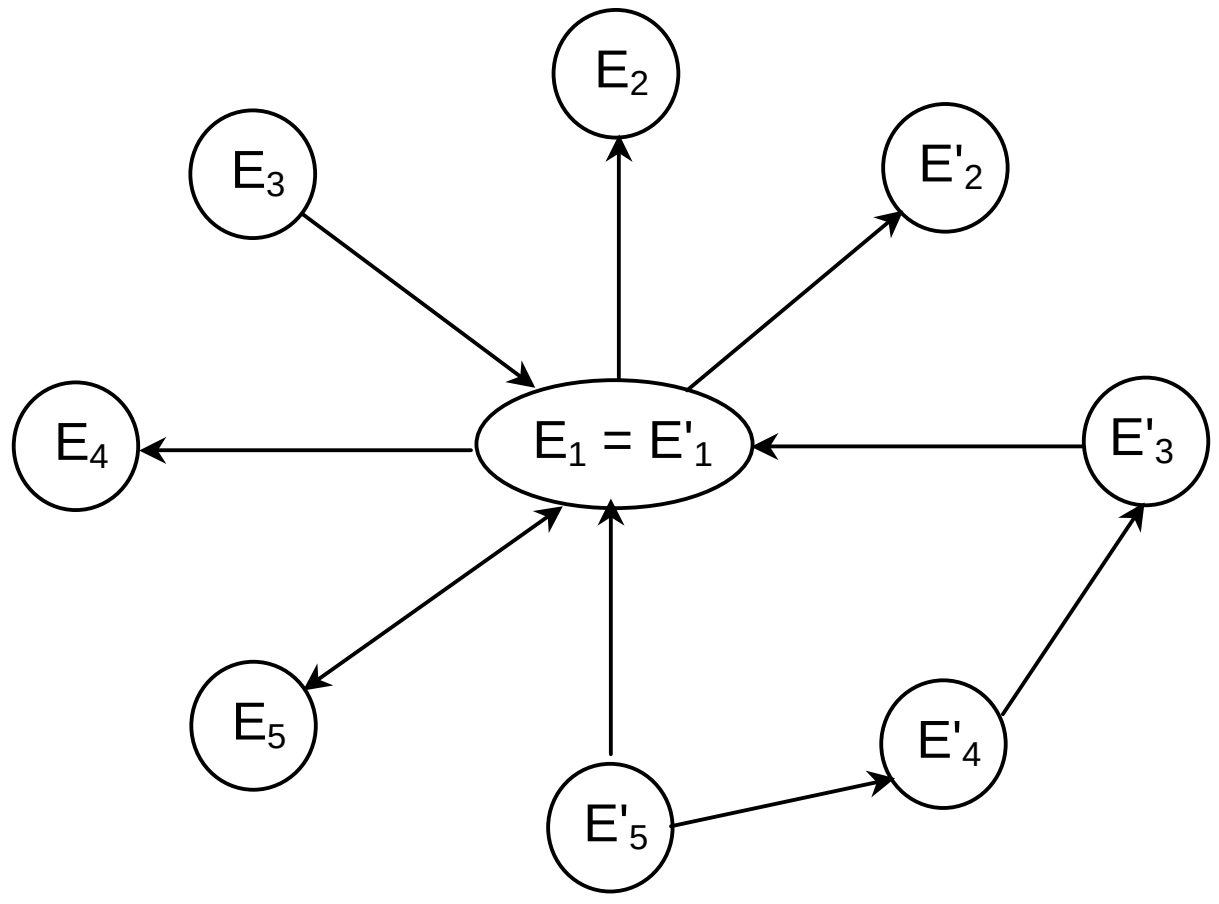

FIGURE: 3.2.3

Clearly this is a connected directed bigraph.

Now the associated connection bimatrix is given by B = $B_1 \cup B_2$

$$\begin{array}{c} \\ E_1 \\ E_2 \\ E_3 \\ E_4 \\ E_5 \end{array}\begin{array}{c} \begin{array}{ccccc} E_1 & E_2 & E_3 & E_4 & E_5 \end{array} \\ \begin{bmatrix} 0 & 1 & 0 & 1 & 1 \\ 0 & 0 & 0 & 0 & 0 \\ 1 & 0 & 0 & 0 & 0 \\ 0 & 0 & 0 & 0 & 0 \\ 1 & 0 & 0 & 0 & 0 \end{bmatrix} \end{array} \cup \begin{array}{c} \\ E_1' \\ E_2' \\ E_3' \\ E_4' \\ E_5' \end{array}\begin{array}{c} \begin{array}{ccccc} E_1' & E_2' & E_3' & E_4' & E_5' \end{array} \\ \begin{bmatrix} 0 & 1 & 0 & 0 & 0 \\ 0 & 0 & 0 & 0 & 0 \\ 1 & 0 & 0 & 0 & 0 \\ 0 & 0 & 1 & 0 & 0 \\ 1 & 0 & 0 & 1 & 0 \end{bmatrix} \end{array}$$

Suppose one is interested in studying the effect of the row bivector

$$\begin{aligned} X &= X_1 \cup X_2 \\ &= (1\ 0\ 0\ 0\ 0) \cup (1\ 0\ 0\ 0\ 0). \end{aligned}$$

The effect of X on B is given by

$$\begin{aligned} XB &= X (B_1 \cup B_2) \\ &= X_1 B_1 \cup X_2 B_2 \\ &= (0\ 1\ 0\ 1\ 1) \cup (0\ 1\ 0\ 0\ 0) \end{aligned}$$

after updating and thresholding we get

$$\begin{aligned} Y &= (1\ 1\ 0\ 1\ 1) \cup (1\ 1\ 0\ 0\ 0) \\ &= Y_1 \cup Y_2. \end{aligned}$$

The effect of Y on B is given by

$$\begin{aligned} YB &= Y_1 B_1 \cup Y_2 B_2. \\ &= (1\ 1\ 0\ 1\ 1) \cup (0\ 1\ 0\ 0\ 0). \end{aligned}$$

After updating and thresholding we get the resultant vector as

$$Z = (1\ 1\ 0\ 1\ 1) \cup (1\ 1\ 0\ 0\ 0);$$

which is a fixed bipoint. Good business according to the first expert has no impact on Customer satisfaction but it has influence on good investment, establishment and good Marketing strategies. But at the same time we see according to the second expert good business has to do only with the appropriate locality and nothing to do with selling quality products or updation or knowledge about the polices of the government. As they are the opinion given by the expert we have no right to change or modify these effects.

Next we go for the opinion of the third expert. He gives the attributes as

| | | |
|---|---|---|
| $E''_1$ | – | Good business |
| $E''_2$ | – | Geographical situation |
| $E''_3$ | – | Rendering good service |
| $E''_4$ | – | Previous experience of the owner |
| $E''_5$ | – | Demand and supply. |

Now using the first and third experts opinions we have the following directed bigraph.

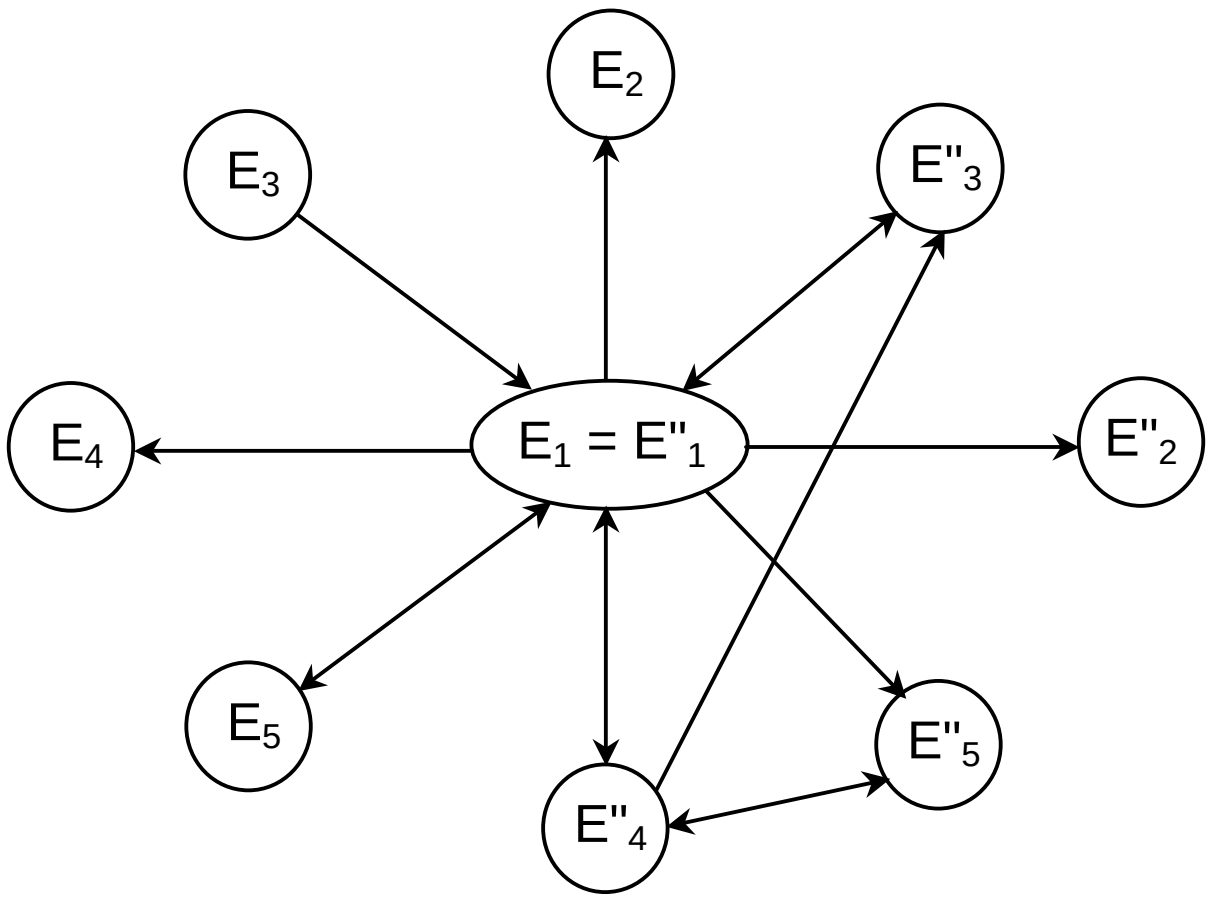

FIGURE: 3.2.4

The related connection bimatrix.

$$B = B_1 \cup B_2.$$

$$\begin{array}{c} \\ E_1 \\ E_2 \\ E_3 \\ E_4 \\ E_5 \end{array}\begin{array}{c} \begin{array}{ccccc} E_1 & E_2 & E_3 & E_4 & E_5 \end{array} \\ \begin{bmatrix} 0 & 1 & 0 & 1 & 1 \\ 0 & 0 & 0 & 0 & 0 \\ 1 & 0 & 0 & 0 & 0 \\ 0 & 0 & 0 & 0 & 0 \\ 1 & 0 & 0 & 0 & 0 \end{bmatrix} \end{array} \cup \begin{array}{c} \\ E''_1 \\ E''_2 \\ E''_3 \\ E''_4 \\ E''_5 \end{array}\begin{array}{c} \begin{array}{ccccc} E''_1 & E''_2 & E''_3 & E''_4 & E''_5 \end{array} \\ \begin{bmatrix} 0 & 1 & 1 & 1 & 1 \\ 0 & 0 & 0 & 0 & 0 \\ 1 & 0 & 0 & 0 & 0 \\ 1 & 0 & 0 & 0 & 1 \\ 0 & 0 & 0 & 1 & 0 \end{bmatrix} \end{array}$$

Now we study the effect of the state bivector on the dynamical bisystem.

Let

$$X = (0\ 1\ 0\ 0\ 1) \cup (1\ 0\ 1\ 0\ 0)$$

i.e. the first expert considers good investment and good marketing strategies in the on state and good business and rendering good service to be in the on states i.e. the bivector given by them has $E_2, E_5, E''_1$ and $E''_5$ to be in the on state and all other nodes are in the off state. We study the effect of

$$\begin{aligned} X &= X_1 \cup X_2 \text{ on } B. \\ XB &= (X_1 \cup X_2)(B_1 \cup B_2) \\ &= X_1 B_1 \cup X_2 B_2. \\ &= (1\ 0\ 0\ 0\ 0) \cup (1\ 1\ 1\ 1\ 1). \end{aligned}$$

after updating we get the resultant

$$\begin{aligned} Y &= (1\ 1\ 0\ 0\ 1) \cup (1\ 1\ 1\ 1\ 1) \\ &= Y_1 \cup Y_2. \end{aligned}$$

the effect of Y on B gives

$$\begin{aligned} YB &= Y_1 B_1 \cup Y_2 B \\ &= (1\ 1\ 0\ 1\ 1) \cup (2\ 1\ 1\ 2\ 2) \end{aligned}$$

after thresholding and updating we get

$$Z = (1\ 1\ 0\ 1\ 1) \cup (1\ 1\ 1\ 1\ 1)$$
$$= Y.$$

Thus the bihidden pattern is a fixed bipoint.
Thus good investment and good marketing strategies has no impact on other nodes; where as the on state of good business and rendering good service makes on all other states.

***Example 3.2.2:*** Now we study using the same model the problems faced by primary school children in relation with parents. We after using a linguistic questionnaire obtained the following unsupervised data. In the unsupervised data 10 important attributes were considered. Here the model is used in an entirely different way for the same set of attributes are used but only the experts differ. We get the directed bigraph which is different and disjoint as all the nodes are same but opinion is given by two different persons. The concepts associated with the parents in relation with their children who are just in their primary level are as follows:

| | | |
|---|---|---|
| $P_1$ | – | Well educated children |
| $P_2$ | – | Educated parents |
| $P_3$ | – | Uneducated parents |
| $P_4$ | – | Moderate fees |
| $P_5$ | – | High fees |
| $P_6$ | – | Rich parents |
| $P_7$ | – | Middle class parents |
| $P_8$ | – | Poor parents |
| $P_9$ | – | Language problem |
| $P_{10}$ | – | Dropout from school. |

The directed bigraph given by two experts is as follows:

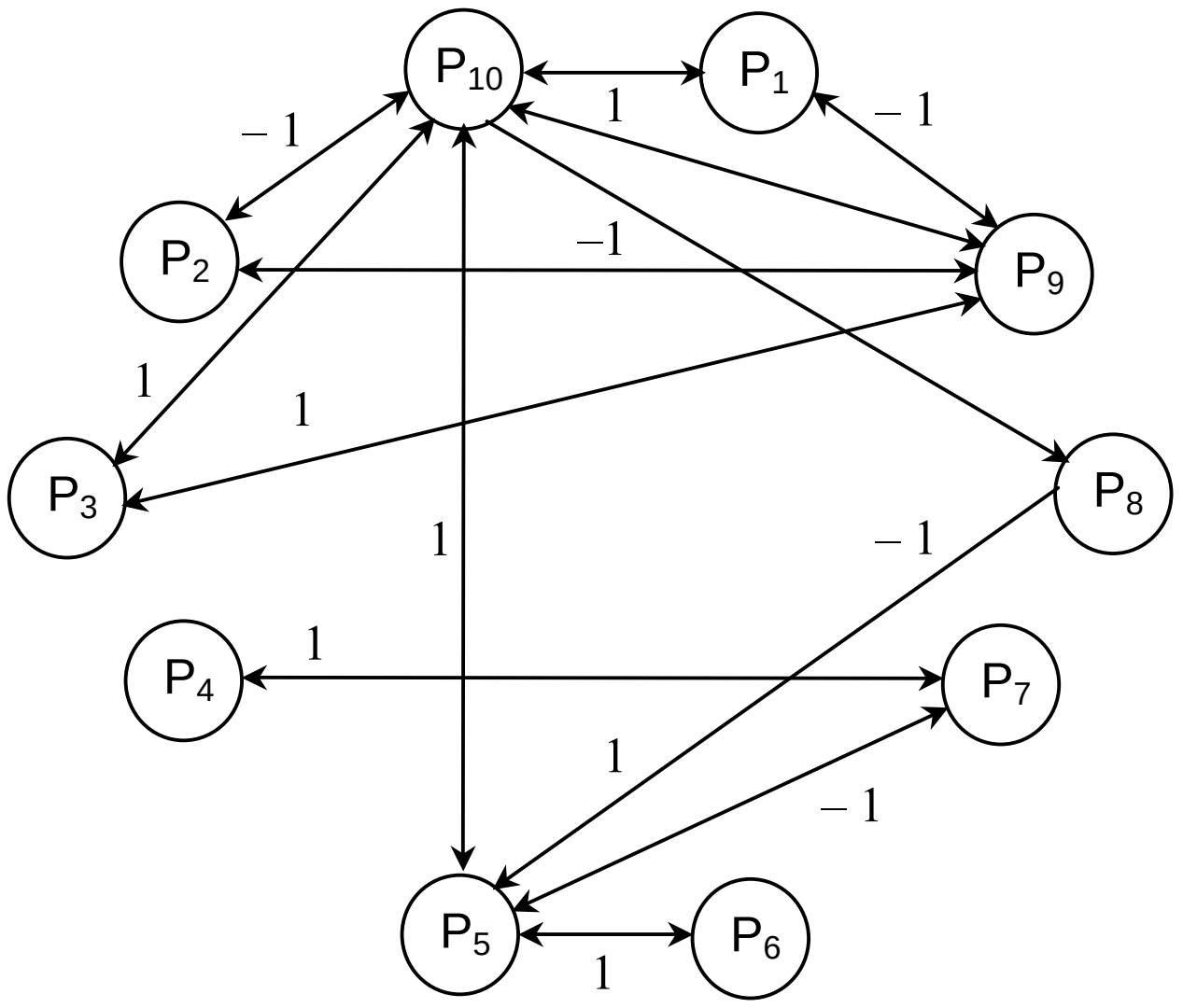

FIGURE: 3.2.5

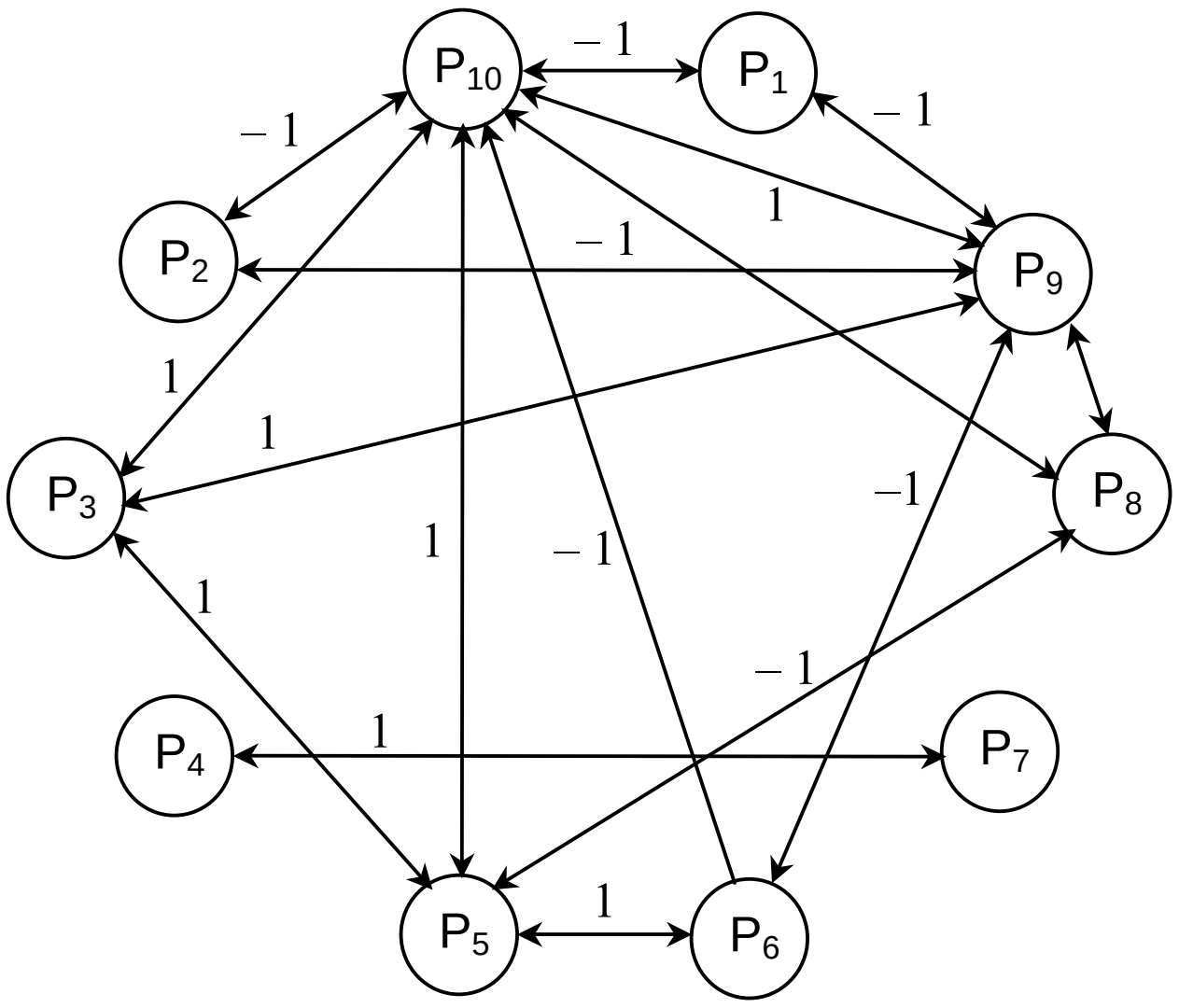

FIGURE: 3.2.6

The directed bigraph given by the experts. The related connection bimatrix of the directed bigraph is a 10 × 10 square bimatrix given by

$$B = B_1 \cup B_2$$

$$\begin{array}{c} \\ P_1 \\ P_2 \\ P_3 \\ P_4 \\ P_5 \\ P_6 \\ P_7 \\ P_8 \\ P_9 \\ P_{10} \end{array}
\begin{array}{c} \begin{array}{cccccccccc} P_1 & P_2 & P_3 & P_4 & P_5 & P_6 & P_7 & P_8 & P_9 & P_{10} \end{array} \\
\left[\begin{array}{cccccccccc}
0 & 0 & 0 & 0 & 0 & 0 & 0 & 0 & -1 & 1 \\
0 & 0 & 0 & 0 & 0 & 0 & 0 & 0 & -1 & -1 \\
0 & 0 & 0 & 0 & 0 & 0 & 0 & 0 & 1 & 1 \\
0 & 0 & 0 & 0 & 0 & 0 & 1 & 0 & 0 & 0 \\
0 & 0 & 0 & 0 & 0 & 1 & -1 & 0 & 0 & 1 \\
0 & 0 & 0 & 0 & 1 & 0 & 0 & 0 & 0 & 0 \\
0 & 0 & 0 & 1 & -1 & 0 & 0 & 0 & 0 & 0 \\
0 & 0 & 0 & 0 & -1 & 0 & 0 & 0 & 0 & 0 \\
-1 & -1 & 1 & 0 & 0 & 0 & 0 & 0 & 0 & 1 \\
-1 & -1 & 1 & 0 & 1 & 0 & 0 & 1 & 1 & 0
\end{array}\right] \end{array}$$

$$\cup \quad \begin{array}{c} \\ P_1 \\ P_2 \\ P_3 \\ P_4 \\ P_5 \\ P_6 \\ P_7 \\ P_8 \\ P_9 \\ P_{10} \end{array}
\begin{array}{c} \begin{array}{cccccccccc} P_1 & P_2 & P_3 & P_4 & P_5 & P_6 & P_7 & P_8 & P_9 & P_{10} \end{array} \\
\left[\begin{array}{cccccccccc}
0 & 0 & 0 & 0 & 0 & 0 & 0 & 0 & -1 & -1 \\
0 & 0 & 0 & 0 & 0 & 0 & 0 & 0 & -1 & -1 \\
0 & 0 & 0 & 0 & 1 & 0 & 0 & 0 & 1 & 1 \\
0 & 0 & 0 & 0 & 0 & 0 & 1 & 0 & 0 & 0 \\
0 & 0 & 1 & 0 & 0 & 1 & 0 & -1 & 0 & 1 \\
0 & 0 & 0 & 0 & 1 & 0 & 0 & 0 & -1 & -1 \\
0 & 0 & 0 & 1 & 0 & 0 & 0 & 0 & 0 & 0 \\
0 & 0 & 0 & 0 & -1 & 0 & 0 & 0 & 1 & 1 \\
-1 & -1 & 1 & 0 & 0 & -1 & 0 & 1 & 0 & 1 \\
-1 & -1 & 1 & 0 & 1 & 0 & 0 & 1 & 1 & 0
\end{array}\right] \end{array}$$

Clearly B is bimatrix which have some rows and columns in common but they are at the same time distinct .

Now suppose we input the vector

$$\begin{aligned} A &= (0\ 0\ 1\ 0\ 0\ 0\ 0\ 0\ 0\ 0) \cup (0\ 0\ 0\ 0\ 0\ 0\ 0\ 1\ 0\ 0) \\ &= A_1 \cup A_2. \end{aligned}$$

That is in the row bivector only the nodes uneducated parents and poor parents are in the on state and all other nodes are in the off state. Now we study the effect of the birow vector on the model B.

$$\begin{aligned} AB &= (A_1 \cup A_2)\ (B_1 \cup B_2) \\ &= A_1 B_1 \cup A_2 B_2 \\ &= (0\ 0\ 0\ 0\ 0\ 0\ 0\ 0\ 1\ 1) \cup (0\ 0\ 0\ 0\ {-1}\ 0\ 0\ 0\ 1\ 1) \end{aligned}$$

After thresholding and updating we get the following resultant bivector say $C = C_1 \cup C_2$ where

$$C = (\ 0\ 0\ 1\ 0\ 0\ 0\ 0\ 0\ 1\ 1) \cup (0\ 0\ 0\ 0\ 0\ 0\ 0\ 1\ 1\ 1)$$

Now we study the effect of the row bivector on the dynamical system B

$$\begin{aligned} CB &= (C_1 \cup C_2)\ (B_1 \cup B_2) \\ &= C_1 B_1 \cup C_2 B_2. \\ &= ({-2}\ {-2}\ 2\ 0\ 1\ 0\ 0\ 1\ 2\ 2) \\ &\quad \cup ({-2}\ {-2}\ 2\ 0\ 0\ {-1}\ 0\ 2\ 2\ 2) \end{aligned}$$

after thresholding and updating we get the resultant bivector to be

$$\begin{aligned} D &= (0\ 0\ 1\ 0\ 1\ 0\ 0\ 1\ 1\ 1) \cup (0\ 0\ 1\ 0\ 1\ 0\ 0\ 1\ 1\ 1) \\ &= D_1 \cup D_2. \end{aligned}$$

The effect D on the dynamical system B is given by

$$\begin{aligned} DB &= (D_1 \cup D_2)\ (B_1 \cup B_2) \\ &= D_1 B_2 \cup D_2 B_2. \\ &= ({-2}\ {-2}\ 2\ 0\ 0\ 1\ {-1}\ 1\ 2\ 4) \\ &\quad \cup (0\ 0\ 1\ 0\ 1\ 0\ 0\ 1\ 1\ 1). \end{aligned}$$

After thresholding and updating we get the resultant bivector $G = G_1 \cup G_2$

$$= (0\ 0\ 1\ 0\ 0\ 1\ 0\ 1\ 1\ 1) \cup (0\ 0\ 1\ 0\ 1\ 0\ 0\ 1\ 1\ 1)$$

$$\begin{aligned} GB &= (G_1 \cup G_2)\ (B_1 \cup B_2) \\ &= G_1\ B_1 \cup G_2\ B_2. \\ &= (-2\ -2\ 2\ 0\ 1\ 0\ 0\ 1\ 2\ 3) \cup (-2\ -2\ 3\ 0\ 1\ 0\ 0\ 1\ 3\ 4). \end{aligned}$$

After thresholding and updating we get the resultant row bivector as

$$D = (0\ 0\ 1\ 0\ 1\ 0\ 0\ 1\ 1\ 1) \cup (0\ 0\ 1\ 0\ 1\ 0\ 0\ 1\ 1\ 1).$$

Now as we proceed on the bihidden pattern of the dynamical system is a limit cycle for the first expert and a fixed point in case of the second expert given by the following

$$\{(0\ 0\ 1\ 0\ 0\ 1\ 0\ 1\ 1\ 1), (0\ 0\ 1\ 0\ 1\ 0\ 0\ 1\ 1\ 1)\} \text{ or I}$$
$$\{(0\ 0\ 1\ 0\ 1\ 0\ 0\ 1\ 1\ 1), (0\ 0\ 1\ 0\ 1\ 0\ 0\ 1\ 1\ 1)\} \text{ or II.}$$

Thus the uneducated parents and poor parents at one time give the same resultant II i.e. uneducated parents invariable mean poor parents and vice versa. In case of both uneducated parents and poor parents we see the children suffer with language problem and they are invariably the dropout in primary school one of the causes in both cases may be due to high fee structure in primary school.

We read the resultant bivector given by *I* we see that rich parents come to on state so the limit cycle *I* is meaningless and hence it is dropped.

***Example 3.2.3:*** This model is mainly given to show that the simultaneous process of opinion can be got using bimatrix given by the two experts with same sets of attributes yielding to stage by stage comparison such study is an impossibility, without a bimatrix.

Now we illustrate one more new model which is little different from the models discussed earlier.

Suppose we are interested in studying the relation between the health hazards the agricultural labourer suffers due to chemical pollution. Let us taken the 11 attributes ($P_1$, $P_2$,…, $P_{11}$) after having interviews and discussion with more than 50 agricultural labourers and their families. These attributes were chosen by the experts.

| | | |
|---|---|---|
| $P_1$ | – | Swollen limbs |
| $P_2$ | – | Ulcer / skin ailments in legs and hands |
| $P_3$ | – | Manuring the fields with chemical fertilizers |
| $P_4$ | – | Vomiting |
| $P_5$ | – | Mouth and stomach ulcer |
| $P_6$ | – | Pollution by drinking water |
| $P_7$ | – | Indigestion |
| $P_8$ | – | Loss of appetite |
| $P_9$ | – | Headache |
| $P_{10}$ | – | Spraying of pesticides |
| $P_{11}$ | – | Blurred vision. |

Suppose two expects choose the attributes in the following way. The first experts lakes the attributes {$P_3$….$P_9$} and the second expert choose the attributes {$P_9$ $P_{10}$ $P_{11}$ $P_1$ $P_2$ $P_3$} Let their views be given as a directed bigraph.

The directed bigraph given by the experts is given in the figure.

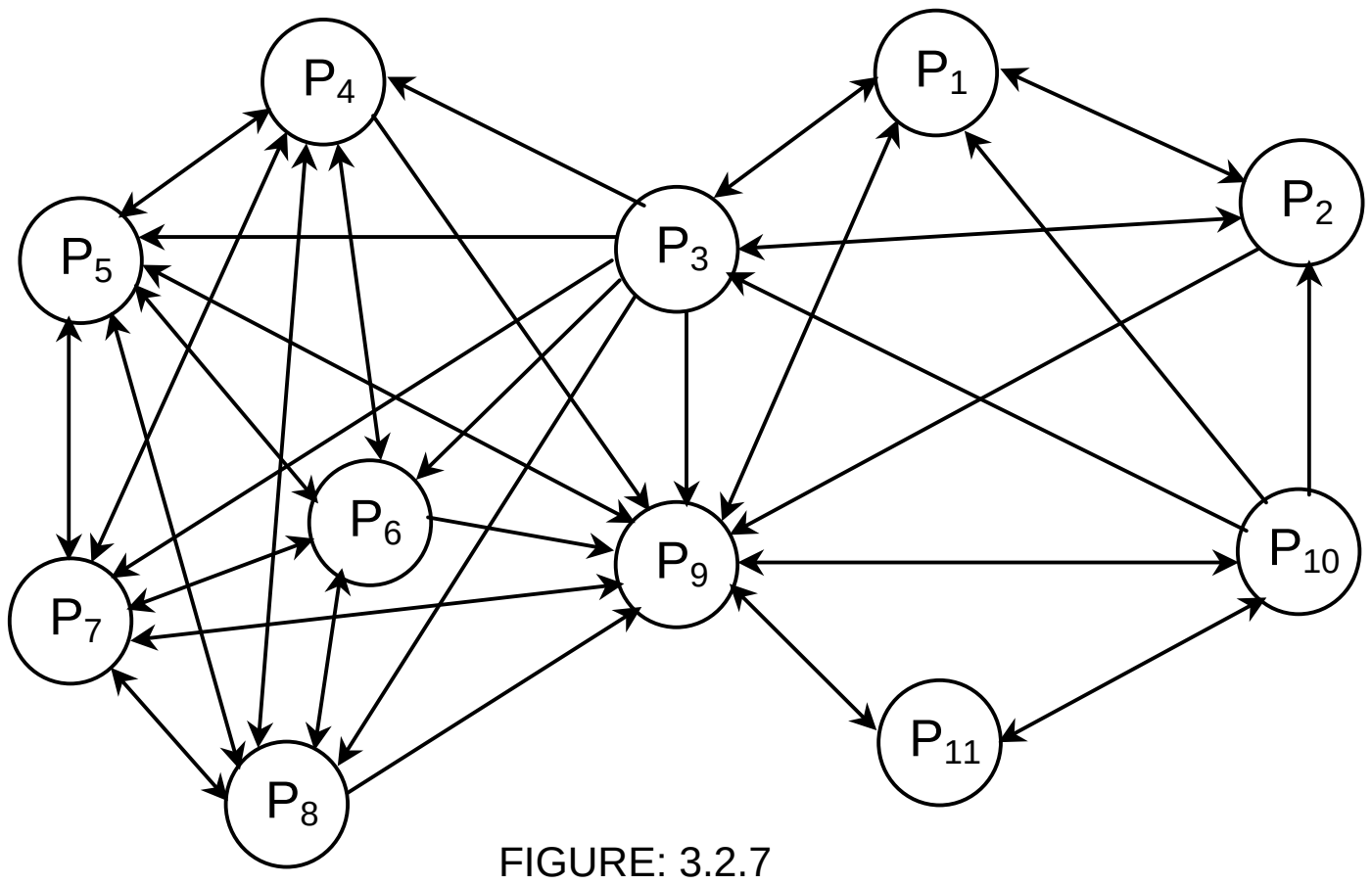

FIGURE: 3.2.7

The directed bigraph is connected and is glued by an edge.

The related connection matrix is given by $B = B_1 \cup B_2$

$$
\begin{array}{c} \\ P_3 \\ P_4 \\ P_5 \\ P_6 \\ P_7 \\ P_8 \\ P_9 \end{array}
\begin{array}{c} \begin{array}{ccccccc} P_3 & P_4 & P_5 & P_6 & P_7 & P_8 & P_9 \end{array} \\
\left[\begin{array}{ccccccc}
0 & 1 & 1 & 1 & 1 & 1 & 1 \\
0 & 0 & 1 & -1 & 1 & 1 & 1 \\
0 & 1 & 0 & -1 & 1 & 1 & 1 \\
0 & 1 & 1 & 0 & 1 & 1 & 1 \\
0 & 1 & 1 & 1 & 0 & 1 & 1 \\
0 & 1 & 1 & -1 & 1 & 0 & 1 \\
0 & 0 & 1 & 0 & 1 & 0 & 0
\end{array}\right] \end{array}
\cup
\begin{array}{c} \\ P_9 \\ P_{10} \\ P_{11} \\ P_1 \\ P_2 \\ P_3 \end{array}
\begin{array}{c} \begin{array}{cccccc} P_9 & P_{10} & P_{11} & P_1 & P_2 & P_3 \end{array} \\
\left[\begin{array}{cccccc}
0 & 1 & 1 & 0 & 0 & 0 \\
1 & 0 & 1 & 1 & 1 & -1 \\
1 & 1 & 0 & 0 & 0 & 0 \\
1 & 0 & 0 & 0 & 1 & -1 \\
1 & 0 & 0 & 1 & 0 & -1 \\
1 & 0 & 0 & 1 & 1 & 0
\end{array}\right] \end{array}
$$

Now the matrices are distinct so B is a mixed square bimatrix. Further none of the columns or rows can be identical as the bimatrix is mixed.

Now using an experts opinion we find the effect of a state bivector on the dynamical system.

Suppose

$$
\begin{aligned}
X &= (1\ 0\ 0\ 0\ 0\ 0\ 0) \cup (0\ 0\ 0\ 0\ 0\ 1) \\
&= X_1 \cup X_2.
\end{aligned}
$$

i.e. the attribute $P_3$ alone is in the on state in both the state vectors and all other nodes are in the off state. To find the effect of X on B.

$$
\begin{aligned}
XB &= (X_1 \cup X_2)\ (B_1 \cup B_2) \\
&= X_1 B_1 \cup X_2 B_2 \\
&= (0\ 1\ 1\ 1\ 1\ 1\ 1) \cup (1\ 0\ 0\ 1\ 1\ 0)
\end{aligned}
$$

after updating we get the resultant state bivector as

$$
\begin{aligned}
Y &= (1\ 1\ 1\ 1\ 1\ 1\ 1) \cup (1\ 0\ 0\ 1\ 1\ 1) \\
&= Y_1 \cup Y_2.
\end{aligned}
$$

Now we study the effect of $Y_B$ on the dynamical system

$$\begin{aligned} YB &= (Y_1 \cup Y_2)\,(B_1 \cup B_2) \\ &= Y_1 B_1 \cup Y_2 B_2 \\ &= (0\ 5\ 6\ {-1}\ 6\ 5\ 6) \cup (3\ 1\ 1\ 2\ 2\ {-2}). \end{aligned}$$

By thresholding and updating the resulting vector we get

$$\begin{aligned} Z &= (1\ 1\ 1\ 0\ 1\ 1\ 1) \cup (1\ 1\ 1\ 1\ 1\ 1) \\ &= Z_1 \cup Z_2. \end{aligned}$$

Now we study the effect of Z on B

$$\begin{aligned} ZB &= (Z_1 \cup Z_2)\,(B_1 \cup B_2) \\ &= Z_1 B_1 \cup Z_2 B_2 \\ &= (0\ 4\ 5\ {-2}\ 5\ 4\ 5) \cup (5\ 2\ 2\ 3\ 3\ {-3}). \end{aligned}$$

By thresholding and updating the resultant vector we get

$$\begin{aligned} U &= (1\ 1\ 1\ 0\ 1\ 1\ 1) \cup (1\ 1\ 1\ 1\ 1\ 1) \\ &= U_1 \cup U_2 \\ &= Z. \end{aligned}$$

Thus we arrive at a fixed bipoint as the bihidden pattern. Thus manuring the field with of chemical fertilizers makes all the coordinates on except $P_4$ which implies they may not suffer the symptom of vomiting because of manuring the field with chemical fertilizers.

Now we study the effect of the attribute Headache $P_9$ on the system.
Thus let

$$\begin{aligned} V &= (0\ 0\ 0\ 0\ 0\ 0\ 1) \cup (1\ 0\ 0\ 0\ 0\ 0) \\ &= V_1 \cup V_2. \end{aligned}$$

be the state vector. The effect of V on the dynamical system B is given by

$$VB = (0\ 0\ 1\ 0\ 1\ 0\ 0) \cup (0\ 1\ 1\ 0\ 0\ 0).$$

After updating the resultant vector we get

$$\begin{aligned} S &= (0\ 0\ 1\ 0\ 1\ 0\ 1) \cup (1\ 1\ 1\ 0\ 0\ 0) \\ &= S_1 \cup S_2. \end{aligned}$$

The effect of S on the dynamical system B is given by

$$\begin{aligned} SB &= (S_1 \cup S_2)\,(B_1 \cup B_2) \\ &= S_1 B_1 \cup S_2 B_2. \\ &= (0\ 2\ 2\ 0\ 2\ 2\ 2) \cup (2\ 2\ 2\ 1\ 1\ -1). \end{aligned}$$

After thresholding and updating we get the resultant as

$$\begin{aligned} R &= (0\ 1\ 1\ 0\ 1\ 1\ 1) \cup (1\ 1\ 1\ 1\ 1\ 0) \\ &= R_1 \cup R_2. \end{aligned}$$

Now we study the effect of R on the dynamical system B is given by

$$\begin{aligned} RB &= (R_1 \cup R_2)\,(B_1 \cup B_2) \\ &= R_1 B_1 \cup R_2 B_2 \\ &= (0\ 3\ 4\ -2\ 4\ 3\ 4) \cup (4\ 2\ 2\ 2\ 2\ -0). \end{aligned}$$

After thresholding and updating we get the resultant as

$$\begin{aligned} P &= (0\ 1\ 1\ 0\ 1\ 1\ 1) \cup (1\ 1\ 1\ 1\ 1\ 0) \\ &= P_1 \cup P_2 \\ &= R. \end{aligned}$$

Thus the bihidden pattern of the state bivector V is a fixed bipoint. Thus when suffers from headache i.e. when only the note $P_9$ is in the on state we see all nodes become on except $P_3$ and $P_6$ in the first state vector and $P_3$ is off or in the zero state. In the second state vector also as the second expert has not chosen to include the node $P_6$ in his analysis but has chosen both $P_9$ and $P_3$ and even in his opinion $P_3$ continues to be in the off state. Thus in view of both experts $P_3$ remains in the off state.

We have seen three conditions in the study of FCBMs

1. When the directed bigraph is disjoint.
2. When the directed bigraph is connected by a vertex i.e. bigraphs glued by a vertex.
3. When the directed bigraph is connected by an edge i.e. the bigraphs are glued by an edge.

Now we proceed on to study directed bigraphs which have a non trivial subgraph i.e., the subgraph is not a edge or a vertex but different has a common subgraph.

***Example 3.2.4:*** Now we study the problems faced by the agriculture labourers in the context of pollution and health hazards faced by them due to spray of pesticides insecticides and manuring the plants by chemical fertilizers. At the first stage we take arbitrary attributes say $\{a_1, a_2,\ldots, a_{10}\}$ associated with the coolies where $a_1, \ldots, a_{10}$ are defined below.

| | | |
|---|---|---|
| $a_1$ | – | Loss of appetite |
| $a_2$ | – | Headache |
| $a_3$ | – | Spraying of pesticides |
| $a_4$ | – | Indigestion |
| $a_5$ | – | Giddiness / fainting |
| $a_6$ | – | Mouth and stomach ulcer |
| $a_7$ | – | Breathlessness |
| $a_8$ | – | Skin ailments just after spray |
| $a_9$ | – | Diarrhea |
| $a_{10}$ | – | Consuming nearby vegetables / greens just after spray of pesticides. |

The directed bigraph associated with the two experts with related attributes $\{a_1, a_2,\ldots, a_7\}$ and $\{a_4, a_5, a_6, a_7, a_8, a_9, a_{10}\}$ respectively.

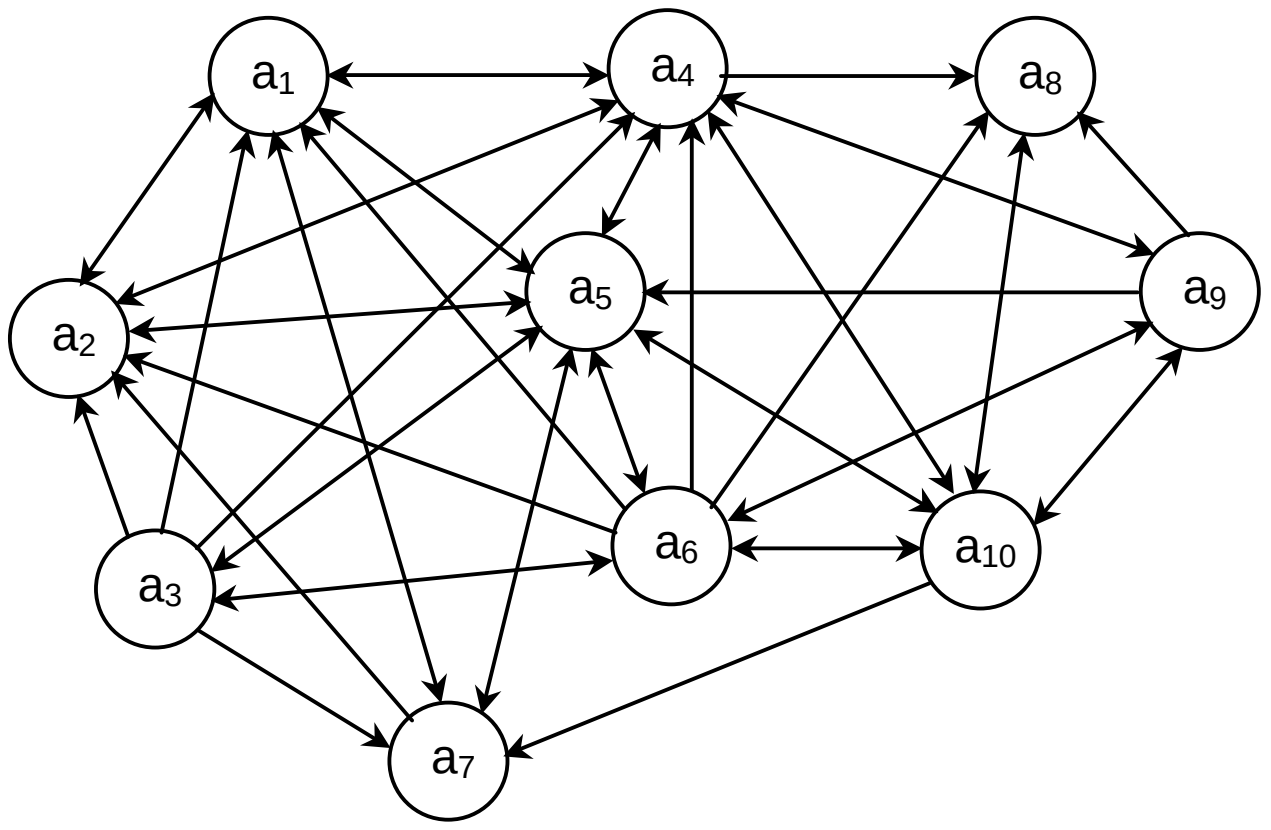

FIGURE: 3.2.7

Now we give the associated bimatrix $B = B_1 \cup B_2$ using the bigraph

$$
B_1 \cup B_2 =
\begin{array}{c}
 \\ a_1 \\ a_2 \\ a_3 \\ a_4 \\ a_5 \\ a_6 \\ a_7
\end{array}
\begin{array}{c}
\begin{array}{ccccccc} a_1 & a_2 & a_3 & a_4 & a_5 & a_6 & a_7 \end{array} \\
\begin{bmatrix}
0 & 1 & 0 & 1 & 1 & 0 & 1 \\
1 & 0 & 0 & 1 & 1 & 0 & 0 \\
1 & 1 & 0 & 1 & 1 & 1 & 1 \\
1 & 1 & 0 & 0 & 1 & 0 & 0 \\
1 & 1 & -1 & 1 & 0 & 1 & 1 \\
1 & 1 & -1 & 1 & 1 & 0 & 0 \\
1 & 1 & 0 & 0 & 1 & 0 & 0
\end{bmatrix}
\end{array}
\cup
\begin{array}{c}
 \\ a_4 \\ a_5 \\ a_6 \\ a_7 \\ a_8 \\ a_9 \\ a_{10}
\end{array}
\begin{array}{c}
\begin{array}{ccccccc} a_4 & a_5 & a_6 & a_7 & a_8 & a_9 & a_{10} \end{array} \\
\begin{bmatrix}
0 & 1 & 0 & 0 & 1 & 1 & -1 \\
1 & 0 & 1 & 1 & 0 & 0 & -1 \\
1 & 1 & 0 & 0 & 1 & 1 & -1 \\
0 & 1 & 0 & 0 & 0 & 0 & 0 \\
0 & 0 & 0 & 0 & 0 & 0 & -1 \\
1 & 1 & 1 & 0 & 1 & 0 & -1 \\
1 & 1 & 1 & 1 & 1 & 1 & 0
\end{bmatrix}
\end{array}
$$

The individual graphs associated with the matrices $B_1$ and $B_2$ are as follows.

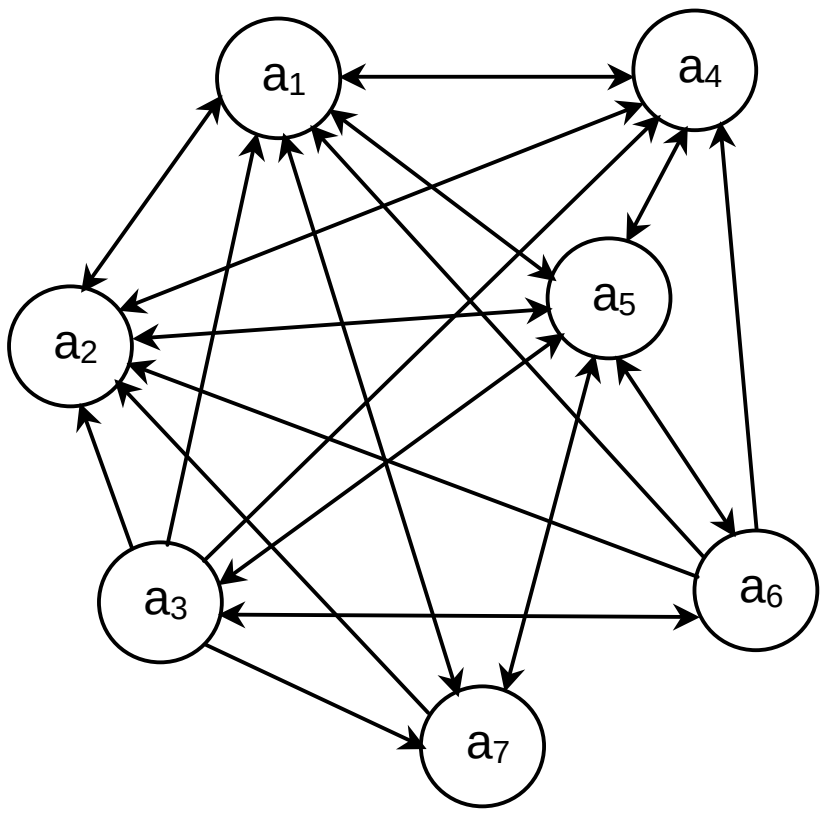

FIGURE: 3.2.7a

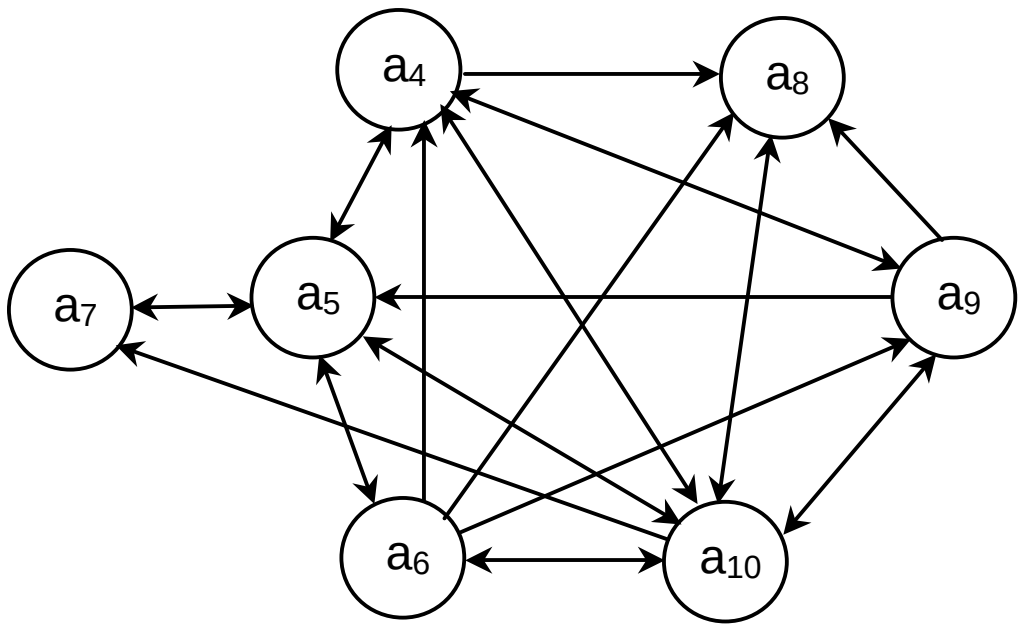

FIGURE: 3.2.7b

Now consider the on state of the row bivector $a_4$ alone and all other nodes are in the off state.

Let X = (0 0 0 1 0 0 0) ∪ (1 0 0 0 0 0 0) be the bivector the effect of X on the dynamical system B is given by

$$\begin{aligned} XB &= (X_1 \cup X_2)(B_1 \cup B_2) \\ &= X_1 B_1 \cup X_2 B_2 \\ &= (1\ 1\ 0\ 0\ 1\ 0\ 0) \cup (0\ 1\ 0\ 0\ 1\ 1\ -1). \end{aligned}$$

After thresholding and updating we get

$$Y = (1\ 1\ 0\ 1\ 1\ 0\ 0) \cup (1\ 1\ 0\ 0\ 1\ 1\ 0) = Y_1 \cup Y_2.$$

Now we study the effect of Y on the dynamical system B.

$$YB = (Y_1 \cup Y_2)\ (B_1 \cup B_2) = Y_1\ B_1 \cup Y_2\ B_2. = (3\ 3\ –1\ 3\ 3\ 1\ 2) \cup (2\ 2\ 2\ 1\ 2\ 1\ –\ 4).$$

After updating and thresholding the resultant vector we get

$$Z = (1\ 1\ 0\ 1\ 1\ 1\ 1) \cup (1\ 1\ 1\ 1\ 1\ 1\ 0) = Z_1 \cup Z_2.$$

$$ZB = (Z_1 \cup Z_2)\ (B_1 \cup B_2) = Z_1\ B_1 \cup Z_2\ B_2 = (5\ 5\ –2\ 4\ 5\ 1\ 2) \cup (3\ 5\ 2\ 1\ 3\ 2\ –5).$$

After updating and thresholding the resultant vector we get

$$= R = R_1 \cup R_2 = (1\ 1\ 0\ 1\ 1\ 1\ 1) \cup (1\ 1\ 1\ 1\ 1\ 1\ 0) = Z.$$

Thus the bihidden pattern happens to be a fixed bipoint when indigestion i.e. $a_4$ alone is in the on state in both the state vectors i.e. in the row bivector, in the resultant all nodes become on, except $a_3$ and $a_{10}$ which implies spraying of pesticides and consuming near by vegetables has no impact on indigestion, but all other medical problems related with indigestion viz. loss of appetite, headache, giddiness / fainting mouth and stomach ulcer breathlessness, skin ailments and diarrhea come to on state which shows the system is well modeled . Thus one can work with any of the on state of the bivectors and arrive at a conclusion using the dynamical system B.

## 3.3 Fuzzy Cognitive Trimaps and their Applications

Now we proceed on to define the notion of Fuzzy Cognitive Trimaps (FCTMs). Here we take either the opinion of three experts or of a single expert on three sets of attributes. The directed graph of a FCTMs may be a disjoint trigraph or a strongly disjoint trigraph or a connected trigraph by vertices or a edge connected trigraph depending on the model under investigation. The FCTMs directed trigraph may help one to analyze the problem simultaneously using three experts. Also the use of this model helps in simplification of the problem or helps one to divide into blocks of lesser size so that it will not be difficult to draw graph or in calculations he will illustrate them with models before which we give a brief description of the model.

**DEFINITION 3.3.1:** *A Fuzzy Cognitive Trimap (FCTM) is a directed trigraph which has its vertices or nodes as concepts or policies and edges are the causal relations between the nodes. The causal relations are denoted by $e_{ij}$. If $e_{ij} > 0$ then the causal relation between the edge $C_i$ and $C_j$ is such that, the increase (or decrease) in the node $C_i$ increases (or decreases) in the node $C_j$, $e_{ij} = 0$ implies there is no casual relation between the nodes $C_i$ and $C_j$. If $e_{ij} < 0$ means the causal relation between $C_i$ and $C_j$ is such that the increase (or decrease) in the node $C_i$ means a decrease (or increase) in the node $C_j$. If $e_{ij}$ takes values from the set {-1, 0, 1} then it implies that the FCTM is a simple one.*

The main advantage of using simple FCTMs is that it will be easy for an expert to work with. He need not be a person with good mathematical background. For an expert may be a doctor or a lay man who can give the association as full values 1, 0 or –1. Thus in such cases the simple FCTMs are best suited. Now once the expert gives his opinion about an FCTM that is made into the directed trigraph. This directed trigraph is converted into the related connection trimatrix of the trigraph. Using the trimatrix which we call as the dynamical trisystem of the model and work for the solutions. The solution is given by the stability of the system. The system is stable if it gives either the fixed tripoint or a limit tricycle or a fixed point and a limit bicycle or a limit cycle and a fixed bipoint. This is known as the trihidden pattern of the system.

Just we describe, the working of the model. Any state vector will be a trirow matrix say

$$\left\{\left(a_1^1,\ldots,a_n^1\right)\cup\left(a_1^2,\ldots,a_m^2\right)\cup\ldots\cup\left(a_1^3,\ldots,a_r^3\right)\right\}$$

where $a_j^i$ = 0 if the node is in the off state

= 1 if the node is in the on state

$i = 1, 2, 3$

$1 \leq j \leq n$ when $i = 1$

$1 \leq j \leq m$ when $i = 2$

$1 \leq j \leq r$ when $i = 3$.

Suppose the trimatrix associated with the trigraph is denoted by $A = A_1 \cup A_2 \cup A_3$ where $A_1$ $A_2$ and $A_3$ are distinct square matrices of same order or of different order. Let $X = X_1 \cup X_2 \cup X_3$ be the trirow matrix. Now to study the effect of X on A we multiply X with A i.e.

$$\begin{aligned} XA &= (X_1 \cup X_2 \cup X_3)(A_1 \cup A_2 \cup A_3) \\ &= X_1 A_1 \cup X_2 A_2 \cup X_3 A_3, \end{aligned}$$

the resultant is once again a trirow vector. Since the dynamical system can realize only the ON or OFF state of the nodes and if the trirow vector by choosing an arbitrary constant K such that if $x_i \in X_1 A_1 \cup X_2 A_2 \cup X_3 A_3$ then

$x_i = 0$ if $x_i < K$.

$x_i = 1$ if $x_i \geq K$.

This process of making the entries in the trimatrix to be either 0 or 1 is called thresholding, by updating we mean we make the entry which we started in the on state to be on i.e. if some $x_k$ entry was started with $x_k = 1$ we see to it that $x_k = 1$ in the resultant i.e., in the trihidden pattern if $x_k = 0$ we make $x_k = 1$. Suppose Y is the resultant trirow vector after we update and threshold i.e. if $Y = Y_1 \cup Y_2 \cup Y_3$ then we find

$$YA = (Y_1 \cup Y_2 \cup Y_3)\,(A_1 \cup A_2 \cup A_3)$$
$$= Y_1 A_1 \cup Y_2 A_2 \cup Y_3 A_3$$

now this resultant vector is once again updated and thresholded. We proceed on till we get a fixed tripoint or a limit tricycle or so on. Unlike in usual FCMs we have the following possibilities in the trihidden pattern. The resultant trirow matrix X is such that each of $X_1$, $X_2$ and $X_3$ are fixed points or one of $X_1$, $X_2$ or $X_3$ is a fixed point and the other two and just limit cycle or one of $X_1$, $X_2$ or $X_3$ is a limit cycle and other two are fixed points or all the resultant vectors are limit cycles.

Thus we have four possibilities to occur. Using the resultant which is known as the trihidden pattern of the system we analyze the problem or interpret our resultant trivector. Now we illustrate the models by a real world problem.

Suppose we are interested in studying the health hazards faced by the agriculture labourers due to pollution by using chemicals in pesticides and fertilizers. We after collecting data from the labourers use experts and get their opinion over the problem.

First we analyze 3 experts opinion $E_1$ $E_2$ and $E_3$ all the three of them use the same set of 12 attributes viz.

- $C_1$ – Consuming nearby green vegetables [after spraying of pesticides and insecticides]
- $C_2$ – Indigestion and loss of appetite
- $C_3$ – Spraying pesticides / using of chemical fertilizers
- $C_4$ – Headache giddiness / fainting
- $C_5$ – Exposure to insecticides and pesticides after spraying in the agricultural field.
- $C_6$ – Mouth and stomach ulcer
- $C_7$ – Breathlessness
- $C_8$ – Skin ailments
- $C_9$ – Diarrhea
- $C_{10}$ – Lack of precaution and treatment

$C_{11}$ – Degradation of the capacity of the land for further cultivation

$C_{12}$ – Financial loss due to less yield and extra expenditure meted out to health care.

The health hazards / symptoms felt by the agricultural labourer after the spray or using fertilizers / pesticides.

Three experts are used to find a solution to the problem. Since over lap of opinion is not possible and all the three experts work only on the same 12 attributes the directed trigraph given by them will be a strongly disjoint trigraph. However as there are twelve vertices or nodes and several casual relations we do not give the directed trigraph. We use the trigraph to obtain the related trimatrix $T = T_1 \cup T_2 \cup T_3$. As each of the matrices is a $12 \times 12$ matrix i.e. trimatrix is a $12 \times 12$ square matrix we just give them as $T = T_1 \cup T_2 \cup T_3$ each $T_i$ alone is described from the inspection of these matrices we see all the three of them are distinct. The connection matrix $T_1$ given by the first expert.

$$T_1 = \begin{bmatrix} 0 & 1 & 0 & 1 & 0 & 1 & 0 & 1 & 1 & 0 & 0 & 1 \\ -1 & 0 & 0 & 1 & 0 & 1 & 0 & 1 & 0 & 0 & 0 & 1 \\ -1 & 1 & 0 & 0 & 0 & 0 & 0 & 0 & -1 & 0 & 1 \\ -1 & 1 & 0 & 0 & 0 & 0 & 1 & 1 & 0 & 0 & 0 & 1 \\ 1 & 1 & 1 & 1 & 0 & 1 & 1 & 1 & 1 & 1 & 1 & 1 \\ -1 & 1 & -1 & 1 & -1 & 0 & 0 & 1 & 0 & -1 & -1 & 1 \\ -1 & 0 & 0 & 1 & -1 & 0 & 0 & 0 & 0 & -1 & 0 & 1 \\ -1 & 1 & -1 & 1 & -1 & 1 & 0 & 0 & 0 & -1 & 0 & 1 \\ 1 & 1 & -1 & 1 & -1 & 0 & 0 & 0 & 0 & -1 & 0 & 1 \\ -1 & 1 & 1 & 1 & -1 & 1 & 1 & 1 & 1 & 0 & 1 & 1 \\ -1 & 0 & 1 & 1 & 1 & 0 & 0 & 0 & 0 & 0 & 0 & 1 \\ -1 & 0 & 1 & 1 & 1 & 0 & 0 & 0 & 0 & -1 & 0 & 0 \end{bmatrix}$$

Now we give the connection matrix $T_2$ as given by the second expert

$$T_2 = \begin{bmatrix} 0 & 1 & 0 & 1 & 0 & 1 & 0 & 0 & 1 & 0 & 0 & 1 \\ -1 & 0 & 0 & 1 & 0 & 1 & 0 & 1 & 0 & 0 & 0 & 1 \\ -1 & 1 & 0 & 0 & 0 & 0 & 0 & 0 & 0 & -1 & 0 & 0 \\ 0 & 1 & 0 & 0 & 0 & 0 & 1 & 1 & 0 & 0 & 0 & 1 \\ 1 & 1 & 1 & 0 & 0 & 1 & 1 & 1 & 1 & 1 & 0 & 1 \\ -1 & 0 & 0 & 1 & -1 & 0 & 0 & 0 & 0 & 0 & 0 & 1 \\ 0 & 1 & -1 & 1 & 0 & 0 & 0 & 1 & 0 & -1 & -1 & 1 \\ -1 & 1 & -1 & 1 & 1 & -1 & 1 & 0 & 0 & 0 & -1 & 0 \\ 1 & 1 & -1 & 1 & -1 & 0 & 0 & 0 & 0 & 1 & 0 & 1 \\ -1 & 1 & 1 & 1 & -1 & 1 & 1 & 1 & 1 & 0 & 1 & 1 \\ -1 & 0 & 1 & 1 & 1 & 0 & 0 & 0 & 0 & 0 & 0 & 1 \\ -1 & 1 & 1 & 1 & 1 & 0 & 0 & 0 & 0 & 0 & 0 & 0 \end{bmatrix}$$

The connection matrix given by the third expert say $T_3$ is given.

$$T_3 = \begin{bmatrix} 0 & 1 & 1 & 0 & 1 & 0 & 0 & 1 & 0 & 0 & 0 & 1 \\ -1 & 0 & 0 & 1 & 0 & 1 & 0 & 1 & 1 & 0 & 0 & 1 \\ -1 & 1 & 0 & 0 & 0 & 0 & 0 & 0 & 0 & -1 & 0 & 0 \\ 0 & 1 & 0 & 0 & 1 & 0 & 0 & 1 & 1 & 0 & 0 & 1 \\ 1 & 1 & 1 & 0 & 0 & 1 & 1 & 1 & 0 & 1 & 1 & 1 \\ -1 & 0 & 0 & 1 & -1 & 0 & 0 & 0 & 0 & 0 & 1 & 1 \\ 0 & 1 & -1 & 1 & 0 & 0 & 0 & 1 & 0 & -1 & -1 & 1 \\ -1 & 1 & 1 & 1 & 0 & 0 & -1 & 0 & 0 & 0 & -1 & 1 \\ 1 & 1 & -1 & 1 & 1 & 0 & 0 & 1 & 0 & 1 & 0 & 0 \\ -1 & 1 & 0 & 1 & -1 & 0 & -1 & 1 & 0 & 0 & 0 & -1 \\ -1 & 0 & 1 & 1 & 0 & 0 & 1 & 0 & 0 & 0 & 0 & 1 \\ 0 & 1 & 1 & 1 & 1 & 0 & 0 & 0 & 0 & 1 & 0 & 0 \end{bmatrix}$$

Now the trimatrix $T = T_1 \cup T_2 \cup T_3$ is the trimodel associated with the FCTM. Now we study the effect of any state trirow vector on $T = T_1 \cup T_2 \cup T_3$ .

Let

$$
\begin{aligned}
X &= X_1 \cup X_2 \cup X_3 \\
&= (0\ 1\ 0\ 0\ 0\ 0\ 0\ 0\ 0\ 0\ 0\ 0\ 0) \cup \\
&\quad (1\ 0\ 0\ 0\ 0\ 0\ 0\ 0\ 0\ 0\ 0\ 0) \cup \\
&\quad (0\ 0\ 0\ 1\ 0\ 0\ 0\ 0\ 0\ 0\ 0\ 0).
\end{aligned}
$$

The node $C_2$ in the first component of the trivector, node $C_1$ in the second component of the trivector and the node $C_4$ in the third component of trivector are in the on state and all other nodes in the trivector are in off state. The effect of X on the dynamical system T is given by

$$
\begin{aligned}
XT &= (X_1 \cup X_2 \cup X_3)\ (T_1 \cup T_2 \cup T_3) \\
&= (-1\ 0\ 0\ 1\ 0\ 1\ 0\ 1\ 0\ 0\ 0\ 1) \cup \\
&\quad (0\ 1\ 0\ 1\ 0\ 1\ 0\ 0\ 1\ 0\ 0\ 1) \cup \\
&\quad (0\ 1\ 0\ 0\ 1\ 0\ 0\ 1\ 10\ 0\ 1).
\end{aligned}
$$

After thresholding and updating we get

$$
\begin{aligned}
Y &= Y_1 \cup Y_2 \cup Y_3 \\
&= (0\ 1\ 0\ 10\ 1\ 0\ 1\ 0\ 0\ 0\ 1) \cup \\
&\quad (1\ 1\ 0\ 1\ 0\ 1\ 0\ 0\ 1\ 0\ 0\ 1) \cup \\
&\quad (0\ 1\ 0\ 1\ 1\ 0\ 0\ 1\ 1\ 0\ 0\ 1).
\end{aligned}
$$

The effect of Y on he dynamical system T is given by

$$
\begin{aligned}
YT &= (Y_1 \cup Y_2 \cup Y_3)\ (T_1 \cup T_2 \cup T_3) \\
&= Y_1\ T_1 \cup Y_2\ T_2 \cup Y_3\ T_3 \\
&= (-5, 4, -1, 4, -1, 1, 1\ 3, 0\ -3, -1, 4) \cup \\
&\quad (-4\ 5\ 0\ 6\ -1\ 3\ 1\ 2\ 2\ 0\ -1\ 5) \cup \\
&\quad (0\ 4\ 2\ 3\ 2\ 2\ 0\ 3\ 1\ 0\ 3\ 3\ ).
\end{aligned}
$$

After thresholding and updating we get the resultant as

$$
\begin{aligned}
Z &= (0\ 1\ 0\ 1\ 0\ 1\ 1\ 1\ 0\ 0\ 0\ 1) \cup \\
&\quad (1\ 1\ 0\ 1\ 0\ 1\ 1\ 1\ 1\ 1\ 0\ 1) \cup \\
&\quad (0\ 1\ 1\ 1\ 1\ 1\ 0\ 1\ 1\ 1\ 1\ 1) \\
&= Z_1 \cup Z_2 \cup Z_3.
\end{aligned}
$$

$(-6\ 3\ -1\ 5\ -2\ 2\ 1\ 3\ 0\ -4\ 15) \cup$
$(-3\ 6\ -2\ 7\ 0\ 1\ 2\ 2\ 1\ 0\ -2\ 6) \cup$
$(-3\ 6\ 3\ 6\ 2\ 2\ -1\ 5\ 2\ 2\ 1\ 6).$

After thresholding and updating the resultant vector we get

$$P = (0\ 1\ 0\ 1\ 0\ 1\ 1\ 1\ 0\ 0\ 0\ 1) \cup (1\ 1\ 0\ 1\ 0\ 1\ 1\ 1\ 1\ 0\ 0\ 1) \cup (0\ 1\ 1\ 1\ 1\ 1\ 0\ 1\ 1\ 1\ 1\ 1) = P_1 \cup P_2 \cup P_3.$$

Now we study the effect of P on T.

$$PT = P_1 T_1 \cup P_2 T_2 \cup P_3 T_3 = (0\ 1\ 0\ 1\ 0\ 1\ 1\ 1\ 0\ 0\ 0\ 1) \cup (1\ 1\ 0\ 1\ 0\ 1\ 1\ 1\ 1\ 0\ 0\ 1) \cup (0\ 1\ 1\ 1\ 1\ 1\ 0\ 1\ 1\ 1\ 1\ 1).$$

The first two components of the row trivector are fixed tripoint the resultant of the 3rd component after updating and thresholding we get as a fixed point. Thus the hidden pattern for this given vector is a fixed tripoint. From this we see when indigestion or loss of appetite is in the on state $C_1$, $C_3$, $C_5$, $C_9$, $C_{10}$ and $C_{11}$ remain unaffected only $C_2$, $C_4$, $C_6$, $C_7$, $C_8$ and $C_{12}$ come to on state which implies that Indigestion is associated with Headache giddiness / fainting, mouth and stomach ulcer, breathless ness, skin ailment, and financial loss due to taking treatment for the same and consuming the nearly green vegetables $C_3$ $C_5$ $C_{10}$ and $C_{11}$ remain unaffected as they are in the off state and nodes $C_2$, $C_4$, $C_6$, $C_7$, $C_8$, $C_9$ and $C_{12}$ come to on states signifying the impact of eating the near by vegetables. Finally when $C_4$ i.e. Headache with giddiness alone in the on state and all other in the off state we see $C_1$ and $C_7$ alone are in off state and all other nodes come to on state. Thus we can consider the effect of any other state vector. This model is nothing but using the notion of disjoint square trimatrix of same rank.

Now we proceed on to study FCTM model using disconnected directed trigraphs.

***Example 3.3.2:*** Now suppose we are interested to study the relation between the health hazards the agricultural labourers suffer due to chemical pollution by using chemical fertilizers and insecticides. Suppose $P_1, \ldots, P_{11}$ are the 11 nodes / concepts given by the expert as given in page 94 of chapter 3.

Now three experts share their opinion on the 11 attributes and they give their opinion which is given by the trigraph.

The first expert uses the nodes $P_1$, $P_2$, $P_3$, $P_4$ as the vertex of the trigraph, the second expert uses $P_5$, $P_6$, $P_7$, $P_8$ as the vertex of the trigraph and the third expert uses $P_9$, $P_{10}$ and $P_{11}$ as the nodes / vertices. The directed trigraph is a disjoint trigraph given by the following figure.

G =

$G_1 \cup G_2 \cup G_3$

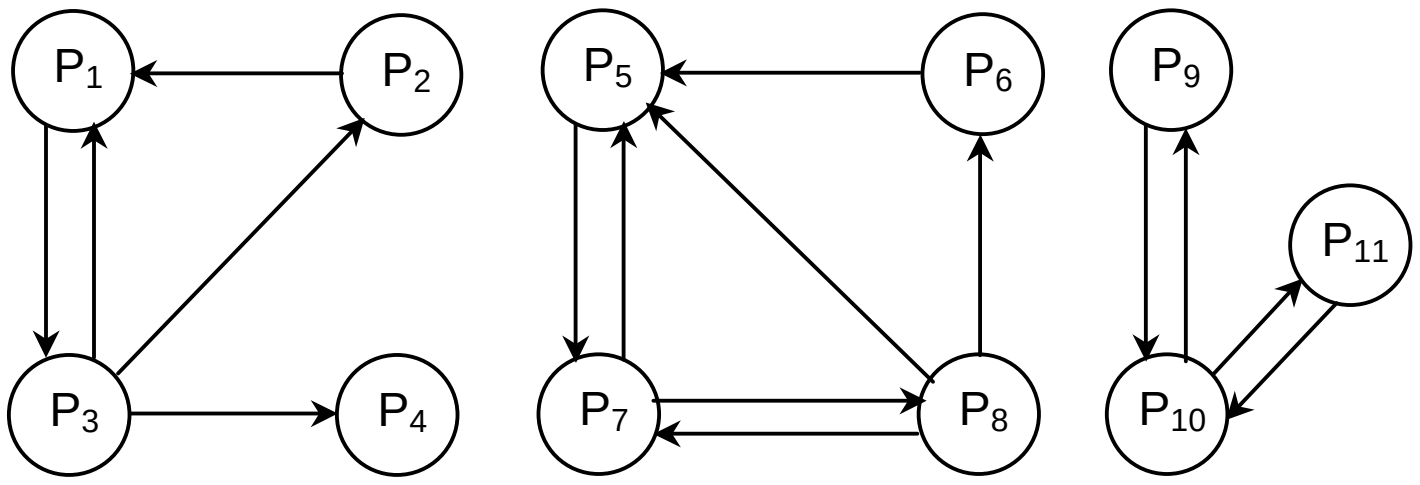

FIGURE: 3.3.1

The related trimatrix.

$$T = T_1 \cup T_2 \cup T_3$$

$$= \begin{array}{c} \\ P_1 \\ P_2 \\ P_3 \\ P_4 \end{array} \begin{array}{c} P_1 \; P_2 \; P_3 \; P_4 \\ \begin{bmatrix} 0 & 0 & 1 & 0 \\ 1 & 0 & 0 & 0 \\ 1 & 1 & 0 & 1 \\ 0 & 0 & 0 & 0 \end{bmatrix} \end{array} \cup \begin{array}{c} \\ P_5 \\ P_6 \\ P_7 \\ P_8 \end{array} \begin{array}{c} P_5 \; P_6 \; P_7 \; P_8 \\ \begin{bmatrix} 0 & 0 & 1 & 0 \\ 1 & 0 & 0 & 0 \\ 1 & 0 & 0 & 1 \\ 1 & 1 & 1 & 0 \end{bmatrix} \end{array} \cup \begin{array}{c} \\ P_9 \\ P_{10} \\ P_{11} \end{array} \begin{array}{c} P_9 \; P_{10} \; P_{11} \\ \begin{bmatrix} 0 & 1 & 0 \\ 1 & 0 & 1 \\ 0 & 1 & 0 \end{bmatrix} \end{array}.$$

Clearly the trimatrix is a square mixed trimatrix.

This model is different from the previous model for in the previous model all the square matrices were of same size and the three experts gave their opinion on the same set of nodes. Here the nodes are different and the experts are also different.

Now we illustrate how our model works. Consider a trirow vector.

$$X = (1\ 0\ 0\ 0) \cup (0\ 1\ 0\ 0) \cup (0\ 0\ 1)$$
$$= X_1 \cup X_2 \cup X_2,$$

where swollen limbs is in the on state for the first expert, pollution after drinking water alone is in the on state for the second expert system and blurred vision is in the on state for the third expert and other the nodes are in the off state. We study the effect of X on the dynamical system T.

$$XT = (X_1 \cup X_2 \cup X_3)\ (T_1 \cup T_2 \cup T_3)$$
$$= X_1\ T_1 \cup X_2\ T_2 \cup X_3\ T_3$$
$$= (0\ 0\ 1\ 0\ ) \cup (1\ 0\ 0\ 0) \cup (0\ 1\ 0).$$

After updating we get the resultant vector as

$$Y = (1\ 0\ 1\ 0) \cup (1\ 1\ 0\ 0) \cup (0\ 1\ 1)$$
$$= Y_1 \cup Y_2 \cup Y_3.$$

Now the effect of Y on the dynamical system T is given by

$$YT = (Y_1 \cup Y_2 \cup Y_3)\ (T_1 \cup T_2 \cup T_3)$$
$$= Y_1\ T_1 \cup Y_2\ T_2 \cup Y_3\ T_3$$
$$= (1\ 1\ 1\ 1) \cup (1\ 0\ 1\ 0) \cup (1\ 1\ 1)$$

after updating we get the resultant vector as

$$Z = (1\ 1\ 1\ 1) \cup (1\ 1\ 1\ 0) \cup (1\ 1\ 1)$$
$$= Z_1 \cup Z_2 \cup Z_3.$$

The effect of Z on the system T is given by

$$
\begin{aligned}
ZT &= (Z_1 \cup Z_2 \cup Z_3)\,(T_1 \cup T_2 \cup T_3) \\
&= Z_1\,T_1 \cup Z_2\,T_2 \cup Z_3\,T_3 \\
&= (2\ 1\ 1\ 1) \cup (2\ 0\ 1\ 1) \cup (1\ 2\ 1).
\end{aligned}
$$

After thresholding and updating we get

$$
\begin{aligned}
R &= (1\ 1\ 1\ 1) \cup (1\ 1\ 1\ 1) \cup (1\ 1\ 1) \\
&= R_1 \cup R_2 \cup R_3.
\end{aligned}
$$

Now in the trirow vector last and first row vectors have become fixed tripoint.
Now the effect of R or T is given by

$$
\begin{aligned}
RT &= (R_1 \cup R_2 \cup R_3)\,(T_1 \cup T_2 \cup T_3) \\
&= R_1\,T_1 \cup R_2\,T_2 \cup R_3\,T_3 \\
&= (1\ 1\ 1\ 1) \cup (3\ 1\ 2\ 1) \cup (1\ 1\ 1).
\end{aligned}
$$

After thresholding and updating we get the resultant as

$$
S = (1\ 1\ 1\ 1) \cup (1\ 1\ 1\ 1) \cup (1\ 1\ 1)
$$

the hidden tripattern is a fixed tripoint.

Thus all the states come to on state proving that the 3 nodes are vital for it makes all other nodes to on state.

Similar study can be carried out for this model in another example where the trigraph happens to be a disconnected directed trigraph.

Now we proceed on to give an example of a model in which the directed trigraph is connected by vertices.

Now the same problem is studied by three experts with one vertex in common for the first expert gives opinion on $\{P_1, P_2, P_3, P_4, P_5\}$ second expert on $\{P_5, P_6, P_7, P_8\}$ and the third expert on $\{P_8, P_9, P_{10}, P_{11}, P_4\}$.

The directed trigraph based on the opinion of these experts are given by the following figure.

$$G = G_1 \cup G_2 \cup G_3.$$

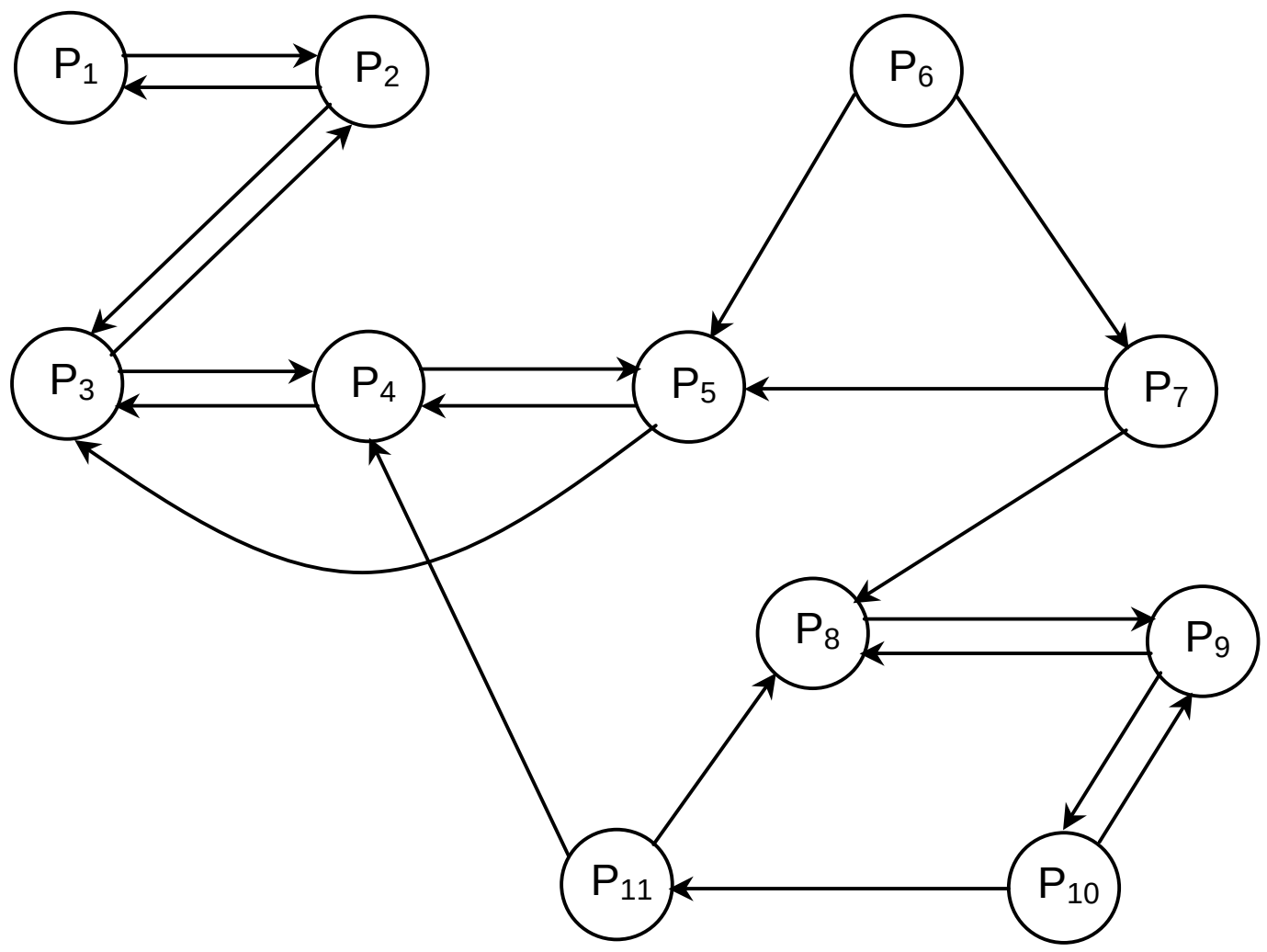

FIGURE: 3.3.2

The directed trigraph is connected by vertices for the graphs $G_1$ and $G_2$ connected by the vertex $P_5$, the graphs $G_2$ and $G_3$ connected by the vertex $P_8$ and the graphs $G_1$ and $G_3$ connected by the vertex $P_4$ such connected trigraphs are called as *cyclically connected trigraphs*.

Now we just give the separated directed graph of each model.

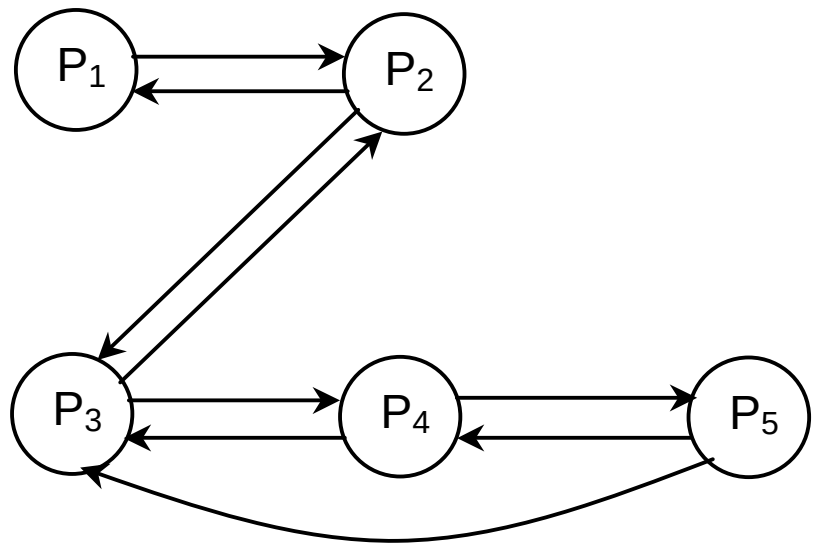

$G_1$

FIGURE: 3.3.2a

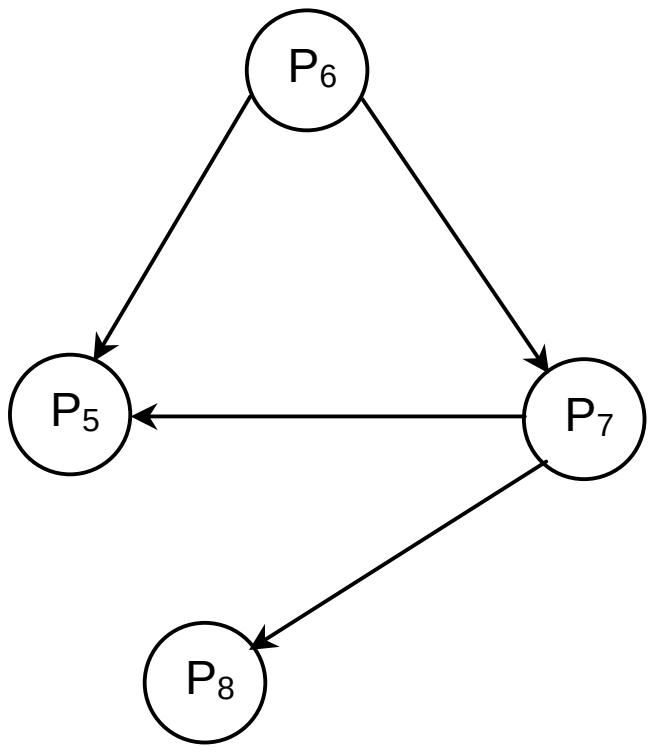

$G_2$

FIGURE: 3.3.2b

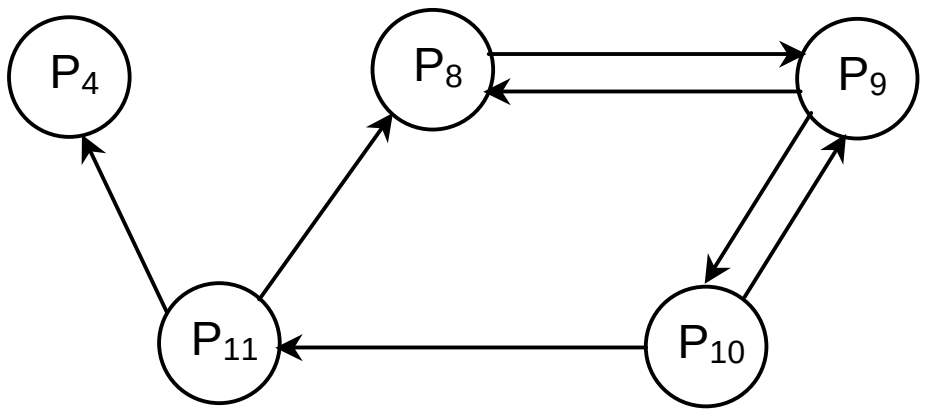

$G_3$

FIGURE: 3.3.2c

The connection trimatrix of the directed trigraph $T = T_1 \cup T_2 \cup T_3$ is as follows:

$$\begin{array}{c} \\ P_1 \\ P_2 \\ P_3 \\ P_4 \\ P_5 \end{array}\begin{array}{c} \begin{array}{ccccc} P_1 & P_2 & P_3 & P_4 & P_5 \end{array} \\ \begin{bmatrix} 0 & 1 & 0 & 0 & 0 \\ 1 & 0 & 1 & 0 & 0 \\ 0 & 1 & 0 & 1 & 0 \\ 0 & 0 & 1 & 0 & 1 \\ 0 & 0 & 1 & 1 & 0 \end{bmatrix} \end{array} \cup \begin{array}{c} \\ P_5 \\ P_6 \\ P_7 \\ P_8 \end{array}\begin{array}{c} \begin{array}{cccc} P_5 & P_6 & P_7 & P_8 \end{array} \\ \begin{bmatrix} 0 & 0 & 0 & 0 \\ 1 & 0 & 1 & 0 \\ 1 & 0 & 0 & 1 \\ 0 & 0 & 0 & 0 \end{bmatrix} \end{array}$$

$$\cup \begin{array}{c} \\ P_4 \\ P_8 \\ P_9 \\ P_{10} \\ P_{11} \end{array} \begin{array}{c} \begin{array}{ccccc} P_4 & P_8 & P_9 & P_{10} & P_{11} \end{array} \\ \left[\begin{array}{ccccc} 0 & 0 & 0 & 0 & 0 \\ 0 & 0 & 1 & 0 & 0 \\ 0 & 1 & 0 & 1 & 0 \\ 0 & 0 & 1 & 0 & 1 \\ 1 & 1 & 0 & 0 & 0 \end{array}\right] \end{array}.$$

Suppose we wish the study the effect of the state trirow vector X where

$$\begin{aligned} X &= (0\ 0\ 0\ 1\ 0) \cup (0\ 0\ 0\ 1) \cup (0\ 1\ 0\ 0\ 0) \\ &= X_1 \cup X_2 \cup X_3. \end{aligned}$$

The effect of X on the dynamical system T is given by

$$\begin{aligned} XT &= X_1\, T_1 \cup X_2\, T_2 \cup X_3\, T_3 \\ &= (0\ 0\ 1\ 0\ 1) \cup (0\ 0\ 0\ 0) \cup (0\ 0\ 1\ 0\ 0). \end{aligned}$$

After updating we get

$$\begin{aligned} Y &= (0\ 0\ 1\ 1\ 1) \cup (0\ 0\ 0\ 1) \cup (0\ 1\ 1\ 0\ 0) \\ &= Y_1 \cup Y_2 \cup Y_3. \end{aligned}$$

Now we study the effect Y on the dynamical system T.

$$YT = (0\ 1\ 2\ 2\ 1) \cup (0\ 0\ 0\ 0) \cup (0\ 1\ 1\ 1\ 0).$$

After thresholding and updating we get the resultant trivector as

$$\begin{aligned} Z &= (0\ 1\ 1\ 1\ 1) \cup (0\ 0\ 0\ 1) \cup (0\ 1\ 1\ 1\ 0). \\ &= Z_1 \cup Z_2 \cup Z_3. \end{aligned}$$

The effect of Z on T is given by

$$\begin{aligned} ZT &= Z_1\, T_1 \cup Z_2\, T_2 \cup Z_3\, T_3 \\ &= (1\ 2\ 3\ 2\ 1) \cup (0\ 0\ 0\ 0) \cup (0\ 1\ 2\ 1\ 1) \end{aligned}$$

after thresholding and updating we get the resultant trivector as

$$W = (1\ 1\ 1\ 1\ 1) \cup (0\ 0\ 0\ 1) \cup (0\ 1\ 1\ 1\ 1).$$

W is the trihidden pattern of the system, when the node vomiting alone in the on state, all nodes come to on state whereas, $P_8$ has no effect on other nodes, i.e. loss of appetite has no influence on or mouth ulcer or pollution by drinking water or indigestion but however it has effect on Headache and spraying of pesticides. Study of this form can be carried out using FCTMs which gives a directed vertex connected trigraphs, we have given this illustration mainly for understanding and the working of the FCTM. .

***Example 3.3.3:*** We can also have another type of FCTMs model in which the directed trigraph is just disconnected but is the union of the graph and a connected bigraph.

We illustrate this model in the study of globalization and its impact on farmers, for globalization directly affects the agriculturists when the seeds have to get the patent right, which will increase the cost of seeds; for as usually farmer in India when they grow any crop (like paddy, groundnut etc) they in the time of harvest first select quality seed from the harvested crop and keep it as the seed to be shown for the next time. But with the patent right for the seeds with the globalization they should buy seed to be sown which has visibly two drawbacks. The first one the farmer has to spend extra cash for the purchase, two the seeds which are used in his soil may give a better yield than the new seed purchased by him for the same purpose. So by the globalization he is still uncertain of his yield leaving the biggest uncertainty the monsoon. When we speak about farmers we do not include very rich land lords but only middle class or poor farmers.

Now we after consulting with several experts give the important attributes related with the problems faced by farmer due to globalization. The attributes as given by the expert are as follows:

$G_1$ – Quality of seeds and its cost.

$G_2$ – Usage of fertilizers banned by foreign countries. (Several times the government of India saying they pay loans to farmers, in many cases distribute to them these banned fertilizers freely, thus they have no other option but to use them).

$G_3$ – Patent right owned by foreign countries. (When the seeds are such that the patient right is owned by foreigners the farmer cannot even be very sure of good yield for no one is sure how the yield would be in Indian soil for the Indian soil itself is very distinctly different from region to region).

$G_4$ – Loan facilities offered by financial institutions.

$G_5$ – Suicide of farmers. (when the monsoon fails and government instead of giving a helping hand demand them to pay back loans the farmers have no other option but to commit suicide. For in the past five years, the suicide by farmers have been on the increase. This has its impact on the globalization for at times one is made to research when the seed of our own soil is used it would with stand the monsoon failure and the loss may not be that great).

$G_6$ – Obtaining electricity and irrigational facilities from government. (In the past when the farmers did not get electricity or any irrigational facilities from the government they lead a peaceful and a contended life; only on the advent of getting help from government which is not always certain depending from party to party which rules they suffer more).

$G_7$ – Usage of technological equipment. (This will have more relevance to rich or very rich farmers certainly not for a middle class or a poor farmer).

$G_8$ – Updation of cultivation method (This also has no relevance to poor or middle class farmers).

$G_9$ – Cultivating quality crops. (It is a matter of fact all these fertilizer so called quality seed etc only has brought down the quality of crops in India. Before such uses the Indian crops had its value now even the food value of these crops have drastically reduced that is why several of the farmers and their children suffer from malnutrition and other types of diseases).

$G_{10}$ – Adverse psychological effect on farmers. (Even a decade back we have never heard of farmers suffering from tensions, and hypertension and related problems. The very impact of globalization had made them more like a person employed in a white coloured job who suffers all types of stress and strain in his work place).

Now using these attributes we obtain the experts opinion using the FCTM model. As all these analysis highly involve uncertainty and impreciseness so we are justified in using FCTM.

Now the directed trigraph $T = T_1 \cup T_2 \cup T_3$ given by the expert is as follows.

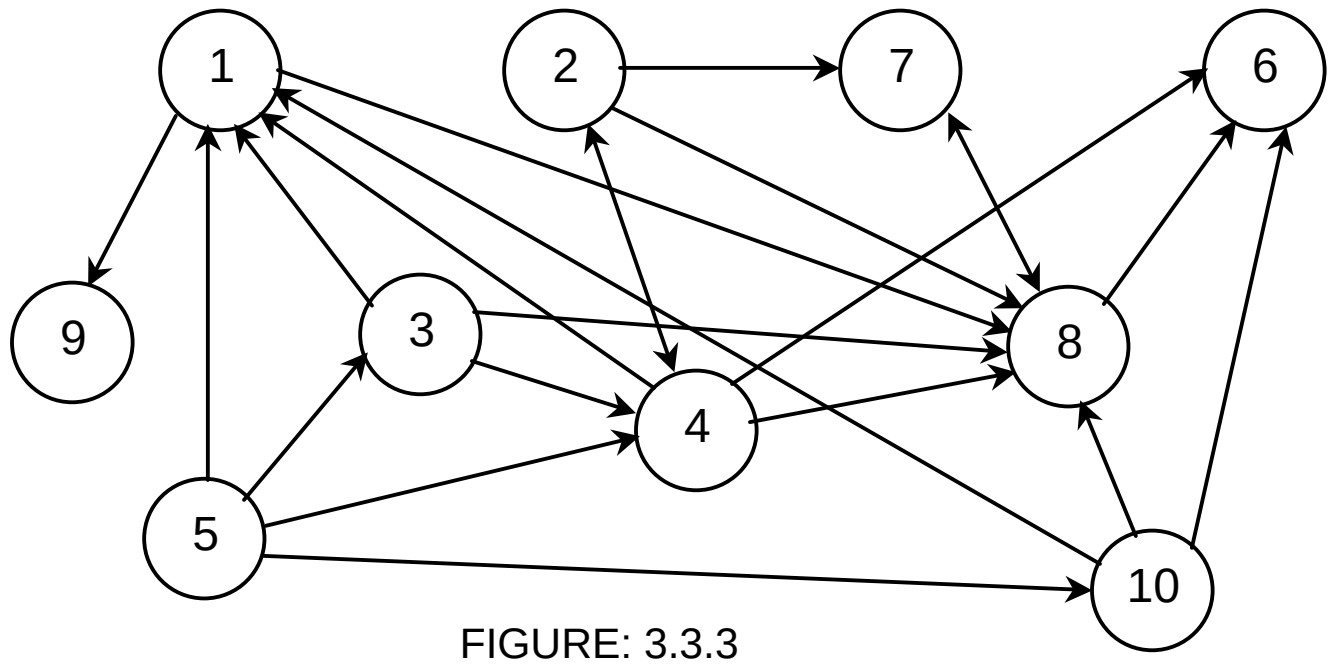

FIGURE: 3.3.3

This trigraph may look confused for the reader so we also give the separated 3 graphs associated with this trigraph

Graph of $T_1$

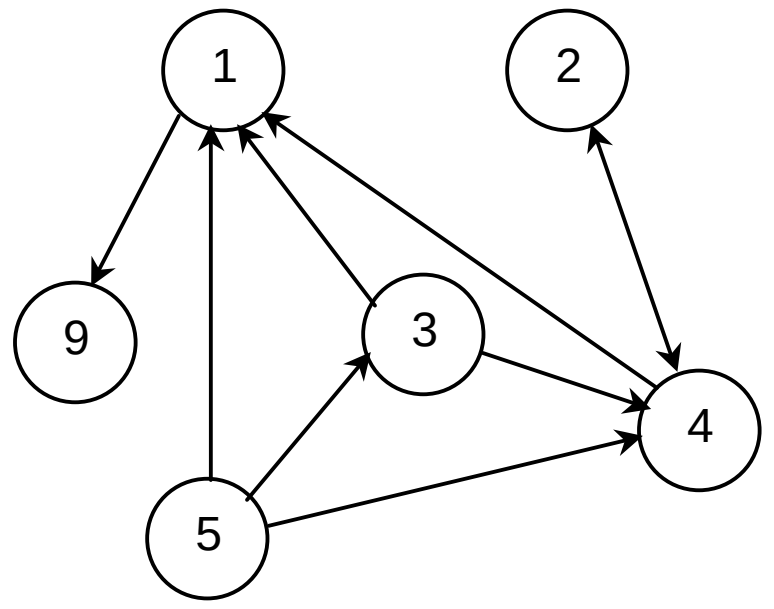

FIGURE: 3.3.3a

Graph of $T_2$

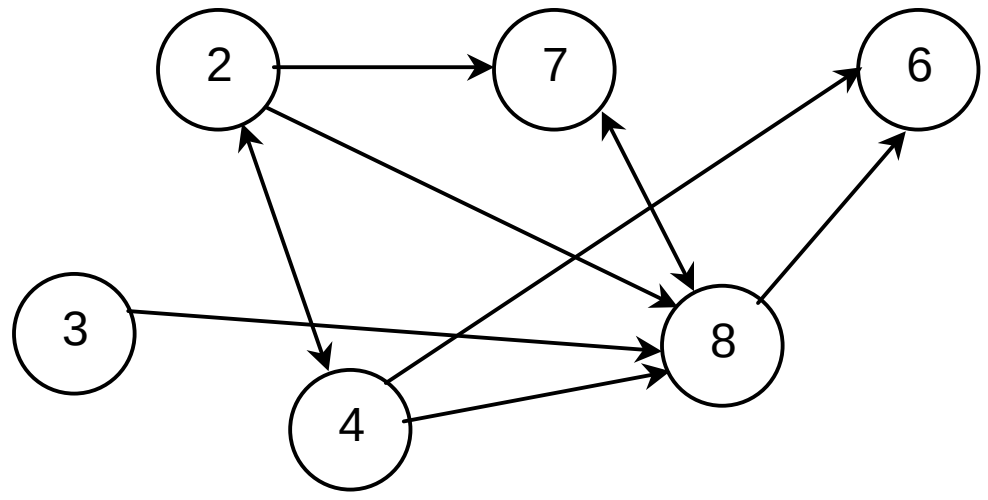

FIGURE: 3.3.3b

Graph of $T_3$

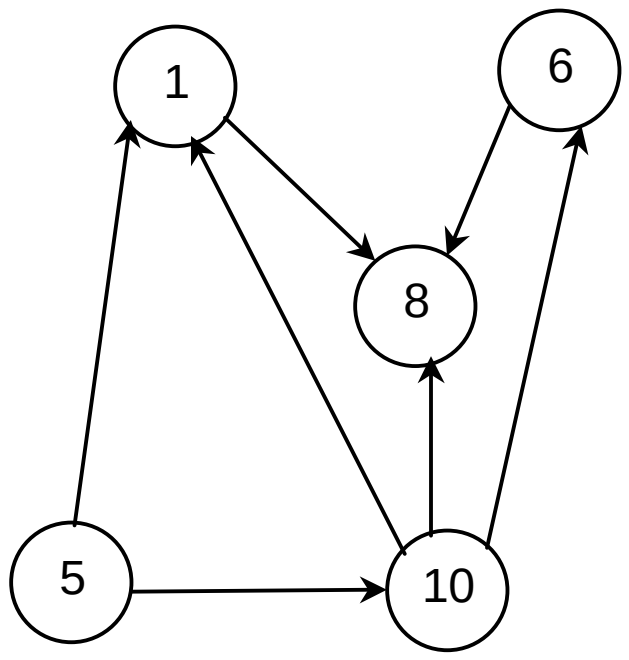

FIGURE: 3.3.3c

Now we give the related connection trimatrix which is a mixed square trimatrix. $A = A_1 \cup A_2 \cup A_3$

$$\begin{array}{c|cccccc} & 1 & 2 & 3 & 4 & 5 & 9 \\ \hline 1 & 0 & 0 & 0 & 0 & 0 & 1 \\ 2 & 0 & 0 & 0 & 1 & 0 & 0 \\ 3 & 1 & 0 & 0 & 1 & 0 & 0 \\ 4 & 0 & 1 & 0 & 0 & 0 & 0 \\ 5 & 1 & 0 & 1 & 1 & 0 & 0 \\ 9 & 0 & 0 & 0 & 0 & 0 & 0 \end{array} \cup \begin{array}{c|cccccc} & 2 & 3 & 4 & 6 & 7 & 8 \\ \hline 2 & 0 & 0 & 1 & 0 & 1 & 1 \\ 3 & 0 & 0 & 0 & 0 & 0 & 1 \\ 4 & 1 & 0 & 0 & 1 & 0 & 1 \\ 6 & 0 & 0 & 0 & 0 & 0 & 0 \\ 7 & 0 & 0 & 0 & 0 & 0 & 1 \\ 8 & 0 & 0 & 0 & 0 & 1 & 0 \end{array}$$

$$\cup \begin{array}{c|ccccc} & 1 & 5 & 6 & 8 & 10 \\ \hline 1 & 0 & 0 & 0 & 1 & 0 \\ 5 & 1 & 0 & 0 & 0 & 1 \\ 6 & 0 & 0 & 0 & 1 & 0 \\ 8 & 0 & 0 & 0 & 0 & 0 \\ 10 & 0 & 1 & 1 & 1 & 0 \end{array}.$$

Now we study the resultant of any state trivector on the dynamical system A.

Let X = (1 0 0 0 0 0) ∪ (1 0 0 0 0 0) ∪ (1 0 0 0 0 0) be the instantaneous state trivector in which the nodes (1) and (2) are in the on state and all other nodes are in the off state.

The effect of X on A is given by

$$\begin{aligned} XA &= (X_1 \cup X_2 \cup X_3)\,(A_1 \cup A_2 \cup A_3) \\ &= X_1 A_1 \cup X_2 A_2 \cup X_3 A_3 \\ &= (0\ 0\ 0\ 0\ 0\ 1) \cup (0\ 0\ 1\ 0\ 1\ 1) \cup (0\ 0\ 0\ 1\ 0). \end{aligned}$$

After updating the resultant state vector we get

$$\begin{aligned} Y &= (10\ 0\ 0\ 0\ 1) \cup (1\ 0\ 1\ 0\ 1\ 1) \cup (1\ 0\ 0\ 1\ 0) \\ &= Y_1 \cup Y_2 \cup Y_3. \end{aligned}$$

The effect of Y on the dynamical system A is given by

$$\begin{aligned} YA &= (Y_1 \cup Y_2 \cup Y_3)\,(A_1 \cup A_2 \cup A_3) \\ &= Y_1 A_1 \cup Y_2 A_2 \cup Y_3 A_3 \end{aligned}$$

$= \quad (0\ 0\ 0\ 0\ 0\ 1) \cup (1\ 0\ 1\ 1\ 2\ 2\ ) \cup (0\ 0\ 0\ 1\ 0)$

After thresholding and updating

$$= \quad (1\ 0\ 0\ 0\ 0\ 1) \cup (1\ 0\ 1\ 1\ 1\ 1) \cup (1\ 0\ 0\ 1\ 0)$$
$$= \quad Z = Z_1 \cup Z_2 \cup Z_3.$$

Now the first and the last trirow vectors are fixed tripoints as they repeat so we only work with the 2nd row vector.

$$ZA \quad = \quad Z_1\ A_1 \cup Z_2\ A_2 \cup Z_3\ A_3$$
$$= \quad Z_1\ A_1 \cup (1\ 0\ 1\ 1\ 2\ 3) \cup Z_3\ A_3$$
$$= \quad (1\ 0\ 0\ 0\ 0\ 1)\ (1\ 0\ 1\ 1\ 2\ 3) \cup (1\ 0\ 0\ 1\ 0)$$

after updating and thresholding we get

$$ZA \quad = \quad (10\ 0\ 0\ 0\ 0\ 1) \cup (1\ 0\ 1\ 1\ 1\ 1) \cup (1\ 0\ 0\ 1\ 0)$$

which is a fixed tripoint of the dynamical system. The on state of $P_1$ and $P_2$ imply the following

- Supply of quality seeds and its related cost
- The usage of fertilizers banned fertilizers lead to untold health problems
- Cultivating of quality crops
- Loan facilities offered by financial institutions / government
- Obtaining irrigational and electrical facilities from government
- Usage of technological equipments and
- Updation of cultivation methods.

Thus we can see the effect of any desired state vector on the dynamical system and derive the conclusions based on the analysis. It is still important to note that this system can give the impact on blocks or as a whole.

## 3.4 Definition and Illustration of Fuzzy Relational Maps (FRMs)

In this section, we introduce the notion of Fuzzy relational maps (FRMs); they are constructed analogous to FCMs described and discussed in the earlier sections. In FCMs we promote the correlations between causal associations among concurrently active units. But in FRMs we divide the very causal associations into two disjoint units, for example, the relation between a teacher and a student or relation between an employee or employer or a relation between doctor and patient and so on. Thus for us to define a FRM we need a domain space and a range space which are disjoint in the sense of concepts. We further assume no intermediate relation exists within the domain elements or node and the range spaces elements. The number of elements in the range space need not in general be equal to the number of elements in the domain space.

Thus throughout this section we assume the elements of the domain space are taken from the real vector space of dimension n and that of the range space are real vectors from the vector space of dimension m (m in general need not be equal to n). We denote by R the set of nodes $R_1,\ldots, R_m$ of the range space, where $R = \{(x_1,\ldots, x_m) \mid x_j = 0 \text{ or } 1 \}$ for $j = 1, 2,\ldots, m$. If $x_i = 1$ it means that the node $R_i$ is in the on state and if $x_i = 0$ it means that the node $R_i$ is in the off state. Similarly D denotes the nodes $D_1, D_2,\ldots, D_n$ of the domain space where $D = \{(x_1,\ldots, x_n) \mid x_j = 0 \text{ or } 1\}$ for $i = 1, 2,\ldots, n$. If $x_i = 1$ it means that the node $D_i$ is in the on state and if $x_i = 0$ it means that the node $D_i$ is in the off state.

Now we proceed on to define a FRM.

**DEFINITION 3.4.1:** *A FRM is a directed graph or a map from D to R with concepts like policies or events etc, as nodes and causalities as edges. It represents causal relations between spaces D and R .*

*Let $D_i$ and $R_j$ denote that the two nodes of an FRM. The directed edge from $D_i$ to $R_j$ denotes the causality of $D_i$ on $R_j$ called relations. Every edge in the FRM is weighted with a number in the set {0, ±1}. Let $e_{ij}$ be the weight of the edge $D_iR_j$, $e_{ij} \in \{0, \pm 1\}$. The weight of the edge $D_i R_j$ is positive if increase in $D_i$ implies increase in $R_j$ or decrease in $D_i$ implies decrease in $R_j$, i.e., causality of $D_i$ on $R_j$ is 1. If $e_{ij}$ = 0, then $D_i$ does not have any effect on $R_j$. We do not discuss the cases when increase in $D_i$ implies decrease in $R_j$ or decrease in $D_i$ implies increase in $R_j$.*

**DEFINITION 3.4.2:** *When the nodes of the FRM are fuzzy sets then they are called fuzzy nodes. FRMs with edge weights {0, ±1} are called simple FRMs.*

**DEFINITION 3.4.3:** *Let $D_1, \ldots, D_n$ be the nodes of the domain space D of an FRM and $R_1, \ldots, R_m$ be the nodes of the range space R of an FRM. Let the matrix E be defined as E = ($e_{ij}$) where $e_{ij}$ is the weight of the directed edge $D_iR_j$ (or $R_jD_i$), E is called the relational matrix of the FRM.*

*Note:* It is pertinent to mention here that unlike the FCMs the FRMs can be a rectangular matrix with rows corresponding to the domain space and columns corresponding to the range space. This is one of the marked differences between FRMs and FCMs.

**DEFINITION 3.4.4:** *Let $D_1, \ldots, D_n$ and $R_1, \ldots, R_m$ denote the nodes of the FRM. Let $A = (a_1, \ldots, a_n)$, $a_i \in \{0, \pm 1\}$. A is called the instantaneous state vector of the domain space and it denotes the on-off position of the nodes at any instant. Similarly let $B = (b_1, \ldots, b_m)$; $b_i \in \{0, \pm 1\}$. B is called instantaneous state vector of the range space and it denotes the on-off position of the nodes at any instant $a_i = 0$ if $a_i$ is off and $a_i = 1$ if $a_i$ is on for i= 1, 2,…, n. Similarly, $b_i = 0$ if $b_i$ is off and $b_i = 1$ if $b_i$ is on, for i= 1, 2,…, m.*

**DEFINITION 3.4.5:** *Let $D_1, \ldots, D_n$ and $R_1, \ldots, R_m$ be the nodes of an FRM. Let $D_iR_j$ (or $R_j D_i$) be the edges of an*

*FRM, $j = 1, 2,\ldots, m$ and $i= 1, 2,\ldots, n$. Let the edges form a directed cycle. An FRM is said to be a cycle if it posses a directed cycle. An FRM is said to be acyclic if it does not posses any directed cycle.*

**DEFINITION 3.4.6:** *An FRM with cycles is said to be an FRM with feedback.*

**DEFINITION 3.4.7:** *When there is a feedback in the FRM, i.e. when the causal relations flow through a cycle in a revolutionary manner, the FRM is called a dynamical system.*

**DEFINITION 3.4.8:** *Let $D_i$ $R_j$ (or $R_j$ $D_i$), $1 \leq j \leq m$, $1 \leq i \leq n$. When $R_i$ (or $D_j$ ) is switched on and if causality flows through edges of the cycle and if it again causes $R_i$ (or $D_j$ ), we say that the dynamical system goes round and round. This is true for any node $R_j$ (or $D_i$ ) for $1 \leq i \leq n$, (or $1 \leq j \leq m$). The equilibrium state of this dynamical system is called the hidden pattern.*

**DEFINITION 3.4.9:** *If the equilibrium state of a dynamical system is a unique state vector, then it is called a fixed point. Consider an FRM with $R_1, R_2,\ldots, R_m$ and $D_1, D_2,\ldots, D_n$ as nodes.*

For example, let us start the dynamical system by switching on $R_1$ (or $D_1$). Let us assume that the FRM settles down with $R_1$ and $R_m$ (or $D_1$ and $D_n$) on, i.e. the state vector remains as (1, 0, …, 0, 1) in R (or 1, 0, 0, … , 0, 1) in D), This state vector is called the fixed point.

**DEFINITION 3.4.10:** *If the FRM settles down with a state vector repeating in the form*

$$A_1 \to A_2 \to A_3 \to \ldots \to A_i \to A_1 \text{ (or } B_1 \to B_2 \to \ldots \to B_i \to B_1)$$

*then this equilibrium is called a limit cycle.*

**METHODS OF DETERMINING THE HIDDEN PATTERN**

Let $R_1, R_2, \ldots, R_m$ and $D_1, D_2, \ldots, D_n$ be the nodes of a FRM with feedback. Let E be the relational matrix. Let us find a hidden pattern when $D_1$ is switched on i.e. when an input is given as vector $A_1 = (1, 0, \ldots, 0)$ in $D_1$, the data should pass through the relational matrix E. This is done by multiplying $A_1$ with the relational matrix E. Let $A_1E = (r_1, r_2, \ldots, r_m)$, after thresholding and updating the resultant vector we get $A_1 E \in R$. Now let $B = A_1E$ we pass on B into $E^T$ and obtain $BE^T$. We update and threshold the vector $BE^T$ so that $BE^T \in D$. This procedure is repeated till we get a limit cycle or a fixed point.

**DEFINITION 3.4.11:** *Finite number of FRMs can be combined together to produce the joint effect of all the FRMs. Let $E_1, \ldots, E_p$ be the relational matrices of the FRMs with nodes $R_1, R_2, \ldots, R_m$ and $D_1, D_2, \ldots, D_n$, then the combined FRM is represented by the relational matrix $E = E_1 + \ldots + E_p$.*

Now we give a simple illustration of a FRM, for more about FRMs please refer [151, 165, 166].

Now we proceed on to give how bigraphs are used in FRMs. To make the book self contained one we have recalled the basic notion of FRMs. Now we define fuzzy relational bimaps.

## 3.5 Fuzzy Relational Bimaps and their Applications

In this section we for the first time introduce the new model fuzzy relational bimaps and illustrate it in the use of real valued problems.

**DEFINITION 3.5.1:** *A Fuzzy Relational bimaps (FRBM) is a directed bigraph or a pair of maps from $D^1$ to $R^1$ and $D^2$ to $R^2$ with concepts like policies or events etc as nodes or* *concepts and causalities as edges. It represents causal relations between spaces $D^i$ and $R^i$ i = 1, 2.*

*The directed bigraph may be disjoint or have a vertex in common or a edge in common or have a subgraph in common. All these depends on the problem model under study.*

*Let $D_i^K$ and $R_j^K$ (K = 1, 2) denote a pair of nodes of an FRBM. The directed edge $D_i^K$ to $R_j^K$ denotes the causality of $D_i^K$ and $R_j^K$ called relations. Every edge in the FRBM is weighted with a number in the set {0, ±1}. Let $e_{ij}^k$ be the weight of the edge $D_i^K$ $R_j^K$ (K = 1, 2) ∈ {0, 1, -1}. The weight of the edge $D_i^K$ $R_j^K$ is positive if increase in $D_i^K$ implies increase in $R_j^K$ or decrease in $D_L^K$ implies decrease in $R_j^K$ (K = 1, 2) i.e. the causality of $D_i^K$ on $R_j^K$ is 1.i.e., $e_{ij}^K$ = 1.*

*If $e_{ij}^k$ = 0 then $D_i^K$ does not have any effect on $R_j^K$. If $e_{ij}^k$ < 0, that is the cases when increases in $D_i^K$ implies decrease in $R_j^K$ or decrease in $D_i^K$ implies increase in $R_j^K$, (K = 1, 2).*

*Let $D_1^K,\ldots,D_{n_k}^K$ and $R_1^K,\ldots,R_{m_k}^K$ denote the nodes of the FRBM, here K = 1, 2. Let $A_1^K = \left(a_1^K,\ldots,a_n^K\right)$ (K = 1, 2) $a_i^K \in$ {0, 1, –1}. $A^1 \cup A^2$ is called the instantaneous state bivector of the domain space and it denotes the on-off position of the nodes at any instant. Similarly let $B^k$ = $\left(b_1^K,\ldots,b_{m_k}^K\right)$, $b_i^K \in$ {0, 1, -1}. $B^K$ is called the instantaneous state vector. Infact $B^1 \cup B^2$ is called the instantaneous state bivector of the range space and denotes the on-off position of the nodes at any instant $a_i^K$ = 0 if $a_i^K$ is off (K = 1,2);*

*$a_i^K$ = 1 if $a_i^K$ is on (K = 1, 2) and i = 1, 2,…, $n_k$. Similarly $b_j^k$ = 0 if $b_j^k$ is off and $b_j^k$ = 1 if $b_j^k$ is on for j = 1, 2,.., $m_k$ K = 1, 2.*

Let $D_i^K \, R_j^K \left(R_j^K \, D_i^K\right)$, K = 1, 2 denote the edges of an FRBM; j = 1, 2, …, $m_k$ and i =1, 2, …, $n_k$. Let the edges of the bigraph form a directed bicycle. An FRBM is said to be cyclic if it possesses a directed cycle. An FRBM is said to be acyclic if it does not posses any directed cycle. An FRBM with cycles is said to an FRBM with a feed back. When there is a feed back in the FRBM i.e. when the causal relations flow through a cycle in a revolutionary manner the FRBM is a dynamical system. When $R_i^K \left(\text{or } D_j^K\right)$ K = 1, 2 in the FRBM is switched on and if causality flows through the edges of the cycle and if it again causes $R_i$ (or $D_j$) we say that the dynamical system goes round and round. This is true for any node $R_j^K$ (or $D_i^K$) , K = 1, 2: 1 ≤ j ≤ $m_k$ and 1 ≤ i ≤ $n_k$. The equilibrium state of this dynamical system is called the bihidden pattern of the FRBM. If the equilibrium state of a dynamical system is a unique state bivector then it is called a fixed bipoint. If the FRBM settles down with a state bivector repeating in the form then their equilibrium is called a limit bicycle.

Now we illustrate the model. FRBMs has two domain spaces and two range spaces which may or may not have common nodes.

***Examples 3.5.1:*** When both the spaces are identical, we illustrate by example in the case of employee-employer relationship model.

The employee- employer relationship is an intricate one. For the employees expect to achieve performance in quality and good production in order to earn profit, on the other hand employees need good pay with all possible allowances. Here we have taken two experts opinion which forms a bigraph. This bigraph is disjoint and both of them are bipartite so a bipartite bigraph.

Suppose the concepts / nodes of domain are taken as $D_1$, …, $D_8$ pertaining to the employee.

- $D_1$ – Pay with allowances and bonus to the employee
- $D_2$ – Only pay to the employee
- $D_3$ – Pay with allowances to the employee
- $D_4$ – Best performance by the employee
- $D_5$ – Average performance by the employee
- $D_6$ – Poor performance by the employee
- $D_7$ – Employee works for more number of hours
- $D_8$ – Employer works for less number of hours.

$D_1$,.., $D_8$ are taken as the employee space i.e. as domain space.

We have taken 5 concepts related with the employer as suggested by the expert. These concepts form the range space which is listed below:

- $R_1$ – Maximum profit to the employer
- $R_2$ – Only profit to the employer
- $R_3$ – Neither profit nor loss to the employer
- $R_4$ – Loss to the employer
- $R_5$ – Heavy loss to the employer.

Now using the above set of nodes as the domain and range space the directed bigraph given by the two experts are

$$G = G_1 \cup G_2.$$

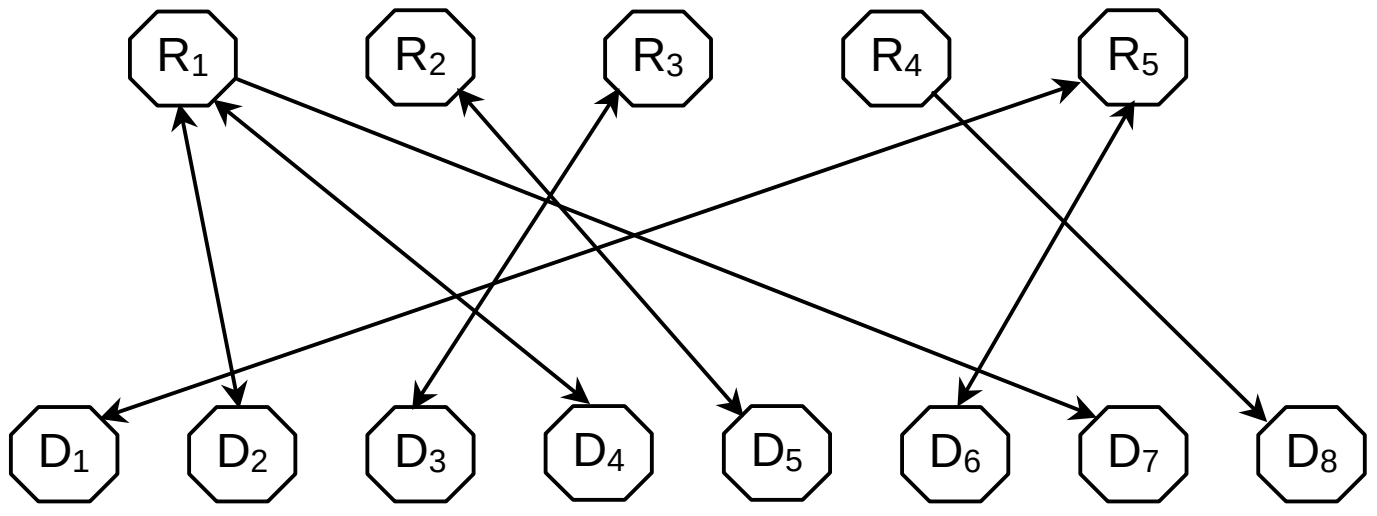

FIGURE: 3.5.1a

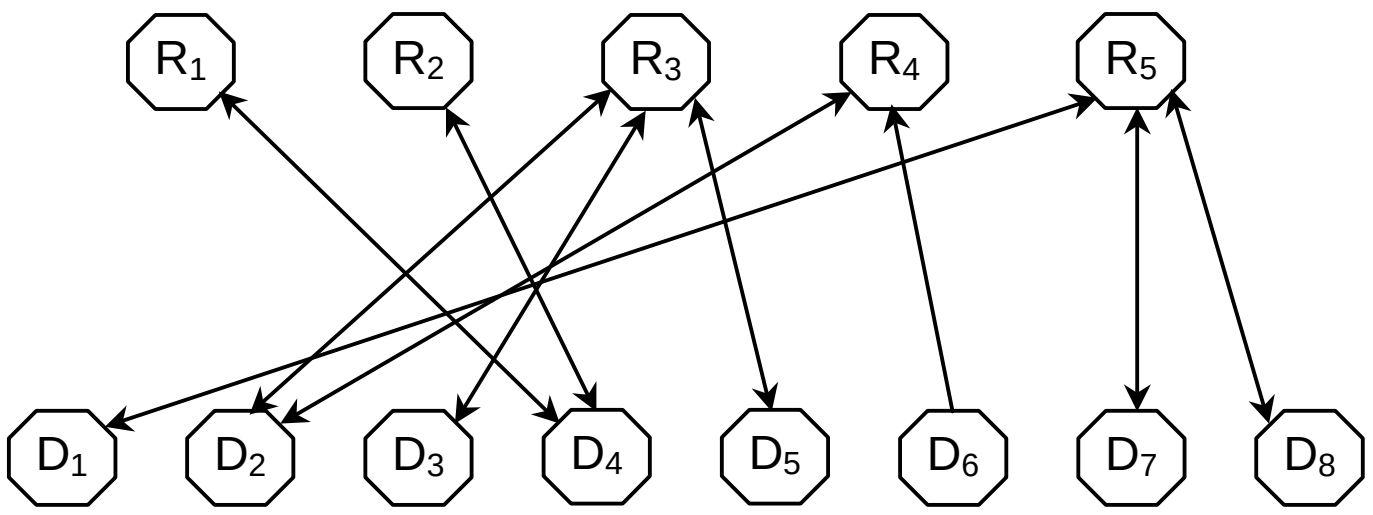

FIGURE: 3.5.1b

It is first noted that the bigraph $G = G_1 \cup G_2$ is a disconnected bigraph both $G_1$ and $G_2$ and distinct and both $G_1$ and $G_2$ are bipartite graphs so $G = G_1 \cup G_2$ is a bipartite bigraph. Now we find the related connection bimatrix, which is a rectangular $8 \times 5$ bimatrix $M = M_1 = M_1 \cup M_2$

$$M = \begin{bmatrix} 0 & 0 & 0 & 0 & 1 \\ 1 & 0 & 0 & 0 & 0 \\ 0 & 0 & 1 & 0 & 0 \\ 1 & 0 & 0 & 0 & 0 \\ 0 & 1 & 0 & 0 & 0 \\ 0 & 0 & 0 & 0 & 1 \\ 1 & 0 & 0 & 0 & 0 \\ 0 & 0 & 0 & 1 & 0 \end{bmatrix} \cup \begin{bmatrix} 0 & 0 & 0 & 0 & 1 \\ 0 & 0 & 1 & 1 & 0 \\ 0 & 0 & 1 & 0 & 0 \\ 1 & 1 & 0 & 0 & 0 \\ 0 & 0 & 1 & 0 & 0 \\ 0 & 0 & 0 & 1 & 0 \\ 0 & 0 & 0 & 0 & 1 \\ 0 & 0 & 0 & 0 & 1 \end{bmatrix}.$$

Let us now input a state bivector

$$\begin{aligned} X &= X_1 \cup X_2. \\ &= (1\ 0\ 0\ 0\ 0\ 0\ 0\ 0) \cup (1\ 0\ 0\ 0\ 0\ 0\ 0\ 0) \end{aligned}$$

i.e. the employer is paid with allowance alone is in the on state and all nodes are in the off state.

$$\begin{aligned} XM &= (X_1 \cup X_2)\ (M_1 \cup M_2) \\ &= X_1 M_1 \cup X_2 M_2 \\ &= (0\ 0\ 0\ 0\ 1) \cup (0\ 0\ 0\ 0\ 1). \end{aligned}$$

Let

$$Y = Y_1 \cup Y_2$$
$$= (0\ 0\ 0\ 0\ 1) \cup (0\ 0\ 0\ 0\ 1)$$

The effect of Y on $M^T$ gives

$$YM^T = (Y_1 \cup Y_2)\,(M_1 \cup M_2)^T$$
$$= Y_1\, M_1^T \cup Y_2\, M_2^T$$
$$= (1\ 0\ 0\ 0\ 0\ 1\ 0\ 0) \cup (1\ 0\ 0\ 0\ 0\ 0\ 1\ 1).$$

After thresholding and updating we get

$$Z = (1\ 0\ 0\ 0\ 0\ 1\ 0\ 0) \cup (1\ 0\ 0\ 0\ 0\ 0\ 1\ 1)$$
$$= Z_1 \cup Z_2.$$

Now we study the influence of the state bivector $Z_2$ on M

$$ZM = (Z_1 \cup Z_2)\,(M_1 \cup M_2)$$
$$= Z_1\, M_1 \cup Z_2\, M_2$$
$$= (0\ 0\ 0\ 0\ 1) \cup (0\ 0\ 0\ 0\ 3).$$

After thresholding we get the resultant as

$$T = (0\ 0\ 0\ 0\ 1) \cup (0\ 0\ 0\ 0\ 1)$$
$$= Y$$

giving a fixed bipoint.

Thus the bihidden pattern of the dynamical system is given by the binary bipair.

$$(\{(1\ 0\ 0\ 0\ 0\ 1\ 0\ 0) \cup (1\ 0\ 0\ 0\ 0\ 1\ 1\ 1)\},$$
$$\{(0\ 0\ 0\ 0\ 1) \cup (0\ 0\ 0\ 1\ 1\ )\}).$$

As in case of the FRMs the resultant state bivector can be interpreted. This makes one clearly see the opinion of the two experts on the given node need not be coincident which is very evident from the resultant binary pair.

Just we have seen the FRBM which is the disconnected bigraph. Now we give an illustration of the model which gives a connected bigraph. Before we proceed to give a FRBM model of special type we define a new type of bipartite bigraph.

**DEFINITION 3.5.2:** *Let $G = G_1 \cup G_2$ be a bigraph. Suppose $G = G_1 \cup G_2$ is a bipartite bigraph such that $G_1$ and $G_2$ are bipartite graphs; suppose $V_1^1, V_2^1$ are the two disjoint subsets by which the vertex set $V^1$ of $G_1$ is divided and if $V_1^2, V_2^2$ are the disjoint subsets of the vertex set $V^2$ of $G_2$.*
*If G is a disconnected bipartite bigraph we see the vertex sets of $G_1$ and $G_2$ have no common element. Suppose $V_2^1 = V_1^2$ (or $\left(V_2^1 = V_2^2\right)$ or $\left(V_1^1 = V_2^2\right)$ or $\left(V_1^1 = V_1^2\right)$ or in the mutually exclusive sense then we say the bipartite bigraph is a strongly biconnected bigraph.*

Now we give an illustration of a strongly biconnected bipartite bigraph before we give an associated model.

***Example 3.5.2:*** Let $G = G_1 \cup G_2$ be a bigraph given by the following figure.

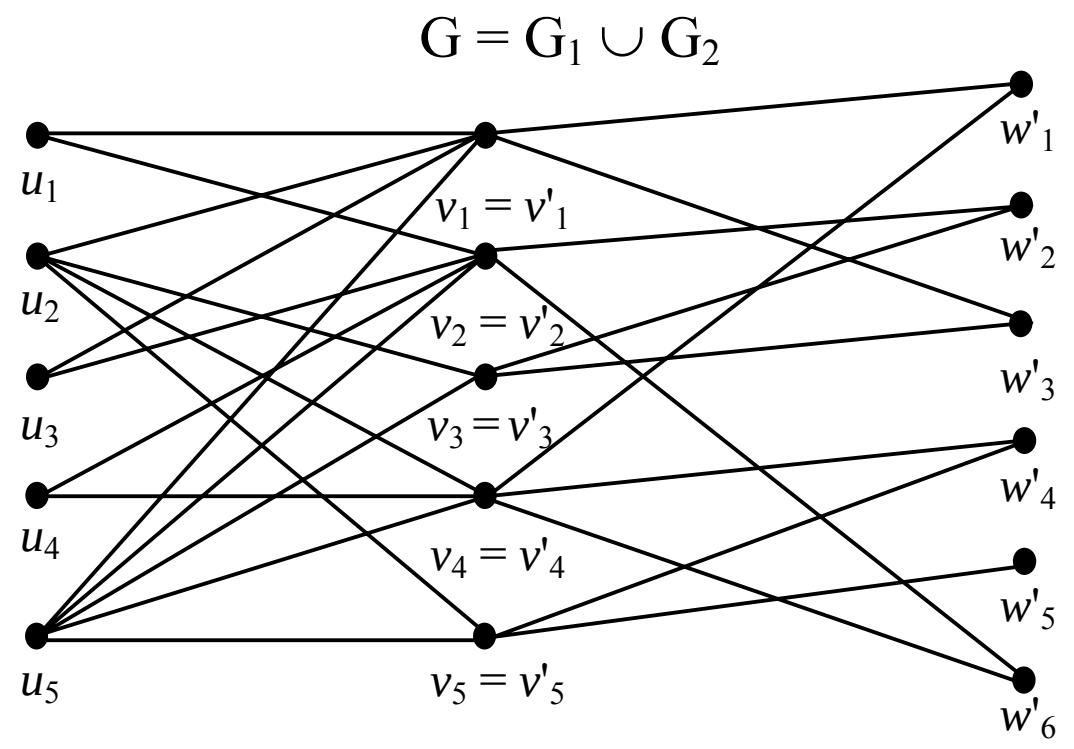

FIGURE: 3.5.2

This is a bipartite bigraph.
The vertex set of $G_1 = \{u_1, u_2, u_3, u_4, u_5, v_1, v_2, v_3, v_4, v_5\}$.

The vertex set of $G_2 = \{v_1, v'_1, v_2 = v'_2, v_3 = v'_3, v_4 = v'_4, v_5 = v'_5, w'_1, w'_2, w'_3, w'_4, w'_5, w'_6\}$.

Clearly both $G_1$ and $G_2$ are bipartite. $G = G_1 \cup G_2$ is a connected bipartite bigraph which is a strongly biconnected bipartite bigraph.

Now we illustrate a model when the resultant bigraph is a strongly biconnected bipartite bigraph.

To model this, the problem of child labour in India is taken. Though child labour problem is one of the major problems of India it is never given that importance because the children are not vote bankers. So the child labour bonded labour among children and misuse of these children who work as labourers are at rampant.

Now when we held discussions with human right activists, NGOs and other people sensitive to child labour problem we could form 3 sets of disjoint nodes one set of attribute connected with the government and its policies, one set of attributes connected with the children working as child labourers and the third set of attributes associated with the public awareness in support of child labour.

We just wish to state that when we say people / public who support child labour we mean people who use children as domestic servants, small scale industries like beedi factories and bangle factories match and crackers factories) who use children as labourers also take work from them for a very long time. Bonded labourers children who are forcefully used as child labourers as they belong to dalits that is due to bias. Also we include the sympathetic public who when come across such atrocities try to voluntarily come forward to help these deprived children.

Now we list the attributes associated with these three sets.

G – The attributes associated with government; {$G_1$, $G_2$, $G_3$ and $G_4$}.

| | | |
|---|---|---|
| $G_1$ | – | Children do not form vote bank |
| $G_2$ | – | Business men industrialists who practice child labour are the main source of vote Bank (In several cases they help with a lot of money as they hoard black money) |
| $G_3$ | – | Free and compulsory education for children |
| $G_4$ | – | No proper punishment by the government for the practice of child labour. |

C – The attributes associated with the children working as child labourers $C_1, C_2, \ldots, C_6$ and $C_7$.

| | | |
|---|---|---|
| $C_1$ | – | Abolition of child labour |
| $C_2$ | – | Uneducated parents |
| $C_3$ | – | School dropouts / never attended any school |
| $C_4$ | – | Social status of child labourer |
| $C_5$ | – | Poverty / source of livelihood |
| $C_6$ | – | Orphans runaways and parents are beggars or father in prison |
| $C_7$ | – | Habits like blue films, drink or smoke etc. |

P – Attributes associated with public awareness in support of child labour $P_1$, $P_2$, $P_3$, $P_4$ and $P_5$.

| | | |
|---|---|---|
| $P_1$ | – | Cheap and long hours of labour exploited from children |
| $P_2$ | – | Children as domestic servants who work for over 15 hours |
| $P_3$ | – | Sympathetic public |
| $P_4$ | – | Motivation by teachers to children to pursue education |
| $P_5$ | – | Perpetuating slavery and caste bias. |

Taking an experts opinion we give the directed bigraph relating the child labour and the government policies and child labour and the role of public in the following figure:

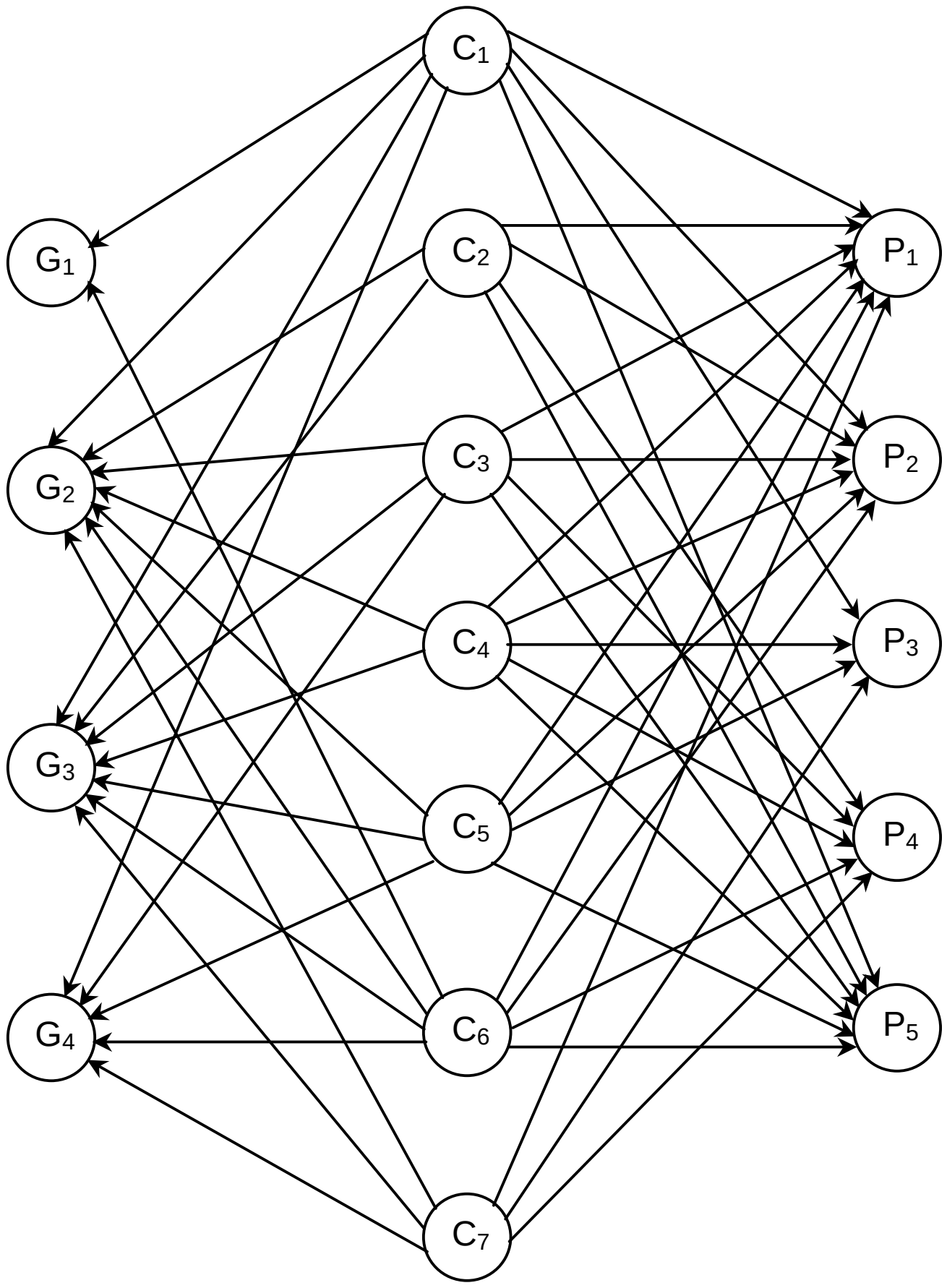

FIGURE 3.5.3

The bigraph which is a connected bipartite bigraph is denoted by $T = G \cup C \cup P$. By Observation we see the bigraph is a bipartite bigraph which is strongly biconnected.

The related connection bimatrix of the bigraph is

$$B = B_1 \cup B_2.$$

$$\begin{bmatrix} 1 & 0 & 0 & 0 & 0 & 1 & 0 \\ -1 & 1 & 1 & 1 & 1 & 1 & 1 \\ 1 & -1 & -1 & -1 & 1 & -1 & -1 \\ -1 & 0 & 1 & 0 & 1 & 1 & 1 \end{bmatrix} \cup \begin{bmatrix} -1 & -1 & 1 & 0 & -1 \\ -1 & 1 & 0 & -1 & 1 \\ 1 & 1 & 0 & -1 & 1 \\ -1 & -1 & 1 & 1 & -1 \\ 1 & 1 & 1 & 0 & 1 \\ 1 & 1 & 0 & -1 & 1 \\ 1 & 0 & -1 & -1 & 0 \end{bmatrix}$$

From these connection matrices of the FRBM one can easily find the related bipartite bigraphs. Suppose we consider the state bivector.

$$\begin{aligned} X &= (1\ 0\ 0\ 0) \cup (0\ 0\ 1\ 0\ 0\ 0\ 0) \\ &= X_1 \cup X_2 \\ XB &= (X_1 \cup X_2)\ (B_1 \cup B_2) \\ &= X_1 B_1 \cup X_2 B_2 \\ &= (1\ 0\ 0\ 0\ 0\ 1\ 0) \cup (1\ 1\ 0\ -1\ 1). \\ &= (1\ 0\ 0\ 0\ 0\ 1\ 0) \cup (1\ 1\ 0\ 0\ 1) \\ Y &= Y_1 \cup Y_2 \end{aligned}$$

$$\begin{aligned} YB^T &= Y_1B_1^T \cup Y_2B_2^T \\ &= (2\ 0\ 0\ 0) \cup (-3\ 1\ 3\ -3\ 3\ 3\ 1). \end{aligned}$$

After thresholding as they are bivector from the respective range spaces we have nothing to update we get the resultant bivector as

$$\begin{aligned} Z &= Z_1 \cup Z_2 \\ &= (1\ 0\ 0\ 0) \cup (0\ 1\ 1\ 0\ 1\ 1\ 1). \end{aligned}$$

The effect of Z on the FRBM is got by applying to B.

$$\begin{aligned} ZB &= Z_1 B_1 \cup Z_2 B_2 \\ &= (\ 1\ 0\ 0\ 0\ 0\ 1\ 0) \cup (3\ 4\ 0\ -4\ 4). \end{aligned}$$

After thresholding we get the resultant we get

$$A = (1\ 0\ 0\ 0\ 0\ 1\ 0) \cup (1\ 1\ 0\ 0\ 1).$$

Thus the bihidden pattern of the dynamical system is given by the fixed point which is a binary bipair.
$(\{(1\ 0\ 0\ 0) \cup (0\ 1\ 1\ 0\ 1\ 1\ 1)\}, \{(1\ 0\ 0\ 0\ 0\ 1\ 0) \cup (1\ 1\ 0\ 0\ 1)\}$
One can as in case of FRM interpret the bihidden pattern of the FRBM.

## 3.6 Fuzzy Relational Trimaps and their Applications

Now we proceed on to show how the notion of trigraphs can be associated with the fuzzy relational trimaps. Thus the connection trimatrix associated with fuzzy relational trimaps (FRTM) will be a trimatrix. In case of FRTMs also we can have disconnected FRTMs which will be a bipartite trigraph. Also FRTMs can also be strongly triconnected bipartite bigraphs. We just give an example of each before we go for the model representations.

***Example 3.6.1*:** Consider the trigraph $G = G_1 \cup G_2 \cup G_3$ given by the following figure.
This is a disconnected trigraph.

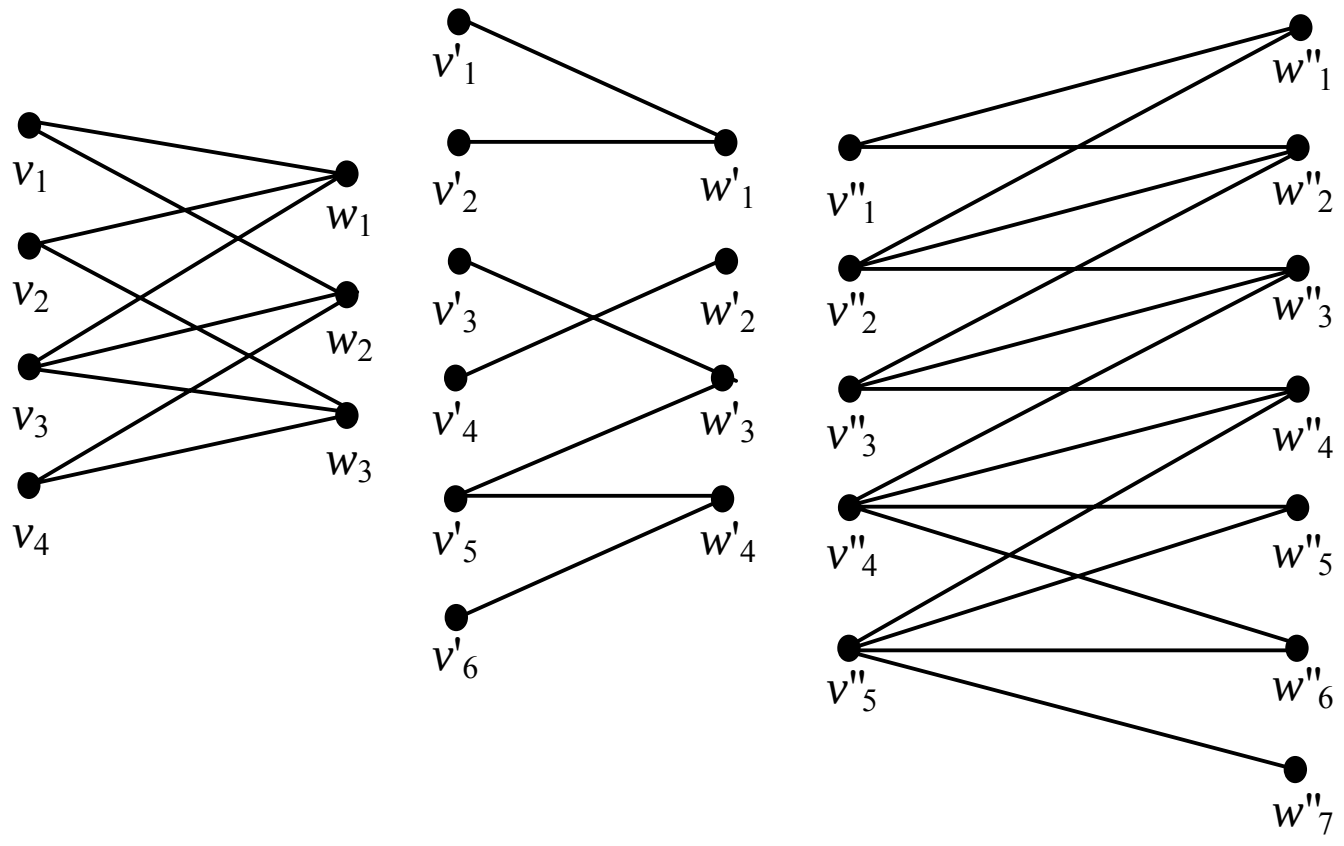

FIGURE: 3.6.1

$$G = G_1 \cup G_2 \cup G_3;$$

each of them is a bipartite graph.

Next we proceed on to give an example of a strongly triconnected bipartite trigraph.

***Example 3.6.2:*** Consider the trigraph which is triconnected bipartite trigraph $G = G_1 \cup G_2 \cup G_3$ given by the following figure 3.6.2.

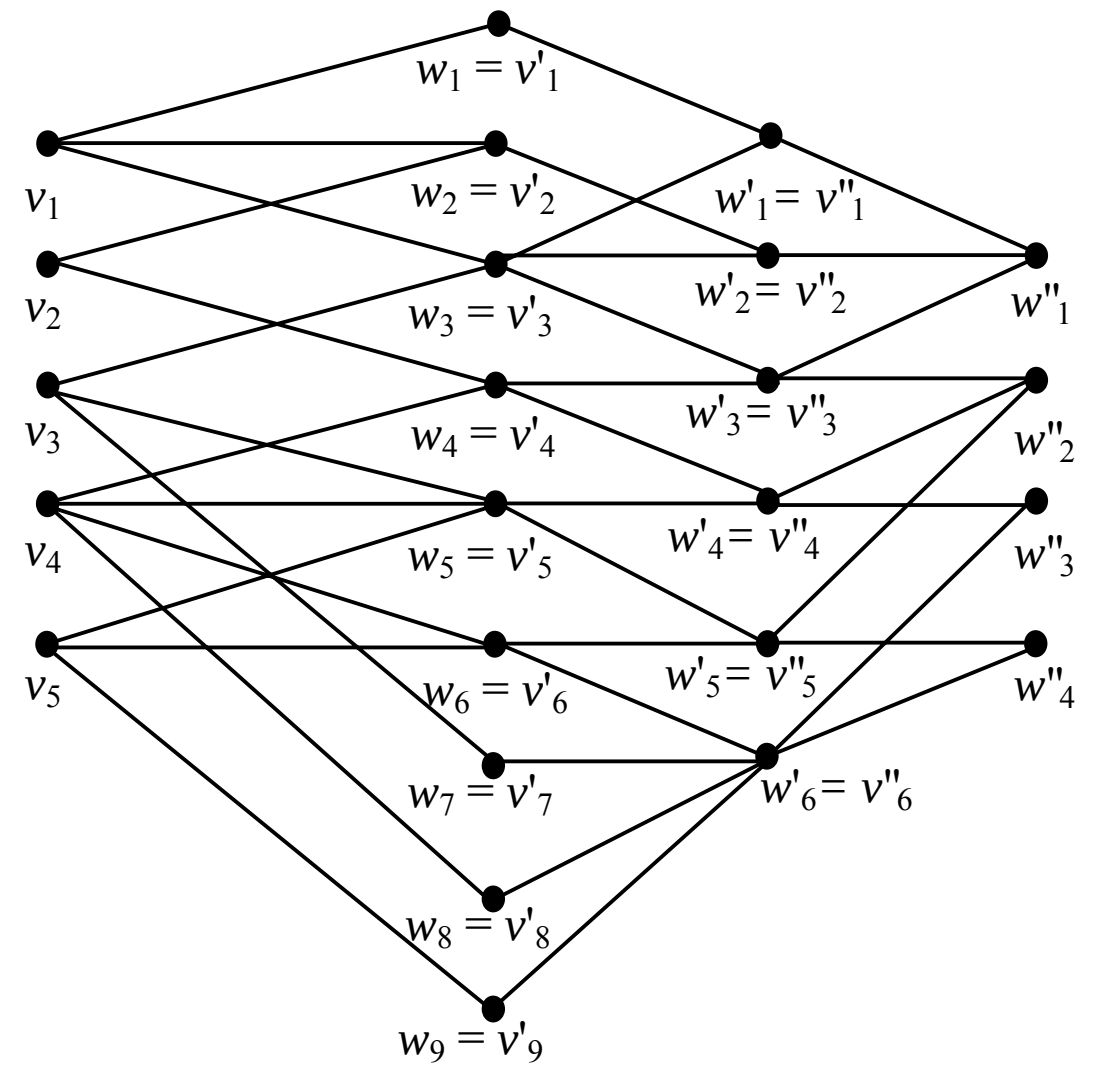

FIGURE: 3.6.2

The graph $G_3$ is given by the following figure.

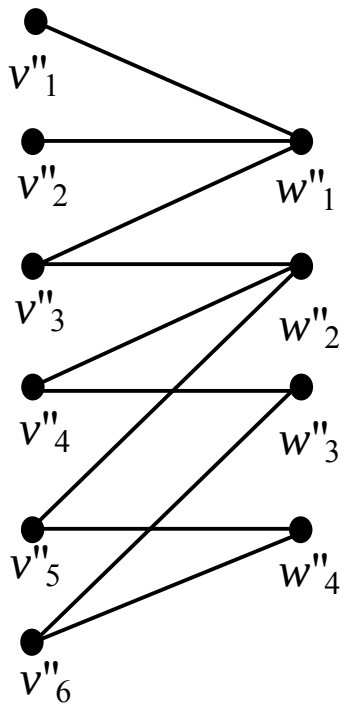

FIGURE: 3.6.2c

The graph of $G_1$ is given by the following figure.

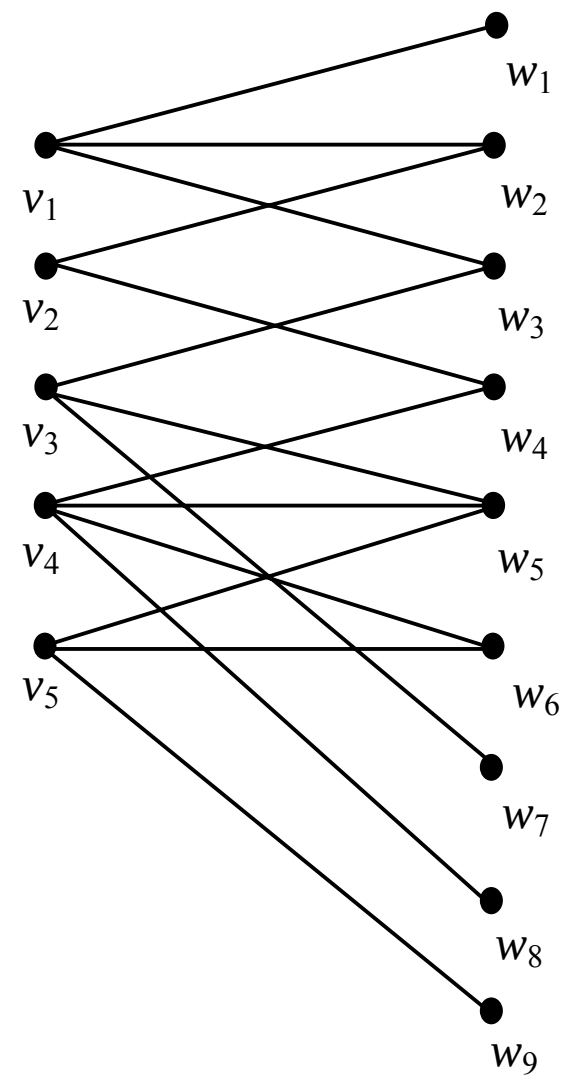

FIGURE: 3.6.2a

The graph of $G_2$ is given by the following figure.

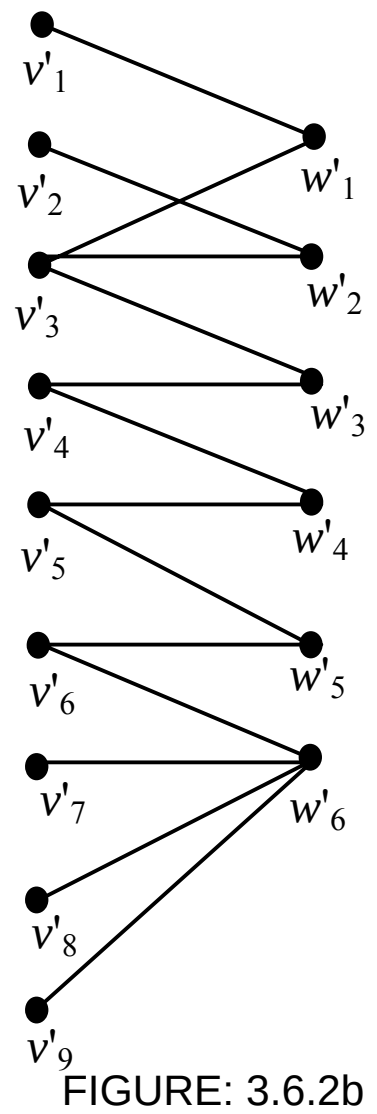

FIGURE: 3.6.2b

The associated trimatrix is given by $M = M_1 \cup M_2 \cup M_3$

$$\begin{bmatrix} 1 & 1 & 1 & 0 & 0 & 0 & 0 & 0 & 0 \\ 0 & 1 & 0 & 1 & 0 & 0 & 0 & 0 & 0 \\ 0 & 0 & 1 & 0 & 1 & 0 & 1 & 0 & 0 \\ 0 & 0 & 0 & 1 & 1 & 1 & 0 & 1 & 0 \\ 0 & 0 & 0 & 0 & 1 & 1 & 0 & 0 & 1 \end{bmatrix}$$

$$\cup \begin{bmatrix} 1 & 0 & 0 & 0 & 0 & 0 \\ 0 & 1 & 0 & 0 & 0 & 0 \\ 1 & 1 & 1 & 0 & 0 & 0 \\ 0 & 0 & 1 & 1 & 0 & 0 \\ 0 & 0 & 0 & 1 & 1 & 0 \\ 0 & 0 & 0 & 0 & 1 & 1 \\ 0 & 0 & 0 & 0 & 0 & 1 \\ 0 & 0 & 0 & 0 & 0 & 1 \\ 0 & 0 & 0 & 0 & 0 & 1 \end{bmatrix}$$

$$\cup \begin{bmatrix} 1 & 0 & 0 & 0 \\ 1 & 0 & 0 & 0 \\ 1 & 1 & 0 & 0 \\ 0 & 1 & 1 & 0 \\ 0 & 1 & 1 & 0 \\ 0 & 0 & 1 & 1 \end{bmatrix}.$$

Now we give an example of the use of disconnected trimatrix in the FRTM model. Already we have studied the employee-employer relationship and this was modeled using the FRBM and related connection bimatrix was a disconnected bipartite bimatrix. We have shown this in page

___. Now we take one more experts opinion so that we have a bipartite trimatrix.

The related bigraph using three experts opinion of the same model is given in the following.

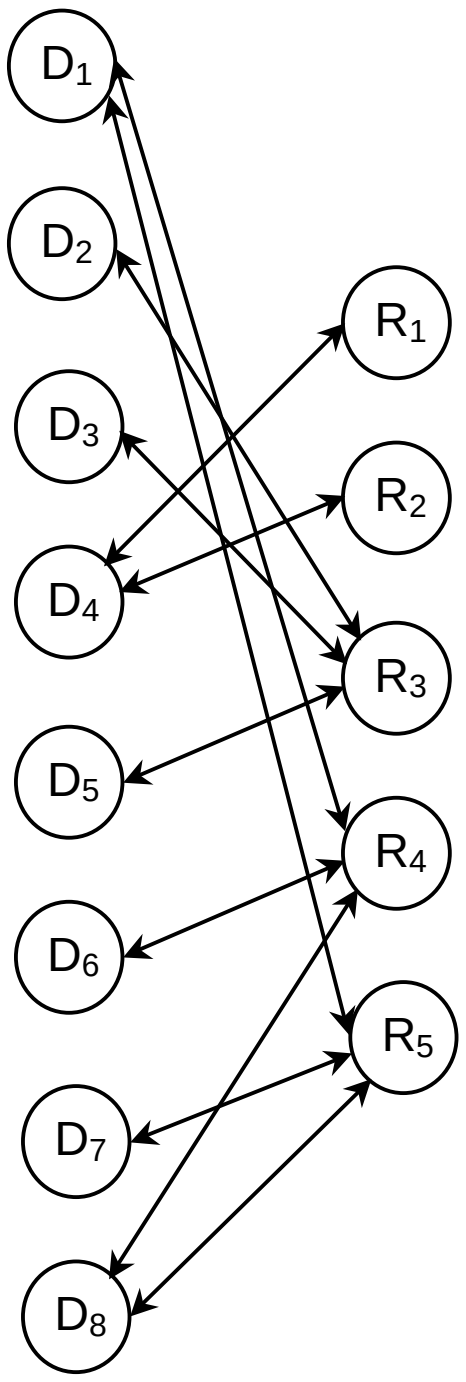

FIGURE: 3.6.3a

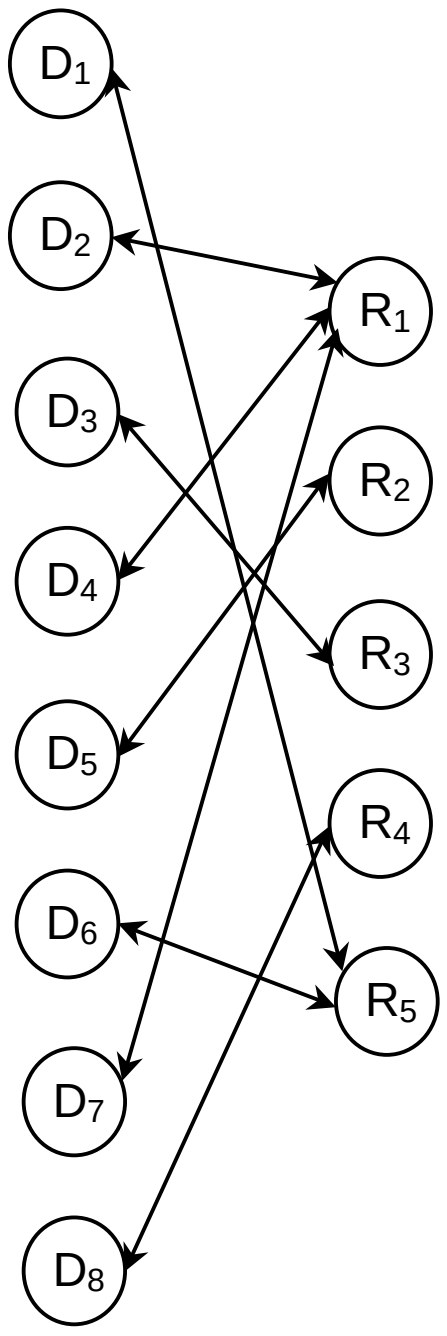

FIGURE: 3.6.3b

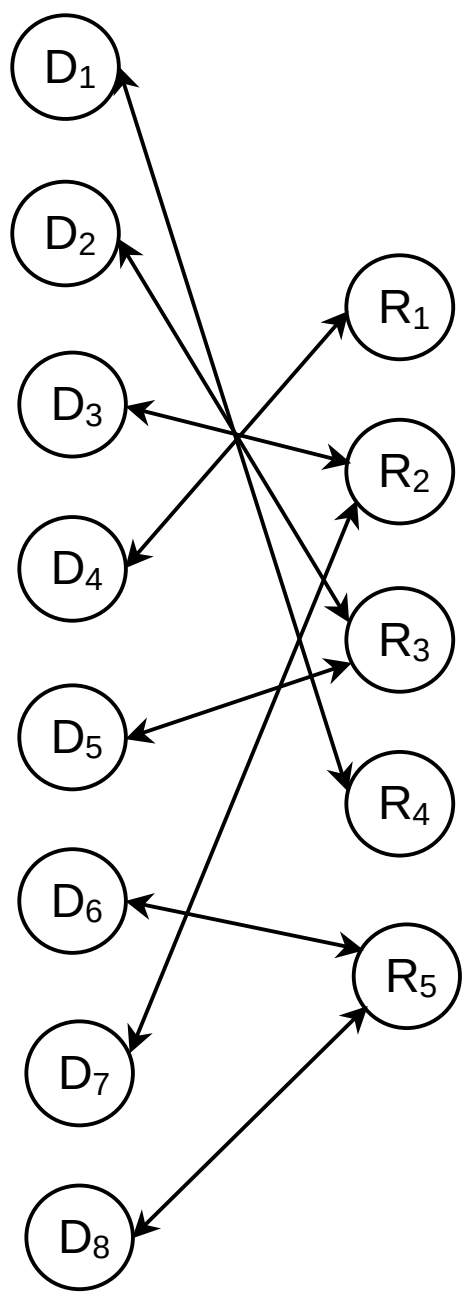

FIGURE: 3.6.3c

FIGURE: 3.6.3

$$T = T_1 \cup T_2 \cup T_3$$

Now we write down the related connection trimatrix

$$A = A_1 \cup A_2 \cup A_3$$

$$A = \begin{bmatrix} 0 & 0 & 0 & 1 & 1 \\ 0 & 0 & 1 & 0 & 0 \\ 0 & 0 & 1 & 0 & 0 \\ 1 & 1 & 0 & 0 & 0 \\ 0 & 0 & 1 & 0 & 0 \\ 0 & 0 & 0 & 1 & 0 \\ 0 & 0 & 0 & 0 & 1 \\ 0 & 0 & 0 & 1 & 1 \end{bmatrix} \cup \begin{bmatrix} 0 & 0 & 0 & 0 & 1 \\ 1 & 0 & 0 & 0 & 0 \\ 0 & 0 & 1 & 0 & 0 \\ 1 & 0 & 0 & 0 & 0 \\ 0 & 1 & 0 & 0 & 0 \\ 0 & 0 & 0 & 0 & 1 \\ 1 & 0 & 0 & 0 & 0 \\ 0 & 0 & 0 & 1 & 0 \end{bmatrix}$$

$$\cup \begin{bmatrix} 0 & 0 & 0 & 1 & 0 \\ 0 & 0 & 1 & 0 & 0 \\ 0 & 1 & 0 & 0 & 0 \\ 1 & 0 & 0 & 0 & 0 \\ 0 & 0 & 1 & 0 & 0 \\ 0 & 0 & 0 & 0 & 1 \\ 0 & 1 & 0 & 0 & 0 \\ 0 & 0 & 0 & 0 & 1 \end{bmatrix}$$

We see this trimatrix is a $8 \times 5$ rectangular trimatrix. Suppose one is interested in studying the effect of $D_8$ i.e. the employee works for less number of hours on the dynamical system A.

That is $D_8$ on the FRTM will simultaneously give the opinion of all the three experts when $D_8$ alone is in the on state.

Let

$$\begin{aligned} X &= (0\ 0\ 0\ 0\ 0\ 0\ 0\ 1) \cup (0\ 0\ 0\ 0\ 0\ 0\ 0\ 1) \\ &\quad \cup (\ 0\ 0\ 0\ 0\ 0\ 0\ 0\ 1) \\ &= X_1 \cup X_2 \cup X_3. \end{aligned}$$

The impact of X on the dynamical system A is given by

$$\begin{aligned} XA &= (X_1 \cup X_2 \cup X_3)\ (A_1 \cup A_2 \cup A_3) \\ &= X_1\ A_1 \cup X_2\ A_2 \cup X_3\ A_3 \\ &= (0\ 0\ 0\ 1\ 1) \cup (0\ 0\ 0\ 1\ 0) \cup (0\ 0\ 0\ 0\ 1) \\ &= Y \\ &= Y_1 \cup Y_2 \cup Y_3. \end{aligned}$$

$$\begin{aligned} YA^T &= (Y_1 \cup Y_2 \cup Y_3)\ (A_1 \cup A_2 \cup A_3)^t \\ &= Y_1\ A_1^t \cup Y_2\ A_2^t \cup Y_3\ A_3^t \end{aligned}$$

$$= (2\ 0\ 0\ 0\ 0\ 1\ 1\ 2) \cup (0\ 0\ 0\ 0\ 0\ 0\ 0\ 1) \cup (0\ 0\ 0\ 0\ 0\ 1\ 0\ 1)$$

After thresholding and updating we get the resultant as

$$\begin{aligned} Z &= (1\ 0\ 0\ 0\ 0\ 1\ 1\ 1) \cup (0\ 0\ 0\ 0\ 0\ 0\ 0\ 1) \cup (0\ 0\ 0\ 0\ 0\ 1\ 0\ 1) \\ &= Z_1 \cup Z_2 \cup Z_3 \end{aligned}$$

$$\begin{aligned} ZA &= (Z_1 \cup Z_2 \cup Z_3)\ (A_1 \cup A_2 \cup A_3) \\ &= Z_1\ A_1 \cup Z_2\ A_2 \cup Z_3\ A_3 \\ &= (0\ 0\ 0\ 3\ 3) \cup (0\ 0\ 0\ 1\ 0) \cup (0\ 0\ 0\ 0\ 2) \end{aligned}$$

After thresholding we get

$$U = (0\ 0\ 0\ 1\ 1) \cup (0\ 0\ 0\ 1\ 0) \cup (0\ 0\ 0\ 0\ 1).$$

Thus the trihidden pattern of the dynamical system is a fixed tripoint given by the bipair of the tristate vectors.

$(\{(1\ 0\ 0\ 0\ 0\ 1\ 1\ 1) \cup (0\ 0\ 0\ 0\ 0\ 0\ 0\ 1) \cup (0\ 0\ 0\ 0\ 0\ 1\ 0\ 1)\}, \{(0\ 0\ 0\ 1\ 1) \cup (0\ 0\ 0\ 1\ 0) \cup (0\ 0\ 0\ 0\ 1)\})$.

Thus we observe that the second experts opinion is some what very inconsistent for the second expert was the union leader but the third expert opinion is best for we see much logic when the employee works for less number of hours, poor performance by the employee and the company suffers a heavy loss. However all the three experts agree on the fact that if the employee works for less number of hours certainly there, will be "Heavy loss to the employer".

Now we give yet one more model of FRTM which gives a resultant matrix which is a strongly triconnected bipartite trigraph. For this we take up the problem of child labour given in pages 131-133 of this book.

The directed trigraph given by the expert is described in the following figure 3.6.4.

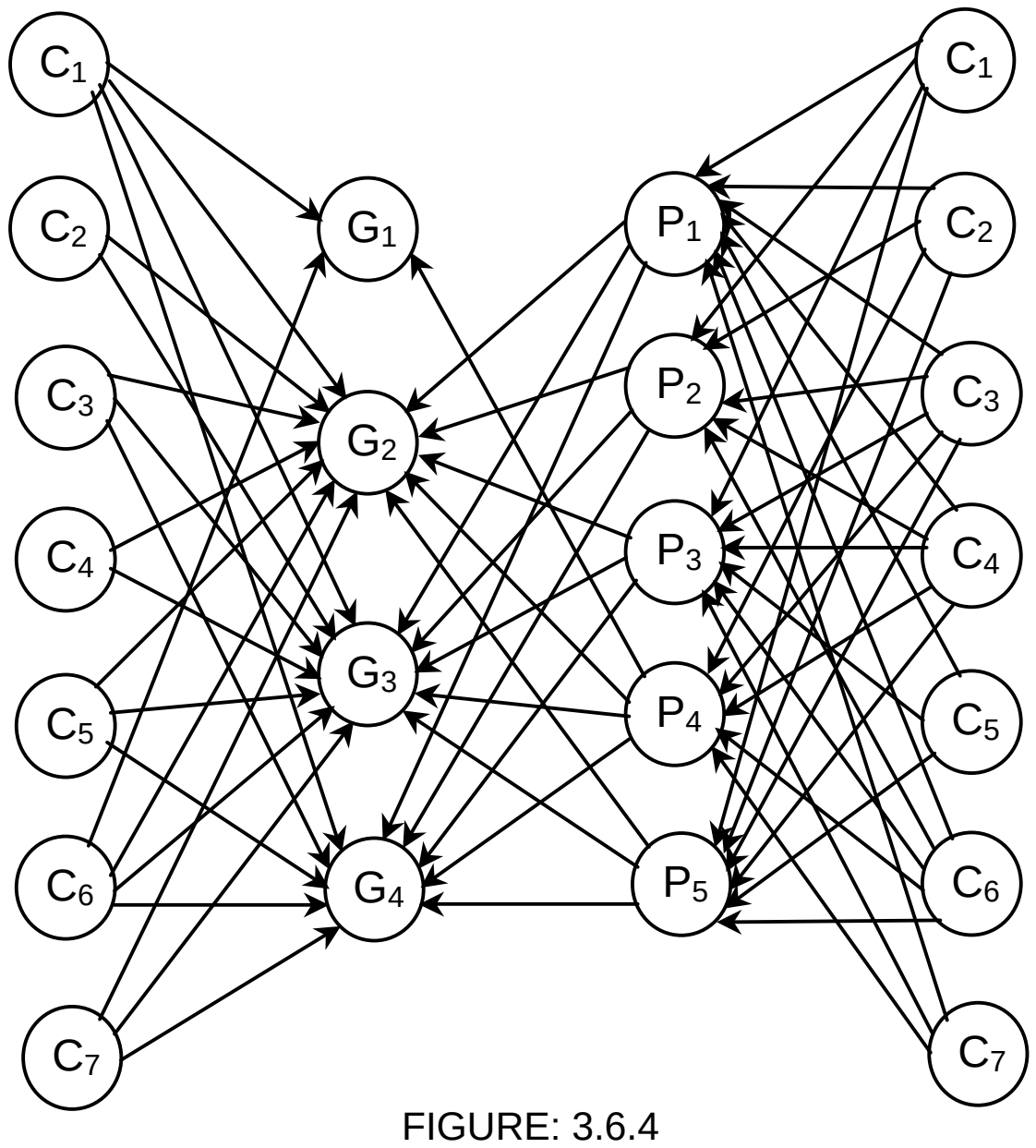

FIGURE: 3.6.4

It can be easily seen that the given trigraph is a strongly triconnected trigraph which is a bipartite trigraph. Suppose $T = T_1 \cup T_2 \cup T_3$ the vertex set of $T_1 = \{(C_1, C_2, C_3, C_4, C_5, C_6, C_7), (G_1, G_2, G_3, G_4)\}$. The vertex set of $T_2 = \{(G_1, G_2, G_3, G_4), (P_1, P_2, P_3, P_4, P_5)\}$. The vertex set of $T_3 = \{(P_1, P_2, P_3, P_4, P_5), (C_1, C_2, C_3, C_4, C_5, C_6, C_7)\}$.

It is easily verified that all the three graphs $T_1$, $T_2$ and $T_3$ are bipartite. Now we give the connection trimatrix associated with the FRTM, $S = S_1 \cup S_2 \cup S_3$.

$$\begin{bmatrix} 1 & -1 & 1 & -1 \\ 0 & 1 & -1 & 0 \\ 0 & 1 & -1 & 1 \\ 0 & 1 & -1 & 0 \\ 0 & 1 & 1 & 1 \\ 1 & 1 & -1 & 1 \\ 0 & 1 & -1 & 1 \end{bmatrix} \cup \begin{bmatrix} 0 & 0 & 0 & 1 & 0 \\ 1 & 1 & -1 & -1 & 1 \\ 1 & -1 & 1 & 1 & 1 \\ 1 & 1 & -1 & 1 & 1 \end{bmatrix} \cup$$

$$\begin{bmatrix} -1 & -1 & 1 & -1 & 1 & 1 & 1 \\ -1 & 1 & 1 & -1 & 0 & 1 & 0 \\ 1 & 0 & -1 & 1 & 1 & -1 & -1 \\ 0 & -1 & -1 & 1 & 0 & -1 & -1 \\ -1 & 1 & 1 & -1 & 1 & 1 & 0 \end{bmatrix}.$$

Thus the trimatrix is a mixed rectangular trimatrix. Now just to show how the FRTM works we consider a state vector

$$X = X_1 \cup X_2 \cup X_3$$
$$X = (0\ 0\ 0\ 0\ 0\ 0\ 1) \cup (1\ 0\ 0\ 0) \cup (0\ 0\ 0\ 0\ 1).$$

With Habits like cinema, smoking, drink etc had forced them to take work alone is in the on state and all other attributes are in the off state in the set of attributes related with children. In the second block, Children don't form vote bank alone is in the on state and all other nodes are in the off state is assumed. In the case of public role perpetuating slavery and caste bias alone is taken as the on state node and all other nodes are assumed to be in the off state. Now we study the effect of X on the FRTM model S.

$$\begin{aligned} XS &= (X_1 \cup X_2 \cup X_3)\,(S_1 \cup S_2 \cup S_3) \\ &= X_1 S_1 \cup X_2 S_2 \cup X_3 S_3 \\ &= (0\ 1\ {-1}\ 1) \cup (0\ 0\ 0\ 1\ 0) \\ &\quad \cup ({-1}\ 1\ 1\ {-1}\ 1\ 1\ 0). \end{aligned}$$

After thresholding we get the resultant trivector as

$$\begin{aligned} R &= (0\ 1\ 0\ 1) \cup (0\ 0\ 0\ 1\ 0) \cup (0\ 1\ 1\ 0\ 1\ 1\ 0) \\ &= R_1 \cup R_2 \cup R_3. \end{aligned}$$

Now we study the effect of R on S.

$$\begin{aligned} RS &= (R_1 \cup R_2 \cup R_3)\,(S_1 \cup S_2\ S_3\ )^T \\ &= R_1\ S_1^T \cup R_2\, S_2^T\ \cup R_3\, S_3^T \\ &= ({-2}\ 1\ 2\ 1\ 2\ 2\ 2) \cup (1\ {-1}\ 1\ 1) \\ &\quad \cup (2\ 3\ {-1}\ {-3}\ 3). \end{aligned}$$

After thresholding and updating we get the resultant W

$$W = (0\ 1\ 1\ 1\ 1\ 1\ 1) \cup (1\ 0\ 1\ 1) \cup (1\ 1\ 0\ 0\ 1)$$
$$= W_1 \cup W_2 \cup W_3.$$

Now we see the effect of W on the FRTM S.

$$WS = W_1 S_1 \cup W_2 S_2 \cup W_3 S_3$$
$$= (1\ 6\ –4\ 3) \cup (2\ 0\ 0\ 3\ 2)$$
$$\cup (–3\ 1\ 3\ –3\ 2\ 3\ 1).$$

After thresholding we get the resultant state trivector

$$P = (1\ 1\ 0\ 1) \cup (1\ 0\ 0\ 1\ 1) \cup (0\ 1\ 1\ 0\ 1\ 1\ 1)$$
$$= P_1 \cup P_2 \cup P_3.$$

Now we study the effect of P on S

$$PS = (P_1 \cup P_2 \cup P_3)\ (S_1 \cup S_2 \cup S_3)^T$$
$$= P_1\ S_1^T \cup P_2 S_2^T \cup P_3 S_3^T$$
$$= (–1\ 1\ 2\ 1\ 2\ 3\ 2) \cup (1\ 1\ 3\ 3)$$
$$\cup (3\ 3\ –2\ –4\ 4).$$

After thresholding we get the resultant trivector as

$$C = (0\ 1\ 1\ 1\ 1\ 1\ 1) \cup (1\ 1\ 1\ 1) \cup (1\ 1\ 0\ 0\ 1)$$
$$= C_1 \cup C_2 \cup C_3.$$

Now we study the influence of C on the dynamical system S.

$$CS = (C_1 \cup C_2 \cup C_3)\ (S_1 \cup S_2 \cup S_3)$$
$$= C_1 S_1 \cup C_2 S_2 \cup C_3 S_3$$
$$= (1\ 5\ –4\ 3) \cup (1\ 0\ 0\ 1\ 1)$$
$$\cup (–3\ 1\ 3\ –3\ 2\ 3\ 1).$$

After thresholding we get the resultant as

$$D = (1\ 1\ 0\ 1) \cup (1\ 1\ 0\ 1\ 1) \cup (0\ 1\ 1\ 0\ 1\ 1\ 1).$$

Thus we see the hidden pattern of the dynamical system is a fixed tripoint given by the binary bipair of row trivectors.

$$(\{(0\ 1\ 1\ 1\ 1\ 1\ 1) \cup (1\ 1\ 1\ 1) \cup (1\ 1\ 0\ 0\ 1)\},$$
$$(1\ 1\ 0\ 1) \cup (1\ 1\ 1\ 1\ 1) \cup (0\ 1\ 1\ 0\ 1\ 1\ 1)\};$$

when the initial vector with which we started was

$$X = (0\ 0\ 0\ 0\ 0\ 0\ 1) \cup (1\ 0\ 0\ 0) \cup (\ 0\ 0\ 0\ 0\ 1).$$

The influence can be interpreted as per the FRMs. Now we start to generalize this notion. We can define any FRn-M by which we mean Fuzzy relational n-maps where $n \geq 1$. Clearly when $n = 1$ we get the Fuzzy Relational Maps (FRM); when $n = 2$ we get the Fuzzy Relational Bimaps (FRBMs) when $n = 3$ we get the Fuzzy Relational trimaps (FRTMs). Thus we can extend this to any finite set say n; $n \geq 1$. That is we have (say) n sets of attributes disjoint then we get a related directed n-graph for the FRn-M. It is still important to say that these sets can be biconnected n-graph which is bipartite or a tri-connected n-graph which is bipartite or any r-connected n-graph $r \leq n$. Just we show by examples each type before we proceed to define the notion of NRMs and NRBMs etc.

***Example 3.6.3:*** Suppose we have the 7-graph given by the following disconnected 7-graph.

$$P = P_1 \cup P_2 \cup \ldots \cup P_7$$

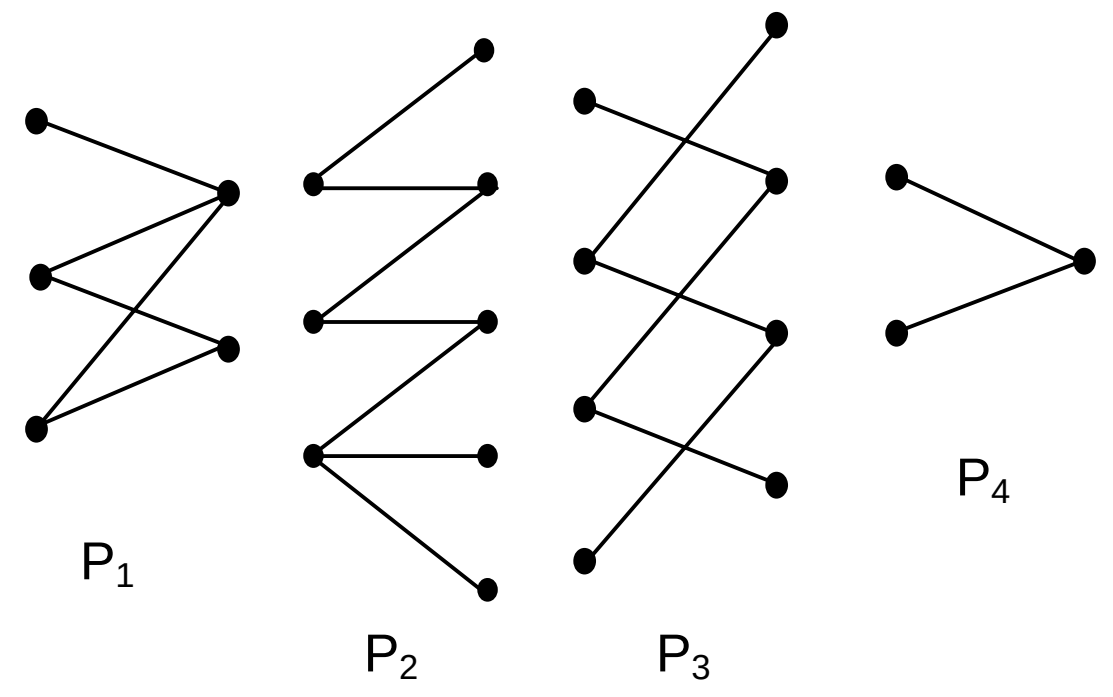

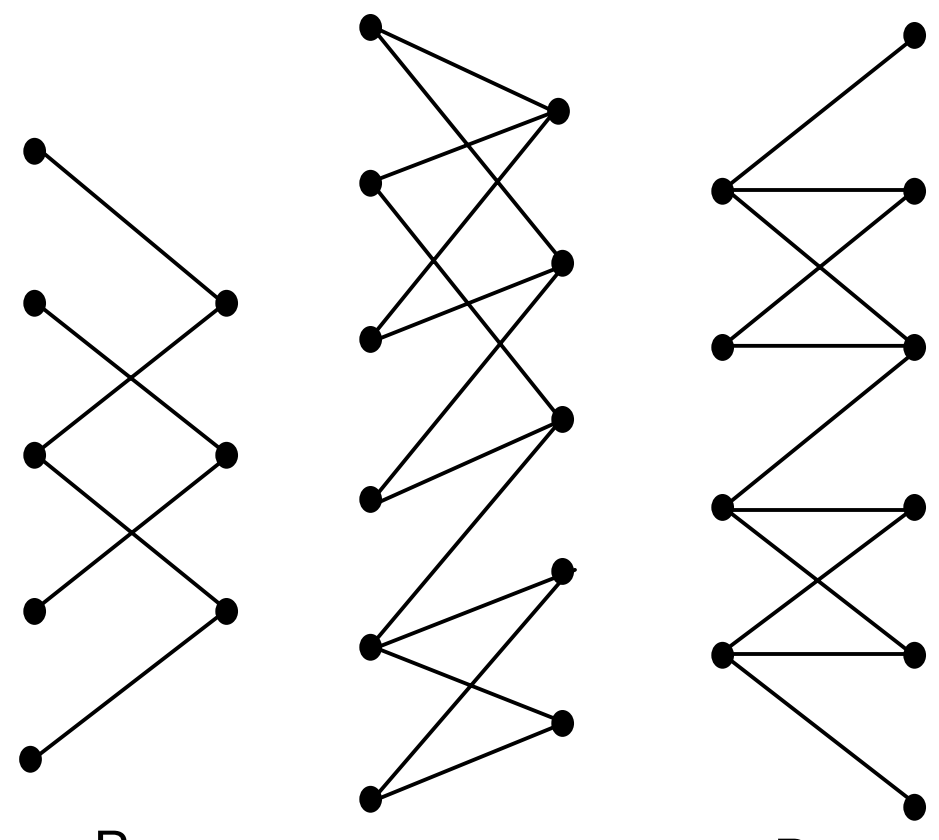

FIGURE 3.6.5

***Example 3.6.4:*** The following gives a strongly 8-connected bipartite 8 graph.

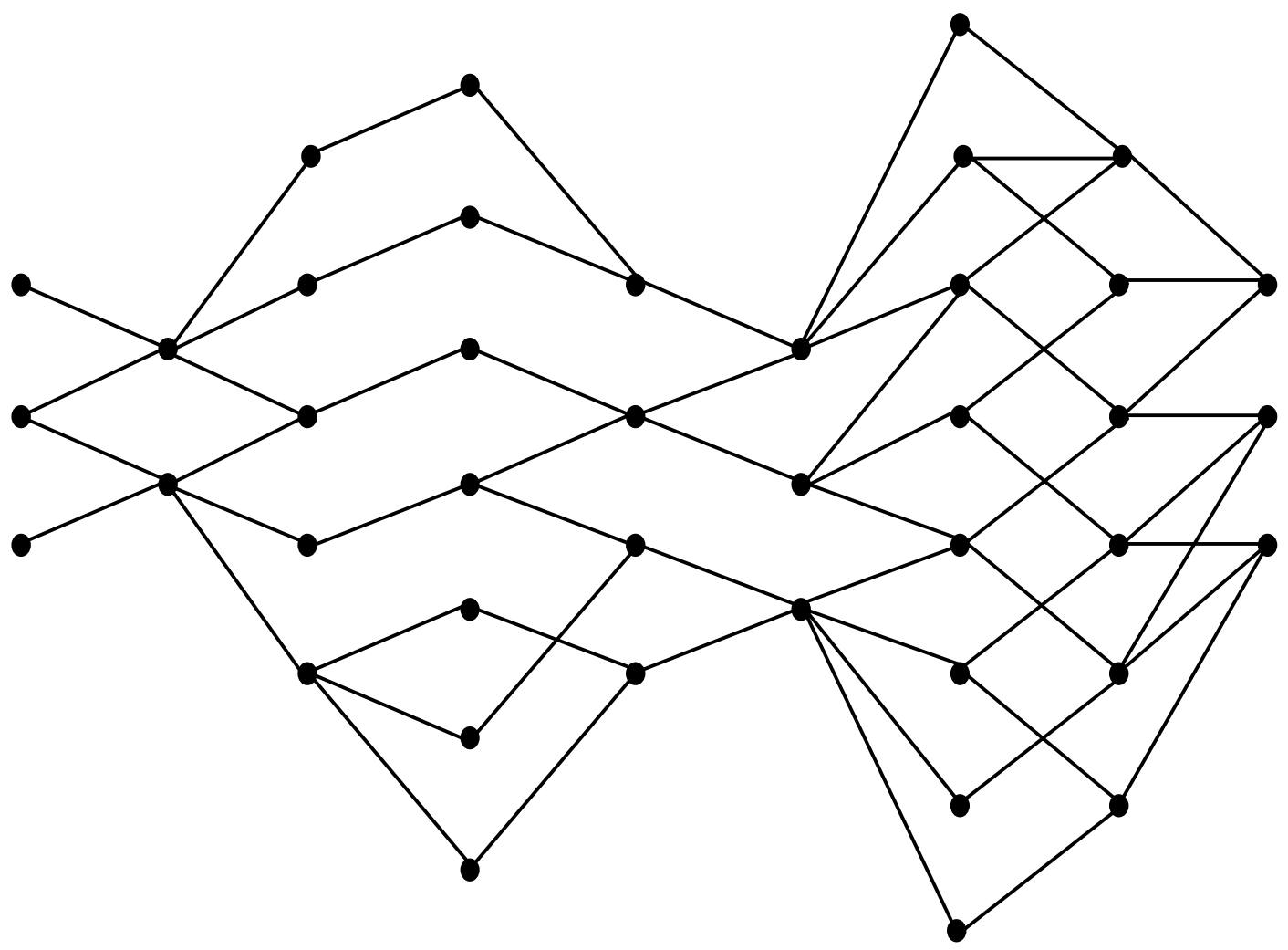

FIGURE 3.6.6

The separate 8-9 bipartite graph follows.

$$G = G_1 \cup G_2 \cup G_3 \cup G_4 \cup G_5 \cup G_6 \cup G_7 \cup G_8$$

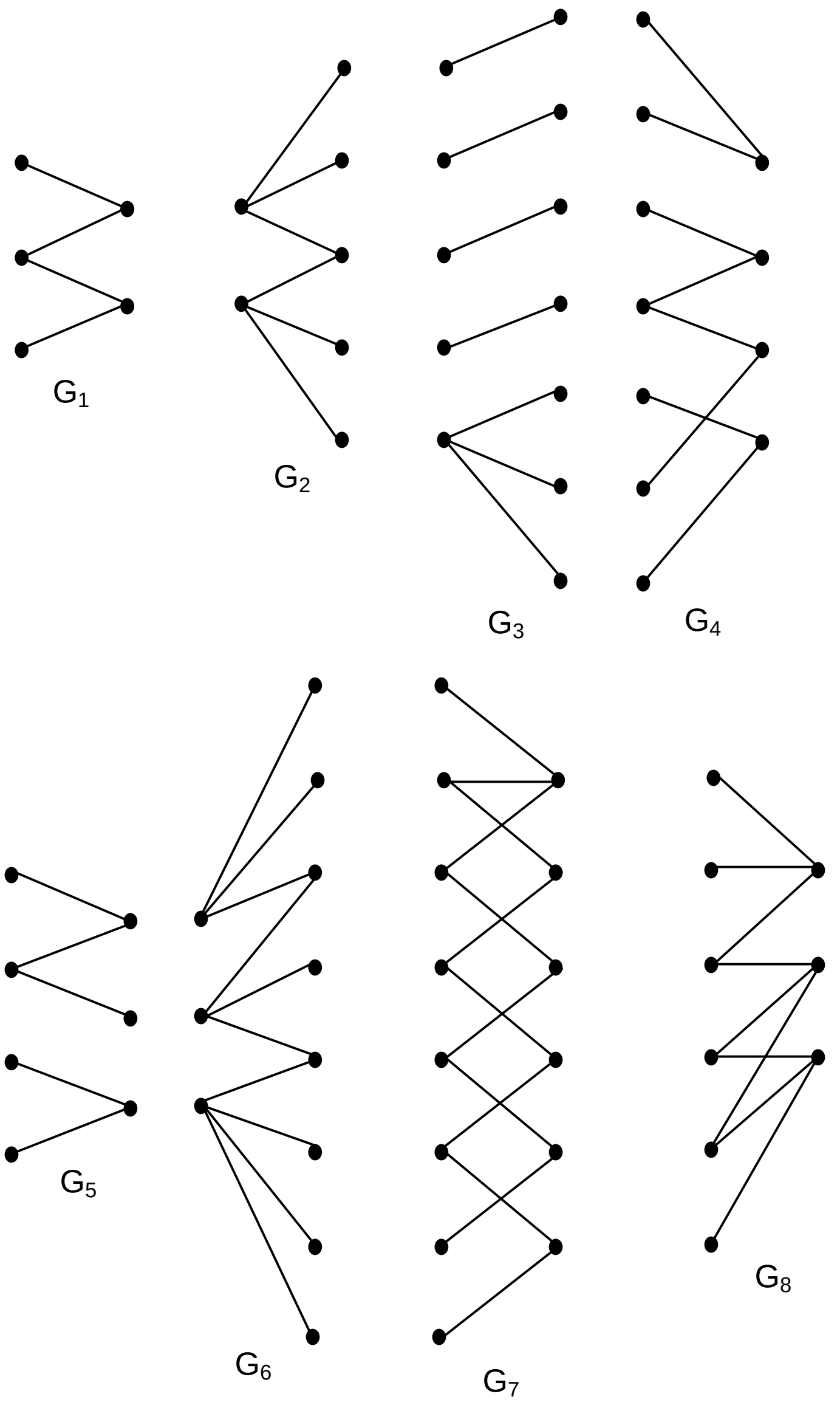

FIGURE: 3.6.7

***Example 3.6.5:*** Now we give an example of a 5-graph which is a biconnected bipartite 5-bigraph. It at least one pair of biconnected bipartite exists then we call such n-graph as biconnected bipartite bigraph.

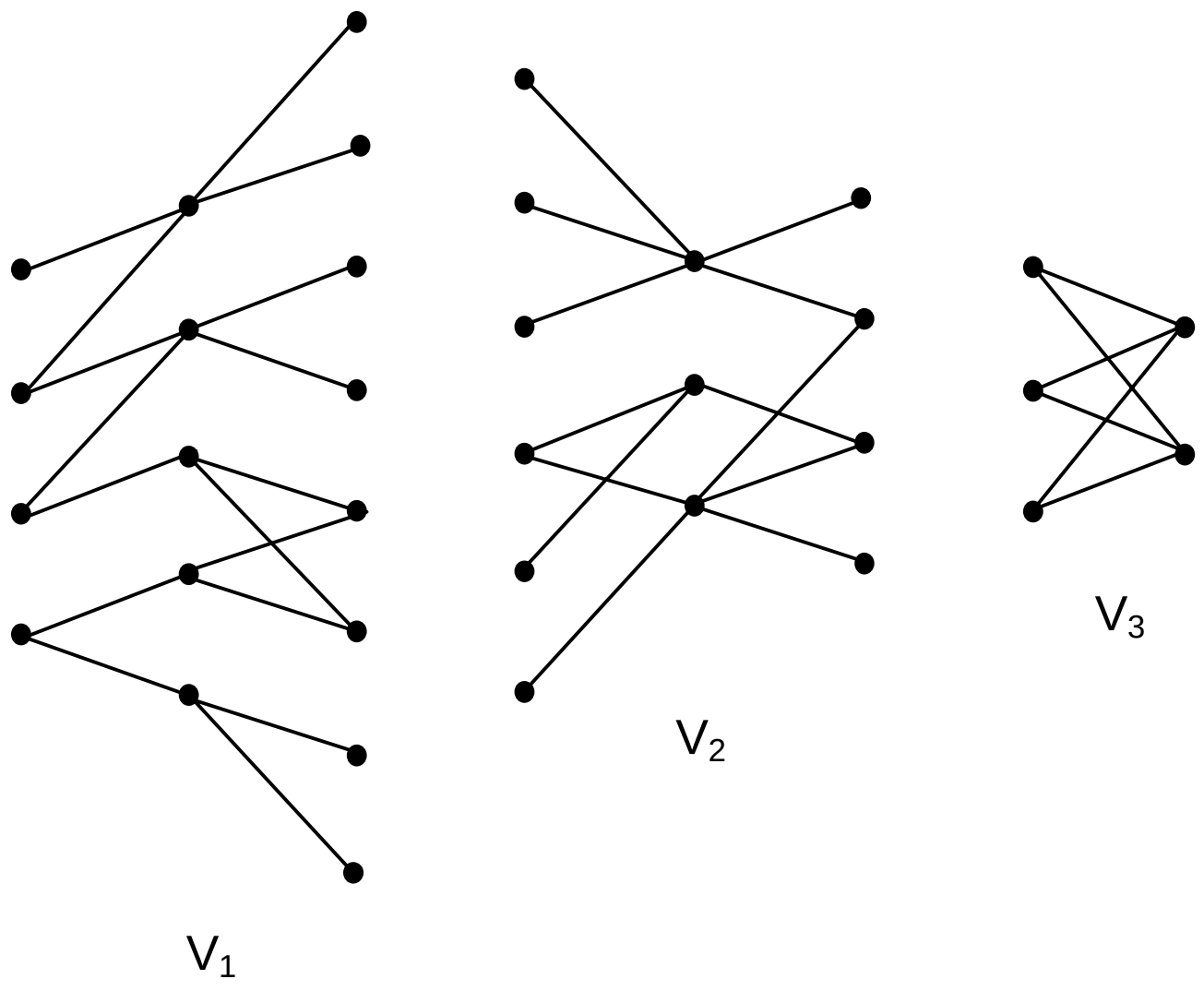

FIGURE 3.6.8

Clearly $V = V_1 \cup V_2 \cup V_3$ is a 5-graph, which is a biconnected bipartite 5-graph.

***Example 3.6.6:*** Now we give an example of a triconnected graph.

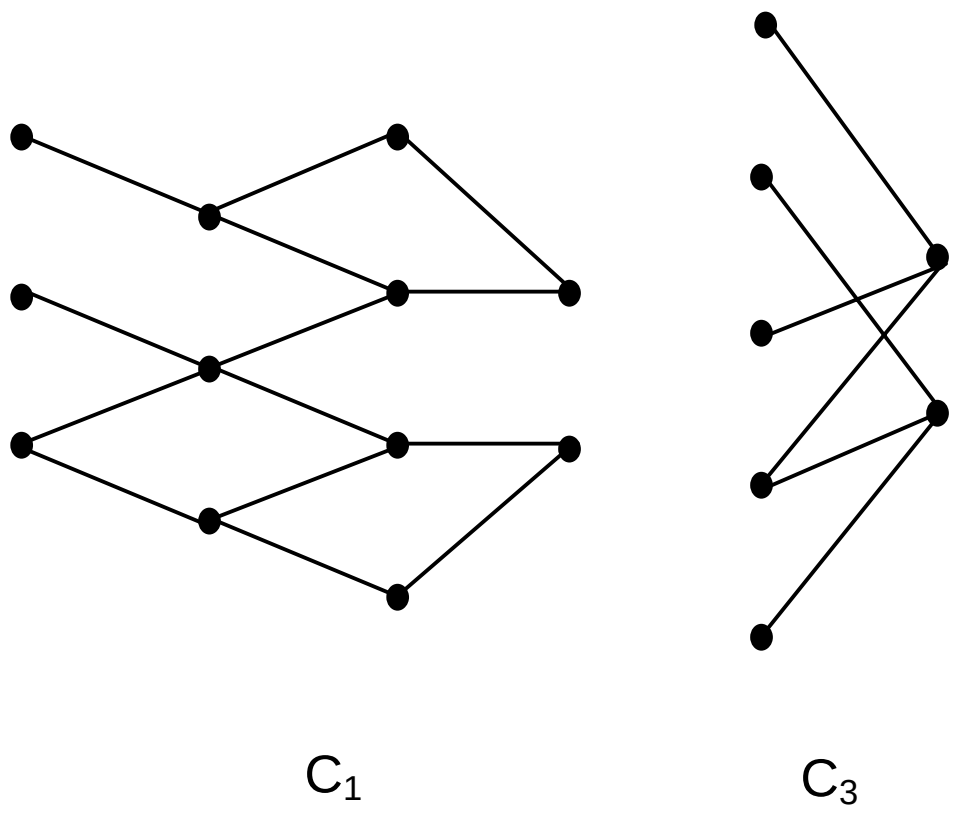

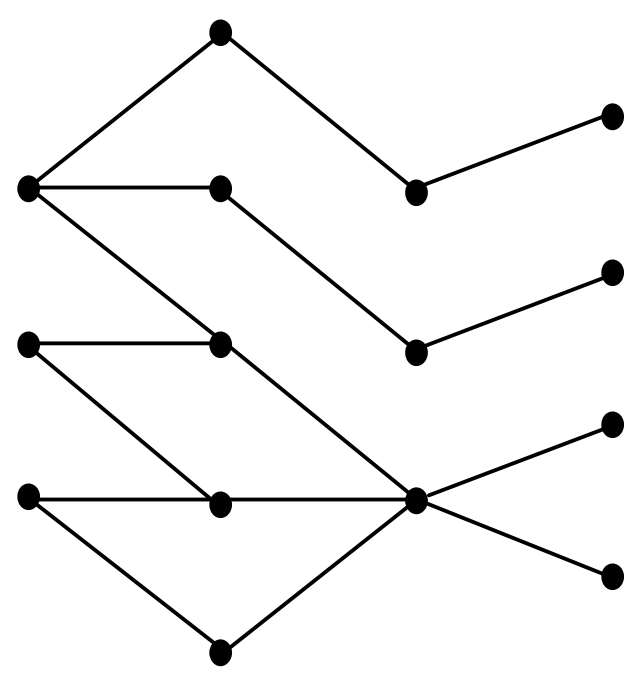

$C_2$

FIGURE 3.6.9

$C = C_1 \cup C_2 \cup C_3$, C is 7-graph which is triconnected bipartite 7-bigraph.

When we have FRMs with more than 3 sets of attributes we seek for 4-graphs or 5 graphs and so on. Thus all these graphs find their application is real world problems which are modeled using fuzzy logic and neutrosophic logic.

Chapter Four

# NEUTROSOPHIC BIGRAPHS, THEIR GENERALIZATIONS AND APPLICATIONS TO NEUTROSOPHIC MODELS

This chapter has six sections. In the first section we just recall the definition of some basic notions of neutrosophic graphs. In section two we for the first time define neutrosophic bigraphs and its properties. Section three introduces neutrosophic cognitive bimaps and its applications. Neutrosophic trimaps are introduced in section four. Section five recalls the notion of neutrosophic relational maps. The final section gives the concept of neutrosophic relational bimaps and their properties.

## 4.1 Some Basics of Neutrosophic Graphs

In this section we recall some basic notion of neutrosophic graphs, illustrate them and give some basic properties. We need the notion of neutrosophic graphs basically to obtain neutrosophic cognitive maps which will be nothing but directed neutrosophic graphs. Similarly neutrosophic relational maps will also be directed neutrosophic graphs.

It is no coincidence that graph theory has been independently discovered many times since it may quite properly be regarded as an area of applied mathematics. The subject finds its place in the work of Euler. Subsequent rediscoveries of graph theory were by Kirchhoff and Cayley. Euler (1707-1782) became the father of graph theory as well as topology when in 1936 he settled a famous unsolved problem in his day called the Konigsberg Bridge Problem.

Psychologist Lewin proposed in 1936 that the life space of an individual be represented by a planar map. His view point led the psychologists at the Research center for Group Dynamics to another psychological interpretation of a graph in which people are represented by points and interpersonal relations by lines. Such relations include love, hate, communication and power. In fact it was precisely this approach which led the author to a personal discovery of graph theory, aided and abetted by psychologists L. Festinger and D. Cartwright. Here it is pertinent to mention that the directed graphs of an FCMs or FRMs are nothing but the psychological inter-relations or feelings of different nodes; but it is unfortunate that in all these studies the concept of indeterminacy was never given any place, so in this chapter for the first time we will be having graphs in which the notion of indeterminacy i.e. when two vertex should be connected or not is dealt with. If graphs are to display human feelings then this point is very vital for in several situations certain relations between concepts may certainly remain as an indeterminate. So this section will purely cater to the properties of such graphs the edges of certain vertices may not find its connection i.e., they are indeterminates, which we will be defining as neutrosophic graphs. The indeterminate edge connection throughout this book will be denoted by dotted lines.

The world of theoretical physics discovered graph theory for its own purposes. In the study of statistical mechanics by Uhlenbeck the points stands for molecules and two adjacent points indicate nearest neighbor interaction of some physical kind, for example magnetic attraction or repulsion. But it is forgotten in all these situations we may have molecule structures which need not attract or repel but remain without action or not able to predict the action for such analysis we can certainly adopt the concept of neutrosophic graphs.

In a similar interpretation by Lee and Yang the points stand for small cubes in Euclidean space where each cube may or may not be occupied by a molecule. Then two points are adjacent whenever both spaces are occupied.

Feynmann proposed the diagram in which the points represent physical particles and the lines represent paths of the particles after collisions. Just at each stage of applying graph theory we may now feel the neutrosophic graph theory may be more suitable for many application.

Now we proceed on to define the neutrosophic graph.

**DEFINITION 4.1.1:** *A neutrosophic graph is a graph in which at least one edge is an indeterminacy denoted by dotted lines.*

**NOTATION***:* The indeterminacy of an edge between two vertices will always be denoted by dotted lines.

***Example 4.1.1:*** The following are neutrosophic graphs:

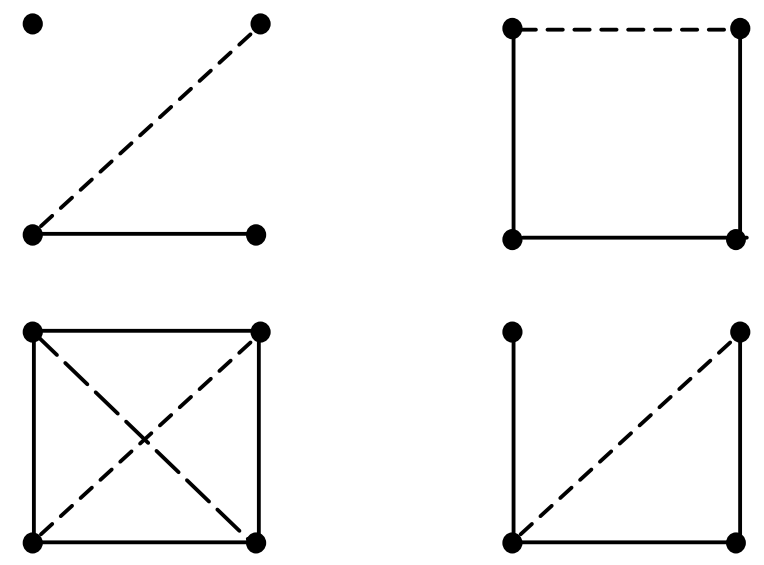

FIGURE: 4.1.1

All graphs in general are not neutrosophic graphs.

***Example 4.1.2:*** The following graphs are not neutrosophic graphs given in figure 4.1.2:

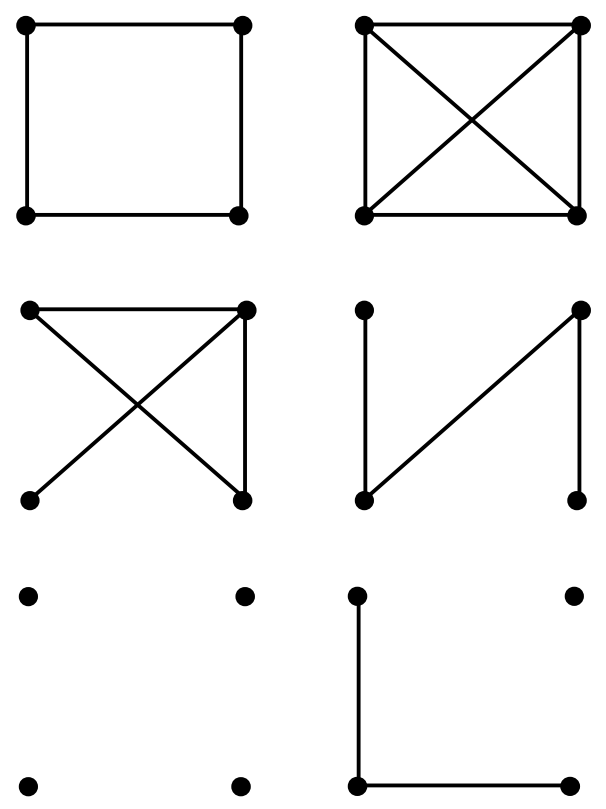

FIGURE: 4.1.2

**DEFINITION 4.1.2:** *A neutrosophic directed graph is a directed graph which has at least one edge to be an indeterminacy.*

**DEFINITION 4.1.3:** *A neutrosophic oriented graph is a neutrosophic directed graph having no symmetric pair of directed indeterminacy lines.*

**DEFINITION 4.1.4:** *A neutrosophic subgraph H of a neutrosophic graph G is a subgraph H which is itself a neutrosophic graph.*

**THEOREM 4.1.1:** *Let G be a neutrosophic graph. All subgraphs of G are not neutrosophic subgraphs of G.*

*Proof:* By an example. Consider the neutrosophic graph given in figure 4.1.3.

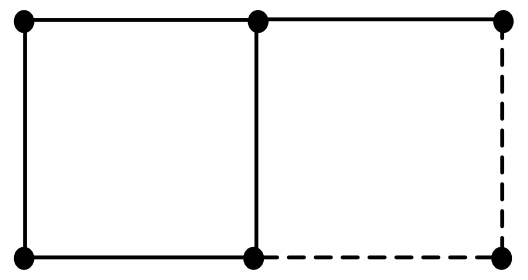

FIGURE: 4.1.3

This has a subgraph given by figure 4.1.4.

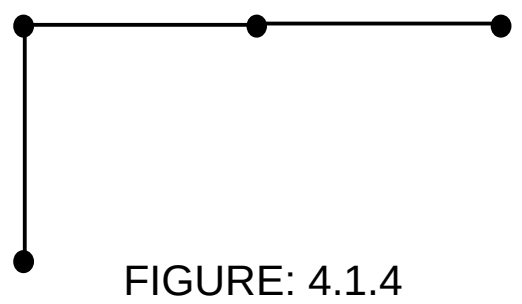

FIGURE: 4.1.4

which is not a neutrosophic subgraph of G.

**THEOREM 4.1.2:** *Let G be a neutrosophic graph. In general the removal of a point from G need not be a neutrosophic subgraph.*

*Proof:* Consider the graph G given in figure 4.1.5.

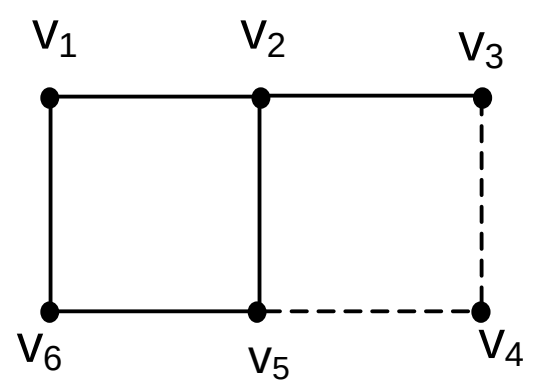

FIGURE: 4.1.5

G \ $v_4$ is only a subgraph of G but is not a neutrosophic subgraph of G.

Thus it is interesting to note that this is a main feature by which a graph differs from a neutrosophic graph.

**DEFINITION 4.1.5:** *Two graphs G and H are neutrosophically isomorphic if*

*i. They are isomorphic.*
*ii. If there exists a one to one correspondence between their point sets which preserve indeterminacy adjacency.*

**DEFINITION 4.1.6:** *A neutrosophic walk of a neutrosophic graph G is a walk of the graph G in which at least one of the lines is an indeterminacy line. The neutrosophic walk is* *neutrosophic closed if $v_0 = v_n$ and is neutrosophic open otherwise.*

*It is a neutrosophic trial if all the lines are distinct and at least one of the lines is a indeterminacy line and a path, if all points are distinct (i.e. this necessarily means all lines are distinct and at least one line is a line of indeterminacy). If the neutrosophic walk is neutrosophic closed then it is a neutrosophic cycle provided its n points are distinct and $n \geq 3$.*

*A neutrosophic graph is neutrosophic connected if it is connected and at least a pair of points are joined by a path. A neutrosophic maximal connected neutrosophic subgraph of G is called a neutrosophic connected component or simple neutrosophic component of G.*

*Thus a neutrosophic graph has at least two neutrosophic components then it is neutrosophic disconnected. Even if one is a component and another is a neutrosophic component still we do not say the graph is neutrosophic disconnected.*

***Example 4.1.3:*** Neutrosophic disconnected graphs are given in figure 4.1.6.

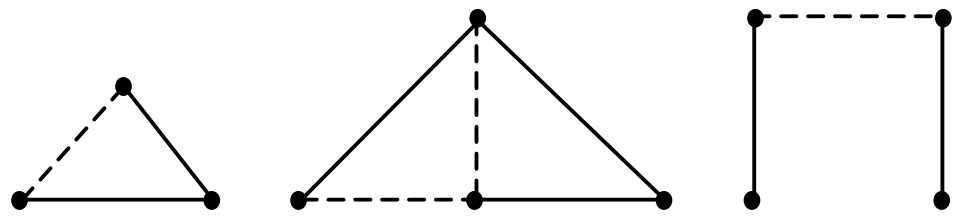

FIGURE: 4.1.6

***Example 4.1.4:*** Graph which is not neutrosophic disconnected is given by figure 4.1.7.

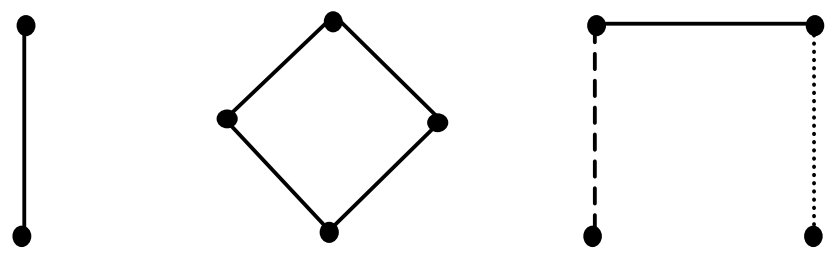

FIGURE: 4.1.7

Several results in this direction can be defined and analyzed.

**DEFINITION 4.1.7:** *A neutrosophic bigraph, G is a bigraph, G whose point set V can be partitioned into two subsets $V_1$ and $V_2$ such that at least a line of G which joins $V_1$ with $V_2$ is a line of indeterminacy.*

This neutrosophic bigraphs will certainly play a role in the study of FRMs and in fact we give a method of conversion of data from FRMs to FCMs.

As both the models FRMs and FCMs work on the adjacency or the connection matrix we just define the neutrosophic adjacency matrix related to a neutrosophic graph G given by figure 4.1.8.

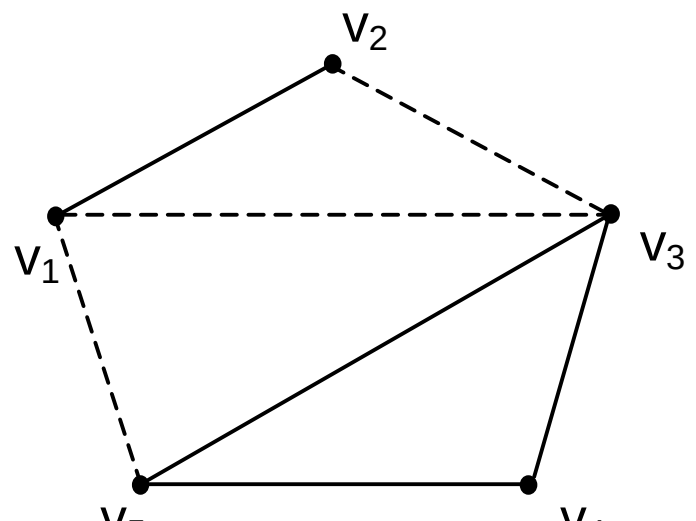

FIGURE: 4.1.8

The neutrosophic adjacency matrix is N(A)

$$N(A) = \begin{bmatrix} 0 & 1 & I & 0 & I \\ 1 & 0 & I & 0 & 0 \\ I & I & 0 & 1 & 1 \\ 0 & 0 & 1 & 0 & 1 \\ I & 0 & 1 & 1 & 0 \end{bmatrix}.$$

Its entries will not only be 0 and 1 but also the indeterminacy *I*.

**DEFINITION 4.1.8:** *Let G be a neutrosophic graph. The adjacency matrix of G with entries from the set (I, 0, 1) is called the neutrosophic adjacency matrix of the graph.*

Now we proceed on to define the notion of neutrosophic bigraphs in a little different way for we have to give it a neutrosophic bimatrix representation.

## 4.2 Neutrosophic Bigraphs and their Properties

In this section we for the first time introduce the notion of neutrosophic bigraphs. The neutrosophic bigraphs mainly differs from the neutrosophic graphs in their matrix representation. In fact the matrix associated with the neutrosophic bigraphs are neutrosophic bimatrices. More so we can extend this notion to any n-graphs when n = 1 we get the usual graphs, and n = 2 we get the bigraph, n = 3 we call them as trigraphs and n = 4 we say the graphs as quad graphs. The main motivation for these structures are they have application in fuzzy and neutrosophic real models.

**DEFINITION 4.2.1:** *A neutrosophic graph $G_N = G_1 \cup G_2$ is said to be a neutrosophic bigraph if both $G_1$ and $G_2$ are neutrosophic graphs that the set of vertices of $G_1$ and $G_2$ are different atleast by one coordinate i.e. $V(G_1) \not\subset V(G_2)$ or $V(G_2) \not\subset V(G_1)$ i.e. $V(G_1) \cap V(G_2) = \phi$ is also possible but is not a condition i.e. the vertex set of $G_1$ is not a proper subset of the vertex set of $G_2$ or vice versa or atleast one edge is different in the graphs $G_1$ and $G_2$. 'or' is not used in the mutually exclusive sense.*

***Example 4.2.1:*** Let $G = G_1 \cup G_2$ be the neutrosophic bigraph given by the following figure 4.2.1.

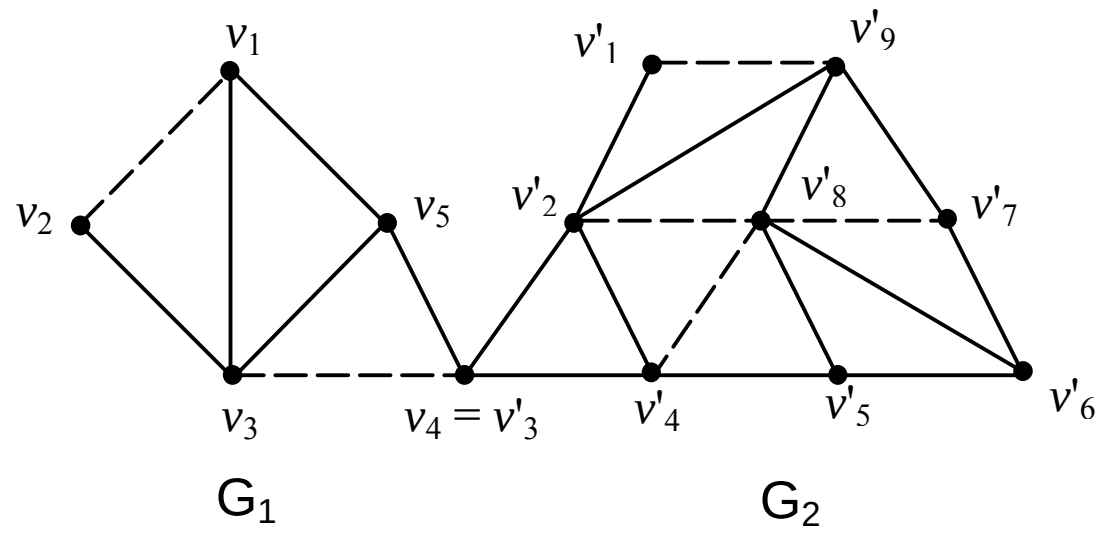

FIGURE: 4.2.1

$$V(G_1) = (v_1, v_2, v_3, v_4, v_5)$$
$$V(G_2) = (v'_1, v'_2, v'_3, v'_4, v'_5, v'_6, v'_7, v'_8, v'_9)$$

dotted edges are the neutrosophic edges. Thus $G_N = G_1 \cup G_2$ is a neutrosophic bigraph.

**DEFINITION 4.2.2:** *A neutrosophic weak bigraph $G = G_1 \cup G_2$ is a bigraph in which at least one of $G_1$ or $G_2$ is a neutrosophic graph and the other need not be a neutrosophic graph.*

***Example 4.2.2:*** Let $G = G_1 \cup G_2$ be a bigraph given by the figure 4.2.2.

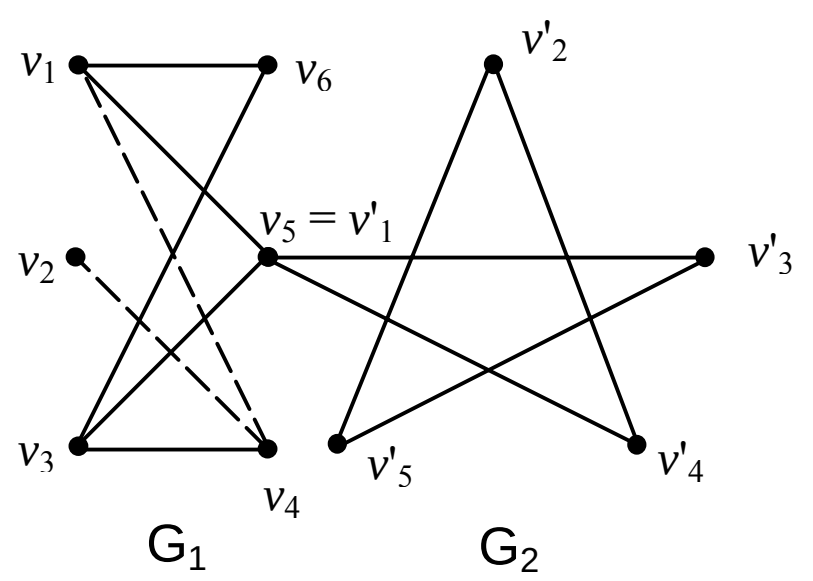

FIGURE: 4.2.2

Clearly $G_1$ is a neutrosophic graph but $G_2$ is not a neutrosophic graph. So $G = G_1 \cup G_2$ is only a weak neutrosophic bigraph.

**THEOREM 4.2.1:** *All neutrosophic bigraph are neutrosophic weak bigraph but a neutrosophic weak bigraph in general is not a neutrosophic bigraph.*

*Proof:* By the very definition we see all neutrosophic bigraphs are weak neutrosophic bigraphs. To show a weak neutrosophic bigraph in general is not a neutrosophic bigraph. Consider the weak neutrosophic bigraph $G = G_1 \cup G_2$ given by the following figure 4.2.3. given by the following example.

***Example 4.2.3:***

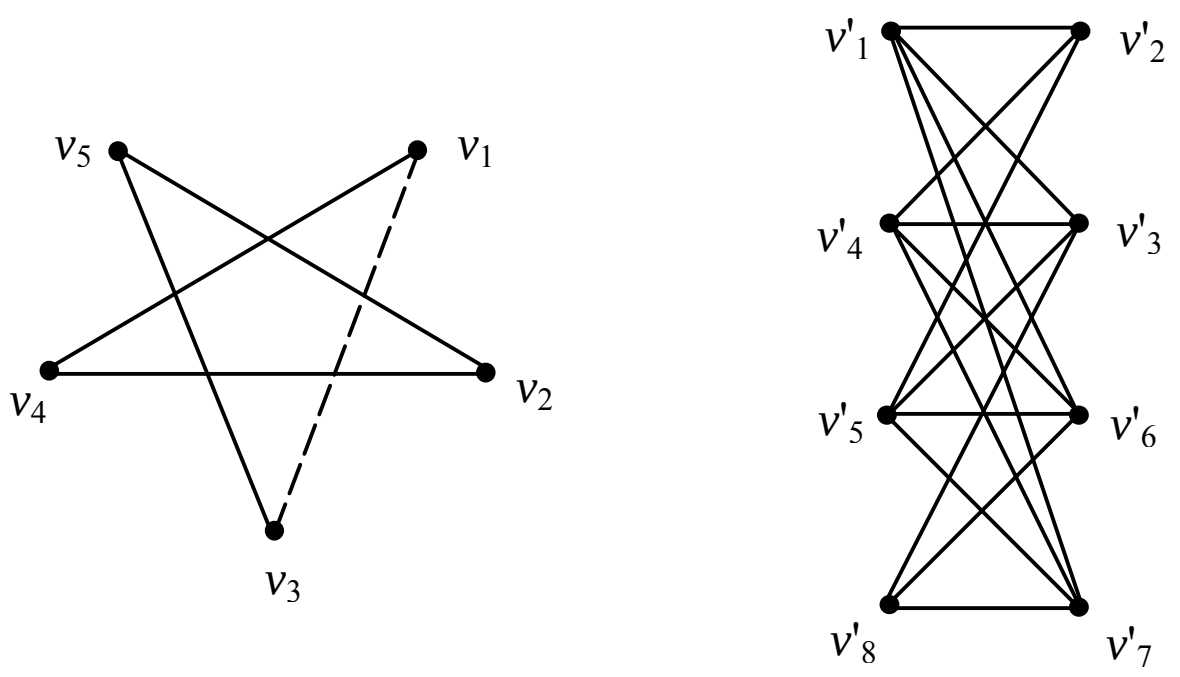

FIGURE: 4.2.3

$$G = G_1 \cup G_2$$

Clearly $G_1$ is a neutrosophic graph but $G_2$ is not a neutrosophic graph. So $G = G_1 \cup G_2$ is not a neutrosophic bigraph but only a weak neutrosophic bigraph.

**DEFINITION 4.2.3:** *A neutrosophic bigraph G is vertex connected or vertex glued neutrosophic bigraph if the bigraph $G = G_1 \cup G_2$ is a vertex glued bigraph.*

***Example 4.2.4:***

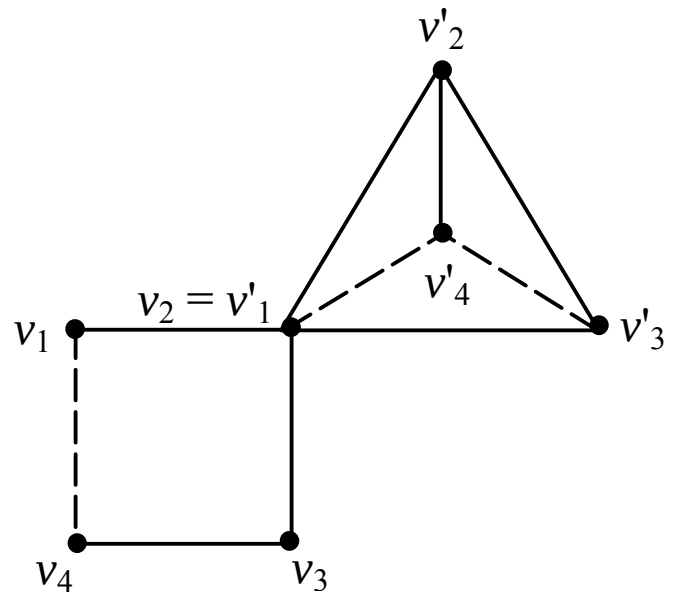

FIGURE: 4.2.4

Clearly G = $G_1 \cup G_2$ where

$$V(G_1) = (v_1, v_2, v_3, v_4)$$
$$V(G_2) = (v'_1, v'_2, v'_3, v'_4)$$

is a vertex glued neutrosophic bigraph given in figure 4.2.4.

**DEFINITION 4.2.4:** *Let G = $G_1 \cup G_2$ be a neutrosophic bigraph. G is said to be a edge glued neutrosophic bigraph if G as a bigraph is a edge glued bigraph.*

***Example 4.2.5:*** Let G = $G_1 \cup G_2$ be a neutrosophic bigraph which is edge glued given by the following figure 4.2.5.

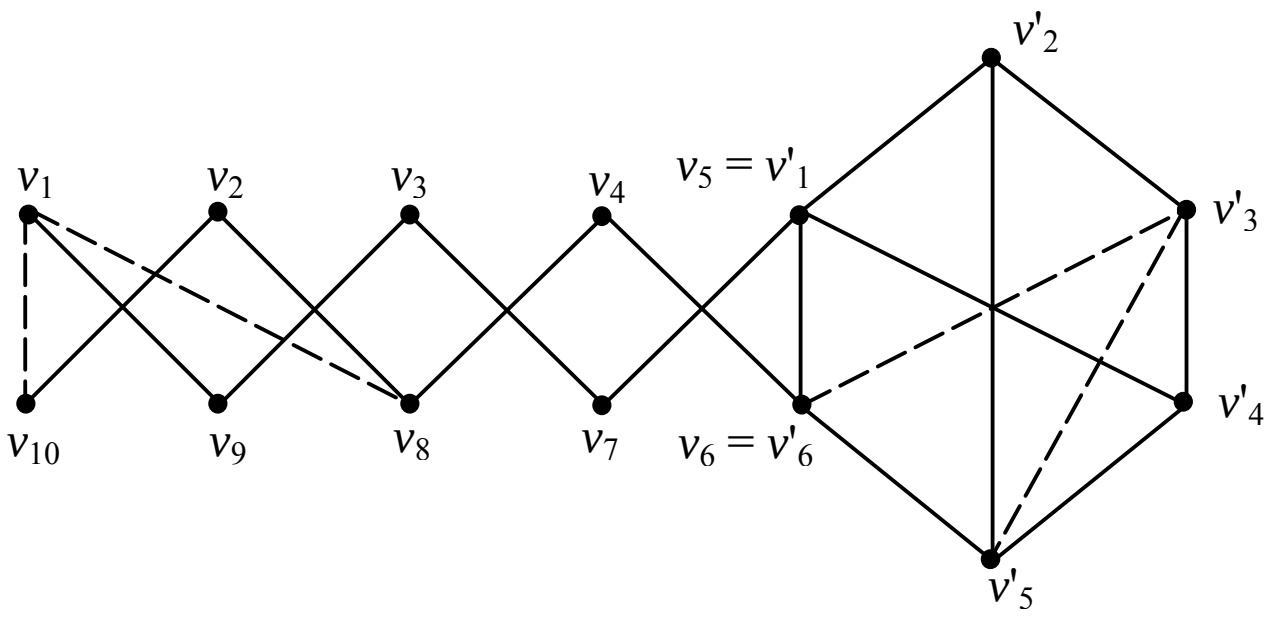

FIGURE: 4.2.5

Thus both $G_1$ and $G_2$ are neutrosophic graphs with the edge joining the vertices $v_5$ and $v_6$ in $G_1$ and $v'_1$ and $v'_6$ in $G_2$. Thus G = $G_1 \cup G_2$ is a edge glued neutrosophic bigraph.

$$V(G_1) = \{ v_1, v_2, \ldots, v_{10}\} \text{ and}$$
$$V(G_1) = \{ v'_1, v'_2, \ldots, v'_6\}$$

Now we proceed on to define the notion of neutrosophically glued neutrosophic bigraph.

**DEFINITION 4.2.5:** *Let G = $G_1 \cup G_2$ be a neutrosophic bigraph. G is said to neutrosophically glued neutrosophic bigraph if the bigraph is edge glued and at least one edge is* *a neutrosophic edge i.e. atleast one edge is joined by dotted lines.*

***Example 4.2.6:*** Let G = $G_1 \cup G_2$ be a neutrosophic bigraph given by the following figure 4.2.6.

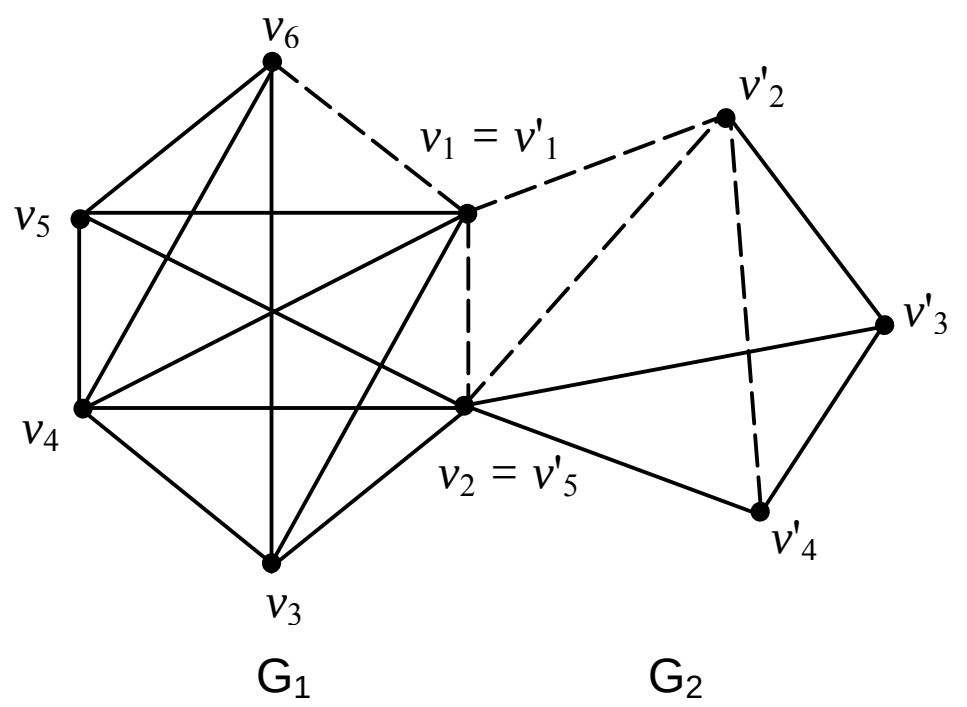

FIGURE: 4.2.6

Clearly the bigraph is a neutrosophic bigraph. We see a neutrosophic edge is common between $G_1$ and $G_2$ for the separate figures of the neutrosophic graphs are given below.

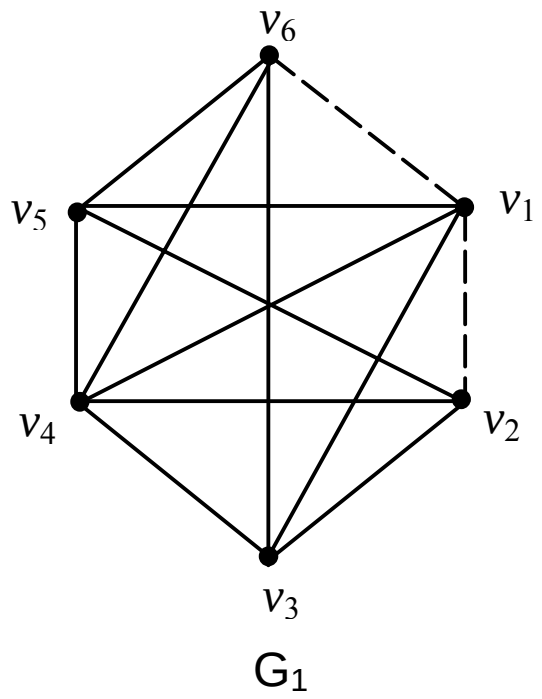

FIGURE: 4.2.6a

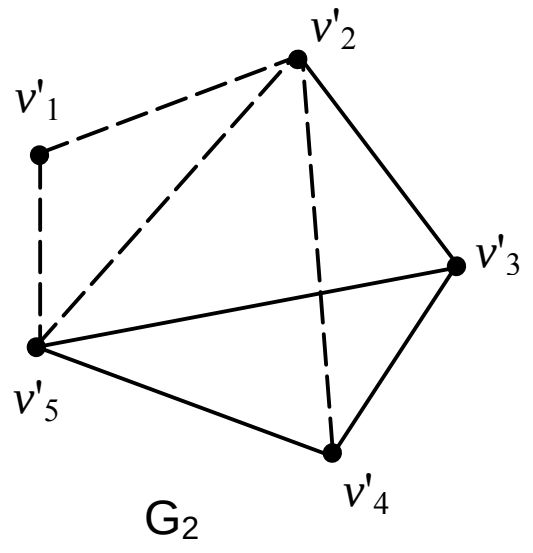

FIGURE: 4.2.6b

Now the notion of vertex glued and edge glued in the case of weak neutrosophic bigraph are also the same. But the concept of neutrosophically glued weak neutrosophic bigraph cannot be defined. This is proved by the following theorem.

*Note:* By the bigraph we mean not a edge or a vertex bigraph it has more than 2 vertices and more than one edge.

**THEOREM 4.2.2:** *Let $G = G_1 \cup G_2$ be a weak neutrosophic graph which is not a neutrosophic bigraph. Then $G = G_1 \cup G_2$ cannot be a neutrosophically glued neutrosophic bigraph.*

*Proof:* Let $G = G_1 \cup G_2$ is a weak neutrosophic bigraph which is not a neutrosophic bigraph; i.e. without loss in generality we assume $G_1$ is a neutrosophic graph and $G_2$ is not a neutrosophic graph. $G = G_1 \cup G_2$ is a weak neutrosophic bigraph only.

Suppose $G = G_1 \cup G_2$ is neutrosophic bigraph then both the graphs $G_1$ and $G_2$ becomes neutrosophic as both $G_1$ and $G_2$ have only one neutrosophic edge in common, which is a very contradiction to our assumption that G is only a weak neutrosophic bigraph.

**DEFINITION 4.2.6:** *Let $G = G_1 \cup G_2$ be a neutrosophic bigraph. G is a said to be a subbigraph connected if the graphs $G_1$ and $G_2$ have a subgraph in common.*

The above definition is true even in case of weak neutrosophic bigraph.

**DEFINITION 4.2.7:** *Let $G = G_1 \cup G_2$ be a neutrosophic bigraph. G is said to be a neutrosophic subbigraph connected if the graph $G_1$ and $G_2$ have neutrosophic subbigraph in common.*

***Example 4.2.7:*** Let $G = G_1 \cup G_2$ be a neutrosophic bigraph given by the following figure.

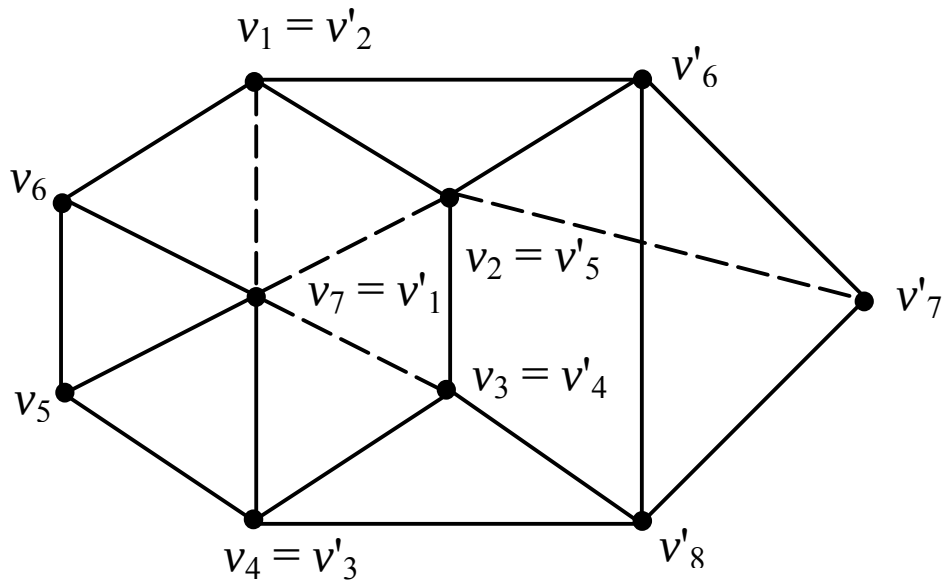

FIGURE: 4.2.7

$$V(G_1) = (v_1, v_2, \ldots, v_7)$$
$$V(G_2) = (v'_1, v'_2, \ldots, v'_8)$$

The separate figure of the two neutrosophic graphs $G_1$ and $G_2$ are as follows.

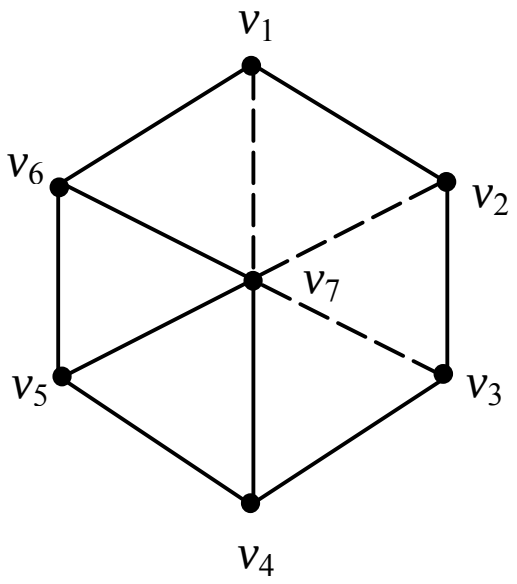

$G_1$

FIGURE: 4.2.7a

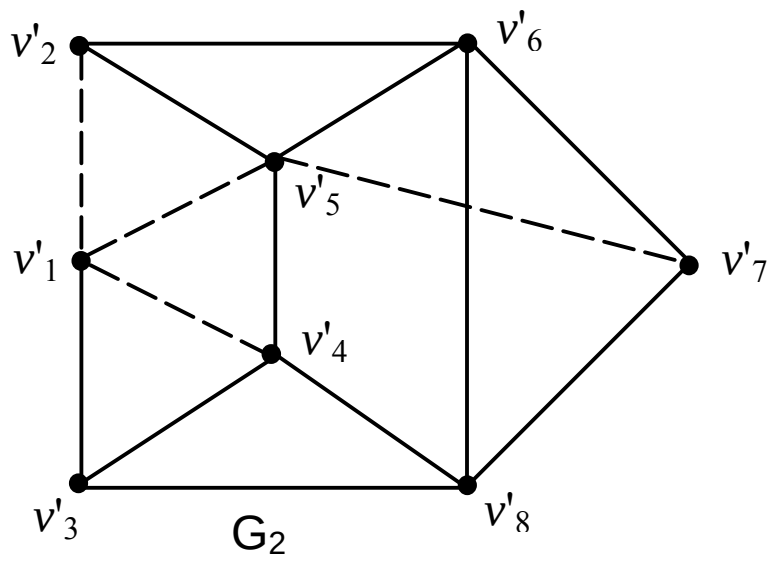

$G_2$

FIGURE: 4.2.7b

The common neutrosophic subbigraph is given by the following figure.

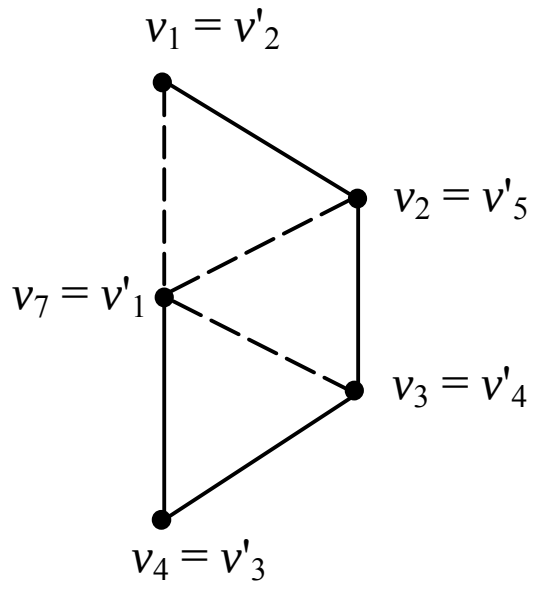

FIGURE: 4.2.7c

Clearly this is a neutrosophic subbigraph of both $G_1$ and $G_2$.

**THEOREM 4.2.3:** *Let $G = G_1 \cup G_2$ be a weak neutrosophic bigraph which is not a neutrosophic bigraph. G cannot be neutrosophic subbigraph connected.*

*Proof*: Given $G = G_1 \cup G_2$ is a weak neutrosophic bigraph which is not a neutrosophic bigraph. That is only one of $G_1$ or $G_2$ is a neutrosophic graph. So by Theorem 4.2.2 the weak neutrosophic bigraph cannot be even neutrosophic edge connected so if, this is to be neutrosophic subbigraph connected $G_1$ and $G_2$ must have a neutrosophic subgraph. which is not possible as only one of $G_1$ or $G_2$ is a neutrosophic graph. Hence the claim.

We know to every neutrosophic graph G there is a neutrosophic matrix associated with it. Likewise with every neutrosophic bigraph, we have neutrosophic bimatrix associated with it.

Further with every weak neutrosophic bigraph we have a weak neutrosophic bimatrix associated with it.

***Example 4.2.8:*** Let G = $G_1 \cup G_2$ be a neutrosophic bigraph given by the following figure 4.2.8.

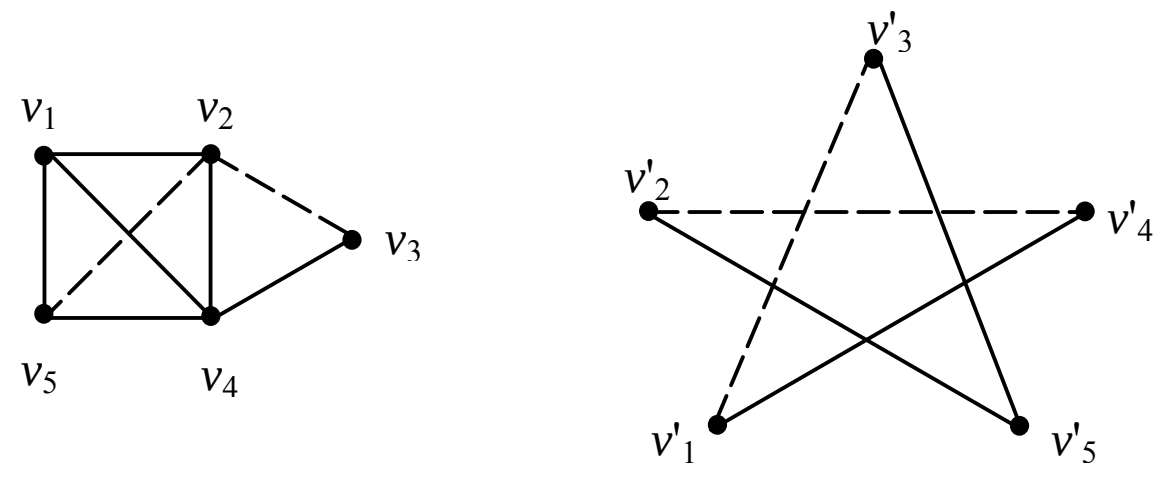

FIGURE 4.2.8

The associated neutrosophic bimatrix is given by M= $M_1 \cup M_2$.

$$M = \begin{array}{c} \\ v_1 \\ v_2 \\ v_3 \\ v_4 \\ v_5 \end{array} \begin{array}{c} \begin{array}{ccccc} v_1 & v_2 & v_3 & v_4 & v_5 \end{array} \\ \begin{bmatrix} 0 & 1 & 0 & 1 & 1 \\ 1 & 0 & I & 1 & I \\ 0 & I & 0 & 1 & 0 \\ 1 & 1 & 1 & 0 & 1 \\ 1 & I & 0 & 1 & 0 \end{bmatrix} \end{array} \cup \begin{array}{c} \\ v'_1 \\ v'_2 \\ v'_3 \\ v'_4 \\ v'_5 \end{array} \begin{array}{c} \begin{array}{ccccc} v'_1 & v'_2 & v'_3 & v'_4 & v'_5 \end{array} \\ \begin{bmatrix} 0 & 0 & I & 1 & 0 \\ 0 & 0 & 0 & I & 1 \\ I & 0 & 0 & 0 & 1 \\ 1 & I & 0 & 0 & 0 \\ 0 & 1 & 1 & 0 & 0 \end{bmatrix} \end{array}$$

We see both $M_1$ and $M_2$ are neutrosophic matrices so M is a neutrosophic bimatrix.

Now we proceed on to give an example of a weak neutrosophic bigraph and give its weak neutrosophic bimatrix representation.

***Example 4.2.9:*** Let G = $G_1 \cup G_2$ be a weak neutrosophic bigraph given by the following figure 4.2.9.

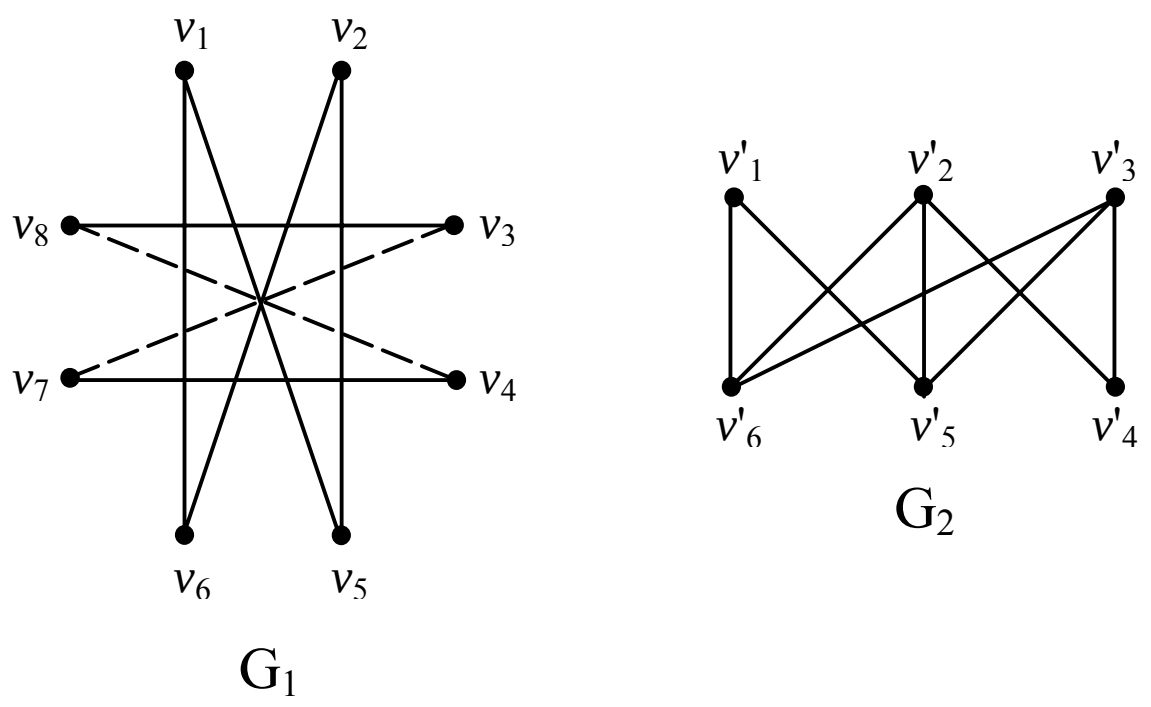

FIGURE: 4.2.9

The related bimatrix M = $M_1 \cup M_2$

$$
\begin{array}{c}
\\ v_1 \\ v_2 \\ v_3 \\ v_4 \\ v_5 \\ v_6 \\ v_7 \\ v_8
\end{array}
\begin{array}{c}
\begin{array}{cccccccc} v_1 & v_2 & v_3 & v_4 & v_5 & v_6 & v_7 & v_8 \end{array} \\
\begin{bmatrix}
0 & 0 & 0 & 0 & 1 & 1 & 0 & 0 \\
0 & 0 & 0 & 0 & 1 & 1 & 0 & 0 \\
0 & 0 & 0 & 0 & 0 & 0 & I & 1 \\
0 & 0 & 0 & 0 & 0 & 0 & 1 & I \\
1 & 1 & 0 & 0 & 0 & 0 & 0 & 0 \\
1 & 1 & 0 & 0 & 0 & 0 & 0 & 0 \\
0 & 0 & I & 1 & 0 & 0 & 0 & 0 \\
0 & 0 & 1 & I & 0 & 0 & 0 & 0
\end{bmatrix}
\end{array}
\cup
\begin{array}{c}
\\ v'_1 \\ v'_2 \\ v'_3 \\ v'_4 \\ v'_5 \\ v'_6
\end{array}
\begin{array}{c}
\begin{array}{cccccc} v'_1 & v'_2 & v'_3 & v'_4 & v'_5 & v'_6 \end{array} \\
\begin{bmatrix}
0 & 0 & 0 & 0 & 1 & 1 \\
0 & 0 & 0 & 1 & 1 & 1 \\
0 & 0 & 0 & 1 & 1 & 1 \\
0 & 1 & 1 & 0 & 0 & 0 \\
1 & 1 & 1 & 0 & 0 & 0 \\
1 & 1 & 1 & 0 & 0 & 0
\end{bmatrix}
\end{array}
$$

Clearly the bigraph has a associated bimatrix which is a weak neutrosophic bimatrix. Almost all properties related with bimatrix can be derived in the case of neutrosophic bimatrix with a possible modifications whenever necessary. This is true even in case of weak neutrosophic bimatrices.

Now we proceed on to define the notion of square, rectangular, mixed square and mixed rectangular neutrosophic bimatrices.

**DEFINITION 4.2.8:** *A bimatrix (weak neutrosophic bimatrix) M = $M_1 \cup M_2$ is said to be a square neutrosophic bimatrix (weak neutrosophic bimatrix) if both $M_1$ and $M_2$ are square neutrosophic bimatrices of same order.*

They are said to be mixed square bimatrices if $M_1$ and $M_2$ are square matrices of different order. A neutrosophic bimatrix $M = M_1 \cup M_2$ is said to be a m × n rectangular neutrosophic bimatrix if both $M_1$ and $M_2$ are m × n neutrosophic rectangular matrices.

A neutrosophic bimatrix $M = M_1 \cup M_2$ is said to be a mixed rectangular neutrosophic bimatrix if both $M_1$ and $M_2$ are rectangular neutrosophic matrices of different orders.

A neutrosophic birow vector

$$M = \left(m_1^1, m_2^1, \ldots, m_t^1\right) \cup \left(m_1^2, m_2^2, \ldots, m_n^2\right), m_i^k,$$

k = 1, 2 are reals with atleast one $m_p^1$ and $m_q^2$ to be *I*. A neutrosophic column bivector

$$C = \begin{bmatrix} a_1^1 \\ \vdots \\ a_r^1 \end{bmatrix} \cup \begin{bmatrix} a_1^2 \\ \vdots \\ a_s^2 \end{bmatrix}$$

with at least one $a_i^1$ and $a_j^2$ to be *I*.

***Example 4.2.10:*** Let B = (1 0 6 0 *I* – 1 1) ∪ (–1 *I* 4 0 *I* 1 –1 4 –2) is a neutrosophic row bivector

$$C = \begin{bmatrix} 0 \\ -1 \\ I \\ 2 \\ 3 \\ -4 \end{bmatrix} \cup \begin{bmatrix} I \\ 1 \\ -4 \\ 0 \\ 2 \\ -5 \end{bmatrix}$$

is a neutrosophic column bivector.

$$B^1 = (1\ 1\ 0\ 6\ –1\ I\ –9) \cup (1\ 2\ 3\ 4\ 5\ –3)$$

is a weak neutrosophic row bivector.

$$C^1 = \begin{bmatrix} 3 \\ -1 \\ 2 \\ 1 \\ 0 \\ 5 \end{bmatrix} \cup \begin{bmatrix} I \\ 1 \\ 0 \\ 5 \\ -3 \end{bmatrix}$$

is a weak neutrosophic column bivector.

Now we define neutrosophic bipartite bigraph. A bipartite bigraph $G = G_1 \cup G_2$ is said to be a neutrosophic bipartite bigraph if both the graphs $G_1$ and $G_2$ are bipartite and neutrosophic i.e., atleast an edge which is an indeterminate.

A bipartite bigraph $G = G_1 \cup G_2$ is said to be a weak neutrosophic bipartite bigraph if and only if one of $G_1$ or $G_2$ is a neutrosophic graph but both of them are bipartite.

Before we go for some more results we just illustrate them by examples.

***Example 4.2.11:*** Let $G = G_1 \cup G_2$ be a neutrosophic bipartite bigraph given by the following figure 4.2.10.

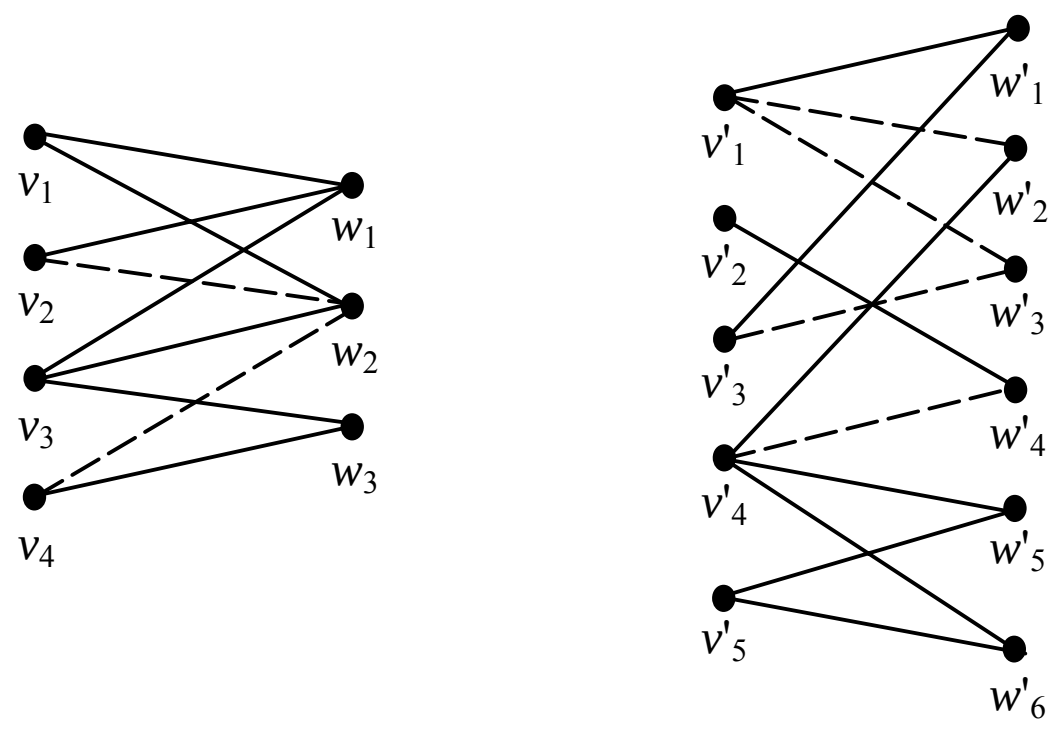

FIGURE: 4.2.10

Clearly both the bipartite graph $G_1$ and $G_2$ are neutrosophic, hence $G = G_1 \cup G_2$ is a neutrosophic bipartite bigraph, which are disconnected. Thus this is an example of a disconnected neutrosophic bipartite bigraph. Now we proceed on to give an example of a weak neutrosophic bipartite bigraph.

***Example 4.2.12:*** Let $G = G_1 \cup G_2$ be a bigraph given by the following figure 4.2.11.

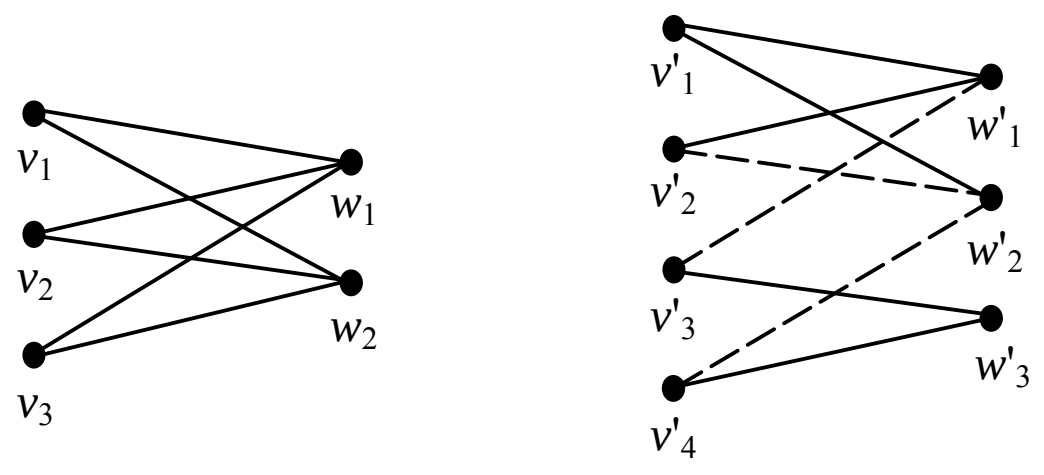

FIGURE: 4.2.11

Clearly $G = G_1 \cup G_2$ is a bipartite bigraph $G_1$ is a bipartite graph which is not neutrosophic but $G_2$ is a bipartite graph which is neutrosophic.

So $G = G_1 \cup G_2$ is a *weak neutrosophic bipartite bigraph*. This bigraph is a disconnected weak neutrosophic bipartite bigraph.

Now we give the following Theorem.

**THEOREM 4.2.4:** *Let $G = G_1 \cup G_2$ be a neutrosophic bipartite bigraph, then G is a weak neutrosophic bipartite bigraph. A weak neutrosophic bipartite bigraph in general need not be a neutrosophic bipartite bigraph.*

*Proof:* By the definition of neutrosophic bipartite bigraph $G = G_1 \cup G_2$ we see both $G_1$ and $G_2$ are neutrosophic bipartite graphs so G is also weak neutrosophic bipartite bigraph. But a weak neutrosophic bipartite bigraph in general need not be a neutrosophic bipartite bigraph.

This is established by the following example.

***Example 4.2.13:*** Consider the neutrosophic bigraph $G = G_1 \cup G_2$ given by the following figure 4.2.12.

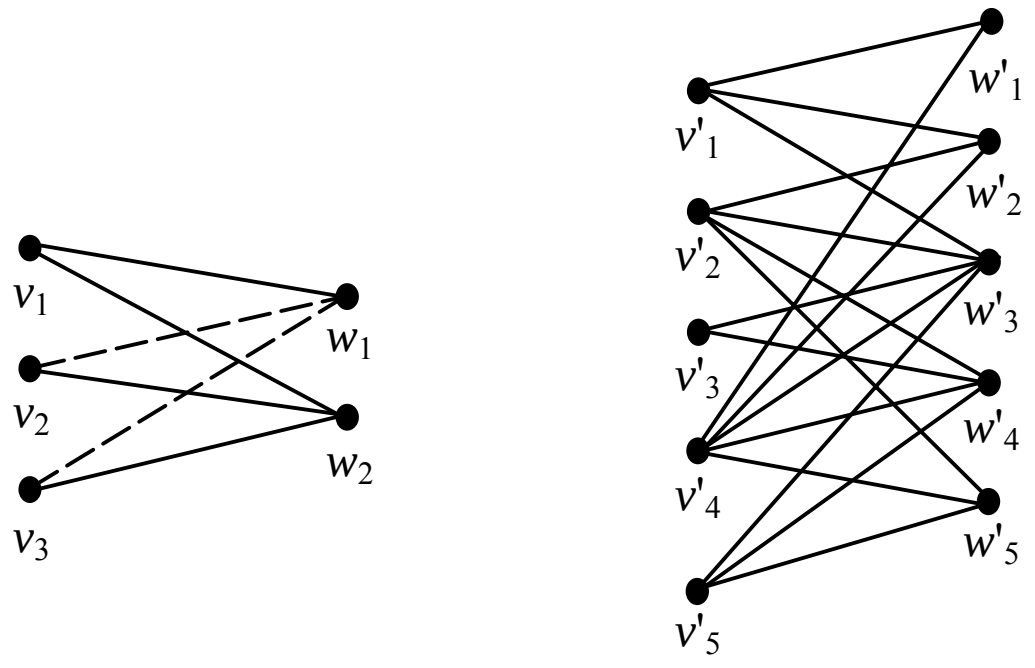

FIGURE: 4.2.12

Clearly $G_1$ is a neutrosophic bipartite graph but $G_2$ only a bipartite graph which is not neutrosophic so $G = G_1 \cup G_2$ is a bipartite bigraph which is not neutrosophic, which is only weak neutrosophic.

Now we have seen only disconnected neutrosophic bipartite bigraphs and weak neutrosophic bipartite bigraphs. Now we define a vertex connected neutrosophic bipartite bigraph and weak neutrosophic bipartite bigraph.

A vertex connected neutrosophic bipartite bigraph $G = G_1 \cup G_2$ is a bipartite bigraph which has only one vertex common between $G_1$ and $G_2$.

***Example 4.2.14:*** Let $G = G_1 \cup G_2$ be a neutrosophic bipartite bigraph (weak neutrosophic bipartite bigraph) given by the following figure.

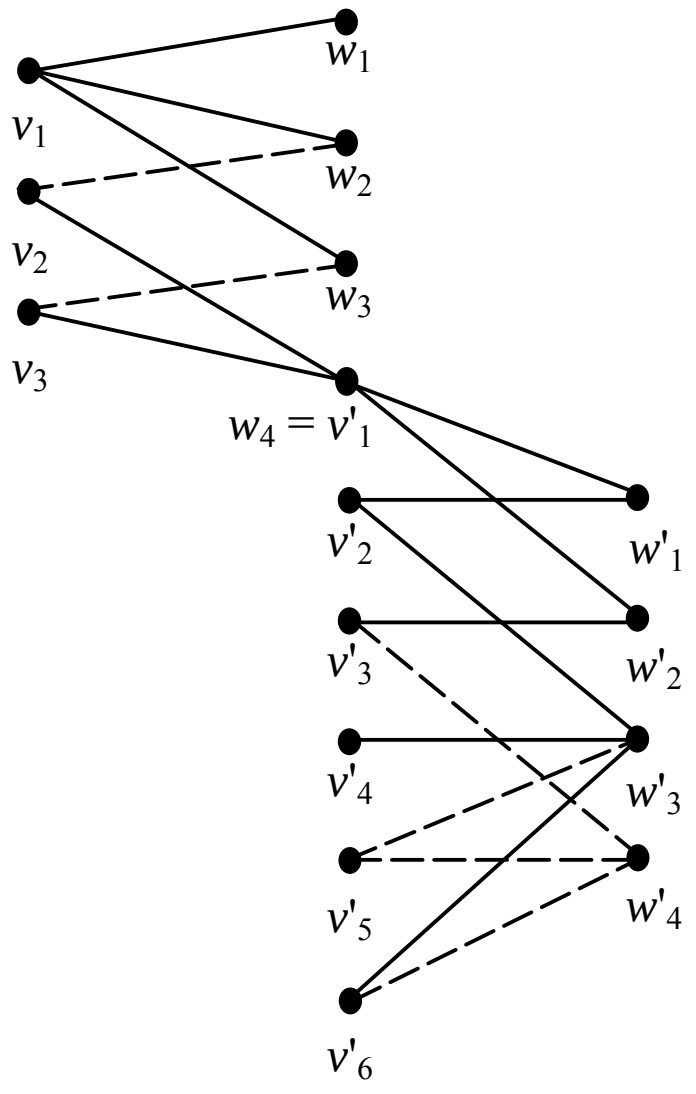

FIGURE: 4.2.13

It is clear from the figure the neutrosophic graphs $G_1$ and $G_2$ have only one vertex $w_4 = v'_1$ in common. Clearly $G = G_1 \cup G_2$ is a vertex connected neutrosophic bipartite bigraph.

***Example 4.2.15:*** $G = G_1 \cup G_2$ be a weak neutrosophic bipartite bigraph given by the following figure.

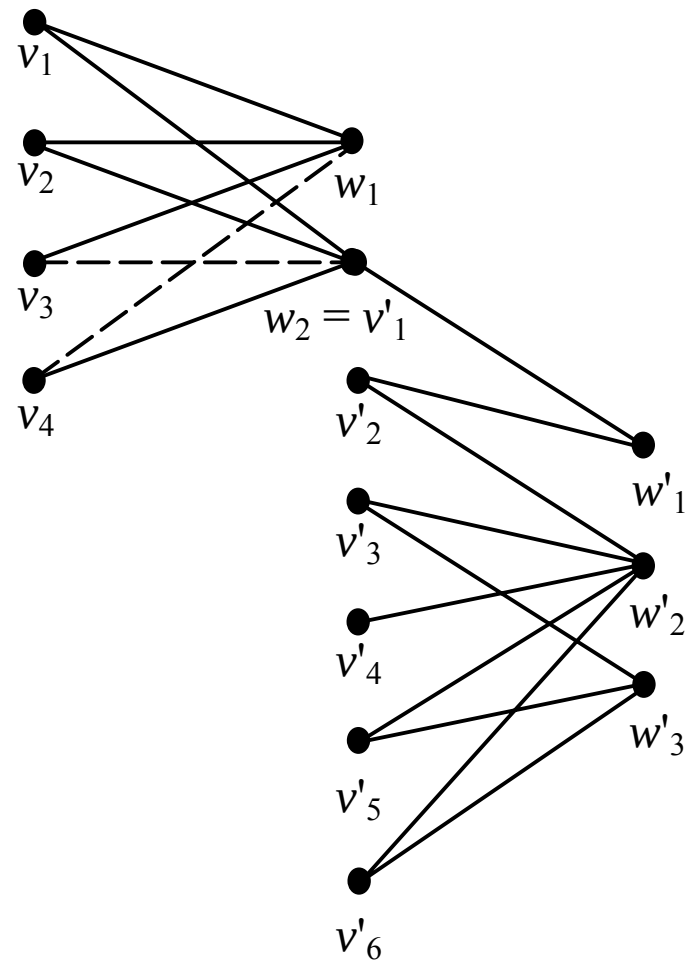

FIGURE: 4.2.14

$$G = G_1 \cup G_2.$$

This graph $G = G_1 \cup G_2$ is a weak neutrosophic bipartite bigraph which is vertex connected they have a vertex $v'_1 = w_2$ to be common.

We proceed on to give few more examples before we proceed to define biconnected neutrosophic bipartite bigraph.

***Example 4.2.16:*** Let $G = G_1 \cup G_2$ be a bigraph given by the following figure.

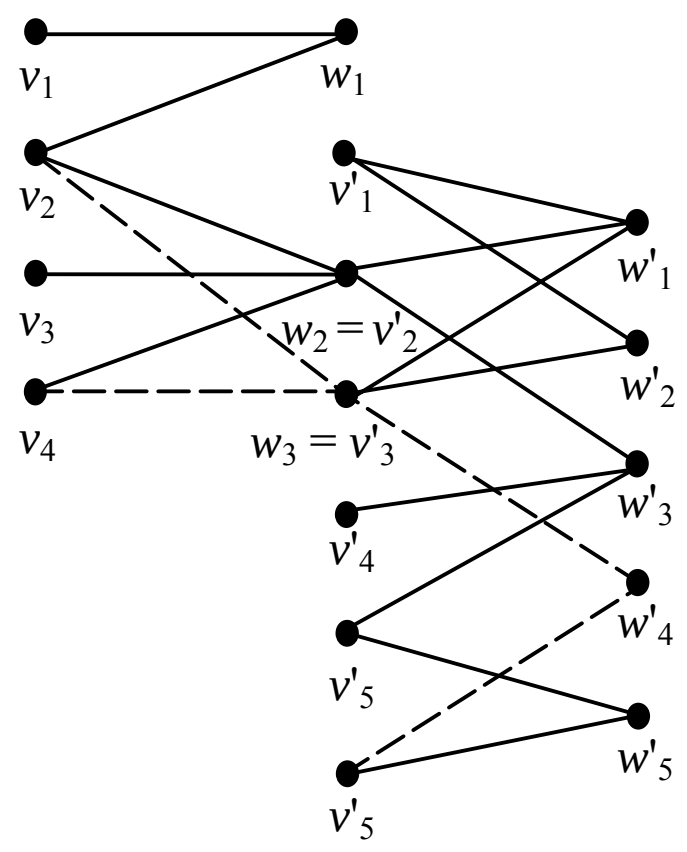

FIGURE: 4.2.15

This bigraph is also a connected neutrosophic bipartite bigraph. The vertices $v'_2 = w_2$ and $w_3 = v'_3$.

Next we give yet another example.

***Example 4.2.17:*** Let $G = G_1 \cup G_2$ be a bigraph given by the following figure.

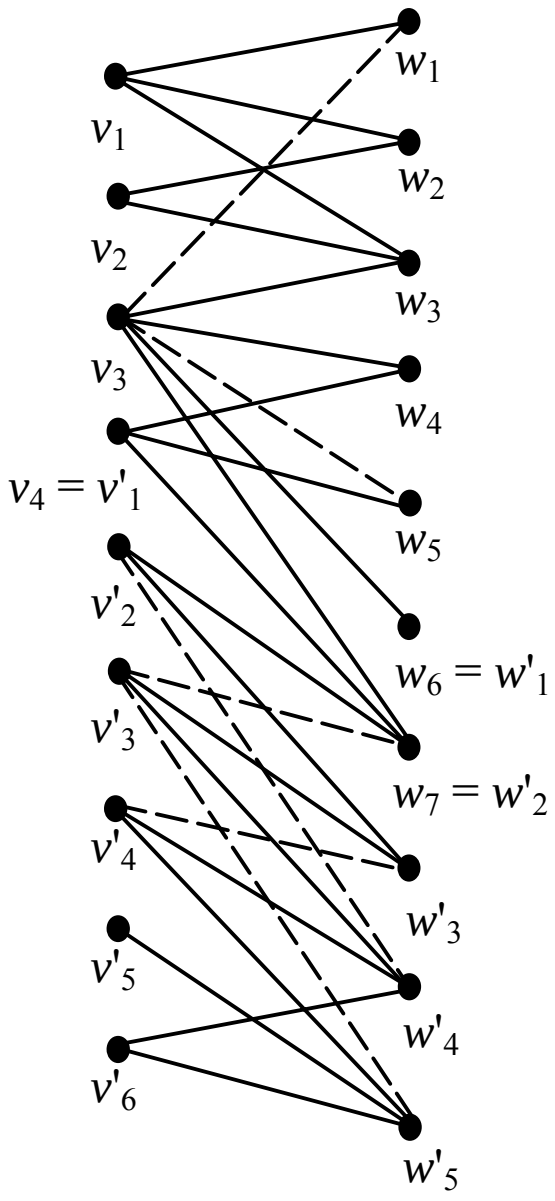

FIGURE: 4.2.16

Clearly $G = G_1 \cup G_2$ is a neutrosophic bipartite bigraph which is edge connected.

$$V(G_1) = (v_1, v_2, v_3, v_4, w_1, w_2, w_3, w_4, w_5, w_6, w_7)$$
$$V(G_2) = (v'_1 (= v_4), v'_2, v'_3, v'_4, v'_5, v'_6, w'_1 (= w_6), w'_2 (= w_7), w'_3, w'_4, w'_5).$$

**DEFINITION 4.2.9:** *Let $G = G_1 \cup G_2$ be a neutrosophic bipartite bigraph, with vertex set of*

$$G_1 = \{(v_1, \ldots, v_n), (w_1, \ldots, w_r)\}$$

*and*

$$G_2 = \{(v'_1 (=w_1), v'_2 (=w_2), \ldots, v'_r (= w_r)), (w'_1, w'_2, \ldots, w'_s)\}.$$

*where $G_1$ and $G_2$ are bipartite graphs, then G is a biconnected neutrosophic bipartite bigraph. i.e. the diagram of a biconnected neutrosophic bipartite bigraph is given in the following figure:*

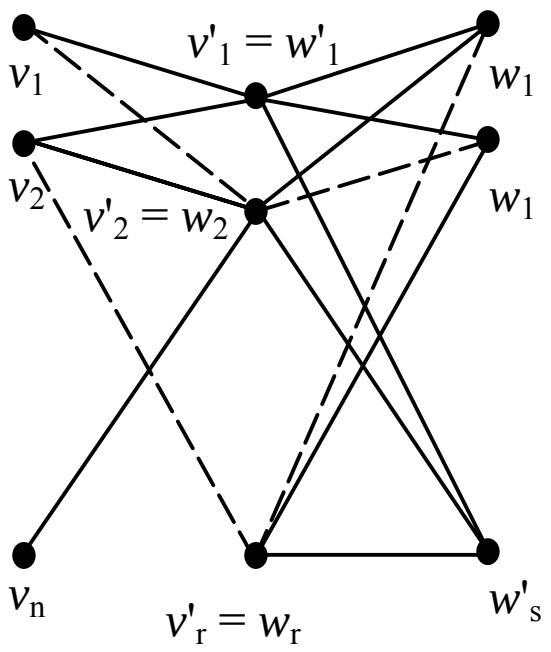

FIGURE: 4.2.17

$$G = G_1 \cup G_2.$$

Now we give a concrete example of a biconnected neutrosophic bipartite bigraph.

***Example 4.2.18:*** Let G = $G_1 \cup G_2$ be a bipartite bigraph given by the following figure .

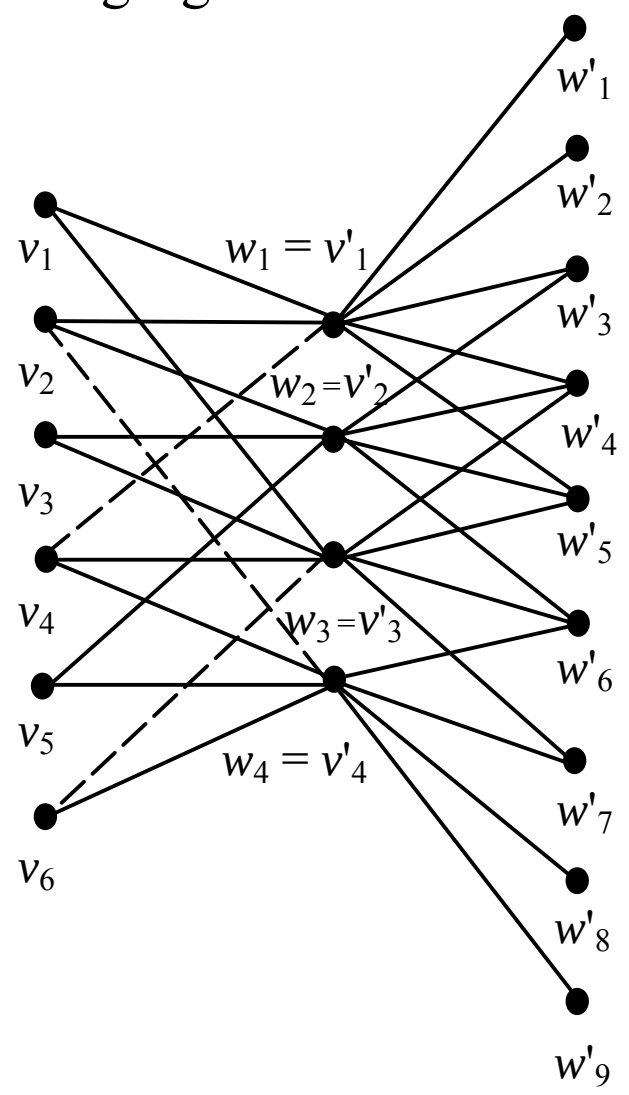

FIGURE: 4.2.18

Clearly G = $G_1 \cup G_2$ is a biconnected neutrosophic bipartite bigraph.

These neutrosophic bipartite bigraph will find its application in FRBMs. The adjacency neutrosophic bimatrix, the incidence neutrosophic bimatrix and weighted bigraph's neutrosophic bimatrix can be obtained for any given neutrosophic bigraph as in case of bigraphs.

**Remark:**

(1) It is important to note that bigraphs cannot always be realized as graphs. Bigraphs are under special conditions understood to be a union of two graphs provided they are subgraph connected.
(2) All bigraphs are graphs.

Now we define the notion of neutrosophic trigraphs and neutrosophic bipartite trigraphs and their generalizations.

**DEFINITION 4.2.10:** *$G = G_1 \cup G_2 \cup G_3$ is a neutrosophic trigraph if $G_1$, $G_2$ and $G_3$ are neutrosophic graphs such that the vertex sets of $G_1$ $G_2$ and $G_3$ are proper subsets of G they are not subsets of the other i.e. V ($G_1$) is not a subset of V($G_2$) or V ($G_3$), V ($G_2$) is not a subset of V ($G_1$) or V ($G_3$) and V ($G_3$) is not a subset of V ($G_1$) or V ($G_2$) or at least there is one edge which is not common with $G_i$ and $g_j$, $1 \le i \le 3$ and $1 \le j \le 3$.*

We now illustrate a trigraph by the following example.

***Example 4.2.19:*** Let $G = G_1 \cup G_2 \cup G_3$ be a neutrosophic trigraph given by the following figure.

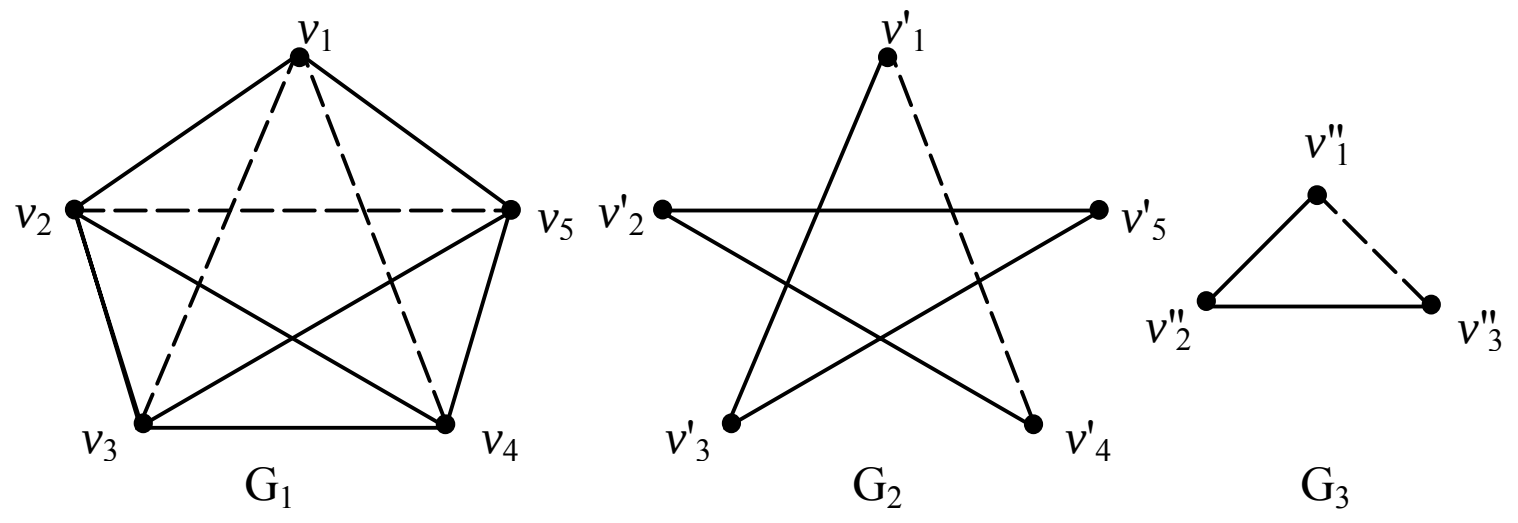

FIGURE: 4.2.19

$$G = G_1 \cup G_2 \cup G_3.$$

Clearly $G_1$, $G_2$ and $G_3$ are neutrosophic graphs i.e. all the three graphs are disjoint.

***Example 4.2.20:*** Let $G = G_1 \cup G_2 \cup G_3$ be a neutrosophic trigraph given by the following figure.

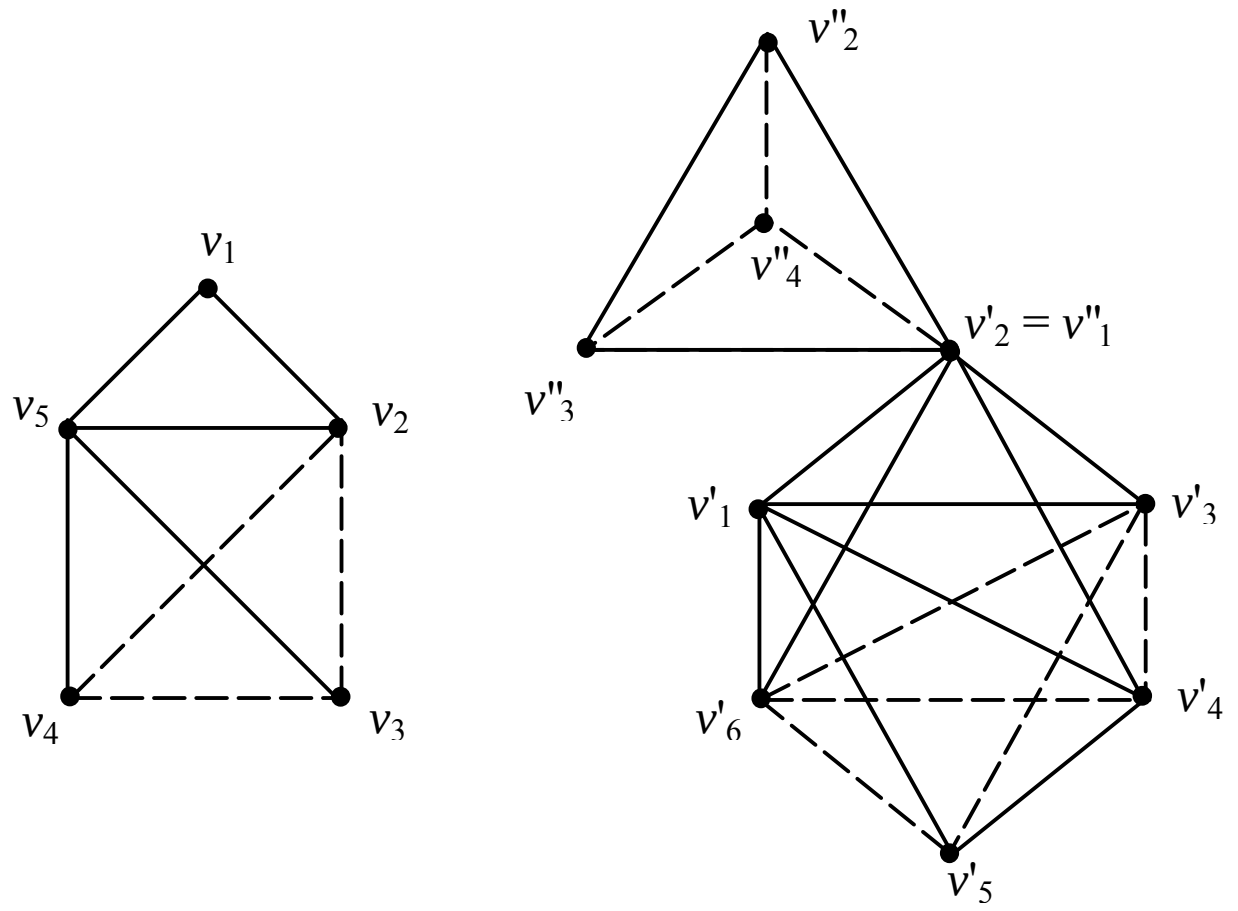

FIGURE: 4.2.20

$$G = G_1 \cup G_2 \cup G_3.$$

Clearly we see $G = G_1 \cup G_2 \cup G_3$ is a disconnected neutrosophic trigraph but it has connected neutrosophic bigraph.

So a neutrosophic trigraph can also be realized as a union of a neutrosophic graph and a neutrosophic bigraph. i.e. if $G = G_1 \cup B$ is trigraph if $G_1$ is a neutrosophic graph and B is a neutrosophic bigraph such that vertex set of $G_1$ is not fully contained in vertex set of B and vice versa.

It is left as an exercise for the reader to show both the definitions are equivalent.

Now we get yet another example of a neutrosophic trigraph.

***Example 4.2.21:*** Let G = $G_1 \cup G_2 \cup G_3$ be a neutrosophic trigraph given by the following figure.

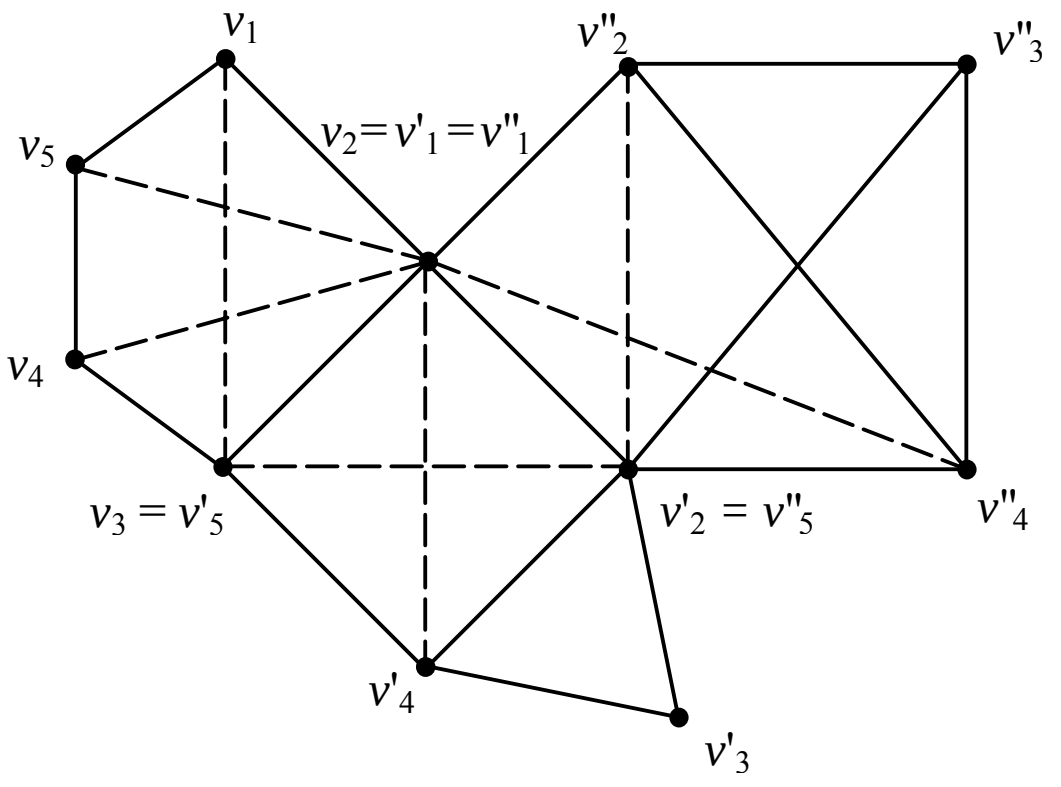

FIGURE: 4.2.21

$$G = G_1 \cup G_2 \cup G_3$$

$$V(G) = (v_1, v_2, v_3, v_4, v_5\} \cup \{v'_1 (= v_2), v'_2, v'_3, v'_4, v'_5 (= v_3)\} \cup \{v''_1 (v'_1 = v_2), v''_2, v''_3, v''_4, v''_5 (= v'_2)\}.$$

The separate neutrosophic graphs related with $G_1$, $G_2$ and $G_3$ are given by the following figure.

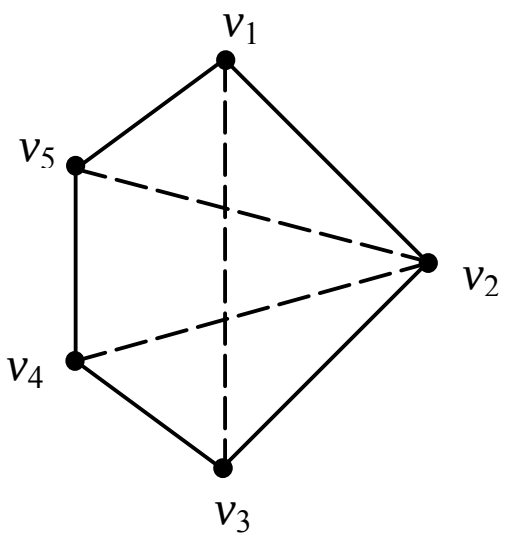

$G_1$

FIGURE 4.2.21a

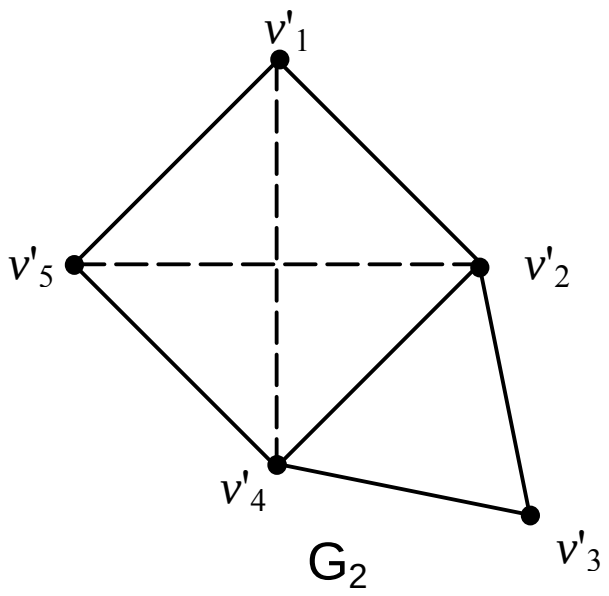

$G_2$

FIGURE 4.2.21b

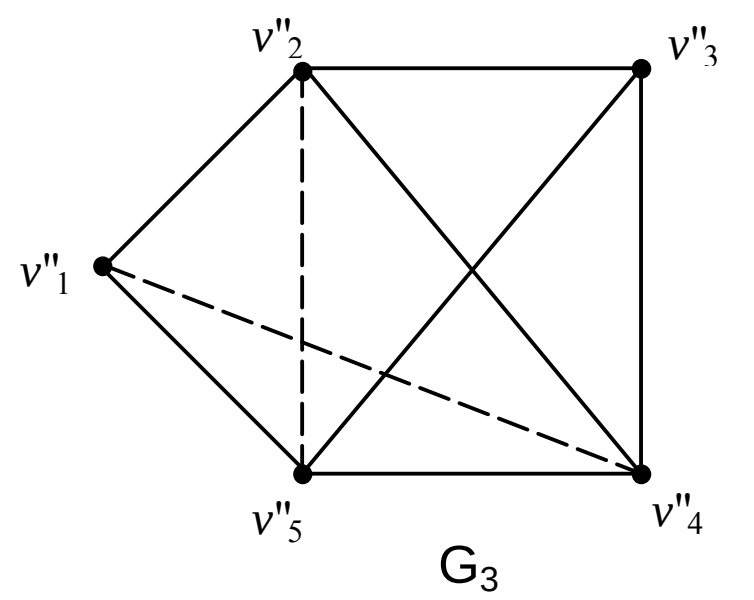

$G_3$

FIGURE 4.2.21c

Clearly $G = G_1 \cup G_2 \cup G_3$ is a connected neutrosophic trigraph. For we have for any neutrosophic trigraph a neutrosophic trimatrix associated with it.

We now define weak neutrosophic trigraph and very weak neutrosophic trigraph.

**DEFINITION 4.2.11:** *Let $G = G_1 \cup G_2 \cup G_3$ be a trigraph if only any two of the trigraphs from $G_1$, $G_2$ and $G_3$ are neutrosophic and one of $G_1$ (or $G_2$ or $G_3$) is not neutrosophic then we say G is a weak neutrosophic trigraph. If only one of $G_1$ or $G_2$ or $G_3$ alone is a neutrosophic graph say $G_1$ and the other two graphs $G_2$ and $G_3$ are not neutrosophic graphs then $G = G_1 \cup G_2 \cup G_3$ is a very weak neutrosophic trigraph.*

We now give examples of such trigraphs.

***Example 4.2.22:*** Let $G = G_1 \cup G_2 \cup G_3$ be a trigraph given by the following figure 4.2.22.

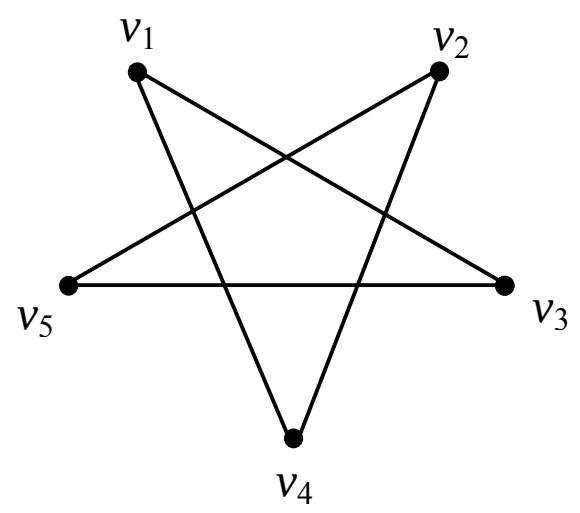

$G_1$

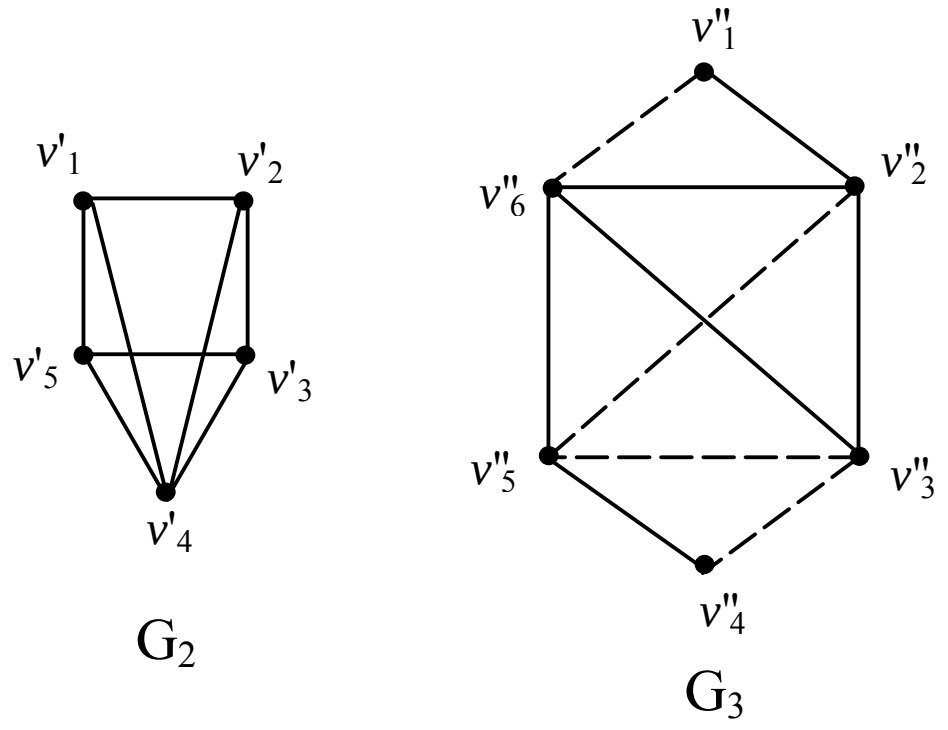

FIGURE 4.2.22

G is a weak neutrosophic trigraph.

We now proceed onto show by an example a trigraph which is a very weak neutrosophic trigraph.

***Example 4.2.23:*** Let $G = G_1 \cup G_2 \cup G_3$ be a trigraph in which only $G_3$ is a neutrosophic graph and $G_1$ and $G_2$ are just graphs given by the following figure 4.2.23:

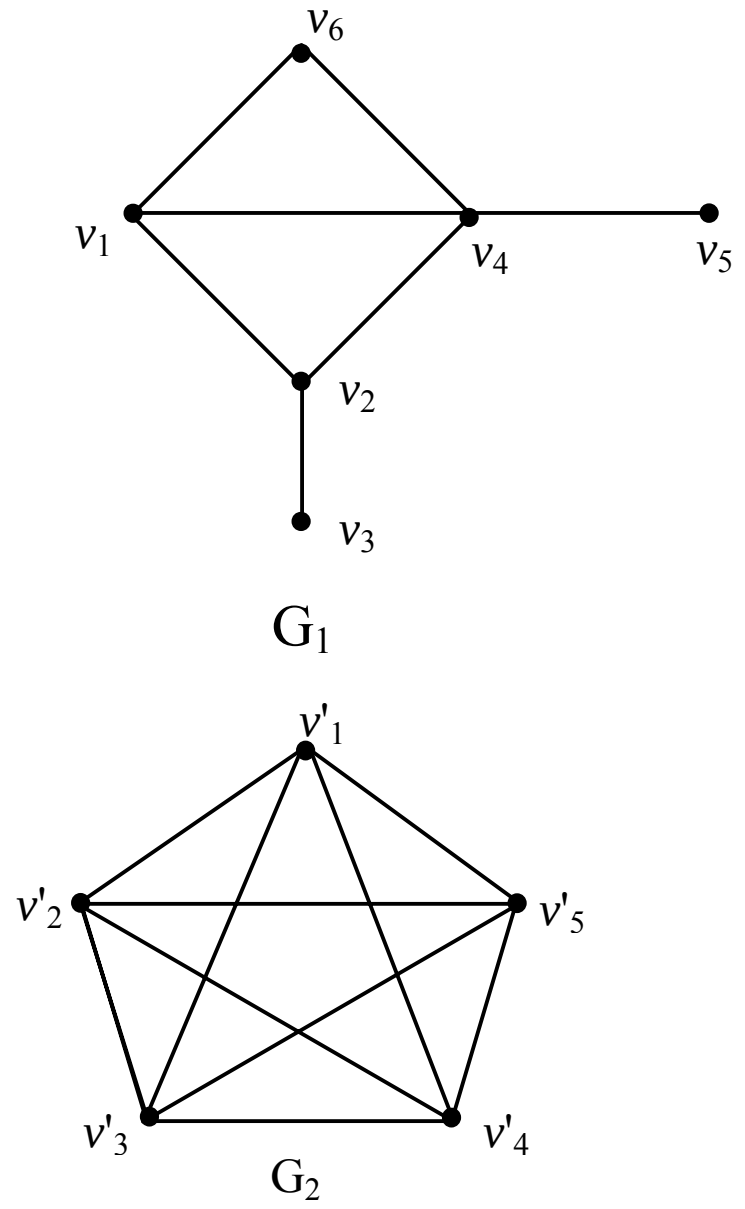

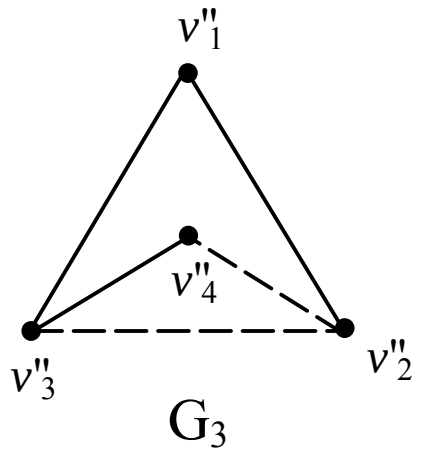

$G_3$

FIGURE: 4.2.23

Clearly G is a very weak neutrosophic trigraph.

The following result is left as an exercise for the reader to prove.

1. A weak neutrosophic trigraph is not in general a neutrosophic trigraph.
2. A very weak neutrosophic trigraph is not in general a weak neutrosophic trigraph. A very weak neutrosophic trigraph is not in general a neutrosophic trigraph.

A trigraph is said to be disconnected if $G = G_1 \cup G_2 \cup G_3$ are such that

$$V(G_i) \cap V(G_j) = \phi,\ i \neq j\ 1 \leq i, j \leq 3.$$

A trigraph is said to be connected if atleast

$$V(G_1) \cap V(G_2) \neq \phi \text{ and}$$
$$V(G_2) \cap V(G_3) \neq \phi$$

It can also happen $V(G_i) \cap V(G_j) \neq \phi$ for $1 \leq i, j \leq 3$.

Now we illustrate a disconnected neutrosophic trigraph.

***Example 4.2.24:*** Let $G = G_1 \cup G_2 \cup G_3$ be the trigraph given by the following figure 4.2.24:

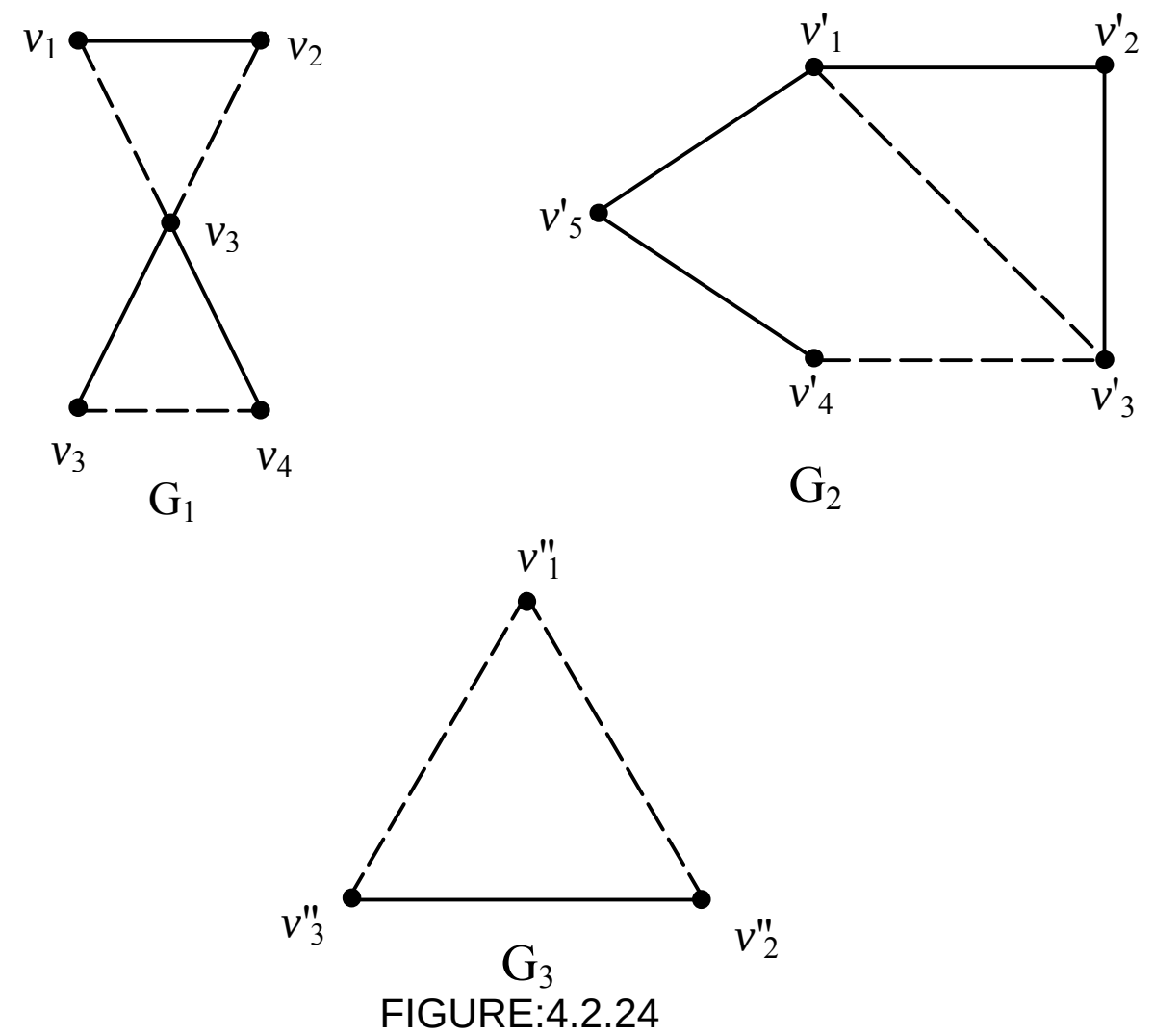

FIGURE:4.2.24

Clearly $G = G_1 \cup G_2 \cup G_3$ is a disconnected neutrosophic trigraph.

***Example 4.2.25:*** Let $G = G_1 \cup G_2 \cup G_3$ be trigraph given by the following figure 4.2.25.

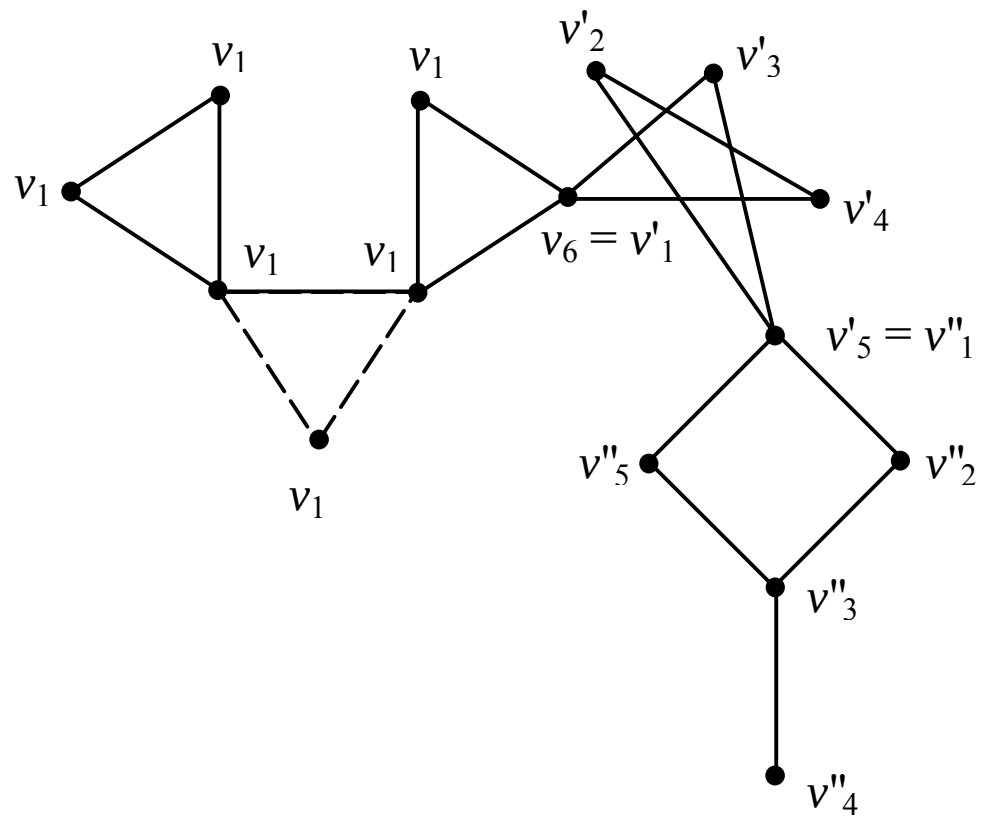

FIGURE: 4.2.25

Clearly $G = G_1 \cup G_2 \cup G_3$ is a neutrosophic trigraph which is connected.

Now we define a weakly connected neutrosophic trigraph.

**DEFINITION 4.2.12:** *Let $G = G_1 \cup G_2 \cup G_3$ be a neutrosophic trigraph if only a pair of graphs $G_1$ and $G_2$ or $G_2$ and $G_3$ or $G_3$ and $G_1$ alone are connected and the other is disconnected we call G a weakly connected neutrosophic trigraph.*

Now we illustrate by an example a weakly connected neutrosophic trigraph.

***Example 4.2.26:*** Let $G = G_1 \cup G_2 \cup G_3$ be a trigraph given by the following figure.

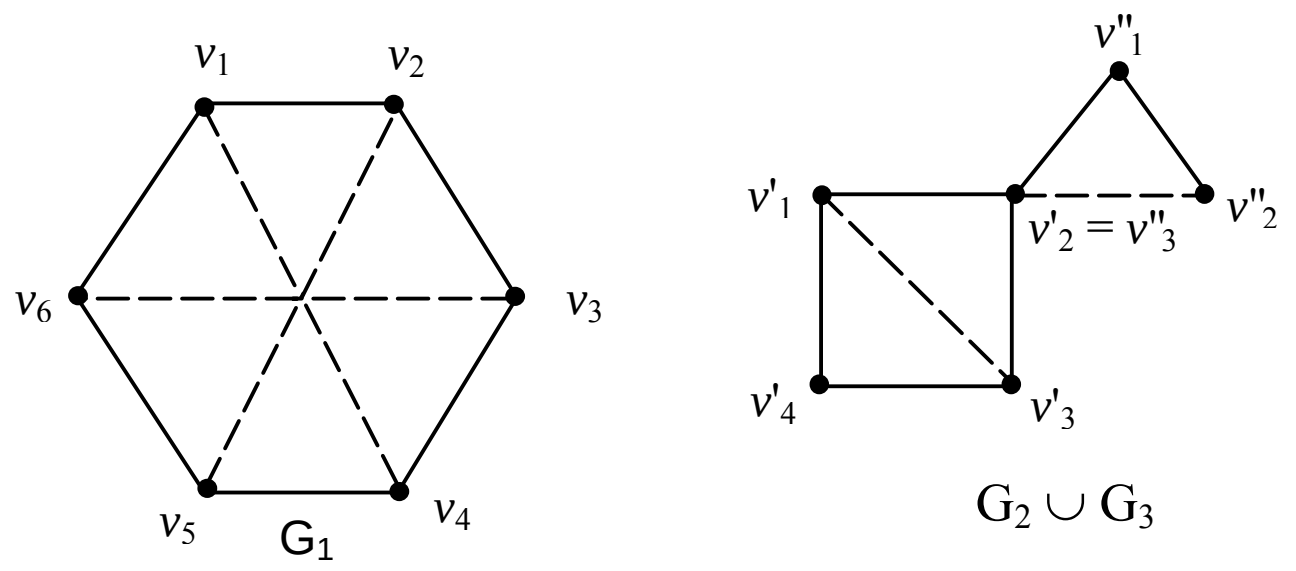

FIGURE 4.2.26

This trigraph is a weakly connected neutrosophic trigraph.

Now we proceed on to define the notion of bipartite trigraph which is neutrosophic.

**DEFINITION 4.2.13:** *Let $G = G_1 \cup G_2 \cup G_3$ be a neutrosophic trigraph if each of $G_1$, $G_2$ and $G_3$ are bipartite neutrosophic graphs then we call G a neutrosophic bipartite trigraph.*

We illustrate this by a simple example.

***Example 4.2.27:*** Consider the neutrosophic trigraph G = $G_1 \cup G_2 \cup G_3$ given by the following figure.

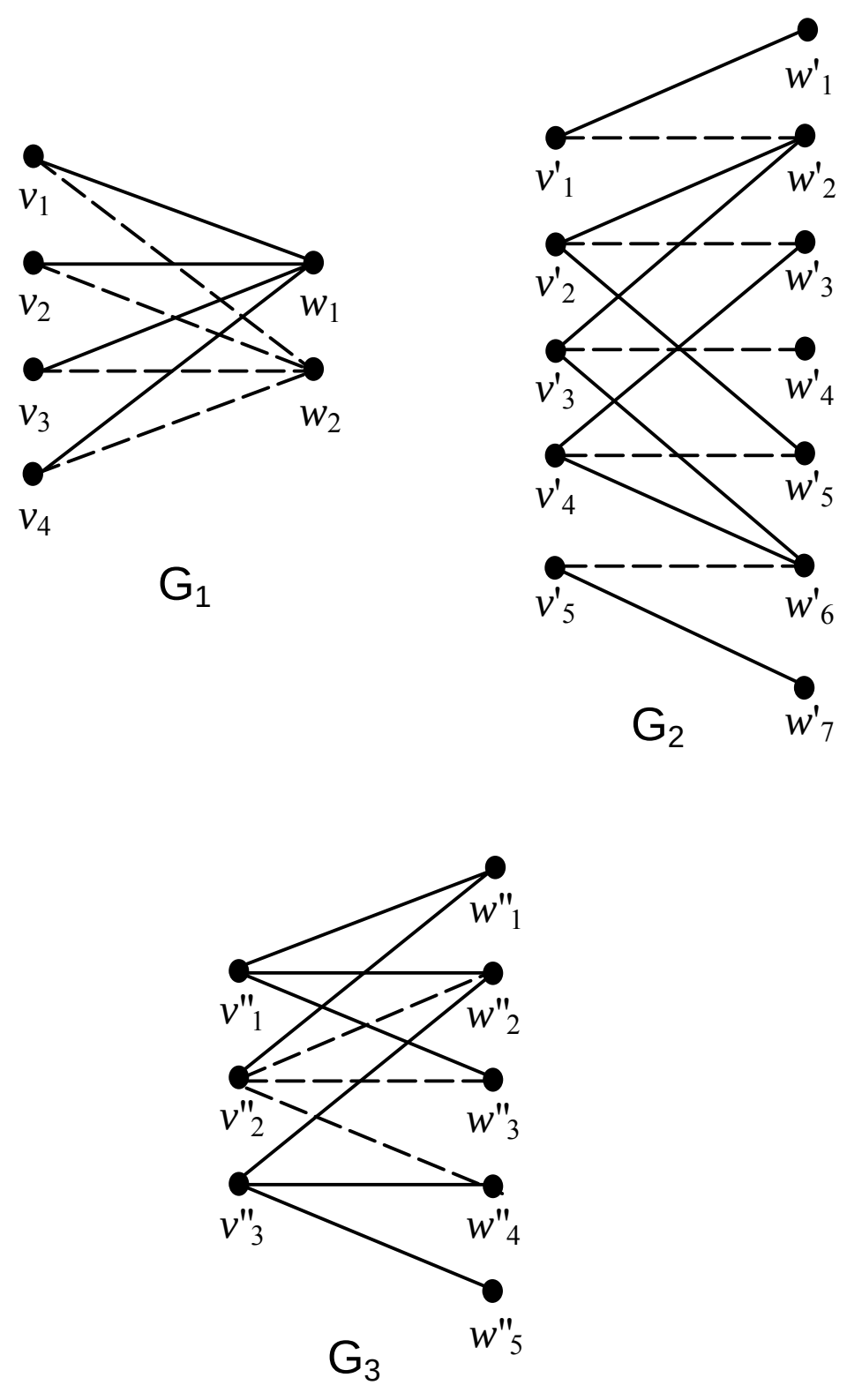

FIGURE: 4.2.27

G = $G_1 \cup G_2 \cup G_3$ is clearly a neutrosophic bipartite trigraph. It is clear from the figure $G_1$, $G_2$ and $G_3$ are neutrosophic bipartite graph.

***Example 4.2.28:*** Let G = $G_1 \cup G_2 \cup G_3$ be a neutrosophic trigraph given by the following figure.

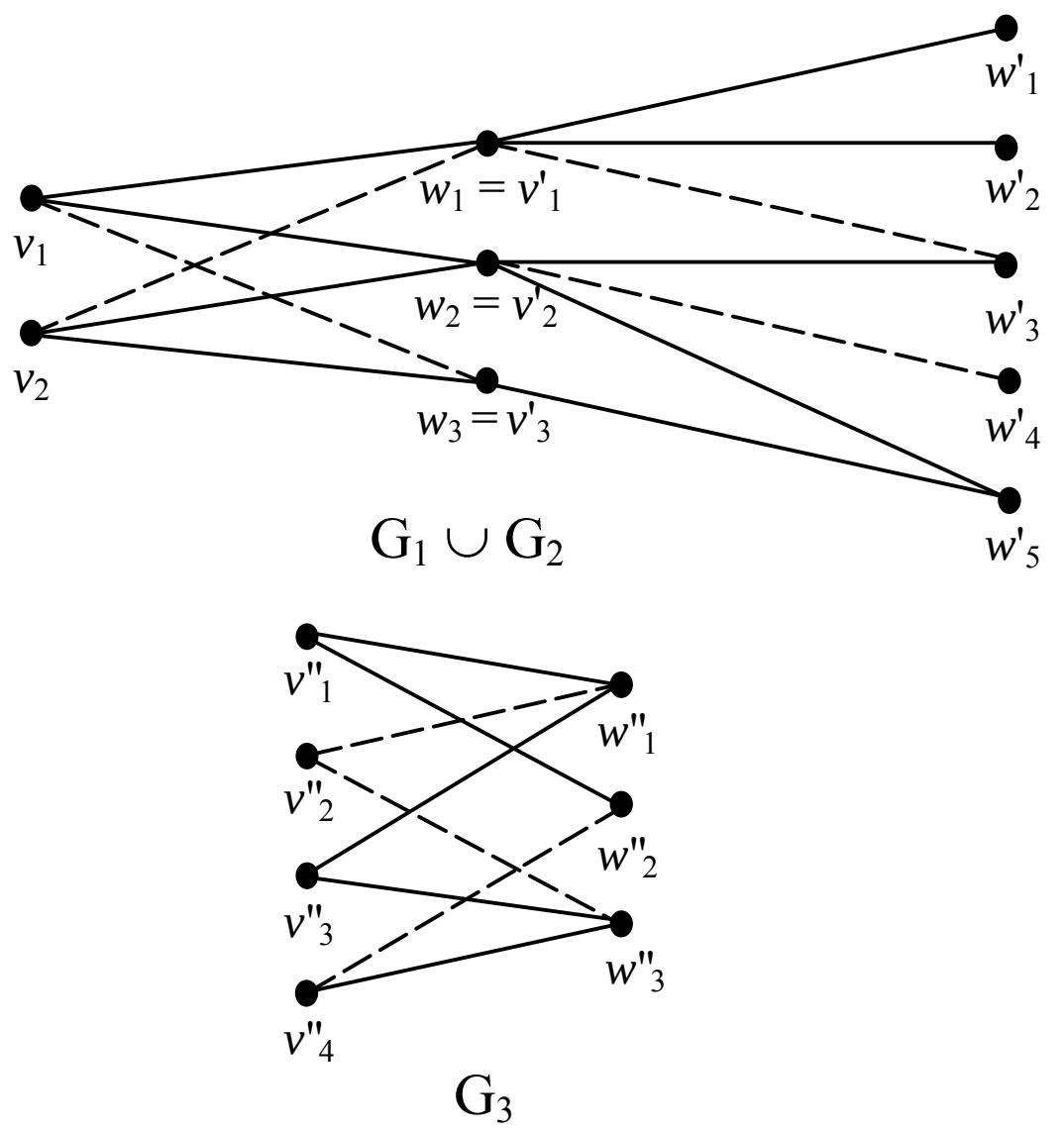

FIGURE: 4.2.28

It is easily seen the trigraph $G = G_1 \cup G_2 \cup G_3$ is a neutrosophic bipartite trigraph which is weakly connected.

Now we proceed on to give one more example of a neutrosophic bigraph.

***Example 4.2.29:*** Let G be a trigraph given by the following figure:

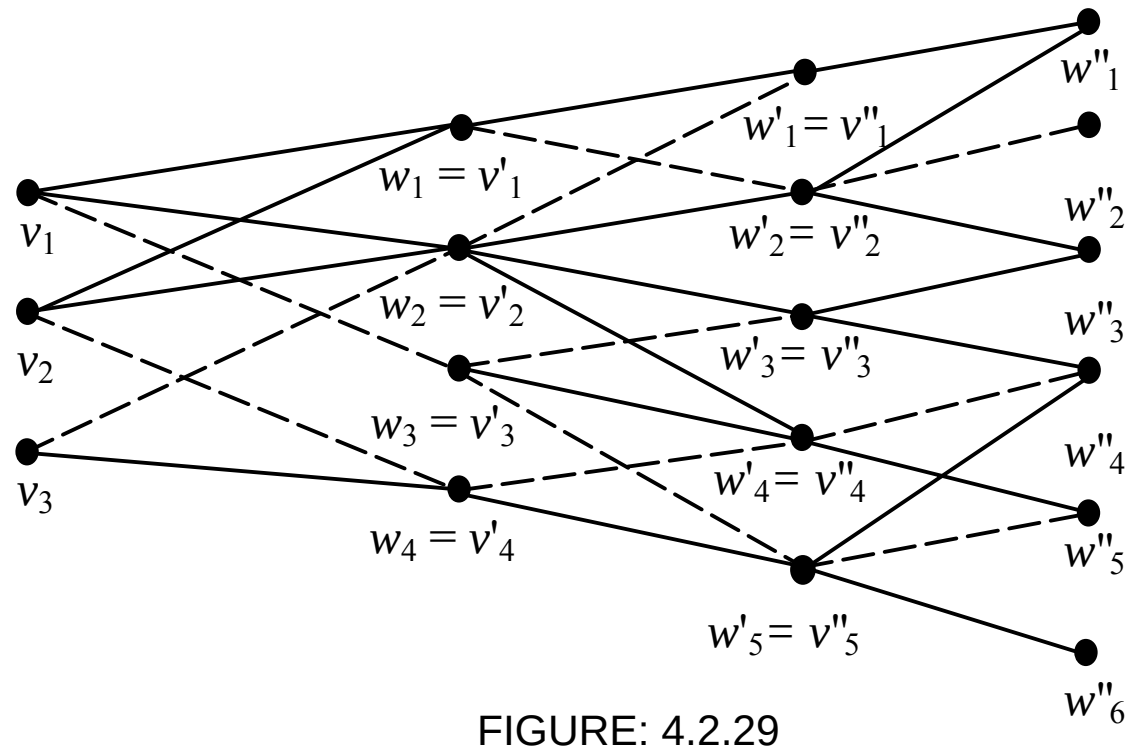

FIGURE: 4.2.29

$$G = G_1 \cup G_2 \cup G_3$$
is a neutrosophic bipartite trigraph which is connected.

We can proceed on to define weak neutrosophic bipartite trigraph as a neutrosophic trigraph $G = G_1 \cup G_2 \cup G_3$ in which one of $G_1$ or $G_2$ or $G_3$ is not a neutrosophic graph but a bipartite graph. So the definitions carry without any problem in case of weak neutrosophic bipartite trigraphs.

Before we go for application we would be defining a few more new concepts in graph theory. They are the concept of neutrosophic tripartite graphs and neutrosophic n-partite graphs analogous to neutrosophic bipartite graphs.

**DEFINITION 4.2.14:** *Let G be a neutrosophic graph, G is said to be a neutrosophic tripartite if its vertex set V can be decomposed into 3 disjoint subsets $V_1$, $V_2$, $V_3$ such that every edge in G join vertex in $V_1$ with a vertex $V_2$ and vertex of $V_2$ with a vertex of $V_3$ or the vertex of $V_3$ with a vertex of $V_1$ and atleast there is a edge connecting $V_1$ to $V_2$ and $V_2$ to $V_3$ or $V_3$ to $V_1$ as a neutrosophic edge.*

***Example 4.2.30:*** Let G be a as graph given by the following figure:

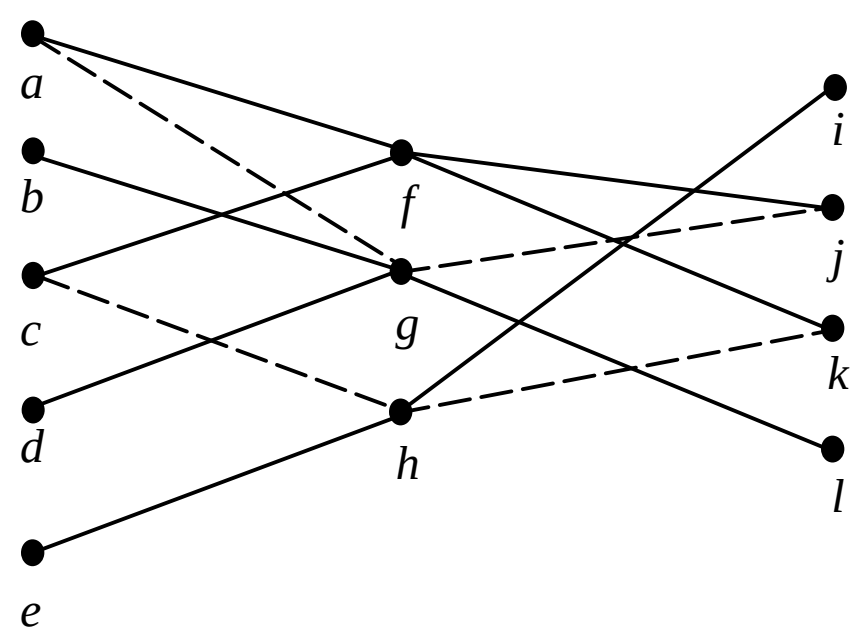

FIGURE: 4.2.30

The vertex set of G is {a, b, c, d, e, f, g, h, *I*, j, k, l} the vertex V (G) is divided into three disjoint classes viz. $V_1$ = {a, b, c, d, e} $V_2$ = {f, g, h} and $V_3$ = {i, j, k, l}.

Thus a neutrosophic tripartite graph can have no self loop. Clearly G is a tripartite graph which is neutrosophic. We see when the edge is itself an indeterminate one, thinking of colouring it etc happens to be an intricate one.

In generalizing this concept a neutrosophic graph G is called p-partite if its vertex set can be decomposed into p-disjoint subsets $V_1, V_2, \ldots, V_p$ such that no edge in G joins the vertex in the same subset. We illustrate a 5-partite neutrosophic graph by the following example.

***Example 4.2.31:*** Let G be a neutrosophic graph given by the figure .

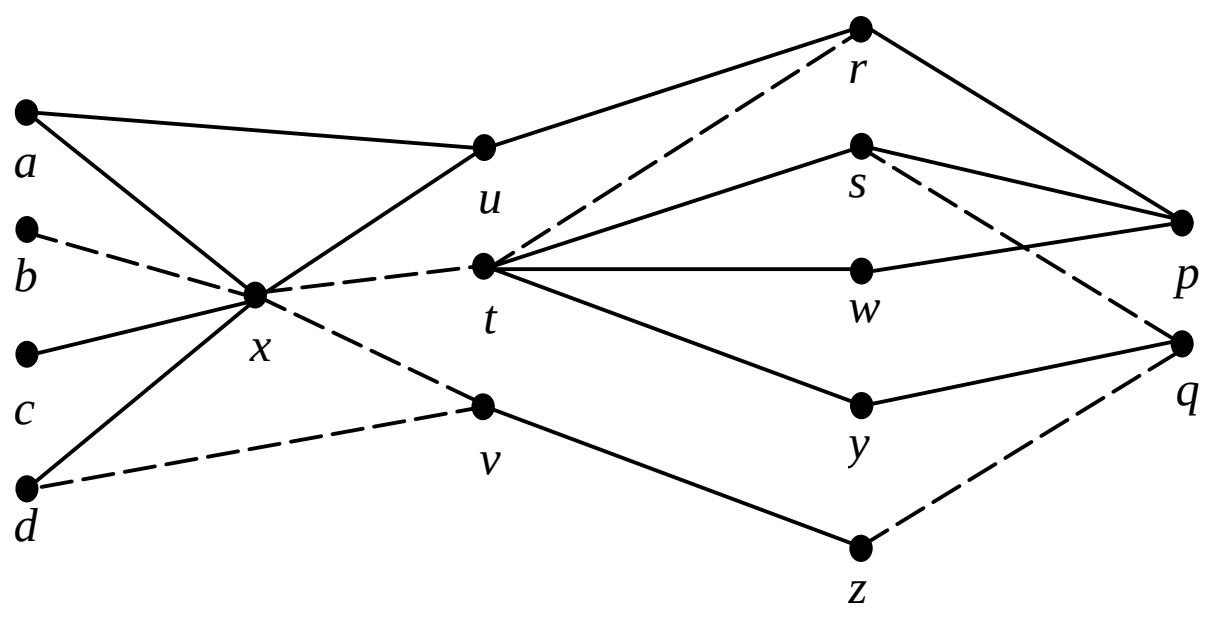

FIGURE: 4.2.31

This has 5 vertex sets $V_1, V_2, \ldots, V_5$ given by $V_1$ = {a, b, c, d}, $V_2$ = {x}, $V_3$ = {u, t, v}, $V_4$ = {r, s, w, y, z} and $V_5$ = {p, q}. Thus G is a 5-partite neutrosophic graph.

The concept of neutrosophic n-partite graph and general p-partite graph will be very useful in the applications of FRn-maps.

Thus we see these models will be drawn with the help of these graphs as the nodes can be classified as the disjoint classes that is no two nodes in the same class are joined by an edge or a neutrosophic edge. All these neutrosophic trigraphs have a neutrosophic trimatrix associated with it.

The matrices may be mixed square trimatrix or a mixed rectangular trimatrix or a just square neutrosophic trimatrix or a rectangular neutrosophic trimatrix.

We now proceed on to give how these neutrosophic bimatrix or neutrosophic trimatrix find its application in neutrosophic models.

## 4.3 Neutrosophic Cognitive Bimaps and their Generalizations

Now for the first time we introduce the Neutrosophic models in which the concept of neutrosophic bigraph and neutrosophic trigraphs find their applications. In a special class of models the notion of neutrosophic bipartite bigraphs and neutrosophic bipartite trigraphs would be used.

We have already defined the notion of neutrosophic cognitive maps in [151] and its generalizations viz. neutrosophic relational maps and further divisions.

The notion of Fuzzy Cognitive Maps (FCMs) which are fuzzy signed directed graphs with feedback are discussed and described. The directed edge $e_{ij}$ from causal concept $C_i$ to concept $C_j$ measures how much $C_i$ causes $C_j$. The time varying concept function $C_i(t)$ measures the non negative occurrence of some fuzzy event, perhaps the strength of a political sentiment, historical trend or opinion about some topics like child labor or school dropouts etc. FCMs model the world as a collection of classes and causal relations between them.

The edge $e_{ij}$ takes values in the fuzzy causal interval [–1, 1] ($e_{ij} = 0$ indicates no causality, $e_{ij} > 0$ indicates causal increase; that $C_j$ increases as $C_i$ increases and $C_j$ decreases as $C_i$ decreases, $e_{ij} < 0$ indicates causal decrease or negative causality $C_j$ decreases as $C_i$ increases or $C_j$ increases as $C_i$ decreases. Simple FCMs have edge value in {–1, 0, 1}. Thus if causality occurs it occurs to maximal positive or negative degree.

*I*t is important to note that $e_{ij}$ measures only absence or presence of influence of the node $C_i$ on $C_j$ but till now any researcher has not contemplated the indeterminacy of any relation between two nodes $C_i$ and $C_j$. When we deal with unsupervised data, there are situations when no relation can be determined between some two nodes. So in this section we try to introduce the indeterminacy in FCMs, and we choose to call this generalized structure as Neutrosophic Cognitive Maps (NCMs). In our view this will certainly give a more appropriate result and also caution the user about the risk of indeterminacy.

Now we proceed on to define the concepts about NCMs.

**DEFINITION 4.3.1:** *A Neutrosophic Cognitive Map (NCM) is a neutrosophic directed graph with concepts like policies, events etc. as nodes and causalities or indeterminates as edges. It represents the causal relationship between concepts.*

Let $C_1, C_2, \ldots, C_n$ denote n nodes, further we assume each node is a neutrosophic vector from neutrosophic vector space V. So a node $C_i$ will be represented by $(x_1, \ldots, x_n)$ where $x_k$'s are zero or one or $I$ where $I$ is the indeterminate and $x_k = 1$ means that the node $C_k$ is in the on state and $x_k = 0$ means the node is in the off state and $x_k = I$ means the nodes state is an indeterminate at that time or in that situation.

Let $C_i$ and $C_j$ denote the two nodes of the NCM. The directed edge from $C_i$ to $C_j$ denotes the causality of $C_i$ on $C_j$ called connections. Every edge in the NCM is weighted with a number in the set $\{-1, 0, 1, I\}$. Let $e_{ij}$ be the weight of the directed edge $C_iC_j$, $e_{ij} \in \{-1, 0, 1, I\}$. $e_{ij} = 0$ if $C_i$ does not have any effect on $C_j$, $e_{ij} = 1$ if increase (or decrease) in $C_i$ causes increase (or decreases) in $C_j$, $e_{ij} = -1$ if increase (or decrease) in $C_i$ causes decrease (or increase) in $C_j$ . $e_{ij} = I$ if the relation or effect of $C_i$ on $C_j$ is an indeterminate.

**DEFINITION 4.3.2:** *NCMs with edge weight from {-1, 0, 1, I} are called simple NCMs.*

**DEFINITION 4.3.3:** *Let $C_1, C_2, \ldots, C_n$ be nodes of a NCM. Let the neutrosophic matrix N(E) be defined as N(E) = ($e_{ij}$) where $e_{ij}$ is the weight of the directed edge $C_i C_j$, where $e_{ij} \in$*

*{0, 1, -1, I}. N(E) is called the neutrosophic adjacency matrix of the NCM.*

**DEFINITION 4.3.4:** *Let $C_1$, $C_2$, …, $C_n$ be the nodes of the NCM. Let A = ($a_1$, $a_2$,…, $a_n$) where $a_i \in$ {0, 1, I}. A is called the instantaneous state neutrosophic vector and it denotes the on – off – indeterminate state position of the node at an instant*

*$a_i$ = 0 if $a_i$ is off (no effect)*
*$a_i$ = 1 if $a_i$ is on (has effect)*
*$a_i$ = I if $a_i$ is indeterminate(effect cannot be determined)*

*for i = 1, 2,…, n.*

**DEFINITION 4.3.5:** *Let $C_1$, $C_2$, …, $C_n$ be the nodes of the FCM. Let $\overrightarrow{C_1C_2}$, $\overrightarrow{C_2C_3}$, $\overrightarrow{C_3C_4}$, … , $\overrightarrow{C_iC_j}$ be the edges of the NCM. Then the edges form a directed cycle. An NCM is said to be cyclic if it possesses a directed cyclic. An NCM is said to be acyclic if it does not possess any directed cycle.*

**DEFINITION 4.3.6:** *An NCM with cycles is said to have a feedback. When there is a feedback in the NCM i.e. when the causal relations flow through a cycle in a revolutionary manner the NCM is called a dynamical system.*

**DEFINITION 4.3.7:** *Let $\overrightarrow{C_1C_2}, \overrightarrow{C_2C_3}, \cdots, \overrightarrow{C_{n-1}C_n}$ be cycle, when $C_i$ is switched on and if the causality flow through the edges of a cycle and if it again causes $C_i$, we say that the dynamical system goes round and round. This is true for any node $C_i$, for i = 1, 2,…, n. the equilibrium state for this dynamical system is called the hidden pattern.*

**DEFINITION 4.3.8:** *If the equilibrium state of a dynamical system is a unique state vector, then it is called a fixed point. Consider the NCM with $C_1$, $C_2$,…, $C_n$ as nodes. For example let us start the dynamical system by switching on $C_1$. Let us assume that the NCM settles down with $C_1$ and $C_n$* *on, i.e. the state vector remain as (1, 0,…, 1) this neutrosophic state vector (1,0,…, 0, 1) is called the fixed point.*

**DEFINITION 4.3.9:** *If the NCM settles with a neutrosophic state vector repeating in the form*

$$A_1 \to A_2 \to \ldots \to A_i \to A_1,$$

*then this equilibrium is called a limit cycle of the NCM.*

**METHODS OF DETERMINING THE HIDDEN PATTERN:**

Let $C_1, C_2, \ldots, C_n$ be the nodes of an NCM, with feedback. Let E be the associated adjacency matrix. Let us find the hidden pattern when $C_1$ is switched on when an input is given as the vector $A_1 = (1, 0, 0, \ldots, 0)$, the data should pass through the neutrosophic matrix N(E), this is done by multiplying $A_1$ by the matrix N(E). Let $A_1N(E) = (a_1, a_2, \ldots, a_n)$ with the threshold operation that is by replacing $a_i$ by 1 if $a_i \geq k$ and $a_i$ by 0 if $a_i < k$ (k – a suitable positive integer) and $a_i$ by *I* if $a_i$ is not a integer. We update the resulting concept, the concept $C_1$ is included in the updated vector by making the first coordinate as 1 in the resulting vector. Suppose $A_1N(E) \to A_2$ then consider $A_2N(E)$ and repeat the same procedure. This procedure is repeated till we get a limit cycle or a fixed point.

**DEFINITION 4.3.10:** *Finite number of NCMs can be combined together to produce the joint effect of all NCMs. If $N(E_1), N(E_2), \ldots, N(E_p)$ be the neutrosophic adjacency matrices of a NCM with nodes $C_1, C_2, \ldots, C_n$ then the combined NCM is got by adding all the neutrosophic adjacency matrices $N(E_1), \ldots, N(E_p)$. We denote the combined NCMs adjacency neutrosophic matrix by $N(E) = N(E_1) + N(E_2) + \ldots + N(E_p)$.*

**NOTATION:** Let $(a_1, a_2, \ldots, a_n)$ and $(a'_1, a'_2, \ldots, a'_n)$ be two neutrosophic vectors. We say $(a_1, a_2, \ldots, a_n)$ is equivalent to $(a'_1, a'_2, \ldots, a'_n)$ denoted by $((a_1, a_2, \ldots, a_n) \sim (a'_1, a'_2, \ldots, a'_n)$ if $(a'_1, a'_2, \ldots, a'_n)$ is got after thresholding and updating the vector $(a_1, a_2, \ldots, a_n)$ after passing through the neutrosophic adjacency matrix N(E).

The following are very important:

***Note 1:*** The nodes $C_1, C_2, \ldots, C_n$ are not indeterminate nodes because they indicate the concepts which are well known. But the edges connecting $C_i$ and $C_j$ may be indeterminate i.e. an expert may not be in a position to say that $C_i$ has some causality on $C_j$ either will he be in a position to state that $C_i$ has no relation with $C_j$ in such cases the relation between $C_i$ and $C_j$ which is indeterminate is denoted by *I*.

***Note 2:*** The nodes when sent will have only ones and zeros i.e. on and off states, but after the state vector passes through the neutrosophic adjacency matrix the resultant vector will have entries from $\{0, 1, I\}$ i.e. they can be neutrosophic vectors.

The presence of *I* in any of the coordinate implies the expert cannot say the presence of that node i.e. on state of it after passing through N(E) nor can we say the absence of the node i.e. off state of it the effect on the node after passing through the dynamical system is indeterminate so only it is represented by *I*. Thus only in case of NCMs we can say the effect of any node on other nodes can also be indeterminates. Such possibilities and analysis is totally absent in the case of FCMs.

***Note 3:*** In the neutrosophic matrix N(E), the presence of *I* in the $a_{ij}^{th}$ place shows, that the causality between the two nodes i.e. the effect of $C_i$ on $C_j$ is indeterminate. Such chances of being indeterminate is very possible in case of unsupervised data and that too in the study of FCMs which are derived from the directed graphs.

Thus only NCMs helps in such analysis.

Now we shall represent a few examples to show how in this set up NCMs is preferred to FCMs. At the outset before we proceed to give examples we make it clear that all unsupervised data need not have NCMs to be applied to it. Only data which have the relation between two nodes to be an indeterminate need to be modeled with NCMs if the data has no indeterminacy factor between any pair of nodes one need not go for NCMs; FCMs will do the best job.

Now we proceed on to define the notion of Neutrosophic Cognitive bimaps (NCBMs).

**DEFINITION 4.3.11:** *A neutrosophic cognitive bimap (NCBM) is a neutrosophic directed bigraph with concepts like policies or events etc as nodes and causalities and indeterminate as edges.*

It represents the casual relationship between concepts.

Note the neutrosophic directed bigraph need not be always a disconnected bigraph. It can be a connected directed neutrosophic bigraph or a disconnected directed neutrosophic bigraph.

Let $\left\{C_1^1,...,C_n^1\right\}$ and $\left\{C_1^2,...,C_m^2\right\}$ be a set of n and m nodes, further we assume each node is a neutrosophic vector from a neutrosophic vector space V. So a node $C = C_1 \cup C_2$ will be represented by $\left\{x_1^1,...,x_n^1\right\} \cup \left\{x_1^2,...,x_m^2\right\}$ where $x_i^t$ are zero or one or $I$ ($I$ is the indeterminate) (t = 1, 2) and $x_i^t = 1$ means the node $C_i$ is in the on state, $x_i^t = 0$ means the node is in the off state $x_i^t = I$ means the node is in the indeterminate state at that time or in that situation or circumstance (For example we are studying about the behaviour of a naughty child in the presence of a stern teacher surrounded by his classmates, the node naughtiness at that movement in that situation is an indeterminate for by the nice behaviour of the child in that circumstance the expert cannot make conclusions that the child is not naughty, he can only say indeterminate without fully knowing about the child's nature from his parents or relatives. Thus the coordinate naughtiness cannot be given 1 as in that circumstance at that time he is so good, cannot be given 0 for the expert is not fully aware of the fact that he is naughty so the expert can say only indeterminate for in the presence of that particular teacher and that room he may be behaving good, might be his behaviour would be very different in the play ground with his group of friends…)

Let $C_i^t$ and $C_j^t$ denote the a pair of two nodes of the NCBM (t = 1, 2). The directed edges from $C_i^1$ to $C_j^1$ and $C_i^2$ to $C_j^2$ denotes the causality of $C_i^t$ on $C_j^t$ (t = 1, 2) called connections. Every edge in the NCBM is weighted with a number in the set {-1 0 1 *I*}. Let $e_{ij}^t, (t = 1, 2)$ be the weight of the directed edges $C_i^t C_j^t$ (t = 1, 2) $e_{ij}^t \in \{1, 0, -1, I\}$, $e_{ij}^t$ is 0 if $C_i^t$ does not have any effect on $C_j^t$, $e_{ij}^t = 1$ if increase (or decrease) in $C_i^t$ causes increase (or decrease) on $C_j^t$. $e_{ij}^t = -1$ if increase (or decrease) on $C_i^t$ causes decrease (or increase) in $C_j^t$. $e_{ij}^t = I$ if the relation or effect of $C_i^t$ on $C_j^t$ is an indeterminate.

**DEFINITION 4.3.12:** *NCBMs with edge weight from {-1, 0, 1, I} are called simple NCBMs.*

**DEFINITION 4.3.13:** *Let $C_1^t, C_2^t, ..., C_{K_t}^t$ (t = 1, 2) be nodes of a NCBM. Let the neutrosophic bimatrix N(E) be defined as N(E) = $\left(e_{ij}^t\right)$ (t = 1, 2) where $e_{ij}^t$ is the weight of the directed edge $C_i^t, C_j^t$ where $e_{ij}^t$ ∈ {1, 0, -1, I}, (t = 1, 2)*

*N(E) is called the neutrosophic adjacency bimatrix of the NCBM.*

**DEFINITION 4.3.14:** *Let $C_1^t, C_2^t, C_3^t, \ldots, C_{n_t}^t$ be nodes of the NCBM. Let A = $\left(a_1^t, a_2^t, \ldots, a_{n_t}^t\right)$ where $a_i^t \in$ {0, 1, I}, A is called the instantaneous state neutrosophic bivector and it denotes the on-off-indeterminate state / position of the node at an instant*

$a_i^t = 0$ *if* $a_i^t$ *is off (no effect)*

$a_i^t = 1$ *if* $a_i^t$ *is on (has effect)*

$a_i^t = I$ *if* $a_i^t$ *is indeterminate (effect cannot be determined)*

*(t = 1, 2) i = 1, 2,…, n.*

**DEFINITION 4.3.15:** *Let $C_1^t, C_2^t, \ldots, C_{n_t}^t$ be the nodes of the NCBM.*

*Let $C_1^t C_2^t,\ C_2^t C_3^t, \ldots, C_i^t C_j^t$ = (t = 1, 2) be the edge of the NCBM. Then the edges form a directed cycle An NCBM is said to be cyclic if it possesses a directed cycle.*

An NCBM is said to be acyclic if it does not possess any directed cycle.

**DEFINITION 4.3.16:** *An NCBM with cycles is said to have a feed back. When there is a feed back in the NCBM i.e. when the causal relation flow through a cycle in a revolutionary manner the NCBM is called a dynamical system.*

**DEFINITION 4.3.17:** *Let $C_1^t C_2^t,\ C_2^t C_3^5, \ldots, C_{n-1}^t C_n^t$ be a cycle when $C_i^t$ (t = 1, 2) is switched on and if the causality flow through the edges of a cycle and it again causes $C_i^t$ we say that the dynamical system goes round and round. This is true for any node $C_i^t$ (t = 1, 2) for i = 1, 2,…, $n_t$ the* *equilibrium state for this dynamical system is called the bihidden pattern.*

**DEFINITION 4.3.18:** *If the equilibrium state of a dynamical system is a unique state bivector, then it is called a fixed bipoint. Consider the NCBM with $C_1$, $C_2$,…, $C_{n_t}$ as nodes. For instance let us start the dynamical system by switching on $C_i$. Let us assume that the NCBM settles down with $C_1^t$ and $C_{m_t}^t$ on i.e. the state vector remains as (1, 0,…,0 1 0) $\cup$ (1, 0,…, 0 1 0) this neutrosophic bivector is called the fixed bipoint.*

**DEFINITION 4.3.19:** *If the NCBM settles with a neutrosophic state bivector repeating in the form $A_1^t \to A_2^t \to, \ldots, \to A_i^t \to A_i$ (t = 1, 2) then this equilibrium is called a limit bicycle of the NCBM.*

We just give a brief description of how the hidden pattern is determined.

Let $C_1^t, C_2^t, \ldots, C_{n_t}^t$ (t = 1, 2) be the nodes of the NCBM with feed back. Let $N(E_B) = E_1 \cup E_2$ be the associated adjacency bimatrix. Let us find the bihidden pattern when $C_1^t$ is switched on (t = 1, 2) when an input is given as the bivector $A = A_1 \cup A_2$ = (1 0 0…0) $\cup$ (1 0 …0) the data should pass through the neutrosophic bimatrix $N(E_B)$ this is done by multiplying the bivector A by $N(E_B)$.

Now

$$\begin{aligned} AN(E_B) &= (A_1 \cup A_2)(E_1 \cup E_2) \\ &= A_1 E_1 \cup A_2 E_2 \\ &= \left(a_1^1, \ldots, a_{n_1}^1\right) \cup \left(a_1^2, \ldots, a_{n_2}^2\right), \end{aligned}$$

with the threshold operation that is by replacing $a_i^t$ by 1 (t = 1, 2) if $a_i^t \geq K$ and $a_i^t = 0$ if $a_i^t < K$ (K – a suitable positive integer) and $a_i^t$ by $I$ if $a_i^t$ is not an integer. We update the resulting concept, the concept $C_1$ is included in the updated bivector by making the first coordinate as 1 in the resulting bivector for we started the operation in this case with the first coordinate in the on state. Suppose AN ($E_B$) $\to$ B = $B_1 \cup B_2$ then consider

$$
\begin{aligned}
BN(B_B) &= (B_1 \cup B_2)(E_1 \cup E_2) \\
&= B_1 E_1 \cup B_2 E_2 \text{ and} \\
&= \left(b_1^1,\ldots,b_{n_1}^1\right) \cup \left(b_1^2,\ldots,b_{n_2}^2\right);
\end{aligned}
$$

and repeat the same procedure. This procedure is repeated till we get a limit bicycle or a fixed bipoint.

**DEFINITION 4.3.20:** *Finite number of NCBM can be combined together to produce the joint effect of all NCBMs. If* $N\left(E_B^1\right),\ldots,N\left(E_B^p\right)$ *be the neutrosophic adjacency bimatrices of a NCBM with nodes* $C_1^t, C_2^t,\ldots,C_{n_t}^t$ *then the combined NCBM is got by adding all the neutrosophic adjacency bimatrices* $N\left(E_B^1\right), \ldots, N\left(E_B^p\right)$. *We denote the combined NCBMs adjacency neutrosophic bimatrices by*

$$N(E_B) = N\left(E_B^1\right) + \ldots + N\left(E_B^p\right).$$

***Example 4.3.1:*** We illustrate this model in the child labour problem. The child labour which is very much prevalent in India but given the least importance is modeled using NCBMs.

Here two experts opinion are taken simultaneously having the same set of conceptual nodes $C_1$, $C_2$, …, $C_7$. described below

| | | |
|---|---|---|
| $C_1$ | - | Child labour |
| $C_2$ | - | Role of political leaders |
| $C_3$ | - | Role of a good teacher |
| $C_4$ | - | Poverty |
| $C_5$ | - | Industrialists attitude in practicing child labour |
| $C_6$ | - | Public practicing / encouraging child labour |
| $C_7$ | - | Good NGOs. |

The directed bigraph of the model is given by the following figure.

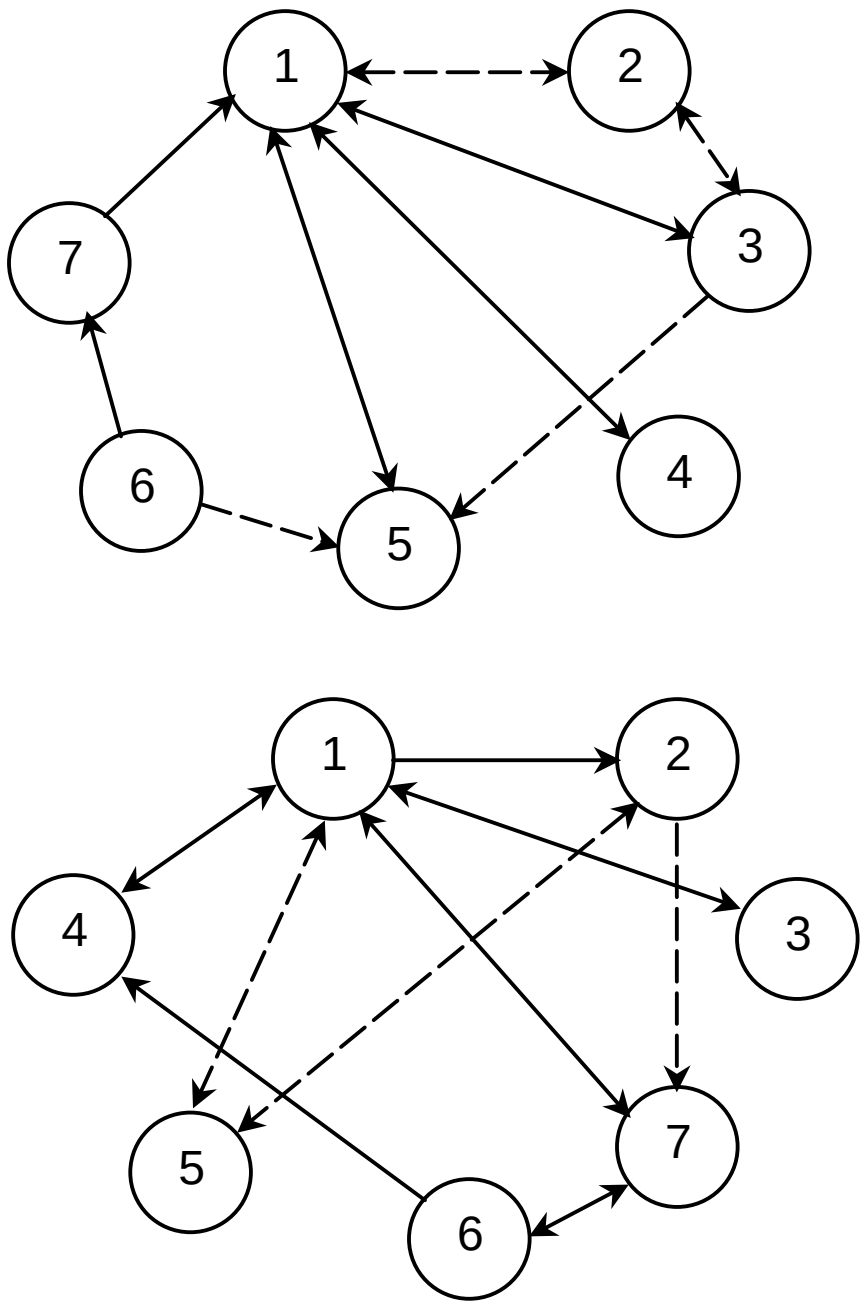

FIGURE: 4.3.1

The corresponding adjacency neutrosophic connection bimatrix N $(E_B)$

$$N(E_B) = N(E_1) \cup N(E_2)$$

$$\begin{bmatrix} 0 & I & -1 & 1 & 1 & 0 & 0 \\ I & 0 & I & 0 & 0 & 0 & 0 \\ 1 & I & 0 & 0 & I & 0 & 0 \\ 1 & 0 & 0 & 0 & 0 & 0 & 0 \\ 1 & 0 & 0 & 0 & 0 & 0 & 0 \\ 0 & 0 & 0 & 0 & I & 0 & -1 \\ -1 & 0 & 0 & 0 & 0 & 0 & 0 \end{bmatrix} \cup \begin{bmatrix} 0 & 1 & -1 & 1 & I & 0 & -1 \\ 0 & 0 & 0 & 0 & I & 0 & I \\ -1 & 0 & 0 & 0 & 0 & 0 & 0 \\ 1 & 0 & 0 & 0 & 0 & 0 & 0 \\ I & I & 0 & 0 & 0 & 0 & 0 \\ 0 & 0 & 0 & 1 & 0 & 0 & -1 \\ -1 & 0 & 0 & 0 & -1 & 0 \end{bmatrix}$$

Suppose the expert wishes to find the effect of the state bivector.

$A = A_1 \cup A_2$
$= (1\ 0\ 0\ 0\ 0\ 0\ 0) \cup (1\ 0\ 0\ 0\ 0\ 0\ 0)$

on the neutrosophic dynamical bisystem $N(E_B)$.

$AN(E_B) = A_1 N(E_1) \cup A_2 N(E_2)$
$= (0\ I - 1\ 1\ 1\ 0\ 0) \cup (0\ 1\ -1\ 1\ I\ 0\ -1).$

After thresholding and updating we get

$B = (1\ I\ 0\ 1\ 1\ 0\ 0) \cup (1\ 1\ 0\ 1\ I\ 0\ 0)$
$= B_1 \cup B_2.$

The same procedure as in case of FCBMs described in the earlier chapters is adopted to arrive at the bihidden pattern.

The important things to be noted in this example is that

1. This is a model to study simultaneously the effect of any state vector on the dynamical system.
2. This study helps to see how at each stage the state vectors behave which is entirely based on the opinion of the expert.
3. The final solution can be compared easily or even stage by stage comparison is possible.
4. At times one experts opinion may give a fixed point and that of the other expert may give a limit cycle where by even the number of times it is processed to arrive at a fixed point or a limit cycle is important.

This will show significance while the problem is an experimental one.

***Example 4.3.2:*** Next we proceed on to give yet another model of the system viz. on the same problem even the nodes chosen by the experts happen to be different except concurring on a node or two. This model describes important factors which influence good business i.e. factor promoting business. Here we use two experts opinion whose choice of nodes / attributes are not identical.

The first expert feels the factors promoting business are

$C_1$ – Good business
$C_2$ – Good investment
$C_3$ – Customer satisfaction
$C_4$ – Establishment
$C_5$ – Marketing strategies

The second expert gives the following nodes.

$D_1$ – Good business
$D_2$ – Geographical situation
$D_3$ – Rendering good service
$D_4$ – Previous experience of the owner
$D_5$ – Demand and supply.

The neutrosophic bigraph related with the model is given in figure.

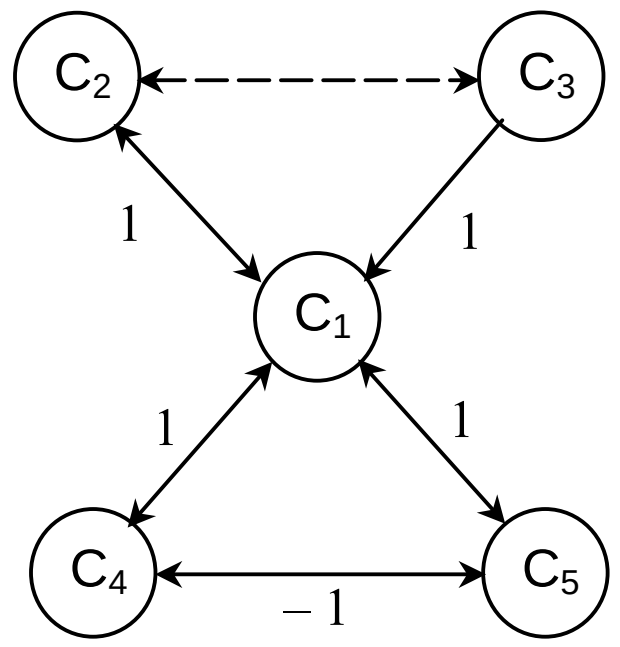

FIGURE: 4.3.2a

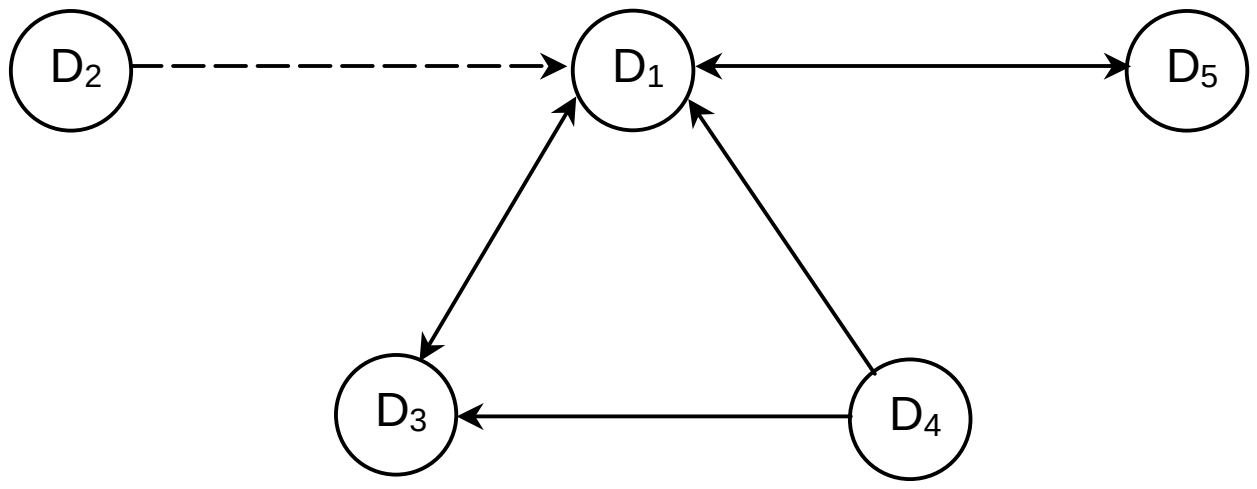

FIGURE: 4.3.2b

Now the related adjacency neutrosophic bimatrix N ($E_B$) is given as

$$\begin{bmatrix} 0 & 1 & 0 & 1 & 1 \\ 1 & 0 & I & 0 & 0 \\ 1 & 0 & 0 & 0 & 0 \\ 1 & 0 & 0 & 0 & -1 \\ 1 & 0 & 0 & -1 & 0 \end{bmatrix} \cup \begin{bmatrix} 0 & 0 & 1 & 0 & 1 \\ I & 0 & 0 & 0 & 0 \\ 1 & 0 & 0 & 0 & 0 \\ 1 & 0 & 1 & 0 & 0 \\ 1 & 0 & 0 & 0 & 0 \end{bmatrix}$$

i.e. $N(E_B) = N(E_1) \cup N(E_2)$

This is clearly a square 5 × 5 neutrosophic bimatrix. Suppose we wish to study the effect of the state bivector A, where

$$\begin{aligned} A &= A_1 \cup A_2 \\ &= (1\ 0\ 0\ 0\ 0) \cup (1\ 0\ 0\ 0\ 0) \\ AN(E_B) &= A_1N(E_1) \cup A_2N(E_2) \\ &= (0\ 1\ 0\ 1\ 1) \cup (0\ 0\ 1\ 0\ 1). \end{aligned}$$

After updating we get

$$\begin{aligned} B &= (1\ 1\ 0\ 1\ 1) \cup (1\ 0\ 1\ 0\ 1) \\ BN(E_B) &= (3\ 1\ I\ 0\ 0) \cup (2\ 0\ 1\ 0\ 1) \\ C &= (1\ 1\ I\ 0\ 0) \cup (1\ 0\ 1\ 0\ 1) \end{aligned}$$

be the resultant vector after updating and thresholding

$$CN(E_B) = (I + 1, 1, I\ 1\ 1) \cup (1\ 0\ 1\ 0\ 1).$$

Let D be the resultant vector after updating and thresholding the resultant is given by

$$F = (1\ 1\ I\ 1\ 1) \cup (1\ 0\ 1\ 0\ 1).$$

Thus the bihidden pattern is a fixed bipoint given by the binary neutrosophic bivector.

$$F = (1\ 1\ I\ 0\ 1) \cup (1\ 0\ 1\ 0\ 1).$$

It is clearly seen that the opinions as well as the nodes are different. For the first expert feels Good business has a direct effect on Good investment and establishment, no effect on Marketing strategies while customer satisfaction is an indeterminate concept to it.

According to the second expert Good business has direct effect on Rendering good service while it has no effect on other factors. Thus we have seen an example of a NCBM which studies the opinion of experts on the same problem but with different nodes / attributes. Now the same neutrosophic bigraph can also be given as a connected neutrosophic bigraph we obtain the following figure.

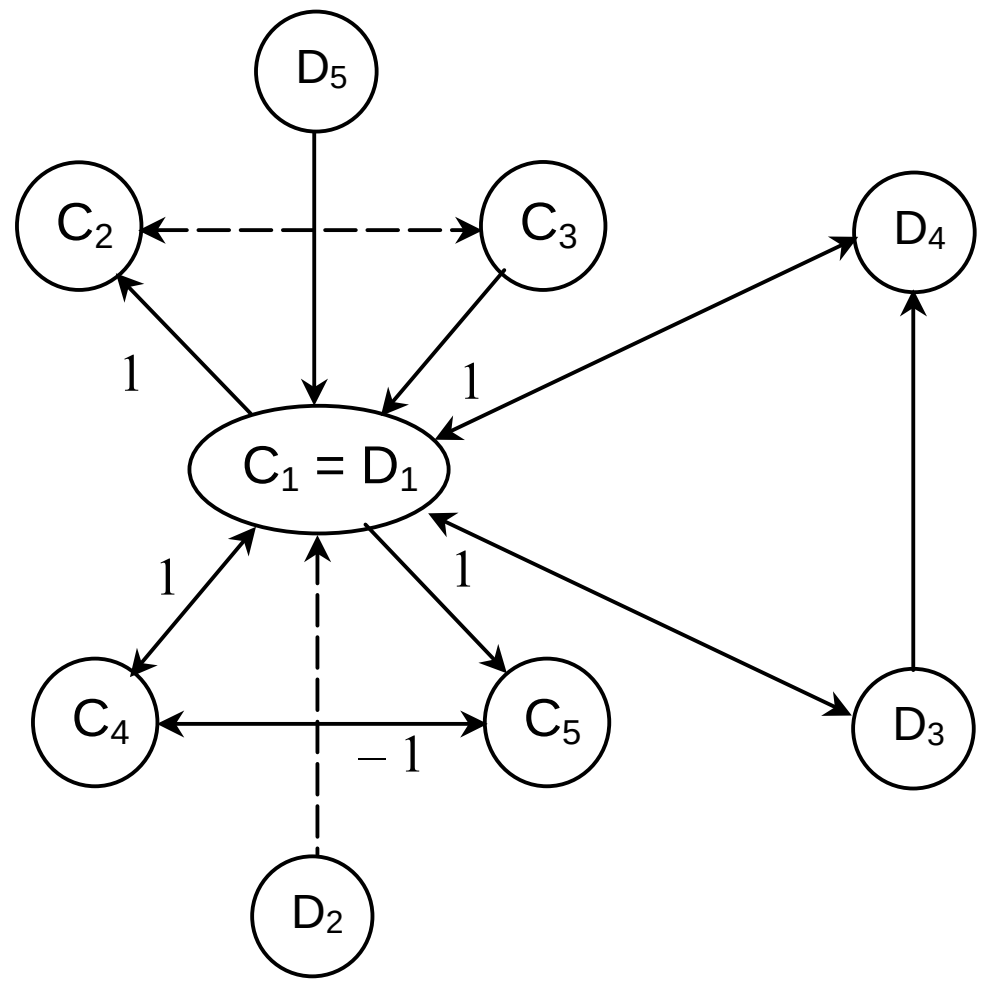

FIGURE: 4.3.3

This is the same bigraph but as "good business" i.e. D = $C_1$ we have drawn the neutrosophic bigraph with a vertex $C_1$ = $D_1$ which is common. This model can be represented as disconnected neutrosophic bigraph or as connected neutrosophic bigraph for the problem of drawing the neutrosophic graphs with given vertices and edges are only in the hands of the expert.

***Example 4.3.3:*** We proceed on to model a problem in which completely different opinions are given by two experts and the opinion also do not over lap. For in the study of health hazards faced by the agriculture coolies due to chemical pollution the two experts give totally disjoint attributes.

| | | |
|---|---|---|
| $P_1$ | – | Loss of appetite |
| $P_2$ | – | Headache |
| $P_3$ | – | Spraying pesticides |
| $P_4$ | – | Indigestion |
| $P_5$ | – | Giddiness / fainting |

The factors given by the second expert are

| | | |
|---|---|---|
| $Q_1$ | – | Swollen limbs |
| $Q_2$ | – | Ulcer / skin ailments in legs and hands |
| $Q_3$ | – | Manuring the field with chemical fertilizers |
| $Q_4$ | – | Vomiting |

The neutrosophic directed bigraph given by them is in figures 4.3.3a and 4.3.3b .

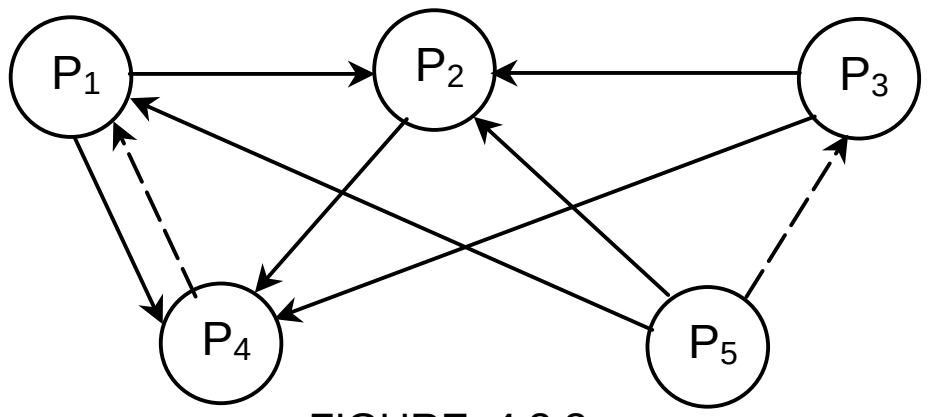

FIGURE: 4.3.3a

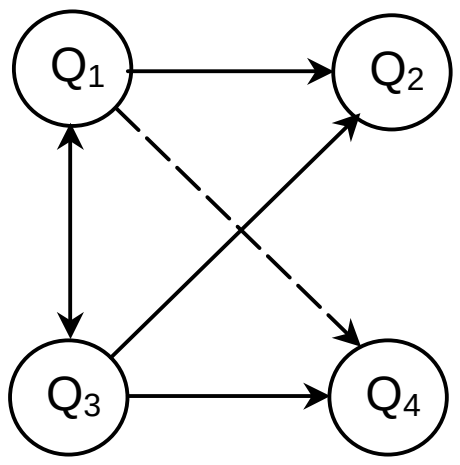

FIGURE: 4.3.3b

The related adjacency neutrosophic connection bimatrix $N(E_B)$ is given by

$$\begin{bmatrix} 0 & 1 & 0 & 1 & 0 \\ 0 & 0 & 0 & 1 & 0 \\ 0 & 1 & 0 & 1 & 0 \\ I & 0 & 0 & 0 & 0 \\ 1 & 1 & I & 0 & 0 \end{bmatrix} \cup \begin{bmatrix} 0 & 1 & 1 & I \\ 1 & 0 & 0 & 0 \\ 1 & 1 & 0 & 1 \\ 0 & 0 & 0 & 0 \end{bmatrix}$$

Using this model one can study the effect of any pair of attributes of the systems which are given as resultant neutrosophic bivector. Thus the related neutrosophic bigraph is a disconnected bigraph.

Now we give an example of a weak neutrosophic bigraph which can be also called as Neutrosophic cognitive weak bimaps. By this one expert may or may not give an indeterminate edge but the other expert will certainly give indeterminate edges so the directed bigraph related with this expert will only be a weak neutrosophic bigraph and the adjacency connection bimatrix would only be a weak neutrosophic bimatrix.

***Example 4.3.4:*** We illustrate this model by the following example. Suppose we study the same type of model to improve business in industries. We take the opinion of two experts one who gives all edges of related concepts as determinate ones and the other uses some neutrosophic theory.

The attributes given by the first expert are

| | | |
|---|---|---|
| $C_1$ | – | Good business |
| $C_2$ | – | Appropriate locality |
| $C_3$ | – | Selling quality products |
| $C_4$ | – | Updation of technologies |
| $C_5$ | – | Knowledge about policies of the government. |

The attributes enlisted by the second expert are

| | | |
|---|---|---|
| $D_1$ | – | Good business |
| $D_2$ | – | Good investment |
| $D_3$ | – | Customer satisfaction |
| $D_4$ | – | Establishment |
| $D_5$ | – | Marketing strategies. |

The related directed weak neutrosophic bigraph is given by the following figure.

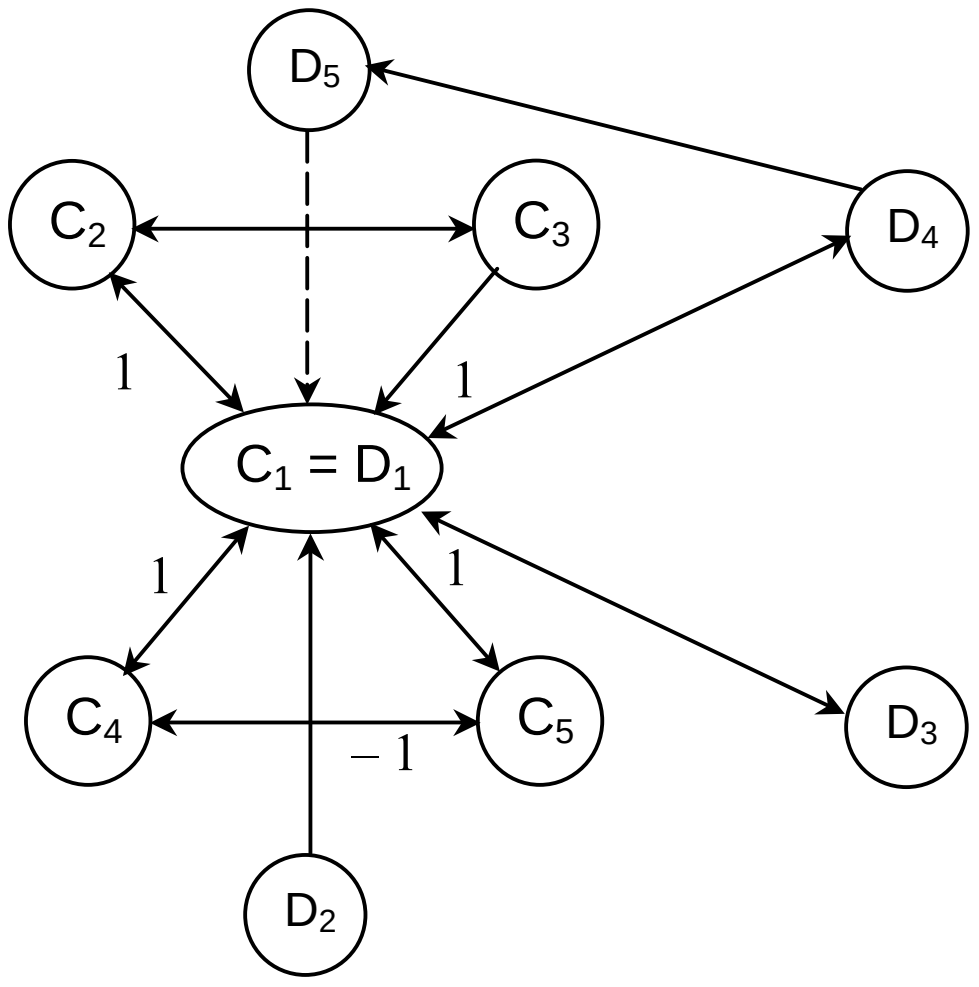

FIGURE: 4.3.4

The related adjacency connection bimatrix which is only a weak neutrosophic bimatrix $N(E_B)$ is given in the following:

$$\begin{bmatrix} 0 & 1 & 0 & 1 & 0 \\ 1 & 0 & 1 & 0 & 0 \\ 1 & 0 & 0 & 0 & 0 \\ 1 & 0 & 0 & 0 & -1 \\ 1 & 0 & 0 & -1 & 0 \end{bmatrix} \cup \begin{bmatrix} 0 & 0 & 1 & 0 & 0 \\ 1 & 0 & 0 & 0 & 0 \\ 1 & 0 & 0 & 0 & 0 \\ 1 & 0 & 0 & 0 & 1 \\ I & 0 & 0 & 0 & 0 \end{bmatrix}$$

$= N(E_1) \cup N(E_2)$
$= N(E_B)$.

Now we can study the effect of any state bivector say

$A = (0\ 0\ 0\ 0\ 1) \cup (0\ 0\ 0\ 0\ 1)$
$= A_1 \cup A_2$
$AN(E_B) = A_1 N(E_1) \cup A_2 N(E_2)$
$= (1\ 0\ 0\ {-1}\ 0) \cup (I\ 0\ 0\ 0\ 0)$.

After thresholding and updating we get the resultant state bivector as

$B = (1\ 0\ 0\ 0\ 1) \cup (I\ 0\ 0\ 0\ 1)$
$= B_1 \cup B_2$.

Now the effect of B on the neutrosophic dynamical system is given by

$BN(E_B) = B_1 N(E_1) \cup B_2 N(E_2)$
$= (1\ 1\ 0\ 0\ 0) \cup (I\ 0\ I\ 0\ 0)$.

After thresholding and updating we get the resultant bivector as

$C = (1\ 1\ 0\ 0\ 1) \cup (I\ 0\ I\ 0\ 1)$
$= C_1 \cup C_2$.

$CN(E_B) = C_1 N(E_1) \cup C_2 N(E_2)$
$= (2\ 1\ 1\ 0\ 0) \cup (2I\ 0\ I\ 0\ 0)$.

After thresholding and updating we get the state bivector F as

$$F = (1\ 1\ 1\ 0\ 1) \cup (I\ 0\ I\ 0\ 1)$$
$$= F_1 \cup F_2\ .$$

Now the effect of F on N ($E_B$) is given by

$$FN(E_B) = F_1\, N(E_1) \cup F_2\, N(E_2)$$
$$= (3\ 1\ 1\ 0\ 0) \cup (2\ I\ 0\ 1\ 0\ 0).$$

After thresholding and updating we get the resultant as

$$S = (1\ 1\ 1\ 0\ 1) \cup (I\ 0\ I\ 0\ 1)$$
$$= S_1 \cup S_2.$$

Thus we see the bihidden pattern of the dynamical system is a fixed bipoint given by the binary bipair {(1 1 1 0 1), (*I* 0 *I* 0 1)}.

Suppose we wish to study the combined effect of several experts on the same set of attributes then we can make use of the concept of combined neutrosophic cognitive bimaps. The following conditions must be a minimum criteria for us to adopt the combined neutrosophic cognitive bimaps (CNCBMs).

***Example 4.3.5:*** Suppose we have a NCBMs on a set of a pair of nodes / concepts given by $A_1, \ldots, A_n$ and $B_1, \ldots, B_m$. Thus when an experts opinion is sought on this problem relating the nodes we would get a mixed square neutrosophic bimatrix $N(E_B) = N(E_1) \cup N(E_2)$ given by

$$\begin{array}{c} \\ A_1 \\ A_2 \\ \vdots \\ A_n \end{array}\begin{array}{c} \begin{array}{ccc} A_1 & A_2 \ \ldots & A_n \end{array} \\ \begin{bmatrix} a^1_{11} & a^1_{12} & a^1_{1n} \\ a^1_{21} & a^1_{22} & a^1_{2n} \\ & & \\ a^1_{n1} & a^1_{n2} & a^1_{nn} \end{bmatrix} \end{array} \cup \begin{array}{c} \\ B_1 \\ B_2 \\ \vdots \\ B_m \end{array}\begin{array}{c} \begin{array}{ccc} B_1 & B_2 \ \cdots & B_m \end{array} \\ \begin{bmatrix} b^1_{11} & b^1_{12} & b^1_{1n} \\ b^1_{21} & b^1_{22} & b^1_{2n} \\ & & \\ b^1_{m2} & \cdots & b^1_{mm} \end{bmatrix} \end{array}$$

with $a^1_{ii} = 0 = b^1_{jj}$, $1 \le i \le n$, $1 \le j \le m$.

Suppose another set of experts give opinion for these attributes and the related neutrosophic square bimatrix is given by

$$N\left(E^1_B\right) = N\left(E^1_1\right) \cup \left(E^1_2\right)$$

$$\begin{array}{c} \\ A_1 \\ A_2 \\ \vdots \\ A_n \end{array}\begin{array}{c} \begin{array}{cccc} A_1 & A_2 & \ldots & A_n \end{array} \\ \begin{bmatrix} a^2_{11} & a^2_{12} & & a^2_{1n} \\ a^2_{21} & a^2_{22} & & a^2_{2n} \\ & & & \\ a^2_{n1} & a^2_{n2} & & a^2_{nn} \end{bmatrix} \end{array} \cup \begin{array}{c} \\ B_1 \\ B_2 \\ \vdots \\ B_m \end{array}\begin{array}{c} \begin{array}{cccc} B_1 & B_2 & \cdots & B_m \end{array} \\ \begin{bmatrix} b^2_{11} & b^2_{12} & & b^2_{1m} \\ b^2_{21} & b^2_{22} & & b^2_{2m} \\ & & & \\ b^2_{m1} & b^2_{m2} & \cdots & b^2_{mm} \end{bmatrix} \end{array}$$

with $a^2_{ii} = 0 = b^2_{jj}$, $1 \le i \le n$, $1 \le j \le m$;

Now we can add $N(E_B)$ and $N\left(E^1_B\right)$ as $N(E_1)$ and $N\left(E^1_1\right)$ have same set of column and row attributes even with order preserved.
Likewise $N(E_2)$ and $N\left(E^1_2\right)$ have the same set of column and row attributes as $B_1, \ldots, B_m$.

Thus the addition is meaningful and in the addition only the columns (rows) of elements pertaining to the same attributes are added. So

$$\begin{aligned} N(E_B) + N\left(E^1_B\right) &= [N(E_1) \cup N(E_2)] + [N\left(E^1_1\right) \cup N\left(E^1_2\right)] \\ &= [N(E_1) + N\left(E^1_1\right)] \cup [N(E_2) + N\left(E^1_2\right)] \end{aligned}$$

$$
= \begin{array}{c} \\ A_1 \\ A_2 \\ \vdots \\ A_n \end{array}
\begin{array}{c} A_1 \quad A_2 \quad \ldots \quad A_n \\ \begin{bmatrix} c_{11} & c_{12} & & c_{1n} \\ c_{21} & c_{22} & & c_{1n} \\ & & & \\ c_{n1} & c_{n2} & & c_{nn} \end{bmatrix} \end{array}
\cup
\begin{array}{c} \\ B_1 \\ B_2 \\ \vdots \\ B_m \end{array}
\begin{array}{c} B_1 \quad B_2 \quad \cdots \quad B_m \\ \begin{bmatrix} d_{11} & d_{12} & \cdots & d_{1m} \\ d_{21} & d_{22} & & d_{2m} \\ & & & \\ d_{m1} & d_{m2} & & d_{mn} \end{bmatrix} \end{array}
$$

$$c_{ii} = d_{jj} = 0 \text{ with } 1 \le i \le n,\ 1 \le j \le m$$
$$\text{and } c_{il} = a^1_{ik} + a^2_{ik},\ 1 \le i, k \le n$$
$$\text{and } d_{jt} = b^1_{jt} + b^2_{jt},\ 1 \le i, j, t \le m.$$

Now using the neutrosophic bimatrix as the adjacency connection bimatrix / dynamical system we can obtain the effect of any state bivector and interpret the related results. Thus FCBMs can be thought of as a generalization of FCMs and likewise NCBMs are the generalization of NCMs.

## 4.4 Neutrosophic Trimaps and their Applications

Now can we still generalize this model? The answer is yes. We can define neutrosophic cognitive Trimaps where they are got by either taking three experts opinion simultaneously or by taking three sets of disjoint or over lapping attributes and studying them separately. The resulting directed graph will be a neutrosophic trigraph may be connected or disconnected depending on the model under investigation.
The related trigraphs of the neutrosophic cognitive trimaps can be connected trigraphs or disconnected trigraphs or weakly connected trigraphs or very weakly connected trigraphs.

We illustrate this model by examples.

***Example 4.4.1:*** Suppose we wish to study the opinion of three experts pertaining to child labour on the same set of conceptual nodes $C_1, C_2, \ldots, C_7$ given in page 197 and 198.

The directed neutrosophic trigraph is given below.

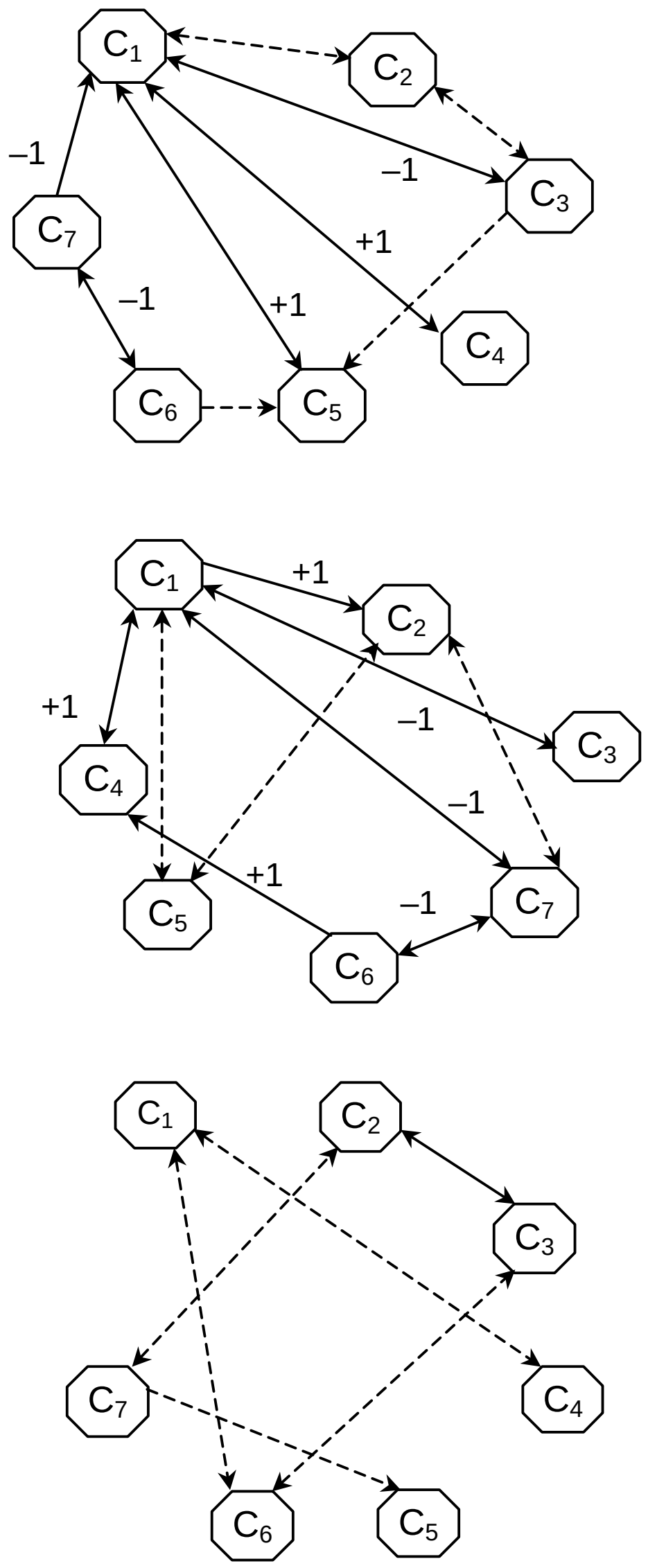

FIGURE: 4.4.1

The related adjacency connection neutrosophic trimatrix $N(E_T)$ is given below.

$$\begin{bmatrix} 0 & I & -1 & 1 & 1 & 0 & 0 \\ I & 0 & I & 0 & 0 & 0 & 0 \\ -1 & I & 0 & 0 & I & 0 & 0 \\ 1 & 0 & 0 & 0 & 0 & 0 & 0 \\ 1 & 0 & 0 & 0 & 0 & 0 & 0 \\ 0 & 0 & 0 & 0 & I & 0 & -1 \\ -1 & 0 & 0 & 0 & 0 & 0 & 0 \end{bmatrix} \cup$$

$$\begin{bmatrix} 0 & 1 & -1 & 1 & I & 0 & -1 \\ 0 & 0 & 0 & 0 & I & 0 & I \\ -1 & 0 & 0 & 0 & 0 & 0 & 0 \\ 1 & 0 & 0 & 0 & 1 & 0 & 0 \\ I & I & 0 & 0 & 0 & 0 & 0 \\ 0 & 0 & 0 & 1 & 0 & 0 & -1 \\ -1 & 0 & 0 & 0 & 0 & -1 & 0 \end{bmatrix} \cup \begin{bmatrix} 0 & 0 & 0 & I & 0 & I & 0 \\ 0 & 0 & 1 & 0 & 0 & 0 & I \\ 0 & 1 & 0 & 0 & 0 & 0 & 0 \\ I & 0 & 0 & 0 & 0 & 0 & 0 \\ 0 & 0 & 0 & 0 & 0 & 0 & I \\ I & 0 & I & 0 & 0 & 0 & 0 \\ 0 & I & 0 & 0 & I & 0 & 0 \end{bmatrix}$$

$$N(E_T) = N(E_1) \cup N(E_2) \cup N(E_3)$$

Just we wish to study the effect of the state trivector

$$\begin{aligned} A &= (1\ 0\ 0\ 0\ 0\ 0\ 0) \cup (1\ 0\ 0\ 0\ 0\ 0\ 0) \cup (1\ 0\ 0\ 0\ 0\ 0\ 0) \\ &= A_1 \cup A_2 \cup A_3 \end{aligned}$$

on the dynamical system

$$\begin{aligned} AN(E_T) = \; & [(0\ I\ {-1}\ 1\ 1\ 0\ 0) \cup (0\ 1\ {-1}\ 1\ I\ 0\ {-1}) \\ & \cup (0\ 0\ 0\ I\ 0\ I\ 0)]. \end{aligned}$$

After thresholding and updating we get the resultant trivector as

$$B = (1\ I\ 0\ 1\ 1\ 0\ 0) \cup (1\ 1\ 0\ 1\ I\ 0\ 0) \cup (1\ 0\ 0\ I\ 0\ I\ 0)$$
$$= B_1 \cup B_2 \cup B_3.$$

Now we study the effect of B on $N(E_T)$

$$BN(E_T) = (B_1 \cup B_2 \cup B_3)(N(E_T))$$
$$= B_1 N(E_1) \cup B_2 N(E_2) \cup B_3 N(E_3)$$
$$= (I+2, I, -1+I, 1, 1, 0\ 0) \cup (1+I, 1+I, -1, -1, 2I+1, 0, -1+I) \cup (I\ I\ 0\ I\ I\ 1\ 0).$$

After thresholding and updating the received trivector we get the resultant trivector as

$$C = (1\ I\ I\ 1\ 1\ 0\ 0) \cup (1\ I\ 0\ 0\ I\ 0\ I) \cup (1\ I\ 0\ I\ I\ 1\ 0)$$
$$= C_1 \cup C_2 \cup C_3.$$

Now we proceed on to work with $CN(E_T) = C_1 N(E_1) \cup C_2 N(E_2) \cup C_3 N(E_3)$ untill we get a fixed tripoint or a limit tricycle. This FCTM is a disjoint FCTM. Using this we can study in parallel the opinion of the expert at each stage. Till date there is no model which can at the same time enable one to study the problem stage by stage i.e., time dependent at each desired time. For the result alone is not the chief criteria in several of the problems. This method / model will enable the analyzer to understand derive conclusion and make further study of the social problem at every stage.

Using the trimatrix or the FCTM model we can also study models in which only a few of the attributes over lap for this three experts opinion are taken and like the previous model all of the three experts do not concur on the same set of attributes. Only a few of the attributes is common. We can study such models.

***Example 4.4.2:*** We will not show the working of the model but just show this by an illustration. In this case the neutrosophic trigraph would only be a connected trigraph and the related neutrosophic trimatrix may be a mixed square neutrosophic trimatrix. We now consider the problem of globalization and its impact on farmers.

Suppose the first expert considers the five attributes;

| | | |
|---|---|---|
| $C_1$ | – | Usage of fertilizers banned by foreign countries |
| $C_2$ | – | Patent right owned by foreign countries |
| $C_3$ | – | Frustration among farmers |
| $C_4$ | – | Loan facilities offered by financial institutions |
| $C_5$ | – | Burden due to loan and impact of patent right. |

The nodes suggested by the second expert are;

| | | |
|---|---|---|
| $D_1$ | – | Patent right owned by foreign countries |
| $D_2$ | – | Suicide of farmers |
| $D_3$ | – | Cultivating quality crops |
| $D_4$ | – | Loan given government and other financial institutions |
| $D_5$ | – | Use of modern technology in cultivation |
| $D_6$ | – | Failure of crops |
| $D_7$ | – | Repaying of loan and related tensions. |

The opinion given by the third expert.

| | | |
|---|---|---|
| $F_1$ | – | Repaying of loan and related tensions |
| $F_2$ | – | Failure of crops |
| $F_3$ | – | Suicide of farmers |
| $F_4$ | – | Loan facilities offered by government |

Now we give the connected neutrosophic trigraph given by the experts in relation to the problem. For the reader to follow however the separate neutrosophic graphs are also given. Now clearly this is a connected trigraph as there is a clear over lapping of attributes.

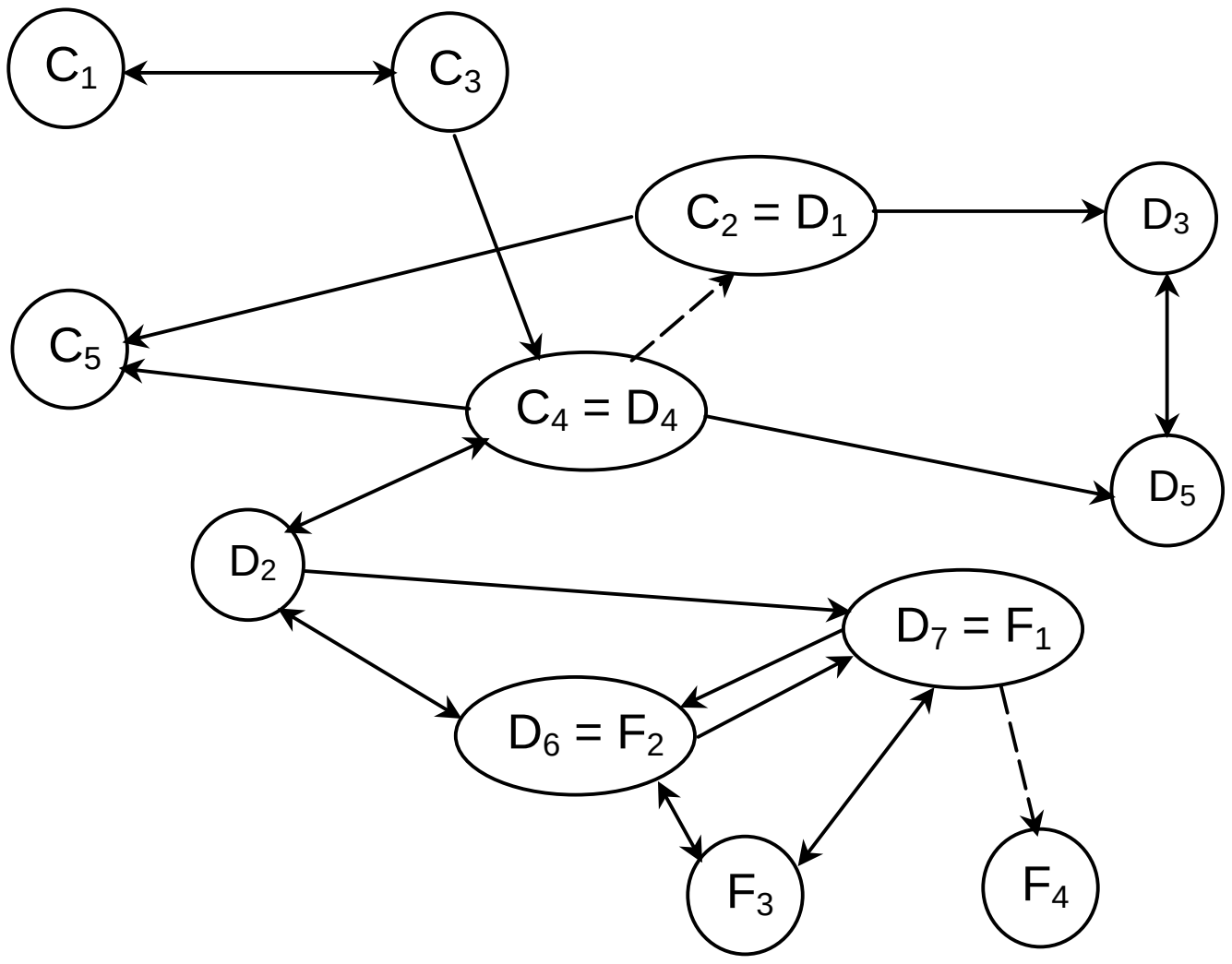

FIGURE: 4.4.2

Now this is also a neutrosophic edge connected trigraph. Now we give the three neutrosophic graphs.

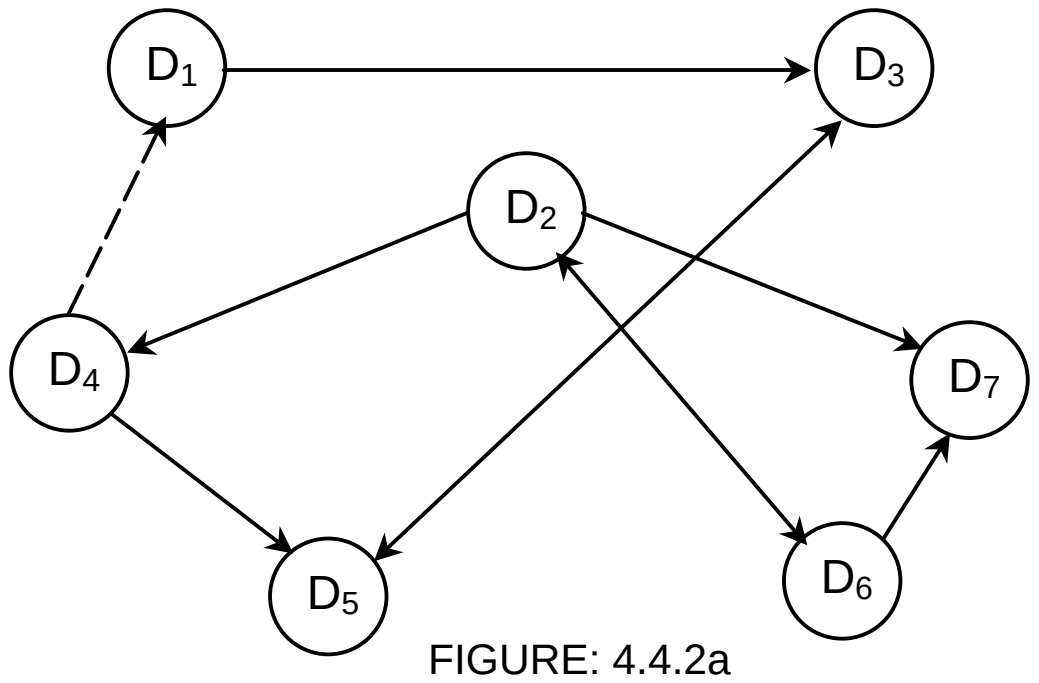

FIGURE: 4.4.2a

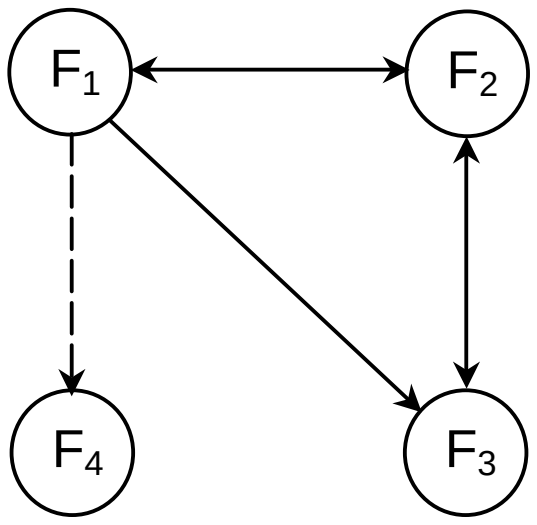

FIGURE: 4.4.2b

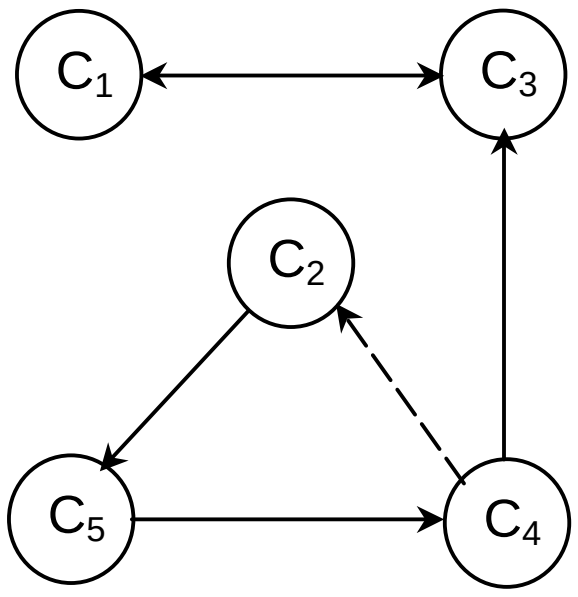

FIGURE: 4.4.2b

Now we proceed on to give the neutrosophic trimatrix which is a mixed square matrix $N(E_T)$

$$\begin{array}{c|ccccc|} & C_1 & C_2 & C_3 & C_4 & C_5 \\ \hline C_1 & 0 & 0 & 1 & 0 & 0 \\ C_2 & 0 & 0 & 0 & 0 & 1 \\ C_3 & 1 & 0 & 0 & 0 & 0 \\ C_4 & 0 & I & 0 & 0 & 0 \\ C_5 & 0 & 0 & 0 & 1 & 0 \end{array} \cup \begin{array}{c|ccccccc|} & D_1 & D_2 & D_3 & D_4 & D_5 & D_6 & D_7 \\ \hline D_1 & 0 & 0 & 1 & 0 & 0 & 0 & 0 \\ D_2 & 0 & 0 & 0 & 1 & 0 & 1 & 1 \\ D_3 & 0 & 0 & 0 & 0 & 1 & 0 & 0 \\ D_4 & I & 0 & 0 & 0 & 1 & 0 & 0 \\ D_5 & 0 & 0 & 1 & 0 & 0 & 0 & 0 \\ D_6 & 0 & 1 & 0 & 0 & 0 & 0 & 1 \\ D_7 & 0 & 0 & 0 & 0 & 0 & 0 & 0 \end{array} \cup$$

$$\begin{array}{c|cccc|} & F_1 & F_2 & F_3 & F_4 \\ \hline F_1 & 0 & 1 & 1 & I \\ F_2 & -1 & 0 & 1 & 0 \\ F_3 & 0 & 1 & 0 & 0 \\ F_4 & 0 & 0 & 0 & 0 \end{array}$$

This is clearly a mixed square neutrosophic trimatrix. Now we just show how this model works. Suppose an expert is interested in studying the initial state trivector

$$\begin{aligned} A &= (1\ 0\ 0\ 0\ 0) \cup (0\ 1\ 0\ 0\ 0\ 0\ 0) \cup (0\ 0\ 1\ 0) \\ &= A_1 \cup A_2 \cup A_2. \end{aligned}$$

To study the effect of A on the dynamical system

$N(E_T) = N(E_1) \cup N(E_2) \cup N(E_3)$

$AN(E_T) = (A_1 \cup A_2 \cup A_3) (N(E_1) \cup N(E_2) \cup N(E_3))$
$= A_1 N(E_1) \cup A_2 N(E_2) \cup A_3 N(E_3)$
$= (0\ 0\ 1\ 0\ 0) \cup (0\ 0\ 0\ 1\ 0\ 1\ 1) \cup (0\ 1\ 0\ 0)$.

After updating the resultant vector we get the resultant vector as

$B = B_1 \cup B_2 \cup B_3$
$= (1\ 0\ 1\ 0\ 0) \cup (0\ 1\ 0\ 1\ 0\ 1\ 1) \cup (0\ 1\ 1\ 0)$

The effect of B on N $(E_T)$

$BN(E_T) = (B_1 \cup B_2 \cup B_3) N(E_T)$
$= B_1 N(E_1) \cup B_2 N(E_2) \cup B_3 N(E_3)$
$= (1\ 0\ 1\ 0\ 0) \cup (I\ 1\ 0\ 1\ 1\ 1\ 2)$
$\cup (-1\ 1\ 1\ 0)$.

After thresholding and updating we get the resultant as C

$C = (1\ 0\ 1\ 0\ 0) \cup (I\ 1\ 0\ 1\ 1\ 1\ 1) \cup (0\ 1\ 1\ 0)$
$= C_1 \cup C_2 \cup C_3$.

We see $C_1$ and $C_3$ are the same so we can with out loss of generality work only with $C_2$.
However the effect of (on N $(E_1)$ gives

$CN(E_I) = (1\ 01\ 0\ 0) \cup (I, 1, 1 + I, 1\ 1\ 1\ 2)$
$\cup (0\ 1\ 1\ 0)$.

After updating and thresholding the resultant vector we get

$D = (1\ 0\ 1\ 0\ 0) \cup (I\ 1\ I\ 1\ 1\ 1\ 1) \cup (0\ 1\ 1\ 0)$
$= D_1 \cup D_2 \cup D_3$.

The effect of D on N $(E_T)$ is given by

$DN(E_T) \quad = \quad (1\ 0\ 1\ 0\ 0) \cup (I\ 1\ 1{+}1\ 1\ 1\ 2) \cup (0\ 1\ 1\ 0).$

After thresholding and updating we get

$F \quad = \quad (1\ 0\ 1\ 0\ 0) \cup (I\ 1\ I\ 1\ 1\ 1\ 1) \cup (0\ 1\ 1\ 0).$

As in case of FCM or NCMs or NCBMs we can read the trihidden pattern which is a fixed tripoint as follows.

1. According to the first expert the usage of fertilizers banned by foreign countries are used alone is in the on state then the node the loan facilities offered by financial institutions alone come to on state their by making one know; for one of the major loan facilities offered by the government is supply of these fertilizers to the farmers and then recover it from them in the form of money.
2. When suicide of farmers alone is in the on state and all other nodes are in the off state we see in the resultant vector patent right owned by foreign countries and cultivating quality crops are in the indeterminate state and all other nodes come to on state. Thus except the two nodes which are indeterminate all other nodes coming to on state implies all the nodes selected by him are so intricately dependent in influencing the suicide of farmers.
3. When suicide of farmers alone is in the on state by the third expert then only failure of crops alone come to an state the other states remain unaffected. Just like in the case of FCBMs for which the associated bimatrix can be a weak neutrosophic bimatrix, even in case of FCTMs we can have the associated trimatrix can be either a weak neutrosophic trimatrix or a very weak neutrosophic trimatrix. All the working as in case of neutrosophic trimatrix can be carried out in a similar way. We see in case when the FCTM model is such that the directed trigraph is a weak neutrosophic trigraph and then the associated connection adjacency trimatrix is a weak neutrosophic trimatrix i.e. if $M = M_1 \cup M_2 \cup M_3$ is the weak neutrosophic trimatrix. We have atleast two of the matrices to be neutrosophic matrices. Likewise in case of very weak neutrosophic trigraph and the related FCTMs.

## 4.5 Neutrosophic Relational Maps

Now we proceed on to define Neutrosophic Relational bimaps and trimaps. Before we do that just we recall the definition of Neutrosophic Relational Maps (NRMs).

When the nodes or concepts under study happens to be such that they can be divided into two disjoint classes then a study or analysis can be made using Fuzzy Relational Maps (FRMs) which is introduced and described. Here we recall the definition of a new concept called Neutrosophic Relational Maps (NRMs), analyse and study them. We also give examples of them.

**DEFINITION 4.5.1:** *Let D be the domain space and R be the range space with $D_1,\ldots, D_n$ the conceptual nodes of the domain space D and $R_1,\ldots, R_m$ be the conceptual nodes of the range space R such that they form a disjoint class i.e. $D \cap R = \phi$. Suppose there is a FRM relating D and R and if at least a edge relating a $D_i R_j$ is an indeterminate then we call the FRM as the Neutrosophic relational maps. i.e. NRMs.*

***Note:*** In everyday occurrences we see that if we are studying a model built using an unsupervised data we need not always have some edge relating the nodes of a domain space and a range space or there does not exist any relation between two nodes, it can very well happen that for any two nodes one may not be always in a position to say that the existence or nonexistence of a relation, but we may say that

the relation between two nodes is an indeterminate or cannot be decided.

Thus to the best of our knowledge indeterminacy models can be built using neutrosophy. One model already discussed is the Neutrosophic Cognitive Model. The other being the Neutrosophic Relational Maps model, which is a further generalization of Fuzzy Relational Maps model.

It is not essential when a study/ prediction/ investigation is made we are always in a position to find a complete answer. In reality this is not always possible (sometimes or many a times) it is almost all models built using unsupervised data, we may have the factor of indeterminacy to play a role. Such study is possible only by using the Neutrosophic logic or the concept of indeterminacy.

***Example 4.5.1:*** Female infanticide (the practice of killing female children at birth or shortly thereafter) is prevalent in India from the early vedic times, as women were (and still are) considered as a property. As long as a woman is treated as a property/ object the practice of female infanticide will continue in India.

In India, social factors play a major role in female infanticide. Even when the government recognized the girl child as a critical issue for the country's development, India continues to have an adverse ratio of women to men. Other reasons being torture by the in-laws may also result in cruel death of a girl child. This is mainly due to the fact that men are considered superior to women. Also they take into account the fact that men are breadwinners for the family. Even if women work like men, parents think that her efforts is going to end once she is married and enters a new family.

Studies have consistently shown that girl babies in India are denied the same and equal food and medical care that the boy babies receive, leave alone the education. Girl babies die more often than boy babies even though medical research has long ago established that girls are generally biologically stronger as newborns than boys. The birth of a male child is a time for celebration, but the birth of female child is often viewed as a crisis. Thus the female infanticide cannot be attributed to single reason it is highly dependent on the feeling of individuals ranging from social stigma, monetary waste, social status etc.

Suppose we take the conceptual nodes for the unsupervised data relating to the study of female infanticide.

We take the social status of the people as the domain space D

| | | |
|---|---|---|
| $D_1$ | – | Very rich |
| $D_2$ | – | Rich |
| $D_3$ | – | Upper middle class |
| $D_4$ | – | Middle class |
| $D_5$ | – | Lower middle class |
| $D_6$ | – | Poor |
| $D_7$ | – | Very poor. |

The nodes of the range space R are taken as

| | | |
|---|---|---|
| $R_1$ | – | Number of female children - a problem |
| $R_2$ | – | Social stigma of having female children |
| $R_3$ | – | Torture by in-laws for having only female children |
| $R_4$ | – | Economic loss / burden due to female children |
| $R_5$ | – | Insecurity due to having only female children (They will marry and enter different homes thereby leaving their parents, so no one would be able to take care of them in later days.) |

Keeping these as nodes of the range space and the domain space experts opinion were drawn which is given by the following figure 4.5.1:

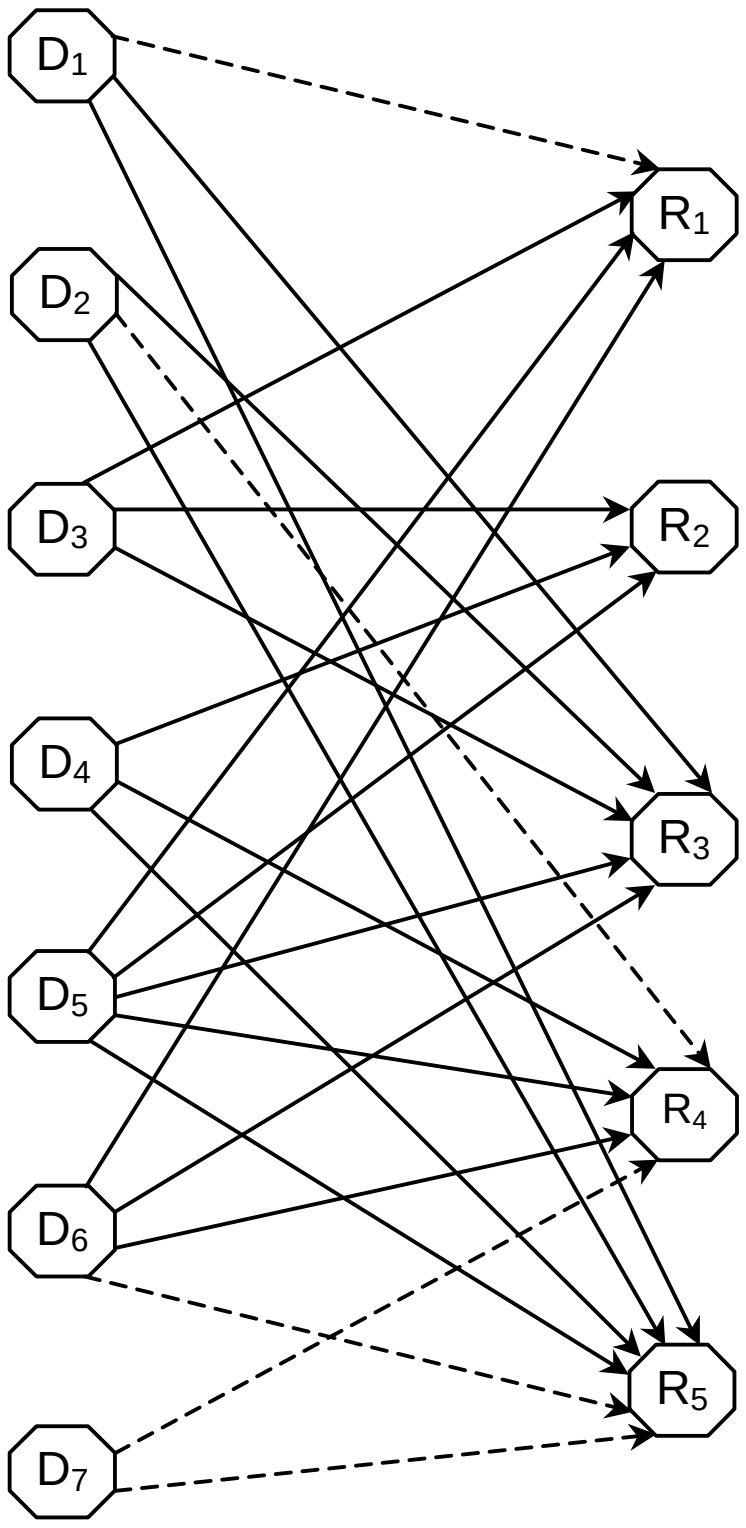

FIGURE: 4.5.1

Figure 4.5.1 is the neutrosophic directed graph of the NRM.

The corresponding neutrosophic relational matrix $N(R)^T$ is given below:

$$N(R)^T = \begin{bmatrix} I & 0 & 1 & 0 & 1 & 1 & 0 \\ 0 & 0 & 1 & 1 & 1 & 0 & 0 \\ 1 & 1 & 1 & 1 & 1 & 1 & 0 \\ 0 & I & 0 & 0 & 1 & 1 & I \\ 1 & 1 & 0 & 1 & 1 & I & I \end{bmatrix}$$

and

$$N(R) = \begin{bmatrix} I & 0 & 1 & 0 & 1 \\ 0 & 0 & 1 & I & 1 \\ 1 & 1 & 1 & 0 & 0 \\ 0 & 1 & 1 & 0 & 1 \\ 1 & 1 & 1 & 1 & 1 \\ 1 & 0 & 1 & 1 & I \\ 0 & 0 & 0 & I & I \end{bmatrix}.$$

Suppose $A_1 = (0\ 1\ 0\ 0\ 0)$ is the instantaneous state vector under consideration i.e., social stigma of having female children.

The effect of $A_1$ on the system N(R) is

$$A_1N(R)^T = (0\ 0\ 1\ 1\ 1\ 0\ 0)$$
$$\to (0\ 0\ 1\ 1\ 1\ 0\ 0) = B_1$$

$$B_1[N(R)] = (2, 3, 3, 1, 2)$$
$$\to (1\ 1\ 1\ 1\ 1) = A_2$$

$$A_2[N(R)]^T = (2 + I, 2 + I, 3, 3, 5, 3 + I, 2I)$$
$$\to (1\ 1\ 1\ 1\ 1\ 1\ I)$$
$$= B_2$$

$$B_2(N(R)) = (I + 3, 3, 5, 2I + 2, 2I + 4)$$
$$\to (1\ 1\ 1\ 1\ 1)$$
$$= A_3 = A_2.$$

Thus this state vector $A_1 = (0, 1, 0, 0, 0)$ gives a fixed point (1 1 1 1 1 1) indicating if one thinks that having female children is a social stigma immaterial of their status they also feel that having number of female children is a problem, it is a economic loss / burden, they also under go torture or bad treatment by in-laws and ultimately it is a insecurity for having only female children, the latter two cases hold were applicable.

On the other hand we derive the following conclusions on the domain space when the range space state vector $A_1$ = (0 1 0 0 0) is sent

$$A_1[N(R)]^T \rightarrow B_1$$
$$B_1[N(R)] \rightarrow A_2$$
$$A_2 [N(R)]^T \rightarrow B_2$$
$$B_2[N(R)] \rightarrow A_3 = A_2 \text{ so}$$
$$A_2 [N(R)]^T \rightarrow B_2 = (1, 1, 1, 1, 1, 1, I)$$

leading to a fixed point. When the state vector $A_1$ = (0, 1, 0, 0, 0) is sent to study i.e. the social stigma node is on uniformly all people from all economic classes are awakened expect the very poor for the resultant vector happens to be a Neutrosophic vector hence one is not in a position to say what is the feeling of the very poor people and the "many female children are a social stigma" as that coordinate remains as an indeterminate one. This is typical of real-life scenarios, for the working classes hardly distinguish much when it comes to the gender of the child.

Several or any other instantaneous vector can be used and its effect on the Neutrosophical Dynamical System can be studied and analysed. This is left as an exercise for the reader. Having seen an example and application or construction of the NRM model we will proceed on to describe the concepts of it in a more mathematical way.

**DESCRIPTION OF A NRM:**

Neutrosophic Cognitive Maps (NCMs) promote the causal relationships between concurrently active units or decides the absence of any relation between two units or the indeterminance of any relation between any two units. But in Neutrosophic Relational Maps (NRMs) we divide the very causal nodes into two disjoint units. Thus for the modeling of a NRM we need a domain space and a range space which are disjoint in the sense of concepts. We further assume no intermediate relations exist within the domain and the range spaces. The number of elements or nodes in the range space need not be equal to the number of elements or nodes in the domain space.

Throughout this section we assume the elements of a domain space are taken from the neutrosophic vector space of dimension n and that of the range space are vectors from neutrosophic vector space of dimension m. (m in general need not be equal to n). We denote by R the set of nodes $R_1, \ldots, R_m$ of the range space, where $R = \{(x_1, \ldots, x_m) \mid x_j = 0$ or 1 for $j = 1, 2, \ldots, m\}$.

If $x_i = 1$ it means that node $R_i$ is in the on state and if $x_i = 0$ it means that the node $R_i$ is in the off state and if $x_i = I$ in the resultant vector it means the effect of the node $x_i$ is indeterminate or whether it will be off or on cannot be predicted by the neutrosophic dynamical system.

It is very important to note that when we send the state vectors they are always taken as the real state vector for we know the node or the concept is in the on state or in the off state but when the state vector passes through the Neutrosophic dynamical system some other node may become indeterminate i.e. due to the presence of a node we may not be able to predict the presence or the absence of the other node i.e., it is indeterminate, denoted by the symbol *I*, thus the resultant vector can be a neutrosophic vector.

**DEFINITION 4.5.2:** *A Neutrosophic Relational Map (NRM) is a Neutrosophic directed graph or a map from D to R with concepts like policies or events etc. as nodes and causalities as edges. (Here by causalities we mean or include the indeterminate causalities also). It represents Neutrosophic Relations and Causal Relations between spaces D and R .*

*Let $D_i$ and $R_j$ denote the nodes of an NRM. The directed edge from $D_i$ to $R_j$ denotes the causality of $D_i$ on $R_j$ called relations. Every edge in the NRM is weighted with a number in the set {0, +1, –1, I}. Let $e_{ij}$ be the weight of the edge $D_i R_j$, $e_{ij} \in \{0, 1, -1, I\}$. The weight of the edge $D_i R_j$ is positive* *if increase in $D_i$ implies increase in $R_j$ or decrease in $D_i$ implies decrease in $R_j$ i.e. causality of $D_i$ on $R_j$ is 1. If $e_{ij} = –1$ then increase (or decrease) in $D_i$ implies decrease (or increase) in $R_j$. If $e_{ij} = 0$ then $D_i$ does not have any effect on $R_j$. If $e_{ij} = I$ it implies we are not in a position to determine the effect of $D_i$ on $R_j$ i.e. the effect of $D_i$ on $R_j$ is an indeterminate so we denote it by I.*

**DEFINITION 4.5.3:** *When the nodes of the NRM take edge values from {0, 1, –1, I} we say the NRMs are simple NRMs.*

**DEFINITION 4.5.4:** *Let $D_1, …, D_n$ be the nodes of the domain space D of an NRM and let $R_1, R_2,…, R_m$ be the nodes of the range space R of the same NRM. Let the matrix N(E) be defined as $N(E) = (e_{ij})$ where $e_{ij}$ is the weight of the directed edge $D_i R_j$ (or $R_j D_i$) and $e_{ij} \in \{0, 1, –1, I\}$. N(E) is called the Neutrosophic Relational Matrix of the NRM.*

The following remark is important and interesting to find its mention in this book.

**Remark**: Unlike NCMs, NRMs can also be rectangular matrices with rows corresponding to the domain space and columns corresponding to the range space. This is one of the marked difference between NRMs and NCMs. Further the number of entries for a particular model which can be treated as disjoint sets when dealt as a NRM has very much less entries than when the same model is treated as a NCM.

Thus in many cases when the unsupervised data under study or consideration can be split as disjoint sets of nodes or concepts; certainly NRMs are a better tool than the NCMs.

**DEFINITION 4.5.5:** *Let $D_1, …, D_n$ and $R_1,…, R_m$ denote the nodes of a NRM. Let $A = (a_1,…, a_n)$, $a_i \in \{0, 1, I\}$ is called the Neutrosophic instantaneous state vector of the domain space and it denotes the on-off position of the nodes at any instant. Similarly let $B = (b_1,…, b_n)$ $b_i \in \{0, 1, I\}$, B is called instantaneous state vector of the range space and it denotes* *the on-off-indeterminate position of the nodes at any instant, $a_i = 0$ if $a_i$ is off and $a_i = 1$ if $a_i$ is on and $a_i = I$ if that time $a_i$ can not be determined for $i = 1, 2, \ldots, n$. Similarly, $b_i = 0$ if $b_i$ is off and $b_i = 1$ if $b_i$ is on, $b_i = I$ if $b_i$ cannot be determined for $i = 1, 2, \ldots, m$.*

**DEFINITION 4.5.6:** *Let $D_1, \ldots, D_n$ and $R_1, R_2, \ldots, R_m$ be the nodes of a NRM. Let $D_i R_j$ (or $R_j D_i$ ) be the edges of an NRM, $j = 1, 2, \ldots, m$ and $i = 1, 2, \ldots, n$. The edges form a directed cycle. An NRM is said to be a cycle if it possess a directed cycle. An NRM is said to be acyclic if it does not possess any directed cycle.*

**DEFINITION 4.5.7:** *A NRM with cycles is said to be a NRM with feedback.*

**DEFINITION 4.5.8:** *When there is a feedback in the NRM i.e. when the causal relations flow through a cycle in a revolutionary manner the NRM is called a Neutrosophic dynamical system.*

**DEFINITION 4.5.9:** *Let $D_i R_j$ (or $R_j D_i$); $1 \leq j \leq m$, $1 \leq i \leq n$, when $R_j$ (or $D_i$ ) is switched on and if causality flows through edges of a cycle and if it again causes $R_j$ (or $D_i$ ) we say that the Neutrosophical dynamical system goes round and round. This is true for any node $R_j$ ( or $D_i$ ) for $1 \leq j \leq m$ (or $1 \leq i \leq n$). The equilibrium state of this Neutrosophical dynamical system is called the Neutrosophic hidden pattern.*

**DEFINITION 4.5.10:** *If the equilibrium state of a Neutrosophical dynamical system is a unique Neutrosophic state vector, then it is called the fixed point. Consider an NRM with $R_1, R_2, \ldots, R_m$ and $D_1, D_2, \ldots, D_n$ as nodes. For example let us start the dynamical system by switching on $R_1$ (or $D_1$). Let us assume that the NRM settles down with $R_1$ and $R_m$ (or $D_1$ and $D_n$) on, or indeterminate on, i.e. the Neutrosophic state vector remains as (1, 0, 0,…, 1) or (1, 0, 0,…, I) (or (1, 0, 0,…, 1) or (1, 0, 0,…, I) in D), this state vector is called the fixed point.*

**DEFINITION 4.5.11:** *If the NRM settles down with a state vector repeating in the form $A_1 \to A_2 \to A_3 \to \ldots \to A_i \to A_1$ (or $B_1 \to B_2 \to \ldots \to B_i \to B_1$) then this equilibrium is called a limit cycle.*

**METHODS OF DETERMINING THE HIDDEN PATTERN IN A NRM**

Let $R_1, R_2, \ldots, R_m$ and $D_1, D_2, \ldots, D_n$ be the nodes of a NRM with feedback. Let N(E) be the Neutrosophic Relational Matrix. Let us find the hidden pattern when $D_1$ is switched on i.e. when an input is given as a vector; $A_1 = (1, 0, \ldots, 0)$ in D; the data should pass through the relational matrix N(E). This is done by multiplying $A_1$ with the Neutrosophic relational matrix N(E). Let $A_1N(E) = (r_1, r_2, \ldots, r_m)$ after thresholding and updating the resultant vector we get $A_1E \in R$, Now let $B = A_1E$ we pass on B into the system $(N(E))^T$ and obtain $B(N(E))^T$. We update and threshold the vector $B(N(E))^T$ so that $B(N(E))^T \in D$.

This procedure is repeated till we get a limit cycle or a fixed point.

**DEFINITION 4.5.12:** *Finite number of NRMs can be combined together to produce the joint effect of all NRMs. Let $N(E_1), N(E_2), \ldots, N(E_r)$ be the Neutrosophic relational matrices of the NRMs with nodes $R_1, \ldots, R_m$ and $D_1, \ldots, D_n$, then the combined NRM is represented by the neutrosophic relational matrix $N(E) = N(E_1) + N(E_2) + \ldots + N(E_r)$.*

## 4.6. Neutrosophic Relational Bimaps and their Applications

Now we define Neutrosophic Relational Bimaps (NRBMs) which is a further generalization of NRMs and illustrate them with examples.

**DEFINITION 4.6.1:** *A Neutrosophic Relational Bimap (NRBM) is a neutrosophic directed bigraph or a bimap*

*from $D^t$ to $R^t$, t = 1,2 with concepts like policies or events etc as nodes and causalities as edges. It represents neutrosophic relations and causal relations between spaces $D^t$ and $R^t$.*

*Let $D_i^t$ and $R_j^t$ denote the nodes of an NRBM (t = 1, 2) The directed edge from $D_i^t$ to $R_j^t$ (t = 1, 2) denotes the causality of $D_i^t$ to $R_j^t$ called relations. Every edge in the NRBM is weighted with a number in the set {0, 1, -1, I} As in case of NRBMs we have edge values taken by $e_{ij}^t$ (t = (1, 2).*

*When the nodes of the NRBM take edge values from {0, -1, 1 I} we say the NRBM is simple NRBM.*

*Let $D_1^t, D_2^t, ..., D_{n_t}^t$ (t = 1, 2) be the nodes of the domain space $D = D_1 \cup D_2$ of the NRBM and $R_1^t, R_2^t, ..., R_{n_t}^t$ be the nodes of the range space $R = R_1 \cup R_2$ of the same NRBM.*
*Let the bimatrix $N(E_B)$ be defined as*

$$N(E_B) = N(E_1) \cup N(E_2)$$
$$= \left(e_{ij}^t\right), t = 1, 2$$
$$= \left(e_{ij}^1\right) \cup \left(e_{ij}^2\right) \text{ where } e_{ij}^t \text{ is}$$

*the weight of the directed edge $D_i^t$ $R_j^t$ (or $R_j^t$ $D_i^t$) and $e_{ij}^t \in$ {0, 1, –1 I}. $N(E_t)$ is called the Neutrosophic Relational bimatrix of the NRBM.*

It is interesting to note in almost all the cases the associated directed bigraphs will always be a bipartite bigraphs. The possible cases would be disconnected bipartite bigraphs and connected bipartite bigraphs.

**DEFINITION 4.6.2:** *Let $D_1^t, ..., D_{n_t}^t$ and $R_1^t, R_2^t, ..., R_{n_t}^t$ denote the nodes of a NRBM (t = 1, 2). Let A = $\left(a_1^t, ..., a_{n_t}^t\right)$, (t = 1,2) $a_i^t \in$ {0, 1, I} is called the neutrosophic* *instantaneous state bivector of the domain spaces $D^t = D^1 \cup D^2$ and it denotes the on – off position of the nodes at any instant. Similarly let B = $\left(b_1^t,...,b_{n_t}^t\right)$ = $B^1 \cup B^2$ = $\left(b_1^1,...,b_{n_1}^1\right) \cup \left(b_1^2,...,b_{n_2}^2\right)$, $b_i^t \in$ {0, 1, I} (t = 1, 2). B is called instantaneous state vector of the range space and position of the nodes at any instant.*
*Let $D_1^t,...,D_{n_t}^1$ and $R_1^t,...,R_{n_t}^t$ (t = 1, 2) be the nodes of a NRBM. Let $D_i^t\ R_j^t$ (or $R_j^t\ D_i^t$) be the edge of an NRBM. j = 1, 2,…, $n_t$ , t = 1, 2 and i =1, 2,…, $n_t$ , t = 1, 2. The edges form a directed cycle. An NRBM is said to be a bicycle if each of the neutrosophic graphs of the neutrosophic bigraph possess a directed cycle. A NRBM is said to be acyclic other wise.*

A NRBM with bicycles is said to be a NRBM with feedback.

When there is a feed back in the NRBM i.e. when the causal relations flow through a bicycle in a revolutionary manner the NRBM is called a Neutrosophic dynamical bisystem. Let $D_i^t.R_j^t$ (or $R_j^t\,D_i^t$) $1 \leq j \leq m_t$ , $1 \leq i \leq n_t$ (t = 1, 2) when $R_j^t$ (or $D_i^t$) is switched on and if causality flows through edges of a cycle and if it again causes $R_i$ (or $D_i$) we say that the neutrosophical dynamical system goes round and round. This is true for any node $R_j^t$ (or $D_i^t$), (t = 1, 2), $1 \leq j \leq m_t$ (or $1 \leq i \leq n_t$). The equilibrium state of this neutrosophical dynamical system is called the Neutrosophic bihidden pattern.

If the equilibrium state of a Neutrosophical dynamical system is a unique Neutrosophic state vector then it is called the fixed point. Consider the NRBM with $R_1^t,...,R_{m_t}^t$ and $D_1^t, D_2^t,...,D_{n_t}^t$ as nodes (t =1, 2). For example if we start the dynamical system by switching on $R_j^t$ (or $D_i^t$), t = 1, 2 on. Let us assume that the NRBM settles down with $R_1$ and $R_{m_t}$ (or $D_1$ and $D_{n_t}$) on, or indeterminate, or off i.e. the neutrosophic state bivector remains as (1 0 0 … 0 1) ∪ (1 0 … 0 1) or ((1 0 0 … *I* ) ∪ (1 0… *I* ) in D) this state bivector is called the fixed bipoint.

It the NRBM settles down with a state bivector repeating in the form.

$$A_1^t \to A_2^t \to \ldots \to A_{it}^t$$
$$\left(\text{or } B_1^t \to B_2^t \to \ldots \to B_{jt}^t\right)$$

(t = 1, 2) then this equilibrium is called the limit bicycle. Let $R_1^t, R_2^t, \ldots, R_{m_t}^t$ and $D_1^t, \ldots, D_{n_t}^t$, (t = 1, 2) be the nodes of the NRBM with feed back. Let N ($E_B$) be the Neutrosophic relational bimatrix. Let us find the bihidden pattern when $D_1^t$ is switched on; i.e., $A_1^t$ = (1 0 …0) ∪ (1 0…0) in $D^t$; the data should pass through the relational bimatrix

$$\begin{aligned} N(E_B) &= N(E_1) \cup N(E_2) \\ A_1^t N(E_B) &= \left(A_1^1 \cup A_1^2\right)\left(N(E_1) \cup N(E_2)\right) \\ &= A_1^1 N(E_1) \cup A_1^2 N(E_2). \end{aligned}$$

Let $A_1^t N(E_B) = \left(r_1^1, r_2^1, \ldots, r_{m_1}^1\right) \cup \left(r_1^2, r_2^2, \ldots, r_{m_2}^2\right)$ after thresholding the resultant vector (we don't need to update for we started only with the on state from the domain space), we get the resultant is in R. Now let $B^t = B^1 \cup B^2 = A_1^t$ (N $E_B$). Now we pass on $B^t$ into $[N(E_B)]^T$ so

$$\begin{aligned} B^t [N(E_B)]^T &= (B^1 \cup B^2) [N(E_1) \cup N(E_2)]^T \\ &= B [N(E_1)]^T \cup B^2 [N(E_2)]^T. \end{aligned}$$

Now we update and threshold the bivector say $C^t = C^1 \cup C^2$ now $C^t \in D^t = D^1 \cup D^2$.

This procedure is repeated till we get a limit bicycle or the fixed bipoint. Now we illustrate this using the model which studies the socio psychological feelings of HIV/AIDS patients, who are mostly uneducated migrant labourers. The disease is very much prevalent among Migrant labourers mainly because they are away from home and secondly they get CSWs for very cheap rates and substantially their job motivitates them to do so. The study of their socio and psychological problem is very important. For from our analysis most men after being affected by HIV/AIDS become sadists. They do not mind infecting, their wife so much so the CSWs. One is not able to understand this temper of these patients. Is it the drugs given to them for treatment makes them desperate / depressed to act so? A research separately in this direction is very essential.

Here we are going to concentrate only on the socio psychological problem not the question raised above.

***Example 4.6.1:*** After discussion with experts we have taken the following attributes related with HIV/AIDS patients, which is taken as the domain space D.

$D_1$ – Feeling of loneliness / aloofness
$D_2$ – Feeling of guilt
$D_3$ – Desperation / fear in public
$D_4$ – Sufferings both mental / physical
$D_5$ – Public disgrace (feeling).

The concepts / nodes related with the public are taken as the nodes of the range space.

$R_1$ – Fear of getting the disease
$R_2$ – No mind to forgive the HIV/AIDS patients sin
$R_3$ – Social stigma to have HIV/AIDS patients as a friend
$R_4$ – No sympathy.

Using these nodes two experts opinion was taken as the directed neutrosophic bigraph which follows:

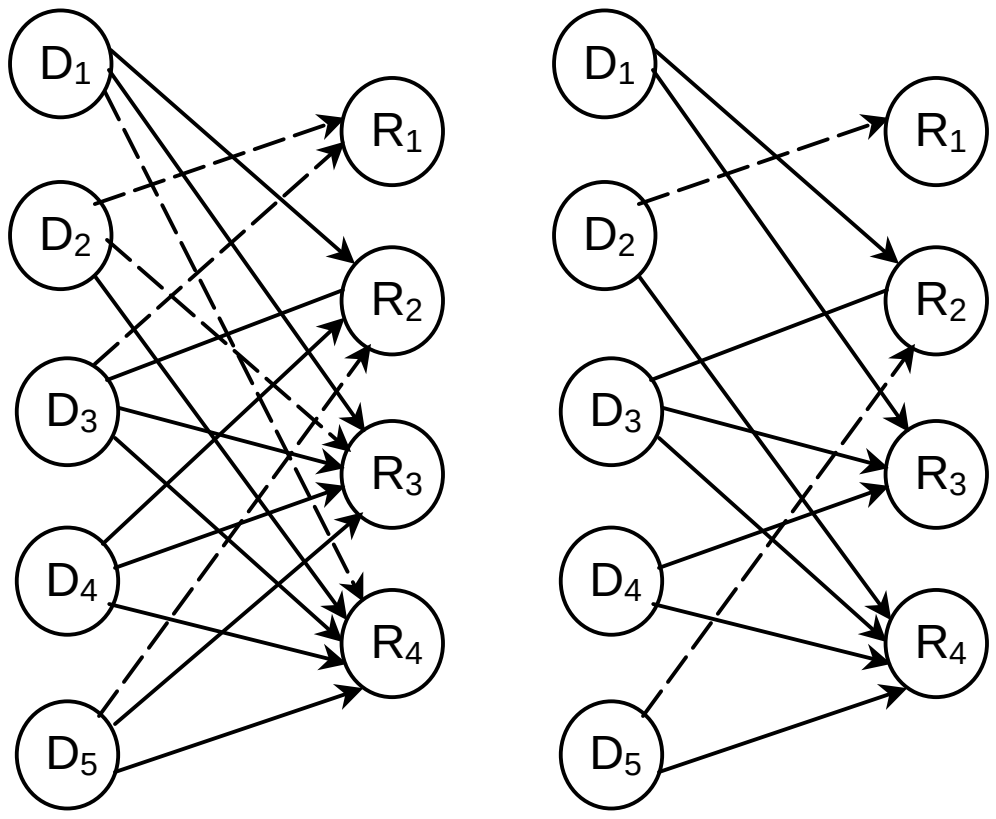

FIGURE: 4.6.1

The related adjacency connection neutrosophic bimatrix $N(E_B) = N(E_1) \cup N(E_2)$

$$\begin{bmatrix} 0 & 1 & 1 & I \\ I & 0 & I & 1 \\ I & 1 & 1 & 1 \\ 0 & 1 & 1 & 1 \\ 0 & I & 1 & 1 \end{bmatrix} \cup \begin{bmatrix} 0 & 1 & 1 & 0 \\ I & 0 & 0 & 1 \\ 0 & 1 & 1 & 1 \\ 0 & 0 & 1 & 1 \\ 0 & I & 0 & 1 \end{bmatrix}$$

Clearly we have a $5 \times 4$ rectangular neutrosophic bimatrix. Suppose we consider the instantaneous state bivector

$$\begin{aligned} A &= (0\ 1\ 0\ 0\ 0) \cup (0\ 1\ 0\ 0\ 0) \\ &= A^1 \cup A^2 \end{aligned}$$

the feeling of guilt alone is in the on state and all other nodes are in the off state. Now the effect of A on the neutrosophic dynamical systems $N(E_B)$ is as follows.

$$AN(E_B) = (A^1 \cup A^2)\ (N\ (E_1) \cup N(E_2)$$

$$= \quad A^1 N(E_1) \cup A^2 N(E_2)$$
$$= \quad (I\ 0\ I\ 1) \cup (I\ 0\ 0\ 1).$$

After thresholding we get the resultant bivector as

$$B = (I\ 0\ I\ 1) \cup (I\ 0\ 0\ 1)$$
$$= B^1 \cup B^2.$$

The effect of B on the Neutrosophic dynamical system $N(E_B)$ as follows:

$$BN(E_B)^T = (B^1 \cup B^2)(N(E_1) \cup N(E_2))^T$$
$$= B^1 N(E_1)^T \cup B^2 N(E_2)^T$$

$$= (2I, 2I+1, 2I+1, I+1, I+1) \cup (0, 1+I, 1, 0, 1).$$

After thresholding and updating using the fact ($2 + I = 1$ and $2I + 1 = I$ as the thresholding function) we get

$$BN(E_B)^T = (I\ 1\ I\ 1\ 1) \cup (0\ 1\ 1\ 0\ 1).$$

Now we proceed on to get a fixed bipoint or a limit bicycle. We see the related bigraph is a neutrosophic bipartite bigraph which is a disconnected one in this case.

Now we would give one more model in which the neutrosophic bipartite bigraph is biconnected. For this we take the problem of women empowerment relative to HIV/AIDS patients and the influence of the public on them.

***Example 4.6.2:*** Since the very concept of women's empowerment and community mobilization in the context of AIDS epidemic is an unsupervised data having no specific statistical values and more so is the opinion of public on these AIDS epidemic affected women. Thus this is modeled using an expert opinion. The related neutrosophic bigraph happens to be a neutrosophic bipartite bigraph.

The attributes related with the women's empowerment and community mobilization are taken as

$W_1$ – Gender balance
$W_2$ – Cost effectiveness
$W_3$ – Large scale operation
$W_4$ – Social service .

The attributes of the range space which is related to HIV/AIDS epidemic are

$A_1$ – Care for the AIDS infected persons
$A_2$ – Prevention of spread of HIV/AIDS epidemic
$A_3$ – Creation of awareness about AIDS
$A_4$ – Medical treatment of AIDS patients
$A_5$ – Social stigma.

The views of public and close kith and kin which make the AIDS / HIV infected women to lead a normal life.

$P_1$ – Public look down upon HIV/AIDS women in major cases as only women with no character or CSW's alone get this disease.
$P_2$ – Close relatives send away the women from the family even if they know for certain only their husbands have infected them.
$P_3$ – These women are left to the mercy of road or in the public hospitals with no one to take care of them.
$P_4$ – No proper women organization ever take up this issue.
$P_5$ – As they have no money and no support they cannot take any proper treatment.

Now relating these three nodes we give the neutrosophic directed bigraph given by figure 4.6.1.

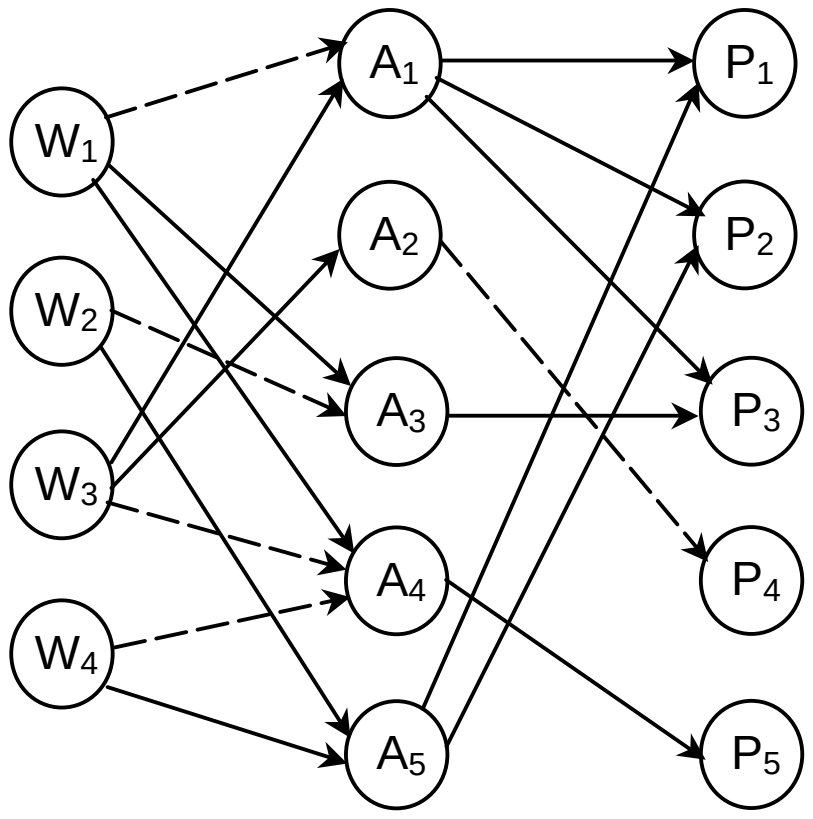

FIGURE: 4.6.2

Clearly this neutrosophic bigraph is biconnected bipartite bigraph.

The neutrosophic bimatrix related with this bigraph is $N(E_B) = N(E_1) \cup N(E_2)$

$$\begin{array}{c} \\ W_1 \\ W_2 \\ W_3 \\ W_4 \end{array}\begin{array}{c} \begin{array}{ccccc} A_1 & A_2 & A_3 & A_4 & A_5 \end{array} \\ \begin{bmatrix} I & 0 & 1 & 1 & 0 \\ 0 & 0 & I & 0 & 1 \\ 1 & 1 & 0 & I & 0 \\ 0 & 0 & 0 & I & 1 \end{bmatrix} \end{array} \cup \begin{array}{c} \\ A_1 \\ A_2 \\ A_3 \\ A_4 \\ A_5 \end{array}\begin{array}{c} \begin{array}{ccccc} p_1 & p_2 & p_3 & p_4 & p_5 \end{array} \\ \begin{bmatrix} -1 & 1 & 1 & 0 & 0 \\ 0 & 0 & 0 & I & 0 \\ 0 & 0 & 1 & 0 & 0 \\ 0 & 0 & 0 & 0 & 1 \\ 1 & 1 & 0 & 0 & 0 \end{bmatrix} \end{array}.$$

Using this neutrosophic bimatrix we can find the effect of any state bivector on the neutrosophical dynamic system $N(E_T)$.

Suppose

$$\begin{aligned} A &= (1\ 0\ 0\ 0) \cup (0\ 0\ 1\ 0\ 0) \\ &= A_1 \cup A_2. \end{aligned}$$

is the state bivector whose effect we have to study on the dynamical system $N(E_B)$.

$$AN(E_B) = (A_1 \cup A_2)\ (N(E_1) \cup N(E_2))$$
$$= A_1N(E_1) \cup A_2N(E_2)$$
$$= (I\ 0\ 1\ 1\ 0) \cup (0\ 0\ 1\ 0\ 0).$$

After updating and thresholding we get the resultant bivector; B = (10110) ∪ (0 0 1 0 0). Now the effect of B on $N(E_T)^T$ is given by

$$BN(E_B)^T = B_1\ N(E_1)^T \cup B_2\ N(E_2)^T$$
$$= (I + 2, I, 2I, I) \cup (1\ 0\ 1\ 0\ 0).$$

After updating and thresholding we get the resultant vector as

$$C = (1\ I\ I\ I) \cup (1\ 0\ 1\ 0\ 0).$$

The effect of C on N $(E_T)$ is given by

$$CN(E_T) = C_1\ N(E_1) \cup C_2\ N(E_2)$$
$$= (2\ I, I, 1 + I\ 1 + 2I\ 2I\ ) \cup (0\ 1\ 2\ 0\ 0).$$

After thresholding and updating we get the resultant neutrosophic state bivector as

$$D = (1\ I\ I\ I\ I) \cup (0\ 1\ 1\ 0\ 0).$$

Now we study the effect of D on N $(E_T)^T$ and so on. Thus we proceed in the some way and arrive at the fixed bipoint or the limit bicycle.

Now this model can further be extended as Neutrosophic Relational trimaps or say any neutrosophic relational n-maps. This model when n = 1 is nothing but the usual neutrosophic relational maps (NRM) when n = 2 we have the neutrosophic relational bimaps (NRBMs) when n = 3 we have the neutrosophic relational trimaps. (NRTMs) and so on. Thus when we want to make a comparative study say r-experts on a same set of attributes, we can adopt the NR r-Ms. (i.e. Neutrosophic relational r-maps). Now we illustrate this by a very simple model when r = 4.

***Example 4.6.3:*** Suppose we are interested in studying the problem of health hazards faced by agricultural labourers and the related types of chemical pollution faced by them. Suppose the experts concur on the four types of pollution.

$C_1$ – Use of chemicals banned by other court
$C_2$ – Pollution of air
$C_3$ – Pollution of soil
$C_4$ – Pollution of food crops.

Now the health problems faced by the agricultural labourers.

$H_1$ – Skin ailment / skin caner
$H_2$ – Indigestion / loss of appetite
$H_3$ – Blurred visions and problems of eye
$H_4$ – Head ache / giddiness / Tension
$H_5$ – Cooked food gets spoiled in a very short duration.

We now give the opinion of the four experts and their related directed neutrosophic 4-graphs, which is bipartite is given by the following figure:

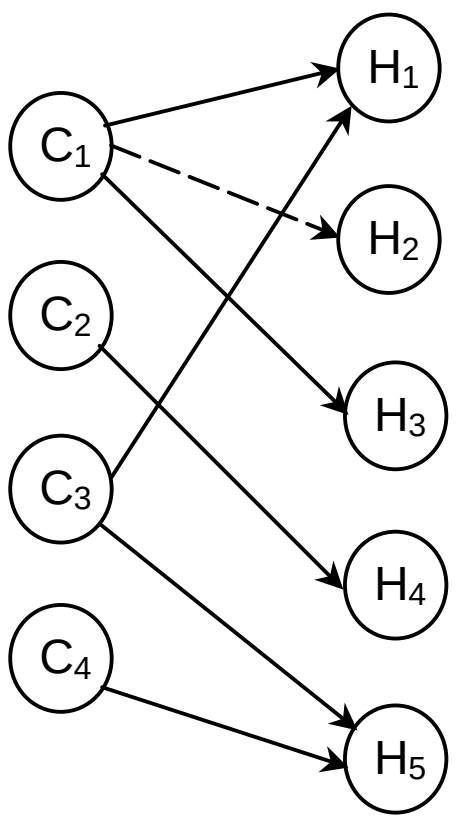

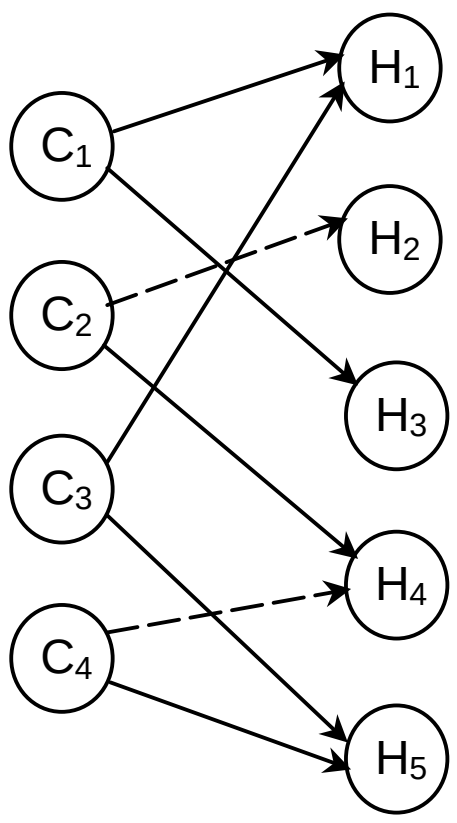

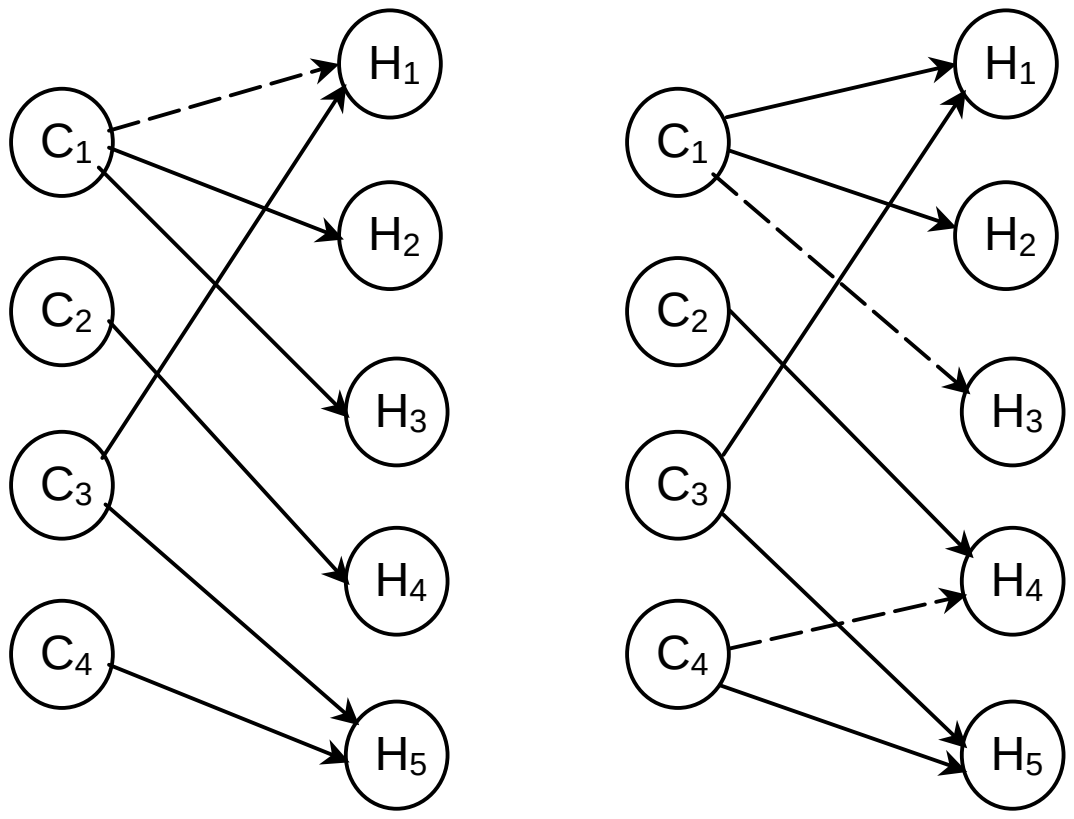

FIGURE: 4.6.4

The related neutrosophic 4-matrix is given below

$$N(E_Q) = N(E_1) \cup N(E_2) \cup N(E_3) \cup N(E_4)$$

$$\begin{bmatrix} 1 & I & 1 & 0 & 0 \\ 0 & 0 & 0 & 1 & 0 \\ 1 & 0 & 0 & 0 & 1 \\ 0 & 0 & 0 & 0 & 1 \end{bmatrix} \cup \begin{bmatrix} 1 & 0 & 1 & 0 & 0 \\ 0 & I & 0 & 1 & 0 \\ 1 & 0 & 0 & 0 & 1 \\ 0 & 0 & 0 & I & 1 \end{bmatrix}$$

$$\cup \begin{bmatrix} I & 1 & 1 & 0 & 0 \\ 0 & 0 & 0 & 1 & 0 \\ 1 & 0 & 0 & 0 & 1 \\ 0 & 0 & 0 & 0 & 1 \end{bmatrix} \cup \begin{bmatrix} 1 & 1 & I & 0 & 0 \\ 0 & 0 & 0 & 1 & 0 \\ 1 & 0 & 0 & 0 & 1 \\ 0 & 0 & 0 & I & I \end{bmatrix}$$

Now we just illustrate how this FR4-M works. Any state quadruple vector say

$$\begin{aligned} A &= A_1 \cup A_2 \cup A_3 \cup A_4 \\ &= (1\ 0\ 0\ 0) \cup (1\ 0\ 0\ 0) \cup \\ &\quad\ (1\ 0\ 0\ 0) \cup (1\ 0\ 0\ 0) \end{aligned}$$

is taken i.e. all nodes are off except the node use of chemicals banned by other countries alone is in the on state. Effect of A on the dynamical system N ($E_Q$) is given by

$$AN(E_Q) = A_1N(E_1) \cup A_2N(E_2) \cup A_3N(E_3) \cup A_4N(E_4)$$

$$= (1\ I\ 1\ 0\ 0) \cup (1\ 0\ 1\ 0\ 0) \cup (I\ 1\ 1\ 0\ 0) \cup (1\ 1\ I\ 0\ 0) = B$$

that is $B = B_1 \cup B_2 \cup B_3 \cup B_4$.

The effect of B on N ($E_Q$) is given by

$$BN(E_Q)^T = B_1N(E_1)^T \cup B_2N(E_2)^T \cup B_3N(E_3)^T \cup B_4N(E_4)^T$$

$$= (2 + I, 0, 1, 0) \cup (2\ 0\ 1\ 0) \cup (2 + I, 0\ I\ 0) \cup (2 + I\ 0\ 1\ 0).$$

After thresholding and updating we get the resultant vector as

$$C = (1\ 0\ 1\ 0) \cup (1\ 0\ 1\ 0) \cup (1\ 0\ I\ 0) \cup (1\ 0\ 1\ 0)$$

$$= C_1 \cup C_2 \cup C_3 \cup C_4.$$

Now we study the effect of C on the dynamical system N ($E_Q$).

$$CN(E_Q) = C_1N(E_1) \cup C_2N(E_2) \cup C_3N(E_3) \cup C_4N(E_4)$$

$$= (2\ I\ 1\ 0\ 1) \cup (2\ 0\ 1\ 0\ 1) \cup (2\ I\ 1\ 1\ 0\ I) \cup (2\ 1\ I\ 0\ 1).$$

After thresholding the resultant n-vector we get the r-vector

$$S = (1\ I\ 1\ 0\ 1) \cup (1\ 0\ 1\ 0\ 1)$$

$\cup (I\ 1\ 1\ 0\ I) \cup (1\ 1\ I\ 0\ 1)$

$= S_1 \cup S_2 \cup S_3 \cup S_4.$

Now we see the effect of S on the dynamical system $N(E_Q)^T$

$$S\,(N(E_Q))^T = S_1N(E_1)^T \cup S_2N(E_2)^T \cup S_3N(E_3)^T \cup S_4(N(E_4))^T$$

$$= (2 + I, 0\ 2, 1) \cup (2\ 0\ 2\ 1) \cup (2 + I\ 0\ 2\ I\ I) \cup (2 + I\ 0\ 2\ 0).$$

Thresholding and updating it we get

$$V = (1\ 0\ 1\ 1) \cup (1\ 0\ 1\ 1) \cup (1\ 0\ I\ I) \cup (1\ 0\ 1\ 1\ 0)$$

$$= V_1 \cup V_2 \cup V_3 \cup V_4.$$

Now we study the effect of V on N $(E_Q)$.

$$SN(E_Q) = S_1N(E_1) \cup S_2N(E_2) \cup S_3N(E_3) \cup S_4N(E_4)$$

$$= (2\ I\ 1\ 0\ 2) \cup (2\ 0\ 1\ I\ 2) \cup (1 + I\ 1\ 1\ 0\ 2\ I) \cup (2\ 1\ I\ 0\ 1).$$

After thresholding we get the resultant quad vector as

$$X = (1\ I\ 1\ 0\ 1) \cup (1\ 0\ 1\ I\ 1) \cup (I\ 1\ 1\ 0\ I) \cup (1\ 1\ I\ 0\ 1).$$

We proceed on till we get a fixed quadpoint or a limit quadcycle of the system.

Thus the neutrosophic disconnected bipartite r-graph can give several other conclusions.

We have taken in this example r = 4 it can also happen the r-graph can be biconnected or even triconnected such model are also possible.

# BIBLIOGRAPHY

This book has included in the references the available papers and books on Fuzzy Cognitive Maps and Neutrosophic Cognitive Maps and other related papers with a main motivation that any innovative reader can apply the bimatrices and their extensions to these models.

# INDEX

**C**

**D**

**E**

**F**

**G**

**H**

**I**

**K**

**L**

**M**

**N**

## O

## P

## Q

## R

## S

**T**

**V**

**W**

**Z**

# About the Authors

**Dr.W.B.Vasantha Kandasamy** is an Associate Professor in the Department of Mathematics, Indian Institute of Technology Madras, Chennai, where she lives with her husband Dr.K.Kandasamy and daughters Meena and Kama. Her current interests include Smarandache algebraic structures, fuzzy theory, coding/ communication theory. In the past decade she has guided nine Ph.D. scholars in the different fields of non-associative algebras, algebraic coding theory, transportation theory, fuzzy groups, and applications of fuzzy theory of the problems faced in chemical industries and cement industries. Currently, six Ph.D. scholars are working under her guidance. She has to her credit 287 research papers of which 209 are individually authored. Apart from this, she and her students have presented around 329 papers in national and international conferences. She teaches both undergraduate and post-graduate students and has guided over 45 M.Sc. and M.Tech. projects. She has worked in collaboration projects with the Indian Space Research Organization and with the Tamil Nadu State AIDS Control Society. This is her 20th book.

She can be contacted at vasantha@iitm.ac.in
You can visit her work on the web at: http://mat.iitm.ac.in/~wbv

---

Dr.Florentin Smarandache is an Associate Professor of Mathematics at the University of New Mexico in USA. He published over 75 books and 100 articles and notes in mathematics, physics, philosophy, psychology, literature, rebus. In mathematics his research is in number theory, non-Euclidean geometry, synthetic geometry, algebraic structures, statistics, neutrosophic logic and set (generalizations of fuzzy logic and set respectively), neutrosophic probability (generalization of classical and imprecise probability). Also, small contributions to nuclear and particle physics, information fusion, neutrosophy (a generalization of dialectics), law of sensations and stimuli, etc.).

He can be contacted at smarand@unm.edu

---

**K. Ilanthenral** is the editor of The Maths Tiger, Quarterly Journal of Maths. She can be contacted at ilanthenral@gmail.com